\documentclass[12pt,reqno]{amsart}
\usepackage[a4paper,margin=30mm,headheight=15pt]{geometry}
\usepackage{amsmath,amssymb,amsthm,mathtools,mathrsfs,tikz-cd,enumitem,booktabs,array,longtable,microtype,placeins}
\DeclareFontFamilySubstitution{TS1}{cmr}{lmr}
\usetikzlibrary{decorations.pathmorphing,arrows.meta}
\definecolor{qIIteal}{HTML}{176568}
\definecolor{qIIclay}{HTML}{A45224}
\definecolor{qIIblue}{HTML}{4C5E87}
\usepackage[hidelinks,bookmarksnumbered]{hyperref}
\setlist[enumerate]{leftmargin=*,itemsep=3pt,topsep=5pt}
\allowdisplaybreaks[2]
\numberwithin{equation}{section}
\newtheorem{theorem}{Theorem}[section]
\newtheorem{proposition}[theorem]{Proposition}
\newtheorem{lemma}[theorem]{Lemma}
\newtheorem{corollary}[theorem]{Corollary}
\theoremstyle{definition}

\newtheorem{conjecture}[theorem]{Conjecture}
\theoremstyle{remark}

\newcommand{\OO}{\mathcal O}\newcommand{\XX}{\mathcal X}
\newcommand{\MM}{\mathcal M}\newcommand{\EE}{\mathfrak E}\newcommand{\GG}{\mathbb G_m}
\newcommand{\CC}{\mathbb C}\newcommand{\PP}{\mathbb P}\renewcommand{\AA}{\mathbb A}\newcommand{\QQ}{\mathbb Q}
\newcommand{\LL}{\mathbf L}\newcommand{\RR}{\mathbf R}\newcommand{\vir}{\mathrm{vir}}

\newcommand{\rel}{\mathrm{rel}}\newcommand{\perf}{\operatorname{Perf}}\newcommand{\Pic}{\operatorname{Pic}}
\newcommand{\Ext}{\operatorname{Ext}}\newcommand{\ext}{\operatorname{ext}}\newcommand{\Hom}{\operatorname{Hom}}
\newcommand{\End}{\operatorname{End}}\newcommand{\Supp}{\operatorname{Supp}}
\newcommand{\Spec}{\operatorname{Spec}}\newcommand{\ch}{\operatorname{ch}}\newcommand{\rk}{\operatorname{rk}}
\newcommand{\depth}{\operatorname{depth}}
\newcommand{\RHom}{\mathbf R\!\operatorname{Hom}}\newcommand{\RG}{\mathbf R\Gamma}
\newcommand{\fib}{\operatorname{fib}}\newcommand{\cofib}{\operatorname{cofib}}
\newcommand{\spc}{\operatorname{sp}}\newcommand{\id}{\operatorname{id}}
\newcommand{\length}{\operatorname{length}}
\newcommand{\DT}{\operatorname{DT}}\newcommand{\dt}{\operatorname{DT}}\newcommand{\Hilb}{\operatorname{Hilb}}
\hypersetup{pdftitle={Surface sheaves on the quintic II. Donaldson--Thomas invariants and the quintic elliptic genus},pdfauthor={Artan Sheshmani}}
\newcommand{\companionref}[1]{\textup{\csname companion@#1\endcsname}}
\expandafter\def\csname companion@eq:Hallloggeneral\endcsname{I.2.4}
\expandafter\def\csname companion@eq:POTdiagram\endcsname{I.14.3}
\expandafter\def\csname companion@cor:noseparatedall\endcsname{I.2.8}
\expandafter\def\csname companion@prop:unlabelledpolarization\endcsname{I.5.2}
\expandafter\def\csname companion@eq:cutdata\endcsname{I.9.5}
\expandafter\def\csname companion@eq:generalEulerpairing\endcsname{I.2.5}
\expandafter\def\csname companion@eq:generalcomponentclasses\endcsname{I.8.3}
\expandafter\def\csname companion@eq:generalpolynomial\endcsname{I.2.1}
\expandafter\def\csname companion@eq:interiorres\endcsname{I.9.8}
\expandafter\def\csname companion@eq:mainfull\endcsname{I.1.6}
\expandafter\def\csname companion@eq:vacuumBzeroeq\endcsname{I.16.14}
\expandafter\def\csname companion@eq:vacuumblowup\endcsname{I.16.1}
\expandafter\def\csname companion@eq:vacuumgraphpolys\endcsname{I.16.13}
\expandafter\def\csname companion@eq:vacuumlocalpolarization\endcsname{I.16.12}
\expandafter\def\csname companion@eq:vacuumordinarypolys\endcsname{I.16.11}
\expandafter\def\csname companion@eq:vacuumpositiveproduct\endcsname{I.16.15}
\expandafter\def\csname companion@eq:vacuumrelativefactors\endcsname{I.16.3}
\expandafter\def\csname companion@eq:vacuumsection\endcsname{I.16.7}
\expandafter\def\csname companion@eq:vacuumstrict\endcsname{I.16.8}
\expandafter\def\csname companion@eq:vacuumtchart\endcsname{I.16.4}
\expandafter\def\csname companion@hyp:compact\endcsname{I.5.3}
\expandafter\def\csname companion@hyp:product\endcsname{I.14.1}
\expandafter\def\csname companion@lem:vacuumgenericfull\endcsname{I.16.6}
\expandafter\def\csname companion@lem:vacuumlinearsystems\endcsname{I.16.9}
\expandafter\def\csname companion@lem:virtualnormalization\endcsname{I.15.2}
\expandafter\def\csname companion@prop:Jacobian\endcsname{I.13.5}
\expandafter\def\csname companion@prop:POTbasechange\endcsname{I.12.4}
\expandafter\def\csname companion@prop:determultiple\endcsname{I.2.7}
\expandafter\def\csname companion@prop:divisorial-hull27\endcsname{I.9.16}
\expandafter\def\csname companion@prop:finitecharges\endcsname{I.5.9}
\expandafter\def\csname companion@prop:generalprofiles\endcsname{I.8.2}
\expandafter\def\csname companion@prop:hull-punctual-DVR\endcsname{I.16.1}
\expandafter\def\csname companion@prop:noordinaryproperextension\endcsname{I.2.12}
\expandafter\def\csname companion@prop:primitive\endcsname{I.2.9}
\expandafter\def\csname companion@prop:redglue\endcsname{I.11.2}
\expandafter\def\csname companion@prop:relativePOT\endcsname{I.13.2}
\expandafter\def\csname companion@prop:relativevdim\endcsname{I.13.3}
\expandafter\def\csname companion@prop:structurecohom\endcsname{I.3.3}
\expandafter\def\csname companion@prop:support-separated27\endcsname{I.5.10}
\expandafter\def\csname companion@rem:multiplicities\endcsname{I.15.4}
\expandafter\def\csname companion@sec:boundaryisomproper\endcsname{I.14.2}
\expandafter\def\csname companion@sec:compactification\endcsname{I.5}
\expandafter\def\csname companion@sec:diagonal\endcsname{I.14}
\expandafter\def\csname companion@sec:geometry\endcsname{I.3}
\expandafter\def\csname companion@sec:hall\endcsname{I.2.4}
\expandafter\def\csname companion@sec:reduction\endcsname{I.11}
\expandafter\def\csname companion@sec:smoothmoduli\endcsname{I.2}
\expandafter\def\csname companion@sec:verifiedcompactification\endcsname{I.16.2}
\expandafter\def\csname companion@sec:zerorelative\endcsname{I.13.3}
\expandafter\def\csname companion@thm:ambientPOT\endcsname{I.12.2}
\expandafter\def\csname companion@thm:cutcriterion\endcsname{I.9.4}
\expandafter\def\csname companion@thm:diagonal\endcsname{I.14.2}
\expandafter\def\csname companion@thm:generalCMhull\endcsname{I.2.4}
\expandafter\def\csname companion@thm:generalpolynomialstability\endcsname{I.9.2}
\expandafter\def\csname companion@thm:main\endcsname{I.1.2}
\expandafter\def\csname companion@thm:matrixexpandedallr\endcsname{I.9.14}
\expandafter\def\csname companion@thm:primitiveprofiles\endcsname{I.8.4}
\expandafter\def\csname companion@thm:verifiedcompactification\endcsname{I.16.5}
\expandafter\def\csname companion@thm:zeroboundaryPOT\endcsname{I.13.8}
\title[Surface sheaves on the quintic II]{Surface sheaves on the quintic II.\\ Donaldson--Thomas invariants and the M5-brane modified elliptic genus in string theory}
\author{Artan Sheshmani}
\address{Artan Sheshmani}
\address{Beijing Institute of Mathematical Sciences and Applications, A6, Room 205 No. 544, Hefangkou Village, Huaibei Town, Huairou District, Beijing 101408, China}
\address{Massachusetts Institute of Technology (MIT), IAiFi Institute, 182 Memorial Drive, Cambridge, MA 02139, USA}
\address{National Research University, Higher School of Economics, Russian Federation, Laboratory of Mirror Symmetry, NRU HSE, 6 Usacheva Street, Moscow, Russia, 119048}
\email{artan@mit.edu}
\date{17 November 2026}
\begin{document}
\begin{abstract}
We compute Donaldson--Thomas invariants of pure two-dimensional sheaves in the hyperplane class of a smooth quintic threefold. We take advantage of two degenerating families to relate the DT geometry of Quintic threefold to geometry of pure two dimensional shaves on a quartic threefold and correction terms induced by invariants of relative Hilbert schemes of points. For this we use degeneration technique and virtual gluing formula in the prequel Part I. The marked virtual intersections, Chow specialization and complete moduli on the smooth quintic give nine coefficients of the M5-brane modified elliptic genus in string theory predicted by Gaiotto, Strominger and Yin, including its entire primitive polar part. Four geometric degrees determine the candidate in the prescribed modular space. We also construct a weight-$26$ period series on the complete hyperplane system and identify the contributions of twisted-cubic hulls and their point modifications on both sides. 

\end{abstract}\maketitle\tableofcontents
\section{Introduction}
Let $X\subset\PP^4$ be a smooth quintic threefold and $H$ its hyperplane class. A pure sheaf of divisor class $H$ is the pushforward of a rank-one torsion-free sheaf on an integral hyperplane section $T$. The section may be singular and varies in $|H|$. Gaiotto, Strominger and Yin predict the virtual degrees of these sheaves, for every curve and point charge, from the modified elliptic genus of a single M5-brane \cite[\S3.3]{GSY}.

The two degenerations share a quartic threefold $Q$ with a quartic K3 section $S$. To compare their contributions one must first know what is carried from $Q$ to $S$ when a supporting surface moves onto an expansion. Part~I \cite{Companion} extends the marked-divisor construction of Li and Wu to perfect surface sheaves. On a fixed union the gluing datum is an isomorphism
\[
       \mathbf L i_1^*F_1\simeq\mathbf L i_2^*F_2
                    \quad\text{in }\operatorname{Perf}(S).
\]
Taking endomorphisms gives the difference of the two restriction maps. Its homotopy fiber recovers the deformations and obstructions of the glued sheaf, including the motions of its support. Under the compactification hypotheses of Theorem~\companionref{thm:main}, marking a join therefore specializes the virtual class of the surface sheaves themselves. When the marked space is the full fiber product of independently proper relative factors over a smooth moduli space of simple boundary objects, compatible obstruction theories give its class by refined diagonal pullback. References prefixed by I are to Part~I. For derived Quot schemes and categorified DT invariants see \cite{KryczkaSheshmani2026}.

Consider the first point defect $I_{p/T}$ on a moving hyperplane section. The proper family of supports in Theorem~\companionref{thm:verifiedcompactification} retains the graph on a ruled component when the hyperplane tends to $S$. A point tending to a join can require one further ruled component. In Theorem~\ref{thm:one-defect-global} we identify the expanded incidence of the support and the point with a proper, separated expanded sheaf space containing the complete generic moduli. The identification includes the infinitesimal deformations. Its positive marked classes are restrictions of K3 diagonal pullbacks (Theorem~\ref{thm:point-open-diagonal}). The class with zero boundary comes from the complementary root. Their degrees give
\[
              -320-360-120=-800 .
\]
The virtual class provided by Part~I is what allows these three incidence calculations to be added as contributions to the same DT invariant.

For conics the same quartic class meets two boundary conditions on $S$. Fixing the support curve gives $5016$; requiring a double point centered on that curve, with arbitrary tangent direction, gives $-26\,560$. Together these intersections recover the contribution of the blowup of $\PP^3$.

Table~\ref{tab:intro-charge-comparison} gives nine complete primitive degrees from marked intersections, Chow specialization and the geometry of the smooth quintic. Gholampour, Sheshmani and Thomas computed $-800$ and $8625$ for the point and line incidences in \cite[Example~1.1]{GST}. Here those incidences exhaust their charges. The nine degrees agree with the predicted elliptic genus and determine its entire polar part. In weight $-3/2$ that part determines the series once its transformation law is proved.

On the complete hyperplane system Section~\ref{sec:quintic-period-comparison-new} constructs a modular period series $\mathcal N_X$ of weight $26$. Equation~\eqref{eq:quintic-cubic-period-sector} identifies the contributions of the twisted-cubic hulls, their positive powers and all their point modifications in both $\mathcal N_X$ and $\eta^{55}h^{\mathrm{DT}}_1$. Equation~\eqref{eq:quintic-two-boundary} gives the first two neutral boundary degrees. The full comparison still requires the virtual multiplicities at excess Hodge intersections and on singular sections.

\subsection{Which invariant is being counted?}
More generally, we fix a pure two-dimensional sheaf $F$ of character
\begin{equation}\label{eq:charge}
       v=\ch(F)=(0,rH,\beta,n),\qquad r\geq 1.
\end{equation}
Purity means that every nonzero subsheaf has two-dimensional support. In particular the divisor counted with multiplicity is
\[
       [F]_{\rm div}=\sum_{T}\length_{\OO_{X,\eta_T}}(F_{\eta_T})[T],\qquad [F]_{\rm div}\sim rH,
\]
where $T$ runs over the irreducible components of the support and $\eta_T$ is its generic point. A divisor in $|rH|$ is cut out on $X$ by a form of degree $r$. For $r=1$ it is integral, its equation annihilates $F$, and $F=i_*G$ for a rank-one torsion-free sheaf $G$ on it (Proposition~\companionref{prop:primitive}). We call this the \emph{primitive divisor class}.

For $r>1$ the support may have several components or nonreduced structure, and the generic rank on a component may exceed one. All ample divisor classes occur since $\Pic(X)=\mathbb Z H$.

Write $d=\int_X H\beta$, $n=\int_X\ch_3(F)$ and $\chi=\chi(X,F)$. As $H^3=5$ and $c_2(X)=10H^2$, Riemann--Roch gives
\[
 P_F(m)=\chi(X,F(mH))=\frac{5r}{2}m^2+dm+\chi,
 \qquad \chi=n+\frac{25r}{6}.
\]
For the reduced Hilbert polynomial $p_F=P_F/\mathrm{lc}(P_F)$, stability requires $p_A(m)<p_F(m)$ for every $0\ne A\subsetneq F$ and $m\gg0$. Semistability allows equality.

When every semistable sheaf of character $v$ is stable, the moduli space $M_X(v)$ is projective. At $[F]$ its tangent and obstruction spaces are $\Ext^1_X(F,F)$ and $\Ext^2_X(F,F)$. Serre duality on the quintic identifies one with the dual of the other. The virtual dimension is zero, and the invariant is
\[
 \Omega_X(v)=\deg[M_X(v)]^{\vir}
     =\sum_{j\in\mathbb Z}j\,e_c(\nu_{M_X(v)}^{-1}(j)).
\]
Here $\nu$ is the Behrend function and $e_c$ denotes compactly supported Euler characteristic \cite[Theorem~4.18]{Behrend27}. In general we use the rational Joyce--Song invariant $\overline\Omega_X(v)$ of the semistable stack \cite[\S5]{JS}. It equals $\Omega_X(v)$ when semistability implies stability.

Sections~\companionref{sec:smoothmoduli} and~\companionref{sec:hall} give the deformation theory and the two definitions.

\subsection{The partition function and its predicted transformation}
Tensoring by $\OO_X(aH)$ gives an isomorphism of moduli problems. Multiplication of the Chern character by $e^{aH}$ sends
\[
          (d,n)\longmapsto(d+5ra,\ n+ad+5ra^2/2).
\]
To count each orbit once, put $q_1=-d$, $q_0=n+25r/12$ and set
\[
 \mu=q_1-5r^2/2\pmod{5r},\qquad x=\frac{q_1^2}{10r}-q_0.
\]
Both $\mu$ and $x$ are invariant under the twist; $\mu$ is integral. We call $\mu=0$ the neutral component. Each coefficient counts one twist orbit (Proposition~\ref{prop:generalspectralflow}). With $q=e^{2\pi i\tau}$ put
\begin{equation}\label{eq:intro-partition}
 h^{\mathrm{DT}}_r=(h^{\mathrm{DT}}_{r,\mu})_{\mu\bmod 5r},\qquad
 h^{\mathrm{DT}}_{r,\mu}(\tau)=\sum_x\overline\Omega_{r,\mu}(-x)q^x .
\end{equation}
The reduced D0 charge is $\widehat q_0=-x$. Proposition~\ref{prop:generalvacuumbound} bounds the polar part, the sum over $x<0$.

For one hyperplane let $Z^{\mathrm{GSY}}=(Z^{\mathrm{GSY}}_\mu)_{\mu\bmod5}$ be the coefficient vector of Gaiotto, Strominger and Yin \cite[\S3.3, (3.10)--(3.11), and Appendix~A]{GSY}. Their modified elliptic genus is $\sum_{\mu=0}^4Z^{\mathrm{GSY}}_\mu(\tau)\Theta_\mu(\bar\tau,y)$, where $y$ records curve charge. The prediction is $h^{\mathrm{DT}}_1=Z^{\mathrm{GSY}}$.

In particular the series should converge on the upper half-plane, have at most a pole at the cusp, and satisfy
\begin{align*}
 h^{\mathrm{DT}}_{1,\mu}(\tau+1)
  &=e^{2\pi i(-55/24+3\mu^2/5)}h^{\mathrm{DT}}_{1,\mu}(\tau),\\
 h^{\mathrm{DT}}_{1,\mu}(-1/\tau)
  &=\frac{(-i\tau)^{-3/2}}{\sqrt5}
          \sum_{\nu\bmod5}e^{-2\pi i\mu\nu/5}h^{\mathrm{DT}}_{1,\nu}(\tau).
\end{align*}
The exponent congruence gives the $T$-transformation; the $S$-transformation couples the five residues \cite[\S3.3, (3.6)--(3.8)]{GSY}. For $r>1$ the predicted completion includes nonholomorphic contributions from smaller divisor multiplicities (Conjecture~\ref{conj:generalSduality}, following \cite[\S\S2.4--2.5]{AFKPmock}).

One symmetry already follows from the sheaves. Theorem~\ref{thm:rank-zero-charge-conjugation} proves, including in classes with strictly semistable sheaves,
\[
 \overline\Omega_X(0,rH,\beta,n)=\overline\Omega_X(0,rH,-\beta,n),
 \qquad h^{\mathrm{DT}}_{r,\mu}=h^{\mathrm{DT}}_{r,-\mu}.
\]
Indeed the shifted dual $\RR\mathcal Hom_X(F,\OO_X)[1]$ has character $(0,rH,-\beta,n)$. It can have a second cohomology sheaf $\mathcal Ext_X^2(F,\OO_X)$ in degree one, supported at finitely many points. Section~\ref{sec:hall-charge-conjugation} keeps track of this finite-length sheaf by counting extensions in the Hall algebra. Its Euler pairing with a surface sheaf vanishes, and the integrated identity gives the equality above.

In primitive class there are consequently three components to compare, represented by $\mu=0,1,2$.

For $r>1$, the exponent also tells us when strictly semistable sheaves can contribute. Put $c_r=5r^3+50r$. Proposition~\ref{prop:generalvacuumbound} gives $x\ge -c_r/24$. If a semistable sheaf has at least two stable factors with the same reduced Hilbert polynomial, Proposition~\ref{prop:step22-jhgap} sharpens this to
\[
                 x+\frac{c_r}{24}\ \ge\ \frac{5r(r-1)}8 .
\]
Thus a nonempty moduli space below this second bound consists entirely of stable sheaves. Section~\ref{sec:hall22} obtains some coefficients in higher multiplicity by wall crossing on the smooth quintic. Table~\ref{tab:intro-charge-comparison} specifies the geometric calculation used for each primitive coefficient.

\subsection{The degeneration theorem of Part I}
Take general forms $f_1,f_4,f_5$ of the indicated degrees in $\PP^4$. For $t\ne0$, the equation $f_1 f_4=t f_5$ defines a smooth quintic after shrinking the parameter curve. A small resolution of the total space has central fiber
\[
        Q\ \cup_S\ B,\qquad
  B=\operatorname{Bl}_{\Gamma}\PP^3,\qquad
  \Gamma=(f_4=f_5=0)\subset\PP^3=V(f_1).
\]
Here $Q=(f_4=0)\subset\PP^4$ is a quartic threefold and $S=Q\cap V(f_1)$ is a quartic K3 surface. Where the primitive profile calculation uses $\Pic(S)=\mathbb Z h$, we take $S$ very general in the polarized quartic family. The counts of lines and conics impose the additional open conditions specified in their respective sections. The two normal bundles of $S$ are dual. Inserting a chain of ruled threefolds over $S$ lets a supporting surface specialize onto a ruled component even when it contracts to $S$ on the original fiber.

For a sheaf $G$ on a smooth root $Y$ and the inclusion $i:S\hookrightarrow Y$, let $s$ cut out $S$. The boundary complex fits in the exact triangle
\[
 G(-S)\xrightarrow{s}G\longrightarrow i_*\mathbf L i^*G
 \longrightarrow G(-S)[1].
\]
In particular $\mathcal H^{-1}(i_*\mathbf L i^*G)=\ker(s)$. This kernel records sections annihilated by $s$ when a support component lies in $S$. The derived evaluation takes values in $\operatorname{Perf}(S)$. A stable rank-one torsion-free sheaf on a plane section gives a point of the relative compactified Jacobian.

For $r\ge2$ the point $\OO_T^{\oplus r}$ retains the nonproper stabilizer $\operatorname{PGL}_r$ after rigidification by scalars (Corollary~\companionref{cor:noseparatedall}). Stable sheaves can specialize to this polystable point (Proposition~\companionref{prop:noordinaryproperextension}). For $r=1$, rank one on an integral support gives stability on the smooth threefold. In higher multiplicity the Hall logarithm supplies the semistable count.

On an expansion Proposition~\companionref{prop:unlabelledpolarization} gives
\[
   \deg(\mathcal L|_\ell)=0
\]
for a contracted ruling $\ell$ and any line bundle $\mathcal L$ descending on the unmarked universal target, even for $r=1$. Part~I instead descends the stable loci defined by ample bundles on the charts. For the marked calculations we verify expansion-length bounds, existence and uniqueness in the charges concerned. Proposition~\companionref{prop:support-separated27} gives uniqueness when the support fixes the target and the two limits are stable for a common polarization. Proposition~\companionref{prop:divisorial-hull27} and the punctual modification analysis give existence under their local hypotheses.

Suppose a joining K3 is marked, and let $g_\gamma$ forget that marking in a given charge splitting $\gamma$. Under the hypotheses of Theorem~\companionref{thm:main}, the identity on the special fiber is
\[
       \spc[M_\eta]^{\vir}=\sum_\gamma (g_\gamma)_*
                                   [\MM^\dagger_\gamma]^{\vir}.
\]
When the marked stratum is the entire fiber product of two proper relative moduli over a smooth separated space $J_\gamma$ of simple boundary objects, with compatible obstruction theories, its summand before pushforward is
\[
   \Delta^!_{J_\gamma}\!
    \left([M_{Q,\gamma}^{\rel}]^{\vir}\times[M_{B,\gamma}^{\rel}]^{\vir}\right).
\]
Applying normalization to $F^\vee\otimes^{\LL}F$ gives the difference of the two restriction maps. Its homotopy fiber controls the intrinsic sheaf on the marked expansion of Li and Wu \cite{LW}. The scalar and determinant reductions of Part~I identify the virtual pullback on the full matching space under the stated hypotheses. A primitive pure curve boundary has $\dim J_\gamma=6$ and relative virtual dimensions $3,3$. A zero boundary contributes a separate class (Theorem~\ref{thm:primitivecompleteformula}).

In the line-defect charge the central limits have two forms. A line on $Q$ with $p=\ell\cap S\in\Gamma$ acquires the exceptional ruling over $p$. A four-secant of $\Gamma$ on $\PP^3$ has strict transform disjoint from $S$. Theorems~\ref{thm:direct-line} and~\ref{thm:expanded-line-compact} give, in either residue $1$ or $4$,
\[
       [q^{-83/120}]h^{\mathrm{DT}}_{1,\mu}
            =3\cdot1\,600+3\cdot1\,275=8\,625 .
\]

\begin{theorem}[Primitive coefficient comparison]\label{thm:intro-comp}
For a general smooth quintic the complete DT degrees in Table~\ref{tab:intro-charge-comparison} agree with the indicated coefficients of $Z^{\mathrm{GSY}}_\mu$. Five values are sums of marked relative virtual degrees. The two-point value is computed by specialization of the virtual zero-cycle in the projective sheaf moduli space. The other entries follow from the geometry and obstruction theory of the complete moduli on the smooth quintic.
\end{theorem}
\begin{table}[!htbp]\centering\small
\caption{Complete geometric DT degrees for the primitive quintic and the coefficients of the quintic elliptic genus in \cite[(3.10)]{GSY}. Each pair of opposite residues occupies one position.}\label{tab:intro-charge-comparison}
\begin{tabular}{@{}rccrrl@{}}\toprule
No. & $\mu$ & $x$ & $[q^x]h^{\mathrm{DT}}_{1,\mu}$ & $[q^x]Z^{\mathrm{GSY}}_\mu$ & DT calculation\\\midrule
1 & $0$ & $-55/24$ & $5$ & $5$ & marked\\
2 & $0$ & $-31/24$ & $-800$ & $-800$ & marked\\
3 & $1,4$ & $-83/120$ & $8\,625$ & $8\,625$ & marked\\
4 & $0$ & $-7/24$ & $58\,500$ & $58\,500$ & specialization\\
5 & $2,3$ & $13/120$ & $-1\,218\,500$ & $-1\,218\,500$ & marked\\
6 & $1,4$ & $37/120$ & $-1\,138\,500$ & $-1\,138\,500$ & marked\\
7 & $0$ & $17/24$ & $5\,817\,125$ & $5\,817\,125$ & smooth quintic\\
8 & $2,3$ & $133/120$ & $441\,969\,250$ & $441\,969\,250$ & smooth quintic\\
9 & $1,4$ & $157/120$ & $3\,777\,474\,000$ & $3\,777\,474\,000$ & smooth quintic\\
\bottomrule\end{tabular}\end{table}

In the two-point row the local algebra of a punctual double is $\CC[\epsilon]/(\epsilon^2)$. The diagonal terms on $Q$ and $B$ and the distinct-point terms in Theorem~\ref{thm:two-points-direct} give
\[
           (-600-900)+(9\,600+28\,800+21\,600)=58\,500 .
\]
Theorems~\ref{thm:marked-conics-18} and~\ref{thm:mixed-all-marked} give the marked conic and mixed degrees. Proposition~\ref{prop:polar-from-specializations} proves that the negative powers in the table exhaust the polar part. Four of the degrees, $5,-800,8\,625,-1\,218\,500$, together with five geometric vanishings select a unique vector $F^{\mathrm{DT}}$ in the predicted modular space (Theorem~\ref{thm:step21-four-seeds}). This gives $F^{\mathrm{DT}}=Z^{\mathrm{GSY}}$ \cite[text preceding (3.10)]{GSY}. The remaining five degrees test $h^{\mathrm{DT}}_1=F^{\mathrm{DT}}$ in further charges.

The neutral term at $q^{17/24}$ distinguishes the modular prediction \cite[(3.10)]{GSY} from the dilute-gas value $5\,797\,625$ in \cite[(3.5)]{GSY}. The three-point incidence includes a correction of one for each of the $2875$ lines on $X$. The other sheaves have hull $\OO_T(L_1-L_2)$ on the section spanned by two skew lines. Thus
\[
           -2\,445\,625+2875\cdot2874=5\,817\,125 .
\]
Theorems~\ref{th:seven-total}--\ref{th:nine-total} determine the complete moduli in the last three rows of the table. In particular $\OO_T(-L_1-L_2)$ occurs once per unordered pair. Corollary~\ref{cor:twoskew-all-lengths} treats their point modifications. Earlier rank-zero calculations by wall crossing appear in \cite[\S5]{AFKPS} and \cite{FTperiod,FeyzOSV}.

For $r>1$ the Hall logarithm \textup{(\companionref{eq:Hallloggeneral})} counts extensions of stable factors. Section~\ref{sec:generalmodularity} evaluates the strictly semistable contribution in a doubled class. Proposition~\ref{prop:step22-shadowunique} gives uniqueness after fixing the nonholomorphic correction and the holomorphic polar part.

\subsection{A pencil of K3 surfaces}\label{intro:pencilgeometry}
The first degeneration leaves us with a relative invariant of the quartic threefold $Q$. To evaluate it, choose a pencil of hyperplane sections of $Q$ containing $S$. Two hyperplanes in the pencil meet in a plane, so their sections of $Q$ meet in a plane quartic curve $C_*\subset S$. Blowing up this common curve gives the morphism
\[
 \kappa:Y=\operatorname{Bl}_{C_*}Q\longrightarrow\PP^1,
 \qquad Y_b=Q\cap H_b .
\]
Here $H_b$ is the hyperplane indexed by $b$. Its section is canonically isomorphic to the fiber $Y_b$ since $C_*$ is Cartier on that section. A general fiber is a smooth quartic K3 surface.

A sheaf is \emph{vertical} if its support lies over finitely many points of $\PP^1$. Let $D=\sum m_iD_i$ be its effective divisorial cycle, here of class $[Y_b]$. For an ample class $A$ on $Y$,
\[
 0=A\cdot Y_b^2=\sum_i m_i(A\cdot D_i\cdot Y_b).
\]
Every summand is nonnegative and is positive if $D_i$ dominates the base. Each $D_i$ is therefore a fiber, since the fibers are integral. Intersecting with $A^2$ gives $\sum m_i=1$.

\begin{figure}[htbp]
\centering
\begin{tikzpicture}[x=1cm,y=1cm,font=\small,>=Stealth]
\definecolor{fibershade}{RGB}{216,235,228}
\draw[black!60] (-3.65,0.37)--(3.95,0.37)--(4.35,1.12)--(-3.25,1.12)--cycle;
\draw[black!60] (-3.65,.37)--(-3.65,1.26)--(-3.25,2.01)--(4.35,2.01)--(4.35,1.12);
\draw[black!60] (-3.65,1.26)--(3.95,1.26)--(4.35,2.01);
\draw[black!60] (3.95,.37)--(3.95,1.26);
\foreach \t in {-2.65,2.75}{
  \draw[black!40] (\t,.37)--(\t,1.26)--(\t+.4,2.01);
}
\filldraw[fill=fibershade,draw=black!85,line width=.9pt]
    (-.4,.37)--(-.4,1.26)--(0,2.01)--(0,1.12)--cycle;
\node[anchor=south west] at (0,1.98) {$Y_b$, carrying $V$};
\node[anchor=east] at (-3.75,1.10) {$Y$};
\draw[black!75,-{Stealth}] (-3.85,-.56)--(4.28,-.56);
\node[anchor=west] at (4.32,-.56) {$\PP^1$};
\draw[->] (-.22,.38)--(-.22,-.46);
\fill (-.22,-.56) circle (1.6pt);\node[below] at (-.22,-.58) {$b$};
\node[anchor=east] at (-.34,-.04) {$\kappa$};
\end{tikzpicture}
\caption{A sheaf of fiber class on the pencil threefold. Its support projects to one point $b$. Its virtual multiplicity is computed on the threefold, allowing the supporting quartic to move in the pencil.}\label{fig:intro-vertical-sheaf}
\end{figure}
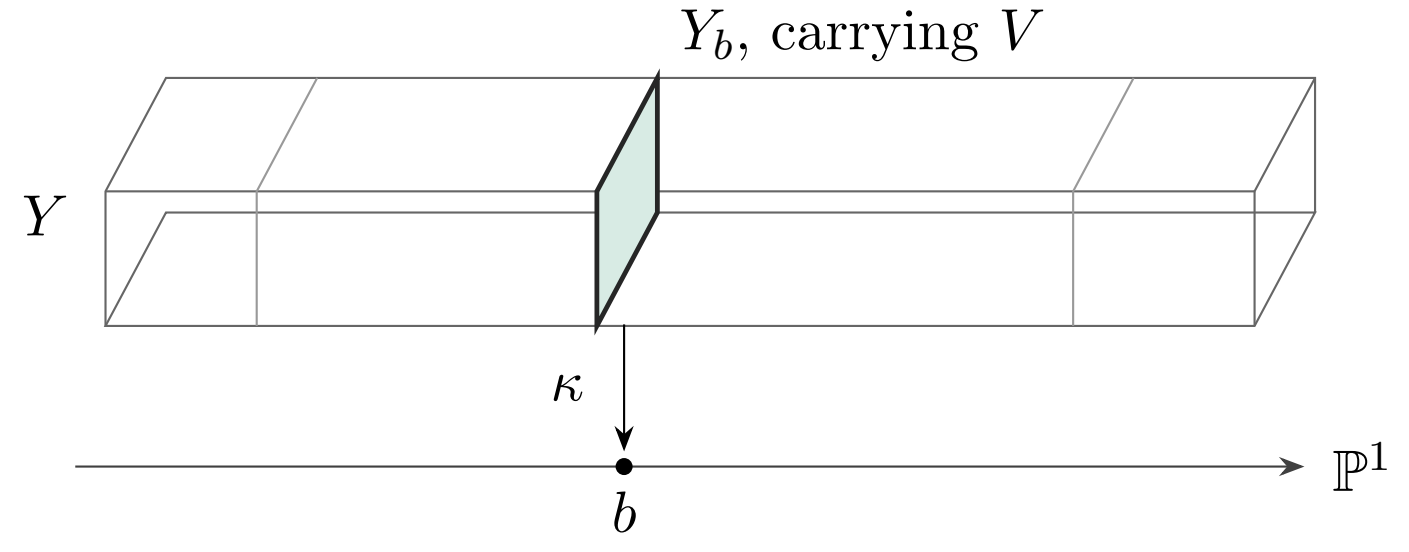

Its equation annihilates the sheaf by purity. Thus the objects counted are
\[
                   i_{b*}V,\qquad i_b:Y_b\hookrightarrow Y,
 \qquad V\text{ torsion-free of rank one on }Y_b .
\]
The point $b$ varies in the moduli problem (Figure~\ref{fig:intro-vertical-sheaf}).

For $P_V(m)=2m^2+dm+c$ one has $c=\chi(V)$ and, on a smooth fiber, $d=h_b\cdot c_1(V)$ for the hyperplane class $h_b$. Let $\dt_Y(d,c)$ be the sum of the absolute threefold DT degrees for this polynomial. As $h_b^2=4$, the twist $V\mapsto V(1)$ sends $(d,c)$ to $(d+4,c+d+2)$. Set
\[
 K^{\rm raw}_{a}(q)=\sum_c \dt_Y(d,c)\,q^{1+d^2/8-c},
                  \qquad a=d\pmod 4 .
\]
The exponent is twist invariant. The neutral pencil component has $4\mid d$. The quintic index $\mu$ is defined in \eqref{eq:intro-partition}.

For a very general pencil the fibers contributing in degree $4\nmid d$ are smooth. Indeed $d=h_b\cdot c_1(M)$ for the line bundle $M=V^{**}$, so $\Pic(Y_b)=\mathbb Zh_b$ would force $4\mid d$. Theorem~\ref{thm:pencil-all-charge-normalization} excludes the nodal fibers. Writing $V=M\otimes I_Z$, the Noether--Lefschetz multiplicity counts the class of $M$ and the K3 Hilbert-scheme factor sums over $Z$. Thus
\[
  K^{\rm raw}_a=K^{\rm per}_a=\widehat\Phi_a/\Delta
  \quad(a=1,2,3),\qquad \Delta=q\prod_{m\ge1}(1-q^m)^{24}.
\]
Here $\widehat\Phi$ is normalized by the degree of the double cover for simultaneous resolution. For $a=0$ choose the representative $d=0$.

For $d=0$, the point ideals $V=I_{Z/Y_b}$ occur also on singular fibers. Corollary~\ref{cor:pencil-ideal-reduction} expresses the discrepancy in terms of their virtual degrees $\mathcal I_N$, $N=\length Z$. By Proposition~\ref{prop:pencil-ideal-three-m27}, $\mathcal I_3=31\,792$ and $(K^{\rm per}-K^{\rm raw})_0=O(q^3)$.

Section~\ref{sec:seconddegfull} degenerates $Y$ to $Q\cup_S W$. Blow up $C_*$ on one section of a $\PP^1$-bundle over $S$ to form $W$. The untouched section glues to $Q$; the strict transform $B_0$ of the other section is the outer divisor. Theorem~\ref{thm:outerbitrarymultiplicity} compares the classical relative and absolute moduli of stable vertical sheaves. Using their virtual degrees also requires the comparison of obstruction theories.

\subsection{Two boundary conditions for conics on the quartic}
Suppose $\Pic(S)=\mathbb Zh$, with $h=c_1(\OO_S(1))$. For a smooth quartic $C_*\in|h|$ the support morphism of the sixfold $J_{-4}$ has fiber
\[
 \rho:J_{-4}=M_S(0,h,-4)\longrightarrow|h|,
 \qquad \rho^{-1}(C_*)=\Pic^{-2}(C_*).
\]
Indeed $g(C_*)=3$ and $\chi(L)=\deg L-2$. Let $f=[\rho^{-1}(C_*)]$. Intersection with $f$ fixes the curve and allows the line bundle to vary.

For general $Q$ and $S$ as above, a plane conic $K\subset Q$ meets $S$ in a subscheme $Z$ of length two. It consists of two points or of a point together with a tangent direction. Write $D$ for the class of the conic parameter surface pushed to the Hilbert scheme $S^{[2]}$ by $K\mapsto K\cap S$. The planes through $Z$ form a pencil. For each such plane its intersection $C$ with $S$ carries the ideal $I_{Z/C}$ of $Z$.

For any surface class $B'$ in $S^{[2]}$ we use the same construction.
\[
 S^{[2]}\xleftarrow{\ p\ } U=\{(Z,C):Z\subset C\in|h|\}
 \xrightarrow{\ a\ } J_{-4},\qquad
 p(Z,C)=Z,\quad a(Z,C)=I_{Z/C}.
\]
The resulting three-cycle $A_{B'}=a_*p^*B'$ is called an Abel cycle. On a smooth $C$ the map sends an effective divisor $Z$ of degree two to $\OO_C(-Z)$.

\begin{theorem}[Two intersections for the conic family]\label{thm:intro-two20}
Fix a smooth plane quartic $C_*\subset S$. Let $\tau$ be the surface in $S^{[2]}$ of nonreduced length-two subschemes supported at a point of $C_*$, with any tangent direction in $S$. For the blowup center $\Gamma\in|5h|$ set $g=[\Gamma^{[2]}]$. Then
\[
 \deg(A_D f)=5016,\qquad \deg(A_D A_\tau)=-26560,\qquad
 \deg(A_D A_g)=-516400 .
\]
\end{theorem}
In the first intersection $Z$ lies on the prescribed curve $C_*$. In the second, only its center must lie on $C_*$. The curve $C$ carrying its ideal varies through the planes containing the tangent subscheme. The incidence $U$ is a divisor in $J_{-4}$ whose normal line bundle has degree $-2$ on each pencil of planes. Its self-intersection gives $-2\int_{S^{[2]}}D\tau=-2(13280)$. The two conditions give independent intersection pairings on the classes $A_{B'}$.

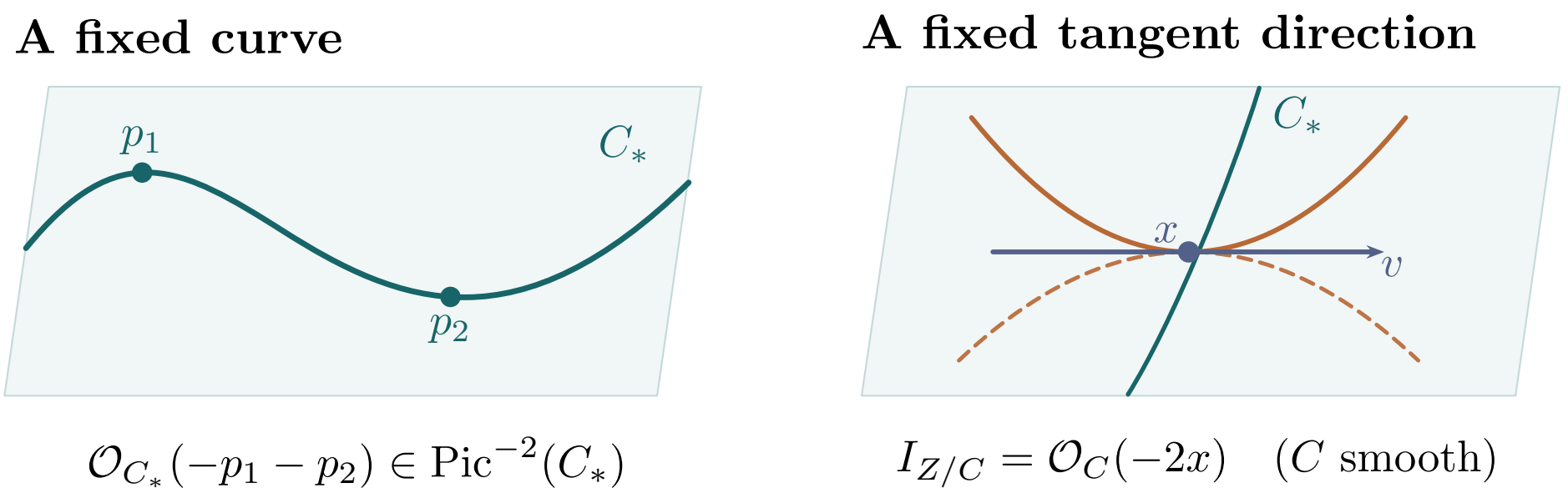
\begin{figure}[htbp]
\centering
\begin{tikzpicture}[x=1cm,y=1cm,font=\footnotesize,line cap=round,line join=round]
\definecolor{ksea}{HTML}{176568}\definecolor{ksrf}{HTML}{B15C25}\definecolor{kamb}{HTML}{53618A}
\node[anchor=west,font=\small\bfseries] at (.06,3.62) {A fixed curve};
\node[anchor=west,font=\small\bfseries] at (6.77,3.62) {A fixed tangent direction};
\path[fill=ksea!6,draw=ksea!25] (.13,.77)--(.48,3.22)--(5.65,3.22)--(5.30,.77)--cycle;
\draw[ksea,line width=1.3pt] (.30,1.94) .. controls (1.04,2.85) and (1.54,2.58) .. coordinate[pos=.46] (abelpa) (2.33,2.08)
 .. controls (3.27,1.48) and (4.2,1.14) .. coordinate[pos=.46] (abelpb) (5.55,2.46);
\fill[ksea] (abelpa) circle (2.3pt) (abelpb) circle (2.3pt);
\node[above,text=ksea] at (abelpa) {$p_1$};
\node[below,text=ksea] at (abelpb) {$p_2$};
\node[text=ksea] at (5.04,2.76) {$C_*$};
\node at (2.92,.24) {$\mathcal O_{C_*}(-p_1-p_2)\in\mathrm{Pic}^{-2}(C_*)$};
\begin{scope}[shift={(9.51,1.91)}]
\path[fill=ksea!6,draw=ksea!25] (-2.59,-1.14)--(-2.23,1.31)--(2.94,1.31)--(2.59,-1.14)--cycle;
\draw[ksea,line width=1.0pt] (-.48,-1.13) .. controls (-.14,-.59) and (.33,.55) .. (.56,1.30);
\node[text=ksea,anchor=west] at (.52,1.08) {$C_*$};
\draw[ksrf!92,line width=1.1pt] plot[domain=-1.72:1.72,samples=50] (\x,{.36*\x*\x});
\draw[ksrf!85,densely dashed,line width=.9pt] plot[domain=-1.82:1.82,samples=50] (\x,{-.26*\x*\x});
\draw[-{Stealth[length=4pt]},kamb,line width=1.1pt] (-1.55,0)--(1.55,0);
\fill[kamb] (0,0) circle (2.4pt);
\node[above left,text=kamb,inner sep=2pt] at (0,0) {$x$};
\node[below right,text=kamb,inner sep=1pt] at (1.48,0) {$v$};
\end{scope}
\node at (9.58,.24) {$I_{Z/C}=\mathcal O_C(-2x)\quad(C\text{ smooth})$};
\end{tikzpicture}
\caption{The two curve conditions in Theorem~\ref{thm:intro-two20}.
On the left the supporting curve is prescribed. Distinct points $p_1,p_2\in C_*$ give the line bundle $\OO_{C_*}(-p_1-p_2)\in\Pic^{-2}(C_*)$. On the right, only $x$ lies on the
prescribed curve $C_*$. In coordinates $(u,t)$ with $v=\partial_u$, the
length-two subscheme $Z$ has ideal $(t,u^2)$. The solid and dashed curves
represent smooth germs of moving plane sections containing $Z$. Each has the
same tangent direction at $x$. The picture shows a direction transverse to
$C_*$, which is one of the directions allowed in $\mathbb P(T_xS)$.}
\label{fig:picardjetII}
\end{figure}

Figure~\ref{fig:picardjetII} shows which curve and point data are fixed in the two intersections.

In particular the first degree alone cannot recover the conic contribution from $\Gamma$, which equals $-50(5016)+10(-26560)$ (Proposition~\ref{prop:mark20-twoprobe} and Corollary~\ref{cor:mark20-conictest}).

To obtain the second intersection starting from a pencil with fixed base curve $C_*$, let $\mathscr T$ parametrize triples $(x,[v],C)$, where $x\in C_*$, $[v]\in\PP(T_xS)$, and $C$ contains the tangent subscheme $Z(x,v)$. The ruled component of the pencil degeneration supplies the Picard variety at the left of the correspondence
\[
\begin{tikzcd}[column sep=large]
 \Pic^{-2}(C_*) & C_* \arrow[l, "\delta"'] &
 \mathscr T \arrow[l,"\pi"'] \arrow[r,"b"] & J_{-4}.
\end{tikzcd}
\]
Here $\delta(x)=\OO_{C_*}(-2x)$, $\pi(x,[v],C)=x$, and $b(x,[v],C)=I_{Z(x,v)/C}$. Pullback of the ruled component's virtual class $[\Pic^{-2}(C_*)]$ gives $[\mathscr T]$, with $b_*[\mathscr T]=A_\tau$ (Theorem~\ref{thm:step21-flag}). For smooth $C$ transverse to $C_*$ the four bundles are $\OO_C(-2x)$ with $x\in C\cap C_*$ and $[v]=T_xC$. Intersecting with $A_D$ and pushing back gives the Chow class
\[
  -6\,640\,\delta_*\big(c_1(\omega_{C_*})\cap[C_*]\big)
  \quad\text{of degree}\quad -6\,640(2g(C_*)-2)=-26\,560
\]
(Proposition~\ref{prop:step21-cap}).

Varying $C_b=F\cap S$ with $K\subset F\subset Q$ gives a smooth proper threefold of flags $(K,F)$. On noncentral targets these flags exhaust the stable fiber sheaves of polynomial $2m^2-2m+1$. Their supports miss the outer fiber and their exceptional restrictions project to $I_{K\cap S\,/\,C_b}$. Theorem~\ref{thm:rev27-moving} identifies the sheaf obstruction theory along the double-point intersection and specializes its degree $-26\,560$ in the flag space. Decomposition on a marked expansion further requires stable reduction and compatible obstruction morphisms there.

\subsection{From the pencil to the quintic series}
Let $\mathcal B$ and $\mathcal C$ pair a quartic boundary class with the classes from $W$ and from $B=\operatorname{Bl}_{\Gamma}\PP^3$, respectively, retaining the output charges. For $\mathcal C$ to factor through $\mathcal B$, one needs
\[
  \mathcal B(\alpha)=\mathcal B(\alpha')\ \Longrightarrow  \ \mathcal C(\alpha)=\mathcal C(\alpha'),
  \qquad\text{equivalently}\qquad \ker\mathcal B\subset\ker\mathcal C .
\]
Write $R_Y,R_X$ for the remaining terms in the pencil and quintic formulas.

\begin{theorem}[Modularity criterion in divisor class $H$]\label{thm:intromod}
Assume the two proper virtual degeneration formulas with identical relative quartic classes, and the virtual comparison at the outer divisor of the pencil. Suppose the relative intersections in every degree give $\mathcal C=T_\Theta\mathcal B$, where $T_\Theta$ contracts the four pencil components with a holomorphic modular matrix $\Theta$ of weight $31/2$ and divides by $\eta^{31}$. Suppose the geometric calculation of $c_Y=K^{\rm per}-K^{\rm raw}$ and the other sectors gives
\[ U=\eta^{55}\bigl(R_X-T_\Theta(R_Y+c_Y)\bigr)\in M_{26}(\rho_5),\]
where $\rho_5$ is the representation on the five quintic components. Then
\begin{equation}\label{eq:intromodularconclusion}
 \eta^{55}h^{\mathrm{DT}}_{1,\mu}=\sum_{a\bmod4}\widehat\Phi_a\Theta_{a,\mu}+U_\mu,
 \qquad \mu\in\mathbb Z/5 .
\end{equation}
In particular the primitive quintic series has the transformation law of \cite[\S3.3, (3.6)--(3.8)]{GSY}.
\end{theorem}
Here $\eta=q^{1/24}\prod_{m\ge1}(1-q^m)$ and $\Delta=\eta^{24}$. The weights in \eqref{eq:intromodularconclusion} give $21/2+31/2-55/2=-3/2$. Theorem~\ref{thm:twodegmodularityfull} proves this conditional implication. Alternatively, the modularity of $\mathcal N_X$ in Section~\ref{sec:quintic-period-comparison-new} reduces the conjecture to the virtual identity $\eta^{55}h^{\mathrm{DT}}_1=\mathcal N_X$. Neither comparison in all charges follows from the finite coefficient calculation.

\subsection{Main results and plan of the paper}

Sections~\ref{sec:firstII}--\ref{sec:boundarycalibrations} construct the two degenerations and follow the boundary maps from their common quartic component. On the quintic side, Theorem~\ref{thm:primitivecompleteformula} retains both the positive and the zero boundary in the primitive specialization under the hypotheses of Part~I. Theorem~\ref{thm:vacuumvirtual} evaluates these two terms on the projective family of expanded supports constructed there, giving $4+1=5$. On the pencil side we compute the deformation complex with both boundary values fixed on the ruled component. At a smooth outer fiber we identify the classical relative and absolute moduli of stable vertical sheaves in arbitrary multiplicity. For the structure sheaf and the spherical rank-three bundle the obstruction theories give the same virtual degree $-1$ on both moduli spaces (Theorems~\ref{thm:K3structurecomparison} and~\ref{thm:rankthreepencilspherical}).

Sections~\ref{sec:generalmodularity} and~\ref{sec:modularity} establish what the charges and the pencil already determine. Shifted duality and the Hall algebra prove charge conjugation for the geometric sheaf series, including strictly semistable classes. Wall crossing gives entire ranges of exact polar coefficients for point and line charges in higher divisor multiplicity (Theorems~\ref{thm:hall22strip} and~\ref{thm:hall22lines}). In class $2(1-e^{-H})$ the five-arrow Kronecker quiver gives a contribution $-155/4$, with a residual stable locus as in \eqref{eq:Kronfull22}. For the very general quartic hyperplane pencil we identify the nonzero degree residues with the normalized Noether--Lefschetz series divided by $\Delta$. In residue zero the discrepancy is expressed entirely in terms of point ideals on the fibers. Their three-point virtual degree is $31\,792$.

The calculations in Sections~\ref{sec:initial-checks}--\ref{sec:smooth-evaluations} establish the nine geometric degrees of Table~\ref{tab:intro-charge-comparison}. For a point, a line and a line with a point defect, the task is to reconstruct the sheaf from the expanded incidence and identify its obstruction theory at the join. The two-point charge instead admits a direct Chow specialization in the proper sheaf moduli, retaining the punctual doubles on both roots. For a conic the quartic boundary carries two independent intersections. The moving-pencil calculation identifies the sheaf obstruction theory along the tangent incidence and gives its degree $-26\,560$ (Theorem~\ref{thm:rev27-moving}). The two intersections recover the blowup contribution $-516\,400$ to the full conic degree. Proposition~\ref{prop:mark20-allcontacts} compares the two curve conditions in cohomology at every contact length.

For the last three complete degrees, Section~\ref{sec:smooth-evaluations} works directly on the smooth quintic. Identifying all the possible hulls brings in pairs of skew lines, twisted cubics and elliptic quartics. Specified two-line hulls contribute a Hilbert scheme of points at every length (Corollary~\ref{cor:twoskew-all-lengths}). Section~\ref{sec:arithfull} constructs the weight-$26$ period vector $\mathcal N_X$ over the entire hyperplane system and the contraction of Weil representations from the four pencil residues to the five quintic residues (Proposition~\ref{prop:intertwinerfull}). For twisted-cubic hulls and their positive powers, the transverse period multiplicity agrees with the contribution of the hull and all its point modifications to $\eta^{55}h^{\mathrm{DT}}_1$. The coefficients of $q$ and $q^2$ in the neutral geometric numerator $\eta^{55}h^{\mathrm{DT}}_{1,0}$ are $-1075$ and $109\,650$, expressed in \eqref{eq:quintic-two-boundary} as Euler integrals over the singular sections. The residual comparison \eqref{eq:quintic-residual-period} specifies the virtual multiplicities still needed for equality with the period vector. For a singular hyperplane section, Part~III \cite{QuinticThreePerfect} extracts the contribution of perfect point modifications at locally free stalks, where the ambient Behrend weight is preserved, and identifies the residual sum to be compared with the periods.

Sections~\ref{sec:twodegoperators} and~\ref{sec:verified-fourier-table} return to the two evaluations of the common quartic class. On a very general quartic K3, eight middle-degree insertions on each side determine the numerical pairing on the primitive compactified Jacobian (Theorem~\ref{thm:mark20-eight}). The factorization of the two evaluations in all charges and the modularity of the remaining terms, including the neutral pencil correction, imply the predicted transformation law by Theorem~\ref{thm:twodegmodularityfull}. Independently of these hypotheses, four marked degrees and the support bound determine $F^{\mathrm{DT}}=Z^{\mathrm{GSY}}$ in the prescribed modular space. Table~\ref{tab:twenty-rev27} gives its first twenty Fourier positions. Under the depth-one quasimodularity hypothesis of Corollary~\ref{cor:step21-depth}, twelve positions suffice to test the full identity. The five complete sheaf degrees beyond the four inputs in Table~\ref{tab:intro-charge-comparison} test that reconstruction, while equality with $h^{\mathrm{DT}}_1$ in every charge remains the geometric modularity question.

\subsection*{Acknowledgements}
The author thanks Amin Gholampour who extensively discussed this project with the author over the past 15 years. The idea of double degeneration was Amin's, however due to difficulties in compactifying the relative branch moduli spaces and also having no control on cohomology theory of compactified Jacobian (and even worse, no morphisms to it either) the project halted for many years, it was the latter obstruction, which led the author to develop ideas of intersection theory on moduli space of perfect complexes on double divisor rather than the compactified Jacobian. The author would like to thank Eyal Markman for teaching him about middle cohomology of compactified Jacobian which helped the author to come up with some of numerical calculations in some specific examples. The author expresses sincere gratitude to Sheldon Katz as well as Richard Thomas for teaching him many facts about geometry of torsion sheaf moduli spaces over the years. The author would like to especially thank John Pardon for stimulating conversations which might provide a different approach to the problem discussed in this article, taking advantage of auxiliary divisors and how they induce maps to compactified Jacobians. Finally especial thanks go to Yukinobu Toda, Maxim Kontsevich, Ludmil Katzarkov, and Vladimir Baranovsky for discussions during the earlier stages of this project. He is grateful to Charles Doran, Tony Pantev, Shing-Tung Yau, Martijn Kool for valuable comments and discussions. The author acknowledges support by grants from Beijing Institute of Mathematical Sciences and Applications (BIMSA), the Beijing NSF BJNSF-IS24005, and the China National Science Foundation (NSFC) NSFC-RFIS program W2432008. He would also like to thank China's National Program of Overseas High Level Talent for generous support. Finally, he would like to thank NSF AI Institute for Artificial Intelligence and Fundamental Interactions at Massachusetts Institute of Technology (MIT) which is funded by the US NSF grant under Cooperative Agreement PHY-2019786.

\section{The first degeneration and its relative moduli}\label{sec:firstII}
\begin{figure}[htbp]
\centering
\begin{tikzpicture}[x=1cm,y=1cm,font=\footnotesize,line cap=round,line join=round]
\definecolor{ksea}{HTML}{176568}\definecolor{ksrf}{HTML}{B15C25}\definecolor{kamb}{HTML}{53618A}
\node[anchor=west,text=kamb] at (0,3.91) {$Q$}; 
\node[anchor=east,text=kamb] at (12.6,3.91) {$B=\mathrm{Bl}_{\Gamma}\mathbb P^3$};
\path[fill=kamb!5] (.15,2.85)--(1.22,3.62)--(6.56,3.62)--(5.59,2.85)--cycle;
\path[fill=kamb!5] (5.59,2.85)--(6.56,3.62)--(12.44,3.62)--(11.46,2.85)--cycle;
\draw[kamb!50,thin] (.15,.77)--(.15,2.85)--(1.22,3.62)--(12.44,3.62)--(12.44,1.55);
\draw[kamb!45,thin] (.15,2.85)--(11.46,2.85)--(12.44,3.62) (11.46,2.85)--(11.46,.77);
\path[fill=ksea!13,draw=ksea!58] (5.59,.77)--(6.56,1.55)--(6.56,3.62)--(5.59,2.85)--cycle;
\path[fill=ksrf!17,draw=ksrf!80,line width=.7pt]
 (.60,1.35) .. controls (2.2,1.98) and (3.5,.96) .. (5.59,1.42)
 .. controls (5.84,1.74) and (6.31,1.69) .. (6.56,2.12)
 .. controls (4.65,1.84) and (3.36,3.05) .. (1.66,2.17) --cycle;
\foreach \t in {.26,.52,.78}
 {\draw[ksrf!35] ({.60+1.06*\t},{1.35+.82*\t}) .. controls
 ({2.2+.36*\t},{1.98+.81*\t}) and ({3.5+.35*\t},{.96+1.15*\t}) .. ({5.59+.97*\t},{1.42+.70*\t});}
\path[fill=ksrf!11,draw=ksrf!80,line width=.7pt] 
 (5.59,1.42) .. controls (7.4,.98) and (9.62,1.24) .. (11.10,1.50)
 --(12.03,2.24) .. controls (9.72,1.98) and (8.45,1.87) .. (6.56,2.12)
 .. controls (6.31,1.69) and (5.84,1.74) .. (5.59,1.42)--cycle;
\draw[ksrf!32] (5.95,1.71) .. controls (7.65,1.37) and (9.66,1.53) .. (11.45,1.78);
\draw[ksrf!32] (6.29,1.84) .. controls (7.94,1.63) and (9.85,1.74) .. (11.72,1.99);
\draw[ksea,line width=1.8pt] (5.59,1.42).. controls (5.84,1.74) and (6.31,1.69)..(6.56,2.12);
\draw[kamb!48,thin] (.15,.77)--(11.46,.77)--(12.44,1.55);
\node[fill=white,inner sep=2pt,text=ksrf] at (3.04,2.19) {$T_Q$};
\node[fill=white,inner sep=2pt,text=ksrf] at (9.30,1.68) {$T_B$};
\node[text=ksea,anchor=south] at (6.00,3.66) {$S$};
\draw[ksea,thin] (4.93,.32) -- (5.34,.32)--(5.93,1.68);
\node[anchor=east,text=ksea] at (4.94,.32) {$C=T_Q\cap S$};
\node[anchor=west,text=ksea] at (6.22,.32) {$C\simeq T_B\cap S$};

\node[anchor=west,font=\footnotesize\bfseries] at (.30,-.45) {The plane section};
\node[anchor=west,font=\footnotesize\bfseries] at (7.03,-.45) {Over one marked point};
\draw[ksea,line width=.9pt] plot[domain=.38:4.88,samples=60]
 (\x,{-1.28+.13*sin((\x-.4)*70)});
\foreach \j in {0,...,19}{
 \pgfmathsetmacro{\tj}{.50+.219*\j}
 \fill[ksrf] (\tj,{-1.28+.13*sin((\tj-.4)*70)}) circle (1.32pt);}
\node[anchor=north,text=ksea] at (2.60,-1.48) {$C\subset\Pi$};
\draw[-{Stealth[length=4pt]},kamb,line width=.65pt] (7.02,-1.26)--(5.14,-1.26);
\node[text=kamb] at (6.08,-.94) {$b$};
\draw[ksea,line width=1.25pt] (7.34,-1.48)..controls (8.26,-1.46) and (8.80,-1.41)..(9.39,-1.21)
 ..controls (10.36,-.90) and (11.02,-.92)..(11.60,-1.04);
\draw[ksrf,line width=1.05pt] (8.89,-.77) .. controls (9.16,-.87) and (9.25,-1.04)..(9.39,-1.21)
 ..controls (9.50,-1.36) and (9.60,-1.60)..(9.60,-1.77);
\fill[ksea] (9.39,-1.21) circle (2.1pt);
\node[anchor=north west,text=ksrf] at (9.61,-1.45) {$E_i$};
\node[anchor=west,text=ksea] at (11.63,-1.04) {$\widetilde C$};
\end{tikzpicture}
\caption{The hyperplane support on the two components of the first degeneration.
Take a general hyperplane $P\subset\mathbb P^4$ and write $\Pi=P\cap\mathbb P^3$.
On the left, $T_Q=P\cap Q$ is a quartic K3 surface. On the right,
$T_B=\mathrm{Bl}_{\Pi\cap\Gamma}\Pi$. On the common K3 surface $S$, the
plane cuts out the smooth quartic $C=\Pi\cap S$. Since $[\Gamma]=5h$ and
$[C]=h$, the plane is blown up at $5h^2=20$ distinct points of $C$.
The blowdown $b:T_B\to\Pi$ maps the strict transform $\widetilde C$ isomorphically to $C$ and contracts each exceptional curve $E_i$. For sheaves locally free on their respective supports near $C$, both boundary values are vector bundles on this curve.}\label{fig:transversehyperplanesII}
\end{figure}
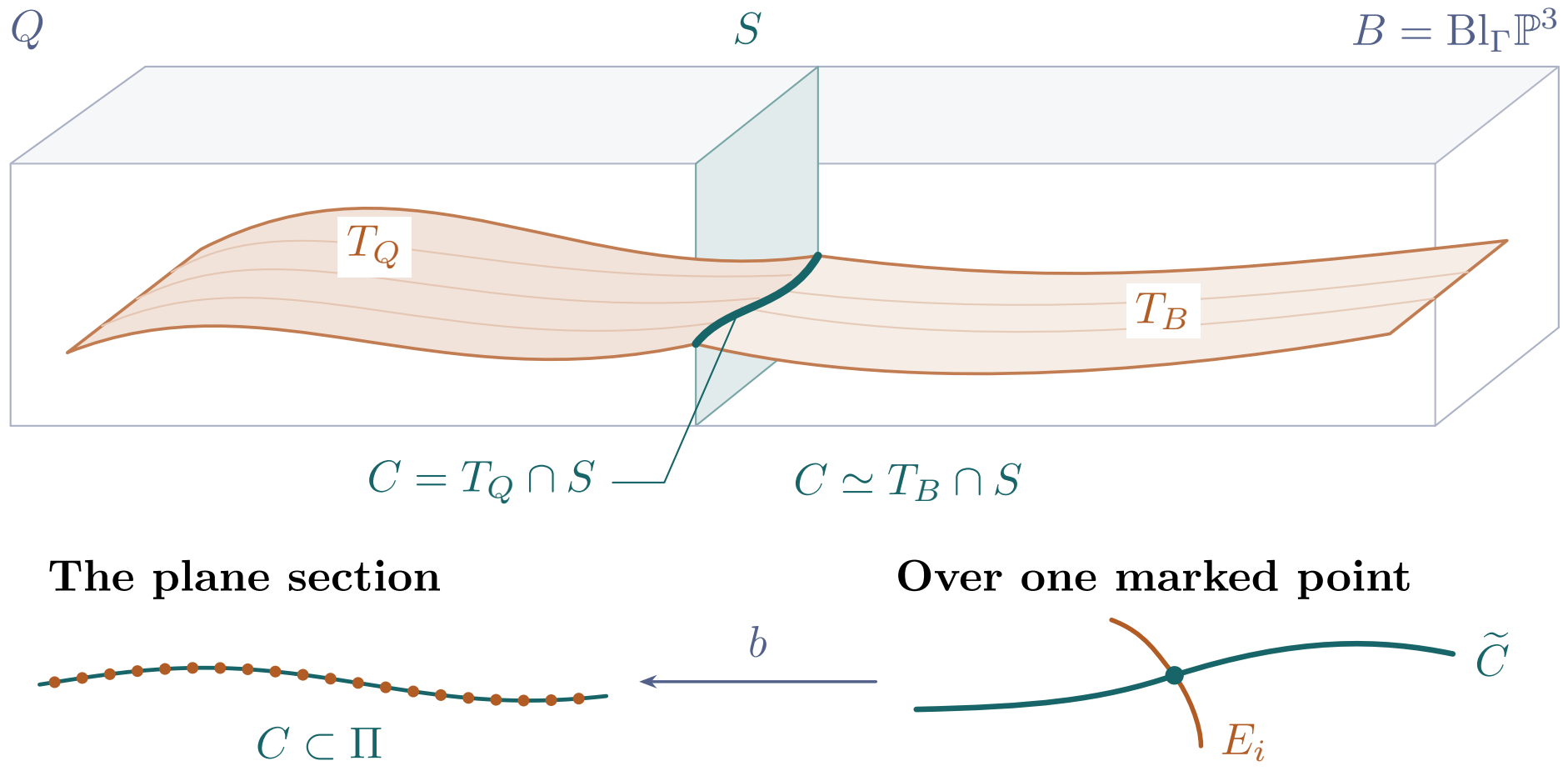

On $Q\cup_S B$ the relative classes record both the hyperplane section and its sheaf. Their intersection takes place on the boundary moduli over $S$.
Figure~\ref{fig:transversehyperplanesII} identifies these two surfaces and the curve on which their restrictions are compared for a general hyperplane.

\subsection{The sheaf and its deformations}
For a stable sheaf $F$ on the smooth quintic, the tangent and obstruction spaces are
\[
           T_{[F]}M_X(v)=\Ext^1_X(F,F),\qquad
           \mathrm{Ob}_{[F]}=\Ext^2_X(F,F)\simeq\Ext^1_X(F,F)^\vee .
\]
Scalar rigidification removes $\mathbb G_m\cdot1_F$. Serre duality gives $\operatorname{vdim}M_X(v)=\ext^1_X(F,F)-\ext^2_X(F,F)=0$. On a smooth proper component $M$ of dimension $d$ the symmetric obstruction theory gives
\[
   [M]^{\vir}=c_d(\Omega_M^1)\cap[M],\qquad
   \deg[M]^{\vir}=(-1)^d e(M)
\]
by \cite[Theorem~4.18]{Behrend27} and the obstruction-bundle formula.

In the primitive case Proposition~\companionref{prop:primitive} gives $F=i_*G$ with $G$ torsion-free of rank one on an integral $T\in|H|$. In class $rH$ all effective supports of that class are allowed.
\subsection{The two component moduli}
Let $H_B$ be the pullback of the hyperplane class and $E_B$ the exceptional divisor. Since $S\sim H_Q$ on $Q$ and $S\sim4H_B-E_B$ on $B$,
\[
 \omega_Q(S)=\OO_Q,\qquad \omega_B(S)=\OO_B,\qquad
 N_{S/Q}=\OO_S(h),\quad N_{S/B}=\OO_S(-h).
\]
Part I, Section~\companionref{sec:geometry}, proves these identities on expansions and provides the relatively ample line bundle used on the two roots.

A relative object consists of an expansion of $(Q,S)$ or $(B,S)$ together with a stable sheaf on it. For a splitting $\gamma$, write $M_{Q,\gamma}^{\rel}$ and $M_{B,\gamma}^{\rel}$ for the chosen relative moduli, with evaluation by derived restriction to the outer copy of $S$. On the simple curve-sheaf boundary the matched tangent complexes are related by the difference of restriction maps
\[
  \RHom_Q(F_Q,F_Q)[1]\oplus\RHom_B(F_B,F_B)[1]
    \xrightarrow{\ \mathrm{res}_Q-\mathrm{res}_B\ }
               \RHom_S(E,E)[1].
\]
For a longer expansion replace each root by its expanded relative target. Before scalar and determinant reduction, the tangent complex of the glued object is the homotopy fiber of this difference. Part~I, Sections~\companionref{sec:reduction}--\companionref{sec:diagonal}, describes the reductions and their compatibility hypotheses. It remains to determine the boundary classes of the proper components and their multiplicities.
\section{A second degeneration through the quartic pencil}\label{sec:seconddegfull}
Degenerating the quartic K3 pencil $Y$ produces a second copy of $(Q,S)$. Its absolute vertical invariants admit a period calculation.
\subsection{The pencil as an incidence variety}
Choose independent sections $s_0,s_1\in H^0(Q,\OO_Q(1))$ such that $S=(s_0=0)$ is smooth and $C_*=(s_0=s_1=0)$ is a smooth plane quartic. Thus $g(C_*)=3$.  Write
\[
 Y=\operatorname{Bl}_{C_*}Q,\qquad b:Y\longrightarrow Q,\qquad
 \kappa:Y\longrightarrow\PP^1.
\]
Let $A=b^*H_Q$ and denote the exceptional divisor by $E_C$.  The class of a fiber of $\kappa$ is denoted $F_{\mathrm{fib}}$.

\begin{proposition}\label{prop:pencilcomplete}
There is an identification
\[
 Y=\{(q,[u:v])\in Q\times\PP^1:u s_0(q)+v s_1(q)=0\}.
\]
Moreover
\begin{gather}
 F_{\mathrm{fib}}=A-E_C,\qquad K_Y=-F_{\mathrm{fib}},\qquad
 \omega_{Y/\PP^1}=\kappa^*\OO_{\PP^1}(1),\label{eq:pencilcanonicalfull}\\
 A^3=4,\quad A^2E_C=0,\quad AE_C^2=-4,\quad E_C^3=-8.\label{eq:pencilintsfull}
\end{gather}
The line bundle $L=A+F_{\mathrm{fib}}$ is ample and restricts to the primitive quartic polarization on every smooth fiber.  The sheaf cohomology of the fibration satisfies
\[
 \kappa_*\OO_Y=\OO_{\PP^1},\qquad R^1\kappa_*\OO_Y=0,
 \qquad R^2\kappa_*\OO_Y=\OO_{\PP^1}(-1).
\]
\end{proposition}
\begin{proof}
The center is a regular sequence. On a trivializing open set its Rees algebra is $\OO_Q[U,V]/(s_1U-s_0V)$. Relative $\operatorname{Proj}$, with $[u:v]=[-V:U]$, gives the incidence variety. Thus
\[
 \OO_Y(F_{\rm fib})=\kappa^*\OO(1),\quad
 F_{\rm fib}=A-E_C,\quad K_Y=b^*K_Q+E_C=-F_{\rm fib}.
\]
Subtracting $\kappa^*K_{\PP^1}$ gives \eqref{eq:pencilcanonicalfull}. The incidence embedding identifies $L$ with the restriction of $\OO_Q(1)\boxtimes\OO_{\PP^1}(1)$, with $L|_{Y_b}=h$.

Write $i:C_*\hookrightarrow Q$. Since $N_{C_*/Q}=\OO_{C_*}(1)^{\oplus2}$, pushforward of the powers of the exceptional divisor gives
\[
 b_*E_C=0,\qquad b_*E_C^2=-i_*[C_*],\qquad
 b_*E_C^3=-i_*c_1(N_{C_*/Q}).
\]
Intersecting with $H_Q$ and using $\deg C_*=4$ yields \eqref{eq:pencilintsfull}.

On every pencil member $T\subset\PP^3$, the sequence $0\to\OO(-4)\to\OO\to\OO_T\to0$ gives $(h^0,h^1,h^2)(\OO_T)=(1,0,1)$. Cohomology and base change give $\kappa_*\OO_Y=\OO$, $R^1\kappa_*\OO_Y=0$. Relative duality gives
\[
 (R^2\kappa_*\OO_Y)^\vee=\kappa_*\omega_{Y/\PP^1}
                   =\kappa_*\kappa^*\OO(1)=\OO(1).
\]
\end{proof}

For reference, the number of nodal fibers is $108$.  Indeed
\[
 c(T_Q)=\frac{(1+H_Q)^5}{1+4H_Q},\qquad c_3(T_Q)=-14H_Q^3,
 \qquad e(Q)=-56.
\]
The blowup adds $e(C_*)=-4$ to $e(Q)$. A nodal quartic has Euler number $24-1=23$. If $N$ members are nodal, then
\[
          -60=e(Y)=24(2-N)+23N,\qquad N=108 .
\]
\subsection{The second semistable family}
Let
\[
 \mathscr D=\operatorname{Bl}_{S\times\{0\}}(Q\times\AA^1_s).
\]
Its central fiber is $Q\cup_{S_\infty}P$, where
\[
 P=\PP_S(\OO_S\oplus\OO_S(h)).
\]
The sections $S_0,S_\infty$ have normal bundles $\OO_S(h),\OO_S(-h)$, respectively; $Q$ meets $P$ along $S_\infty$.  The strict transform of $S\times\AA^1$ is an outer divisor, meeting the central fiber in $S_0$.

Let $\mathscr C$ be the strict transform of $C_*\times\AA^1$ in $\mathscr D$, and form
\[
 \mathscr Y=\operatorname{Bl}_{\mathscr C}\mathscr D.
\]
Denote the strict transform of the outer divisor by $\mathscr B_0$.

\begin{proposition}\label{prop:secondfamilyfull}
After restricting the parameter neighborhood as necessary, $\mathscr Y\to\AA^1_s$ is a projective semistable family.  Its general fiber is $Y$ and its central fiber is
\begin{equation}\label{eq:secondcentralfull}
 \mathscr Y_0=Q\cup_S W,\qquad
 W=\operatorname{Bl}_{C_*\subset S_0}P.
\end{equation}
If $B_0\subset W$ is the strict transform of $S_0$, then
\begin{equation}\label{eq:secondcanonicalfull}
 K_W+S_\infty+B_0=0,\qquad B_0\cap S_\infty=\varnothing,
 \qquad N_{B_0/W}=\OO_S.
\end{equation}
More generally,
\[
 \omega_{\mathscr Y/\AA^1_s}(\mathscr B_0)\simeq\OO_{\mathscr Y}.
\]
\end{proposition}
\begin{proof}
Use local coordinates $(x,u,v)$ on $Q$, with $S=(x=0)$ and $C_*=(x=u=0)$.  Blowing up $(x,s)$ gives charts
\[
 x=sz \quad\text{and}\quad s=xw,\qquad w=z^{-1}.
\]
In the first chart the strict transform of the moving center is $(z,u)=0$.  Its special fiber lies on $S_0=(z=0)$, not on the attaching section at infinity.  Blowing up this center gives the two smooth charts
\[
 z=u a,\qquad\text{and}\qquad u=z b.
\]
The map to the base is $s$ on these two charts and $xw$ on the untouched attaching chart. The central equations $s=0$ and $xw=0$ give \eqref{eq:secondcentralfull} in a smooth total space. Over $s\ne0$ the fiber is $\operatorname{Bl}_{C_*}Q$.

Write $E_P$ for the exceptional divisor of $W\to P$.  Since
\[
 S_0=S_\infty+p^*h,\qquad K_P=-S_0-S_\infty,
 \qquad B_0=b_W^*S_0-E_P,
\]
we obtain
\[
 K_W=b_W^*K_P+E_P=-S_\infty-B_0.
\]
The center is disjoint from $S_\infty$, which proves disjointness.  On $S_0$ it is the Cartier divisor $C_*\in|h|$.  The strict transform $B_0$ is therefore isomorphic to $S$, and
\[
 N_{B_0/W}=N_{S_0/P}\otimes\OO_S(-C_*)
 =\OO_S(h)\otimes\OO_S(-h)=\OO_S.
\]
For each of the two blowups $b$, the center has codimension two and lies with multiplicity one in the outer divisor $D$. With $E$ exceptional, $K_{\rm new}+\widetilde D=(b^*K_{\rm old}+E)+(b^*D-E)=b^*(K_{\rm old}+D)$. Starting with $K_Q+S=0$ proves the relative assertion.
\end{proof}

The two families share the quartic root:
\begin{equation}\label{diag:twodegenerationsfull}
\begin{tikzcd}[column sep=huge,row sep=large]
 X_t \arrow[r,rightsquigarrow] & Q\cup_S B\\
 Y \arrow[r,rightsquigarrow]\arrow[d,"\kappa"'] & Q\cup_S W\\
 \PP^1 &
\end{tikzcd}
\end{equation}
Section~\ref{sec:twodegoperators} specifies when their polarizations and relative moduli on $(Q,S)$ agree.

\subsection{The rigid surface in the complementary component}
Put
\[
 D_W=b_W^*p^*h-E_P.
\]
Let $R_*=\widetilde{p^{-1}C_*}\subset W$.  Write $s_-=R_*\cap S_\infty$, and let $f$ be the numerical class of a ruling fiber on $R_*$.

\begin{proposition}\label{prop:rigidruledfull}
The linear system $|D_W|$ consists of $R_*$ alone.  There are identifications and intersection formulas
\begin{gather}
 R_*\simeq\PP_{C_*}(\OO_{C_*}\oplus\OO_{C_*}(h)),\qquad
 R_*\cap B_0=\varnothing,\label{eq:ruledidentfull}\\
 s_-^2=-4,\quad s_-f=1,\quad f^2=0,\quad
 K_{R_*}=-2s_-,\quad N_{R_*/W}=\OO_{R_*}(-s_-).\label{eq:rulednormalfull}
\end{gather}
In particular $H^j(R_*,N_{R_*/W})=0$ for every $j$.
\end{proposition}
\begin{proof}
The blowup identity $b_{W*}\OO_W(-E_P)=I_{C_*/P}$ gives
\[
 H^0(W,\OO_W(D_W))
 =H^0(P,I_{C_*/P}\otimes p^*\OO_S(h))
 =H^0(S,I_{C_*/S}(h))=\CC.
\]
Here $p_*\OO_P=\OO_S$ and $I_{C_*/S}(h)=\OO_S$. The unique divisor is $R_*\simeq p^{-1}C_*$, since its blowup center is Cartier on $p^{-1}C_*$. The strict transforms of $u=0$ and $z=0$ meet the exceptional line of $\operatorname{Bl}_{(u,z)}$ at $[0:1]$ and $[1:0]$, respectively. Thus $R_*\cap B_0=\varnothing$.

The positive section on $R_*$ is linearly equivalent to $s_-+p_{C_*}^*(h|_{C_*})$ and numerically equivalent to $s_-+4f$. The exceptional divisor restricts to $s_+$. Thus in the Picard group
\[
 D_W|_{R_*}=p_{C_*}^*(h|_{C_*})-s_+=-s_-.
\]
The canonical formula for a ruled surface gives
\[
 K_{R_*}=p_{C_*}^*K_{C_*}-s_--s_+=-2s_-,\qquad K_{C_*}=h|_{C_*}.
\]
On each ruling fiber $f$, $\OO_{R_*}(-s_-)|_f=\OO_{\PP^1}(-1)$ and $H^{0,1}(f,\OO_f(-1))=0$. Base change gives $R(p_{C_*})_*N_{R_*/W}=0$.
\end{proof}

\begin{figure}[htbp]
\centering
\begin{tikzpicture}[x=1cm,y=1cm,font=\footnotesize,line cap=round,line join=round]
\definecolor{ksea}{HTML}{176568}\definecolor{ksrf}{HTML}{B15C25}\definecolor{kamb}{HTML}{53618A}
\node[anchor=west,font=\small\bfseries] at (0,3.50) {Before the blowup};
\node[anchor=west,font=\small\bfseries] at (7.20,3.50) {Above the center};
\path[fill=ksrf!12,draw=ksrf!80] (.56,.65) .. controls (1.44,.31) and (2.98,.36) .. (3.82,.81)
 --(3.82,2.69) .. controls (2.98,2.24) and (1.44,2.19) .. (.56,2.53)--cycle;
\foreach \x/\y in {.92/.54,1.42/.46,1.99/.43,2.56/.46,3.11/.55,3.59/.70}{
 \draw[ksrf!55] (\x,\y)--(\x,{\y+1.88});}
\draw[ksea,line width=1.35pt] (.25,.78)..controls (1.12,.33) and (2.95,.12)..(4.11,1.00);
\draw[kamb,line width=1.35pt] (.25,2.66)..controls (1.12,2.21) and (2.95,2.00)..(4.11,2.88);
\node[anchor=west,text=kamb] at (4.14,2.89) {$S_0$};
\node[anchor=west,text=ksea] at (4.14,.99) {$S_\infty$};
\node[fill=white,inner sep=2pt,text=ksrf] at (2.16,1.62) {$p^{-1}(C_*)$};
\fill[kamb] (2.05,2.315) circle (2.5pt);
\node[anchor=south,text=kamb] at (2.05,2.40) {$x\in C_*$};
\draw[ksea,line width=1.4pt] (8.76,1.20) .. controls (9.33,.90) and (10.29,1.88) .. (10.41,2.27);
\draw[ksea!65,line width=.75pt] (10.41,2.27) .. controls (10.11,2.57) and (8.78,1.76) .. (8.76,1.20);
\draw[ksrf,line width=1.4pt] (7.23,.65)..controls (7.85,.96) and (8.36,1.06)..(8.76,1.20) --(9.18,1.45);
\draw[kamb,line width=1.4pt] (10.04,2.26)--(10.41,2.27)..controls (10.77,2.30) and (11.66,2.61)..(12.27,2.99);
\fill[ksrf] (8.76,1.20) circle (2.5pt);
\fill[kamb] (10.41,2.27) circle (2.5pt);
\node[anchor=north west,text=ksrf] at (7.01,.62) {$R_*$};
\node[anchor=west,text=kamb] at (12.27,2.98) {$B_0$};
\node[text=ksrf,anchor=west] at (8.83,.81) {$[0:1]$};
\node[text=kamb,anchor=south] at (10.43,2.38) {$[1:0]$};
\node[text=ksea,anchor=west] at (10.39,1.39) {$E_x\simeq\mathbb P^1$};
\end{tikzpicture}
\caption{The ruled surface in the second degeneration. On the left the
rulings are fibers over $C_*$. The upper section contains the center of
the blowup, and the lower edge traces the intersection of the ruled surface with the common K3.
On the right we take a surface slice transverse to the center at $x$.
Write $p^{-1}(C_*)=(u=0)$ and $S_0=(z=0)$. The exceptional fiber
has homogeneous coordinates $[u:z]$. The strict transforms $R_*$ and
$B_0$ meet this fiber at its two displayed points.
Along $S_\infty$ the gluing curve remains $R_*\cap S_\infty\simeq C_*$.}
\label{fig:ruledblowupsliceII}
\end{figure}
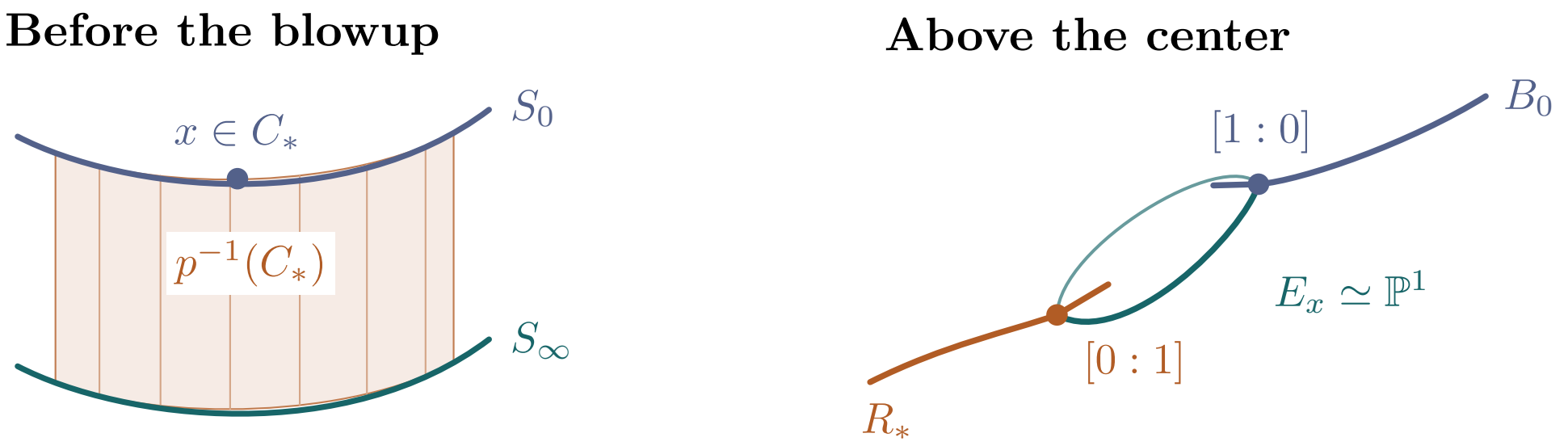

Figure~\ref{fig:ruledblowupsliceII} shows the proper transform of the ruled surface at the blowup center.

\subsection{Compatible polarizations and explicit charge conservation}\label{sec:bothchargesfull}
Start with $3H_Q$ on the deformation to the normal cone, twist by twice the root component, then subtract the exceptional divisor of the second blowup. The resulting line bundle restricts to
\begin{equation}\label{eq:secondpolarizationfull}
 2A+F_{\mathrm{fib}}=3A-E_C\text{ on }Y,\qquad
 H_Q\text{ on }Q,\qquad
 L_W=2b_W^*S_0+D_W\text{ on }W.
\end{equation}
For example, on the ruled component before the second blowup the restriction is $3p^*h+2S_\infty=2S_0+p^*h$, which proves the last expression after subtracting the exceptional divisor.  Both root restrictions on the gluing K3 are $h$.

To check ampleness on $W$, use the elementary transformation
\[
 0\longrightarrow\OO_S^{\oplus2}
 \xrightarrow{\operatorname{diag}(1,s_{C_*})}
 \OO_S\oplus\OO_S(h)
 \longrightarrow\OO_{C_*}(h)\longrightarrow0.
\]
The graph of the induced rational map $P\dashrightarrow S\times\PP^1$ is $W\hookrightarrow P\times_S(S\times\PP^1)$. The two tautological pullbacks are $b_W^*S_0$ and $b_W^*S_0-E_P$. Thus
\[
       L_W=b_W^*S_0+(b_W^*S_0-E_P)+p^*h .
\]
The first two terms are nef and their sum is relatively ample over $S$. Adding the pullback of the ample class $h$ gives ampleness on $W$.

The root and general-fiber restrictions in \eqref{eq:secondpolarizationfull} are ample as well.  Relative ampleness follows near the central fiber.

This general-fiber polarization restricts to $2h$ on a K3 fiber.  The modularity theorem uses the primitive restriction $h$, namely $A+F_{\mathrm{fib}}$.  Rank-one sheaves on integral fiber supports are stable for both polarizations.

For a rank-zero sheaf with character $(0,D,\beta,n)$ on a smooth threefold and polarization $L$, Riemann--Roch is
\begin{equation}\label{eq:RRrankzerofull}
 P(m)=\frac{DL^2}{2}m^2+
 \left(\beta L-\frac{DK_YL}{2}\right)m+
 n-\frac{\beta K_Y}{2}+\frac{D(K_Y^2+c_2(Y))}{12}.
\end{equation}
This formula follows by multiplying $\ch(F)e^{mL}$ by
$1-K_Y/2+(K_Y^2+c_2(Y))/12+\cdots$ and keeping degree six.

Consider first the ordinary nonzero-boundary root sector $D_Q=H_Q$, $D_B=H_B$.  Put
\[
 a=\beta_QH_Q=\beta_BS=\chi(G).
\]
With the original polarization $\mathcal A$, the component and boundary polynomials are
\begin{align}
 P_Q(m)&=2m^2+(a+2)m+n_Q+\frac a2+\frac73,\label{eq:Qpolyfull}\\
 P_B(m)&=8m^2+(2\beta_BH_B+a+2)m+n_B+\frac a2+\frac{11}{6},\label{eq:Bpolyfull}\\
 P_G(m)&=4m+a.\label{eq:boundarypolyfull}
\end{align}
Here $H_Q\operatorname{td}_2(Q)=7/3$, since $c_1(Q)=H_Q$ and $c_2(Q)=6H_Q^2$.  To verify the other Todd number, the plane transform has $\chi(\OO_T)=1$ and class $H_B$.  Applying \eqref{eq:RRrankzerofull} to its structure sheaf gives
\[
 1=H_B\operatorname{td}_2(B)+\frac14K_BH_B^2+\frac16H_B^3
   =H_B\operatorname{td}_2(B)-1+\frac16,
\]
so $H_B\operatorname{td}_2(B)=11/6$.  Also
$H_B(6H_B-E_B)^2=16$ and
$H_B(6H_B-E_B)S=4$, which give the displayed leading and linear terms.  The relation $\beta_BS=4\beta_BH_B-\beta_BE_B=a$ gives the stated form of the linear coefficient.

The normalization sequence says $P_X=P_Q+P_B-P_G$.  Since
\[
 P_X(m)=10m^2+2dm+n+25/6,
\]
we obtain the actual first-degeneration conservation equations
\begin{equation}\label{eq:firstchargesfull}
 d=a+\beta_BH_B,\qquad n=n_Q+n_B.
\end{equation}

For the second degeneration, for sheaves supported on the fibers, write the $W$ character as $(0,D_W,\beta_W,n_W)$.  Its support is $R_*$, so $\beta_W(D_W+S_\infty)=0$.  Put $a=\beta_QH_Q=\beta_WS_\infty$.  The intersection identities are
\[
 D_W^3=-4,\qquad D_W^2S_\infty=4,\qquad
 K_W=-2S_\infty-D_W,\qquad
 D_W\operatorname{td}_2(W)=-\frac13.
\]
The first two follow by restricting $D_W=-s_-$ and $S_\infty=s_-$ to $R_*$.  The last follows from $\chi(\OO_{R_*})=-2$:
\[
 -2=D_W\operatorname{td}_2(W)+\frac14K_WD_W^2+\frac16D_W^3
 =D_W\operatorname{td}_2(W)-1-\frac23.
\]
Therefore
\begin{align}
 P_W(m)&=6m^2+
 \left(\beta_W(2b_W^*S_0-S_\infty)+2\right)m
 +n_W+\frac a2-\frac13,\label{eq:Wpolyfull}\\
 P_Y(m)&=8m^2+2(\widetilde\beta A)m+\widetilde n+2.
\end{align}
The equality $P_Y=P_Q+P_W-P_G$ gives
\begin{equation}\label{eq:secondchargesfull}
 \widetilde\beta A=\beta_Wb_W^*S_0,\qquad
 \widetilde n=n_Q+n_W,\qquad
 \beta_QH_Q=\beta_WS_\infty.
\end{equation}

For a sheaf $F=j_*V$ on a pencil member $j:Y_b\hookrightarrow Y$,
\[
 F\otimes K_Y=j_*\bigl(V\otimes\OO_{\PP^1}(-1)|_b\bigr).
\]
Thus threefold Serre duality pairs its self-Ext groups after choosing a basis of this one-dimensional fiber, as in \cite{GS}.

On $\mathscr Y$ one has $\omega_{\mathscr Y/\AA^1}=\OO(-\mathscr B_0)$. Vanishing derived restriction to $\mathscr B_0$ identifies $F(-\mathscr B_0)\to F$ with an isomorphism. The determinant trivialization at the outer divisor gives the top-trace reduction of Part~I, Section~\companionref{sec:zerorelative}. To compare with the full pencil invariant, the same construction must extend over sheaves specializing to the chosen outer fiber.
\subsection{Sheaves on the second degeneration and its two boundary maps}\label{sec:secondmoduli}
Fix a fiber charge $(0,F_{\mathrm{fib}},\widetilde\beta,\widetilde n)$ on $Y$.

Let $M_Y$ parametrize pure sheaves stable for $A+F_{\rm fib}$ in this charge. The effective-cycle calculation in Section~\ref{intro:pencilgeometry} gives $F=j_*G$ on one integral fiber $j:Y_b\hookrightarrow Y$. Since $K_Y=\kappa^*\OO(-1)$,
\[
        F\otimes K_Y\simeq j_*(G\otimes\OO(-1)|_b).
\]
Trivializing this one-dimensional factor and applying Serre duality gives
\[
 \ext^2_Y(j_*G,j_*G)=\ext^1_Y(j_*G,j_*G),\qquad
 \ext^3_Y(j_*G,j_*G)=1.
\]
The two-term theory with tangent and obstruction given by the middle Ext groups and its virtual dimension zero are those of \cite[Proposition~2.2 and Theorem~2.3]{GS}. Its degrees for fixed Hilbert polynomial are the coefficients of $K^{\mathrm{raw}}$.

Fix the numerical types on $Q$ and $W$ with the polarizations \eqref{eq:secondpolarizationfull}. The root stacks parametrize base-flat pure sheaves and have their derived enhancements inside the perfect-complex stacks. On $W$, set $D_\partial=S_\infty\sqcup B_0$. The evaluation is
\[
 F_W\longmapsto(E_\infty,E_0)
       =(\LL i_{S_\infty}^*F_W,\LL i_{B_0}^*F_W)
          \in\mathbf{Perf}(S_\infty)\times\mathbf{Perf}(B_0).
\]
The root $Q$ matches $E_\infty$. The second entry is the outer boundary value. Fixing $E_0$ alone, with its identification, takes the homotopy fiber of $K_F[1]\to K_0[1]$ in the notation below. In particular $E_0=0$ specifies an outer boundary condition.
\begin{theorem}\label{thm:secondends}
For a perfect sheaf $F$ on $W$ let $K_F=\RHom_W(F,F)$, $K_\infty=\RHom_{S_\infty}(E_\infty,E_\infty)$ and $K_0=\RHom_{B_0}(E_0,E_0)$. Before rigidification the tangent complex with both restrictions fixed, including the chosen identifications, is $\RHom_W(F,F(-D_\partial))[1]$. For arbitrary derived boundary values one has
\[
 2\chi_W(F,F)=\chi_{S_\infty}(E_\infty,E_\infty)+\chi_{B_0}(E_0,E_0).
\]
\end{theorem}
\begin{proof}
Tensor $\OO_W(-D_\partial)\to\OO_W\to i_*\OO_{D_\partial}$ with $F^\vee\otimes^{\LL}F$. Derived adjunction and disjointness give the triangle
\[
 \RHom_W(F,F(-D_\partial))\longrightarrow K_F
       \longrightarrow K_\infty\oplus K_0\longrightarrow .
\]
After shifting by one, its middle arrow is the differential of the evaluation functor. Taking the fiber proves the assertion for deformations with the two boundary values fixed. Since $K_W=-D_\partial$, Serre duality identifies the first term with $K_F^\vee[-3]$.

Its Euler characteristic is $-\chi_W(F,F)$, which proves the formula.
\end{proof}
On $E_0=0$, assume $F$ is simple and $\Ext^3_W(F,F)=0$. Rigidification gives a two-term theory with
\[
   \operatorname{vdim}=1-\chi_W(F,F)
           =1+\tfrac12\langle v(E_\infty),v(E_\infty)\rangle .
\]
In the primitive stable-curve sector this equals three. The top-Ext vanishing is a separate condition for a general derived boundary.

Compactifying this open problem also requires accounting for sheaves which specialize to the selected outer fiber. The comparison with $M_Y$ is calculated in the fiber structure-sheaf charge in Theorem~\ref{thm:K3structurecomparison}.

\subsection{The numerical invariant relative to a smooth outer K3} Let $Y=\mathrm{Bl}_{C_*}Q\to\mathbb P^1$ be a general Lefschetz quartic pencil and $D_0$ a general smooth fiber. For $\ch_1=F$ the primitive fiber class, and every allowed pair of lower charges, the relative invariant of Gholampour and Sheshmani \cite[\S2.3, Lemmas~2.20--2.21]{GS} satisfies
\begin{equation}\label{eq:outer}
 \DT_{\rm GS}(Y/D_0,\ch)=\DT_{\rm GS}(Y,\ch).
\end{equation}
The fibers are integral, the singular fibers are one-nodal, and $K_Y=\kappa^*\OO(-1)$. Primitivity excludes strict semistability. For the cap $D_0\times\PP^1$ let $A_{\rm cap}$ and $R_{\rm cap}$ denote its absolute and one-relative degrees in the fiber charge. The product obstruction calculation gives $A_{\rm cap}=0$, and degeneration of the cap gives $A_{\rm cap}=2R_{\rm cap}$ \cite[Lemmas~2.20--2.21]{GS}. Degeneration along $D_0$ therefore yields
\[
 \DT_{\rm GS}(Y,\ch)=\DT_{\rm GS}(Y/D_0,\ch)+\sum_{\text{cap lifts}}R_{\rm cap}
                     =\DT_{\rm GS}(Y/D_0,\ch);
\]
only the root or a cap can carry the primitive fiber class. This proves \eqref{eq:outer}. At the internal cut $Q\cup_S W$ one must still identify its induced outer theory with this relative theory.

\section{The primitive quintic comparison}\label{sec:quinticapplication}
Fix $v=(0,H,\beta,n)$ and a moduli space satisfying Hypothesis~\companionref{hyp:compact} with complete generic fiber $M_t(v)$. Primitivity gives stability on the smooth fiber. Descent and the valuative criteria on expansions are verified for the structure sheaf, the first point and the line in Theorems~\companionref{thm:verifiedcompactification}, \ref{thm:one-defect-global} and~\ref{thm:expanded-line-compact}.

Suppose its marked sectors satisfy the pure-boundary and relative product conditions of Hypothesis~\companionref{hyp:product}.

\begin{corollary}\label{cor:quinticVFC}
For every such charge and compactification system,
\begin{equation}\label{eq:quinticVFC}
\spc[M_t(v)]^{\vir}
=\sum_{\gamma\in\Lambda_v}(g_\gamma)_*
\Delta_{J_{\chi_\gamma}}^!
\left([M_{1,\gamma}^{\rel}]^{\vir}
\times[M_{2,\gamma}^{\rel}]^{\vir}\right)
\end{equation}
in $A_0(\MM_0(v))_{\QQ}$ whenever every boundary vector is $(0,h,\chi_\gamma)$. The formula includes all bubbles through the relative moduli, and its virtual dimension is zero. More general smooth stable boundary sectors are included by replacing $J_{\chi_\gamma}$ by their actual smooth moduli and retaining their dimensions. For a general derived boundary one uses the full matching formula \textup{(\companionref{eq:mainfull})} under the compactification and obstruction-theory hypotheses of Theorem~\companionref{thm:main}.
\end{corollary}
\begin{proof}
Part~I, Section~\companionref{sec:geometry} and Proposition~\companionref{prop:primitive}, verifies the geometry and identifies the generic moduli. Put $a_{i,\gamma}=[M_{i,\gamma}^{\rel}]^{\vir}$. The specialization and the product hypotheses in Part~I, Theorem~\companionref{thm:main}, give
\[
 \spc[M_t(v)]^{\vir}=\sum_\gamma(g_\gamma)_*
 [\MM^\dagger_\gamma]^{\vir},\qquad
 [\MM^\dagger_\gamma]^{\vir}=\Delta^!_{J_{\chi_\gamma}}
             (a_{1,\gamma}\times a_{2,\gamma}).
\]
For $v_S=(0,h,\chi_\gamma)$, $v_S^2=4$; the relative dimensions and diagonal codimension are $1+v_S^2/2=3$ and $2+v_S^2=6$ (Part~I, Propositions~\companionref{prop:Jacobian} and~\companionref{prop:relativevdim}). Thus each term lies in $A_{3+3-6}=A_0$, as claimed.
\end{proof}

\begin{corollary}\label{cor:degreepairing}
Let $\operatorname{ev}_i:M_{i,\gamma}^{\rel}\to J_{\chi_\gamma}$ be the proper boundary evaluations and put
\begin{equation}\label{eq:boundarycycles}
\alpha_{i,\gamma}=(\operatorname{ev}_i)_*[M_{i,\gamma}^{\rel}]^{\vir}
\in A_3(J_{\chi_\gamma})_{\QQ}.
\end{equation}
Under the hypotheses of Corollary~\ref{cor:quinticVFC},
\begin{equation}\label{eq:numericalpairing}
\dt_{X_t}(v)
=\sum_{\gamma\in\Lambda_v}
\int_{J_{\chi_\gamma}}\alpha_{1,\gamma}\cdot\alpha_{2,\gamma}.
\end{equation}
\end{corollary}
\begin{proof}
Set $Z=M_{1,\gamma}^{\rel}\times_J M_{2,\gamma}^{\rel}$ and write $e:Z\to J$ for the common evaluation. The square
\[
\begin{tikzcd}[row sep=small]
 Z\arrow[r]\arrow[d,"e"'] &M_{1,\gamma}^{\rel}\times M_{2,\gamma}^{\rel}\arrow[d,"{\operatorname{ev}_1\times\operatorname{ev}_2}"]\\
 J\arrow[r,"\Delta_J"']&J\times J
\end{tikzcd}
\]
is Cartesian. For $a_i=[M_{i,\gamma}^{\rel}]^{\vir}$, proper pushforward commutes with refined pullback:
\[
 e_*\Delta_J^!(a_1\times a_2)
 =\Delta_J^!\big((\operatorname{ev}_1)_*a_1\times
                     (\operatorname{ev}_2)_*a_2\big)
 =\alpha_{1,\gamma}\cdot\alpha_{2,\gamma}.
\]
Taking degrees in \eqref{eq:quinticVFC} gives the formula, since specialization in a proper family preserves degree. The projection formula also applies to insertions extending over that family.
\end{proof}

\subsection{The complete list of boundary types in the primitive formula}\label{sec:primitivefull}
For fixed primitive Hilbert polynomial and finite intrinsic expansion stabilizers, Theorem~\companionref{thm:primitiveprofiles} gives only two boundary curve classes, $h$ and $0$.

Suppose a moduli system as in Hypothesis~\companionref{hyp:compact} is given, and suppose every nonzero marked boundary is pure. Partition the finite set of marked numerical sectors into
\[
 \Lambda_v=\Lambda_v^+\sqcup\Lambda_v^0
\]
according as that boundary is nonzero or zero. These subsets are open and closed in the marked moduli: the boundary sheaf is flat and its Hilbert polynomial is locally constant. The nonzero vectors are $(0,h,\chi_\gamma)$.

For $\gamma\in\Lambda_v^0$, the profile classification gives $F_L=0$. Any empty bubble on the left would have a rescaling stabilizer $\mathbb G_m$. The left target is therefore $Y_1$. Write $M_{2,\gamma}^{\rel,0}$ for the right relative moduli with zero outer restriction and the inherited stability and numerical conditions.

\begin{theorem}\label{thm:primitivecompleteformula}
Assume Hypothesis~\companionref{hyp:compact} and purity of all nonzero marked boundaries. Suppose the positive-boundary sectors satisfy the relative product hypothesis of Section~\companionref{sec:diagonal}. Then
\begin{equation}\label{eq:completeprimitive}
\begin{split}
 \spc[M_\eta(v)]^{\vir}
 ={}&\sum_{\gamma\in\Lambda_v^+}(g_\gamma)_*
 \Delta_{J_{\chi_\gamma}}^!
 \left([M_{1,\gamma}^{\rel}]^{\vir}
           \times[M_{2,\gamma}^{\rel}]^{\vir}\right)\\
 &+\sum_{\gamma\in\Lambda_v^0}(g_\gamma)_*
            [M_{2,\gamma}^{\rel,0}]^{\vir}
 \quad\text{in }A_0(\MM_0(v))_{\QQ}.
\end{split}
\end{equation}
The second sum carries the zero-boundary theory of Theorem~\companionref{thm:zeroboundaryPOT}. Both sums have finitely many proper terms. The positive terms use only the primitive boundary vector $(0,h,\chi_\gamma)$.
\end{theorem}
\begin{proof}
Apply Theorem~\companionref{thm:main} in its full-boundary form and the finite decomposition of Proposition~\companionref{prop:finitecharges}. Theorem~\companionref{thm:primitiveprofiles} proves that every nonzero boundary cycle is $h$ and that the nonzero branches have integral multiplicity-one supports. The additional purity gives the smooth Jacobian evaluation. The product hypotheses and Theorem~\companionref{thm:diagonal} identify the corresponding classes with the first sum.

If the marked boundary is zero, the two branch sheaves glue as a direct sum with disjoint supports. Indeed for $E=0$ the cross-Hom groups vanish and
\[
 F=F_L\oplus F_R,\qquad
 \End(F)=\End(F_L)\oplus\End(F_R)=\CC .
\]
The nonzero summand is $F_R$ by the profile theorem. Restriction to the right and attachment of $(Y_1,0)$ are inverse on flat families and their arrows. Theorem~\companionref{thm:zeroboundaryPOT} identifies the determinant and obstruction theory, giving the second sum.

Each marked numerical piece is open and closed in a finite cover of the proper special fiber. Its pushforward includes the stack multiplicity of the normalized boundary (Remark~\companionref{rem:multiplicities}). Unmarked bubbles belong to the relative factors.
\end{proof}

The numerical label also records the lower Chern characters. For the normalization, with $i_i$ the component inclusions and $j_e$ the joint inclusions, the matching sequence gives
\begin{equation}\label{eq:fullchargeconservation}
 [F]=\sum_i(i_i)_*[F_i]-\sum_e(j_e)_*[G_e].
\end{equation}
Thus the two lower Chern characters and the Euler characteristics must be tracked together with the divisorial profiles. Vertical twists preserve the total original Hilbert polynomial but can change $\chi_\gamma$, and are included in these labels.

\begin{corollary}\label{cor:completepairing}
Under the hypotheses of Theorem~\ref{thm:primitivecompleteformula}, define
\[
 \alpha_{i,\gamma}=(\operatorname{ev}_i)_*[M_{i,\gamma}^{\rel}]^{\vir}
 \in A_3(J_{\chi_\gamma})_{\QQ},\qquad
 z_\gamma=\deg[M_{2,\gamma}^{\rel,0}]^{\vir}.
\]
Then
\begin{equation}\label{eq:completepairing}
 \dt_{X_\eta}(v)
 =\sum_{\gamma\in\Lambda_v^+}
       \int_{J_{\chi_\gamma}}\alpha_{1,\gamma}\cdot\alpha_{2,\gamma}
   +\sum_{\gamma\in\Lambda_v^0}z_\gamma.
\end{equation}
\end{corollary}
\begin{proof}
Take degrees in \eqref{eq:completeprimitive}, use proper pushforward, and apply the projection formula for the regular diagonals, exactly as in Corollary~\ref{cor:degreepairing}.
\end{proof}

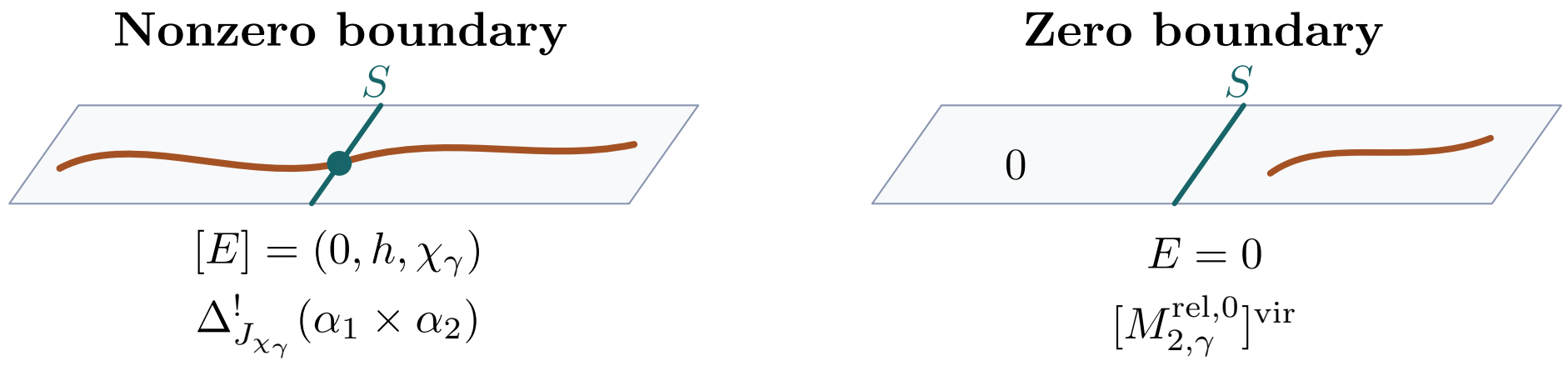
\begin{figure}[htbp]\centering
\begin{tikzpicture}[x=1cm,y=1cm,font=\footnotesize,line cap=round,line join=round]
\foreach \x/\name in {0/{Nonzero boundary},6.85/{Zero boundary}}{
 \begin{scope}[shift={(\x,0)}]
 \node[font=\small\bfseries] at (2.72,2.87) {\name};
 \path[draw=qIIblue!65,fill=qIIblue!4] (.09,1.52)--(.64,2.30)--(3.04,2.30)--(2.49,1.52)--cycle;
 \path[draw=qIIblue!65,fill=qIIblue!4] (2.49,1.52)--(3.04,2.30)--(5.56,2.30)--(5.01,1.52)--cycle;
 \draw[qIIteal,line width=1.15pt] (2.49,1.52)--(3.04,2.30);
 \node[text=qIIteal] at (3.00,2.49) {$S$};
 \end{scope}}
\draw[qIIclay,line width=1.55pt] (.49,1.80) .. controls (1.07,2.11) and (1.96,1.67).. (2.71,1.84)
 .. controls (3.54,2.11) and (4.29,1.83)..(5.05,1.99);
\fill[qIIteal] (2.71,1.84) circle (2.8pt);
\node at (2.70,1.13) {$[E]=(0,h,\chi_\gamma)$};
\node at (2.70,.55) {$\Delta^!_{J_{\chi_\gamma}}(\alpha_1\times\alpha_2)$};
\draw[qIIclay,line width=1.45pt] (10.10,1.76)..controls (10.54,2.08) and (11.22,1.79)..(11.85,2.04);
\node at (8.08,1.83) {$0$};
\node at (9.58,1.13) {$E=0$};
\node at (9.58,.55) {$[M_{2,\gamma}^{\rel,0}]^{\vir}$};
\end{tikzpicture}
\caption{A schematic cut at one marked join. On the left both relative classes meet over the same boundary object, with $\alpha_i=[M_{i,\gamma}^{\rel}]^{\vir}$. On the right the entire left branch is empty. The right relative factor then carries the zero-boundary obstruction theory.}\label{fig:II-occupied-cut}
\end{figure}

Figure~\ref{fig:II-occupied-cut} shows the two types of marked summand in \eqref{eq:completeprimitive}.

\section{The complementary boundary classes and their punctual factors}\label{sec:boundarycalibrations}
The two degenerations in \eqref{diag:twodegenerationsfull} pair the proposed common quartic class with different complementary boundary classes. Their values before punctual modification and their first coefficients fix the normalization required of a correspondence in all degrees.

\subsection{The Abelian fiber and the structure-sheaf section}
Let $J=J_{-2}$ be the six-dimensional compactified Jacobian of pure sheaves of Mukai vector $(0,h,-2)$ on $S$.  Its support map is denoted $\pi:J\to|h|=\PP^3$.  The universal Cartier curve gives a section
\[
 \sigma:\PP^3\longrightarrow J,\qquad C\longmapsto\OO_C.
\]
For a smooth $C_*\in|h|$, let $\iota:\operatorname{Pic}^0(C_*)\hookrightarrow J$ denote its fiber.  Put
\[
 s=\sigma_*[\PP^3],\qquad f=\iota_*[\operatorname{Pic}^0(C_*)]\quad\text{in }A_3(J)_{\QQ}.
\]

\begin{proposition}\label{prop:sfintersectionsfull}
The intersection numbers of these two classes are
\begin{equation}\label{eq:sfintersectionsfull}
 \int_Js f=1,\qquad \int_Jf^2=0,\qquad \int_Js^2=-4.
\end{equation}
In particular $s$ and $f$ are not proportional, even numerically.
\end{proposition}
\begin{proof}
Every $C\in|h|$ is integral and $\OO_C$ is stable. The universal Cartier curve is flat, so $\pi\sigma=\id$ defines a closed section. At a smooth fiber $f=\pi^*[C_*]$ in operational Chow cohomology, capped with $[J]$. The projection formula gives
\[
  \int_Jsf=\int_{\PP^3}\sigma^*\pi^*[C_*]=1,\qquad
  f^2=\pi^*([C_*]^2)\cap[J]=0.
\]
The K3 trace pairing defines the symplectic form $\omega_J$. Since $H^0(\PP^3,\Omega^2)=0$ and $2\dim\PP^3=\dim J$, $\sigma^*\omega_J=0$ identifies $N_\sigma\simeq\Omega^1_{\PP^3}$. Self-intersection and the Euler sequence therefore give
\[
   \int_Js^2=\int_{\PP^3}c_3(N_\sigma)
        =[H^3](1-H)^4=-4.
\]
The pairing matrix in the order $(s,f)$ has determinant $-1$.
\end{proof}

\subsection{Three elementary relative moduli problems}
We next identify the structure-sheaf and degree-zero line-bundle families whose boundary pushforwards are $s$ and $f$.  For each sheaf $F$ on its threefold $W$ we use the deformation and obstruction groups of $R\operatorname{Hom}_W(F,F)$.

\begin{proposition}\label{prop:Bplanefamilyfull}
The structure sheaves of the transforms of planes in $B$ form a smooth proper moduli component $\PP^{3\vee}$.  Their obstruction space for sheaves on $B$ is zero, and their boundary evaluation is $\sigma$.  Its virtual pushforward is $s$.
\end{proposition}
\begin{proof}
Every plane cuts an integral divisor $T\in|H_B|$, since no plane contains the integral curve $\Gamma$.  The sequence
\[
 0\to\OO_B(-H_B)\to\OO_B\to\OO_T\to0
\]
and $R\rho_*\OO_B=\OO_{\PP^3}$ imply $H^0(\OO_T)=\CC$ and $H^{>0}(\OO_T)=0$.  Similarly,
\[
 0\to\OO_B\to\OO_B(H_B)\to\OO_T(H_B)\to0
\]
gives $h^0(\OO_T(H_B))=3$ and $H^{>0}(\OO_T(H_B))=0$.

Applying $\mathcal Hom_B(-,\OO_T)$ to its Cartier resolution gives the complex $[\OO_T\xrightarrow{0}N_{T/B}]$. Hence
\[
 \Ext^1_B(\OO_T,\OO_T)=H^0(T,N_{T/B})=T_{[T]}|H_B|,\qquad
 \Ext^2_B(\OO_T,\OO_T)=0 .
\]
The annihilator recovers $T$. Thus $|H_B|\to M_B$ is proper and injective on geometric points. By the Ext calculation it is \'etale, hence an isomorphism onto an open and closed locus, also at singular plane transforms. Since $C=T\cap S$ is Cartier on $T$,
\[
     \mathbf L i_S^*\OO_T=\OO_C,\qquad
     \operatorname{ev}_*[|H_B|]=\sigma_*[\PP^3]=s .
\]
\end{proof}

\begin{proposition}\label{prop:Wlinefamilyfull}
For $R_*\subset W$ of Proposition~\ref{prop:rigidruledfull}, the sheaves
\[
 i_*(p_{C_*}^*L),\qquad L\in\operatorname{Pic}^0(C_*),
\]
form a smooth proper moduli component with vanishing obstruction space for sheaves on $W$.  Its boundary evaluation at $S_\infty$ is the inclusion of the Jacobian fiber, and its virtual pushforward is $f$.
\end{proposition}
\begin{proof}
Put $G=p^*L$. The two local self-Ext sheaves of $i_*G$ are $i_*\OO_{R_*}$ and $i_*N_{R_*/W}$. By Proposition~\ref{prop:rigidruledfull} the latter is acyclic. Since $Rp_*\OO_{R_*}=\OO_{C_*}$, the local-to-global sequence gives
\[
 \Ext^j_W(i_*p^*L,i_*p^*L)=H^j(C_*,\OO_{C_*})=
 \begin{cases}\CC&j=0,\\\CC^3&j=1,\\0&j\ge2.\end{cases}
\]
The differential of $p^*:\operatorname{Pic}^0(C_*)\to M_W$ is the isomorphism on $H^1$ above. The Jacobian is proper and the sheaf moduli is smooth here, so its image is an entire component. On $S_\infty$ the restriction is $L$; on $B_0$ it is zero since $R_*\cap B_0=\varnothing$.
\end{proof}

\begin{proposition}\label{prop:Qstructurefamilyfull}
The minimal relative expansion compactification of the hyperplane structure-sheaf family on $(Q,S)$ is
\[
 \widehat P=\operatorname{Bl}_{[S]}|H_Q|=\operatorname{Bl}_{[S]}\PP^4.
\]
Write $b_P:\widehat P\to\PP^4$ and $r:\widehat P\to|h|=\PP^3$.  It has the relative obstruction line bundle $b_P^*\OO_{\PP^4}(-1)$ and boundary evaluation $\sigma r$.  Consequently its virtual pushforward is $-s$.
\end{proposition}
\begin{proof}
The family-groupoid identification and the unchanged abstract supporting family are proved in Theorem~\companionref{thm:verifiedcompactification}.  We calculate its obstruction theory here to identify the boundary cycle used by the two degeneration operators.

The universal hyperplane surface over $\PP^4$ has exact sequence
\[
 0\to\OO_Q(-H_Q)\boxtimes\OO_{\PP^4}(-1)
 \to\OO_{Q\times\PP^4}\to\OO_{\mathcal T}\to0.
\]
The quartic threefold has $H^{>0}(\OO_Q)=0$, and duality gives
$H^3(Q,\OO_Q(-H_Q))=\CC$, with the other cohomology groups zero.  Therefore
\[
 R^1\pi_*\OO_{\mathcal T}=0,\qquad
 R^2\pi_*\OO_{\mathcal T}=\OO_{\PP^4}(-1).
\]
The support family on $\widehat P$ is its pullback, as just proved.  Line-bundle deformations vanish since $H^1(\OO_T)=0$.  On an unexpanded support $H^1(N_{T/Q})=0$, so the Cartier self-Ext sequence identifies the obstruction line bundle with $H^2(\OO_T)$.  At an expanded graph the same normal-cohomology vanishing holds, since its normal bundle is $\OO_S(h)$.

The target-smoothing parameter gives an unobstructed direction in the explicit charts above, and the geometric rescaling removes its corresponding scalar direction.  Neither operation changes the degree-two sheaf obstruction line bundle.  The functorial self-Ext edge maps consequently give the global line bundle $b_P^*\OO(-1)$ as the obstruction bundle.

Write $H_P=b_P^*c_1(\OO_{\PP^4}(1))$. The fiber of $r$ maps to a line through $[S]$, so $\int_{r^{-1}(C)}H_P=1$. The obstruction-bundle formula on the smooth fourfold gives
\[
 (\sigma r)_*[\widehat P]^{\vir}
 =-\sigma_*r_*(H_P\cap[\widehat P])=-\sigma_*[\PP^3]=-s.
\]
\end{proof}

Thus the two elementary boundary pairings are
\begin{equation}\label{eq:elementarypairingsfull}
 \int_J(-s)s=4,\qquad \int_J(-s)f=-1.
\end{equation}
These are pairings of the stated proper relative components.  Their use is to calibrate the actual complementary operators: the plane and ruled components are not interchangeable boundary states.

\subsection{Virtual gluing in the structure-sheaf charge}\label{sec:vacuumvirtual}
The compactification of Section~\companionref{sec:verifiedcompactification} verifies the geometric hypotheses of the virtual comparison for this charge.

\begin{lemma}\label{lem:vacuumobstructionbundles}
The obstruction bundles for sheaves on the respective threefolds on the smooth relative spaces of \textup{(\companionref{eq:vacuumrelativefactors})} are
\begin{equation}\label{eq:vacuumobsbundles}
 \operatorname{Ob}_{\mathfrak D_Q}=b_P^*\OO_{\PP^4}(-1),\qquad
 \operatorname{Ob}_{\mathfrak D_B^+}=0,\qquad
 \operatorname{Ob}_{\mathfrak D_B^0}=\OO_{\PP^4}(-1)^{\oplus4}.
\end{equation}
In the last case the theory is the zero-boundary theory of Theorem~\companionref{thm:zeroboundaryPOT}.  These are identifications of the obstruction sheaves of the stated global perfect obstruction theories, including the expanded boundary.
\end{lemma}
\begin{proof}
For a relative Cartier divisor $i:T\hookrightarrow Z$ the resolution
\[
 [\OO_Z(-T)\xrightarrow{s_T}\OO_Z]\longrightarrow i_*\OO_T
\]
gives, after applying $\mathcal Hom_Z(-,i_*\OO_T)$, the two sheaves $i_*\OO_T$ in degree zero and $i_*N_{T/Z}$ in degree one with zero differential.  Both terms in the resolution are globally defined line bundles.  Consequently derived global sections give the global self-Ext identification
\begin{equation}\label{eq:vacuumkoszulglobal}
 \RHom_Z(\OO_T,\OO_T)
 \simeq \RG(T,\OO_T)\oplus\RG(T,N_{T/Z})[-1].
\end{equation}
The construction and its edge maps commute with base change.  It remains to compute $H^1(N)$ and $R^2\pi_*\OO_T$, and to account for the expansion base.

For $\mathfrak D_Q$, the universal supporting family is the pullback of the universal hyperplane surface in $Q\times\PP^4$: this was proved by the invertibility of $(f_1,s)$ on that support.  Its divisor sequence gives
\[
 R^1\pi_*\OO_{\mathcal T}=0,\qquad
 R^2\pi_*\OO_{\mathcal T}=b_P^*\OO_{\PP^4}(-1),
\]
since $H^3(Q,\OO_Q(-H_Q))=\CC$ and all other cohomology of this line bundle vanishes.  On an unexpanded support $H^1(N)=0$ follows from $0\to\OO_Q\to\OO_Q(H_Q)\to N\to0$.  On the expanded graph $N=\OO_S(h)$, again with $H^1(N)=0$.  Thus \eqref{eq:vacuumkoszulglobal} identifies the relative obstruction line bundle with the first bundle in \eqref{eq:vacuumobsbundles}.

For the plane family in $B$, both $H^{>0}(\OO_T)$ and $H^{>0}(N)$ vanish by the two divisor sequences with $\OO_B(-H_B)$ and $\OO_B(H_B)$.  Hence $\operatorname{Ext}^2_B(\OO_T,\OO_T)=0$.

For $\mathfrak D_B^0$, put $D_B=5H_B-E_B$.  The universal abstract support is unchanged by its relative expansion, because the center ideal in \textup{(\companionref{eq:vacuumBzeroeq})} is invertible on that support.  Its sequence on $B\times|D_B|$ is
\[
 0\to\OO_B(-D_B)\boxtimes\OO_{\PP^4}(-1)
 \to\OO_{B\times\PP^4}\to\OO_{\mathcal T}\to0.
\]
Here $H^{>0}(B,\OO_B)=0$ and
\[
 H^3(B,\OO_B(-D_B))
 =H^0(B,\OO_B(K_B+D_B))^*
 =H^0(B,\OO_B(H_B))^*=\CC^4.
\]
Therefore
\begin{equation}\label{eq:vacuumzeroHodge}
 R^1\pi_*\OO_{\mathcal T}=0,\qquad
 R^2\pi_*\OO_{\mathcal T}=\OO_{\PP^4}(-1)^{\oplus4}.
\end{equation}
The supporting surface misses the outer divisor, also on the expanded graph-plus-plane fiber. Relative adjunction and duality give $N_{T/Z}=\omega_T$ and $H^1(N_{T/Z})=H^1(\OO_T)^\vee=0$. The zero-boundary determinant reduction removes the scalar and its top dual; $H^2(\OO_T)$ survives as in \eqref{eq:vacuumzeroHodge}.

On a framed chart for $\mathfrak D_Q$ the exceptional coordinate is its free expansion-smoothing parameter. On $\mathfrak D_B^0$ it is the coefficient $c$ in \textup{(\companionref{eq:vacuumBzeroeq})}. The other support coefficients are independent, so each chart is smooth over the expansion base. Every base tangent lifts, and the tangent sequence gives
\[
 T_{\rm base}\xrightarrow{0}h^1(T_{\rm relative})
               \xrightarrow{\sim}h^1(T_{\rm absolute}).
\]
The rescaling direction has degree $-1$. Thus \eqref{eq:vacuumobsbundles} holds for the absolute theories and descends on overlaps.
\end{proof}

\begin{theorem}\label{thm:vacuumvirtual}
Let $\mathfrak D$ be the projective expanded sheaf moduli of Theorem~\companionref{thm:verifiedcompactification}.  Let $g_+$ and $g_-$ be the two maps from the normalized special fiber to $\mathfrak D_0$.  In its rational Chow group one has
\begin{equation}\label{eq:vacuumVFC}
\begin{split}
 \spc[M_\eta(v_0)]^{\vir}
 ={}&(g_+)_*\Delta_{J_{-2}}^!
 \left(
    [\mathfrak D_Q]^{\vir}\times[\mathfrak D_B^+]^{\vir}
 \right)\\
 &+(g_-)_*\left(
 c_4\bigl(\OO_{\PP^4}(-1)^{\oplus4}\bigr)\cap[\PP^4]
 \right).
\end{split}
\end{equation}
In particular
\begin{equation}\label{eq:vacuumdegreeverified}
 \dt_X(v_0)=4+1=5.
\end{equation}
The same statements hold after tensoring by $\OO_X(aH)$, with the corresponding boundary line-bundle twists and translated charges.
\end{theorem}
\begin{proof}
Theorem~\companionref{thm:verifiedcompactification} proves (K1)--(K4), identifies the complete generic moduli, and identifies both marked relative products as actual proper family groupoids.  Every universal sheaf is perfect by its length-one Cartier resolution.  Thus Theorems~\companionref{thm:ambientPOT} and \companionref{thm:diagonal} and the virtual-normalization argument apply without an additional geometric hypothesis.  The marked normalization is $P_+\sqcup P_-$; its two maps have generic degree one.

The first piece is the full product \textup{(\companionref{eq:vacuumpositiveproduct})}.  The second piece has zero left object and the relative determinant theory of $\mathfrak D_B^0$.  Lemma~\ref{lem:vacuumobstructionbundles} gives its virtual class.  These facts prove \eqref{eq:vacuumVFC}.

We compute the degrees of its two terms.  Write $H_P=b_P^*c_1(\OO_{\PP^4}(1))$ on $\mathfrak D_Q$.  On the smooth moduli with obstruction line bundle in \eqref{eq:vacuumobsbundles}, the intrinsic-cone definition gives
\[
 [\mathfrak D_Q]^{\vir}=-H_P\cap[\mathfrak D_Q].
\]
The projection $\rho:\operatorname{Bl}_{[S]}\PP^4\to\PP^3$ is a projective-line bundle and $H_P$ has degree one on its fibers.  Thus
\[
 (\operatorname{ev}_Q)_*[\mathfrak D_Q]^{\vir}=-s,
 \qquad
 (\operatorname{ev}_B)_*[\mathfrak D_B^+]^{\vir}=s.
\]
By Proposition~\ref{prop:sfintersectionsfull} and compatibility with proper pushforward, the positive-boundary term has degree $\int_J(-s)s=4$.

The second term has degree
\[
 \int_{\PP^4}c_4(\OO(-1)^{\oplus4})
 =\int_{\PP^4}H^4=1.
\]
This proves \eqref{eq:vacuumdegreeverified}.  Independently, the smooth generic moduli $\PP^4$ has the symmetric obstruction bundle on the moduli of sheaves on the quintic $\Omega^1_{\PP^4}$, whose fourth Chern number is $5$, confirming the normalization.  Tensoring is an equivalence of the actual derived sheaf problems and transports their unit, trace and cotangent morphisms, so the assertion for twists follows as well.
\end{proof}

\begin{corollary}\label{prop:quinticleadingfull}
The full leading primitive structure-sheaf coefficient, not merely the contribution from smooth hyperplane supports, is $5$.  Its positive and zero marked-boundary contributions are respectively $4$ and $1$ in the reduced ordered expansion convention.
\end{corollary}
\begin{proof}
The complete generic-moduli identification is Lemma~\companionref{lem:vacuumgenericfull}.  Theorem~\ref{thm:vacuumvirtual} computes the virtual class using all its expanded limits.  Singular supports were never removed in either construction.
\end{proof}

\subsection{Stable reduction at a smooth outer fiber in arbitrary multiplicity}
\label{sec:outer-fiber-multiple}
The outer relative divisor in the pencil is itself a fiber. This special
geometric fact gives a stable reduction theorem for higher-rank sheaves
on the fibers, including limits whose entire supporting surface is the
chosen outer K3. Write $C=\mathbb P^1$, choose a smooth fiber
$D_0=Y_{b_0}$, and fix a numerical character $v$ with
$\ch_1=rF_{\mathrm{fib}}$, $r\ge1$. Write $c=\chi(v)$ and polarize
the fibers by the primitive quartic class $h$.

Fix the ordinary threefold
moduli with polarization $A+F_{\mathrm{fib}}$.

An expansion of $(C,b_0)$ inserts a chain of projective lines at
$b_0$. Its pullback along $Y\to C$ is the expansion of $(Y,D_0)$:
each inserted component is $D_0\times\mathbb P^1$, with the two
attachable sections at $0,\infty$. This follows by flat base change of
the blowup construction for expansions, or directly from the charts
$t=az$, with $t$ a parameter on $C$ at $b_0$.

On such an expansion
consider flat families of pure sheaves of pushed character $v$,
admissible at the joints and the outer divisor, with finite target
stabilizer. Here admissibility includes the usual ordinary matching
of the branch sheaves, with the section equation a nonzerodivisor
on each incident branch sheaf. For stability on a fixed geometric
expansion use a polarization inducing $h$ on each K3 slice.

The pullback of $A$ is relatively ample over the expanded base curve. Tensoring it with a sufficiently ample bundle pulled back from that curve gives the required polarization. We treat the ordinary outer relative problem; the stable groupoid below is independent of this choice.

\begin{theorem}\label{thm:outerbitrarymultiplicity}
The groupoid of these stable sheaves on expansions of $(Y,D_0)$ is
equivalent to the groupoid of ordinary stable vertical sheaves on $Y$
of the specified character. The equivalence commutes with base change
and with scalar rigidification. At most one bubble occurs, and all
intrinsic target stabilizers are trivial.

If every semistable sheaf of
this character on $Y$ is stable, the rigidified expanded moduli is
projective, separated and proper. A sufficient numerical condition is
$\gcd(r,c)=1$.
\end{theorem}
\begin{proof}
Let $c:Y'\to Y$ be contraction of a geometric expansion, $F_b=Y_b$,
and $H_Y=A+F_{\rm fib}$. Put $[G]=\sum_V m_V[V]$, $m_V>0$,
for the fundamental cycle of a pure surface sheaf of the prescribed
character on $Y'$. Then
\[
 0=H_YF_{\rm fib}(rF_{\rm fib})
   =\sum_V m_V H_YF_{\rm fib}\,c_*[V].
\]
Every summand is nonnegative. If a root surface $V$ dominates $C$,
its summand before multiplication by $m_V$ is
$\deg(H_Y|_{V\cap F_b})>0$ for a general $b$. Thus every root
surface is a fiber, since all fibers are integral. This calculation
also applies to the unexpanded target.

\emph{Propagation across the chain.} Number the components from
the root outwards. Let $D_j$ be the join of components $j$ and
$j+1$ and $G_j$ the corresponding branch sheaf. On the root,
near the first join, some power of its equation $t$ annihilates
$G_0$, since its support meets only finitely many fibers. In a
trivialization, admissibility and matching give
\[
  t^N G_0=0,\quad \ker(t|_{G_0})=0
  \ \Longrightarrow\ G_0=0\text{ near }D_0
  \ \Longrightarrow\ G_1|_{D_0}=0 .
\]
On a bubble $D_0\times\PP^1$ with zero incoming restriction, each
irreducible support surface $V$ satisfies
\[
  V\cap(D_0\times\{0\})=\varnothing
  \quad\Longrightarrow\quad
  \mathrm{pr}_2(V)=\{z\}
  \quad\Longrightarrow\quad V=D_0\times\{z\}.
\]
The image of $V$ is a proper closed irreducible subset of the projective line because it misses $0$. The first disjointness follows from Nakayama. The equation
of the outgoing section acts injectively by admissibility; if
$z=\infty$ it would act nilpotently at the generic point of $V$.
Thus $z\ne\infty$ and the outgoing restriction vanishes. Induction,
also at the outer section, places every occupied slice in the
interior.

With no bubble, the same nilpotence argument excludes
$F_{b_0}$ from the unexpanded relative support.

For a stable geometric $G$, $\End(G)=k$ over its algebraically closed
ground field. Disjoint slices give orthogonal idempotents, so exactly
one slice occurs. A uniformizer of its base point acts nilpotently on
$G$: multiplication near the support extends by zero on its complement.
It therefore acts by zero. For the inclusion $j$ of the slice we obtain
\[
 G=j_*V,\quad \rk V=r,\quad
 \{G'\subset G\}=\{j_*V':V'\subset V\},\quad
 P_{j_*V'}=P^h_{V'}.
\]
Purity of $G$ is equivalent to torsion-freeness of $V$ on the
integral slice, also when the slice is nodal. Every empty bubble
has a dilation subgroup $\GG$ fixing $G$. Finite target stabilizer
leaves at most one bubble, and an occupied bubble has trivial
stabilizer because $z\in\GG$ and $\lambda z=z$ implies $\lambda=1$.

\emph{Recovering the section over a base.} Write
$q':Y'_T\to C'_T$ for the fibration over the expanded base curve and
$p:Y'_T\to T$ for its structural morphism. Let $\mathcal G$ be
a flat family as in the statement, with $T$ locally noetherian.
Its support misses the joins by the geometric calculation.
The flat morphism $p$ is smooth near this support by its geometric
fibers: a point on a nodal pencil member is still a point on the
smooth threefold $Y$.

Lifting finite fiberwise resolutions using
flatness makes $\mathcal G$ perfect there. It is zero off the support.
Perfect proper pushforward and base change give
\[
 K=\RR p_*\RR\mathcal Hom(\mathcal G,\mathcal G),\qquad
 K\otimes^{\LL}k(t)=\RHom_{Y'_t}(\mathcal G_t,\mathcal G_t).
\]
For $L=\operatorname{Cone}(\OO_T\xrightarrow{1}K)$, fiberwise simplicity gives $H^{\le0}(L\otimes^{\LL}k(t))=0$.
A minimal finite free representative at $t$ thus begins in degree
one, and so does its pullback. Hence for every base change
$f:T'\to T$,
\[
 H^0(\LL f^*K)=\OO_{T'},\qquad H^{<0}(\LL f^*K)=0 .
\]
These identities include all nilpotents in $\OO_{T'}$.

Locally let $T=\Spec R$. The support image in the geometric expanded
curve over a point of $T$ consists of one point. Take an affine
neighborhood $U=\Spec B\subset C'_T$ of it in the interior.
The closed complement of $(q')^{-1}U$ in $\Supp\mathcal G$ is
proper over $T$; after shrinking $T$ it is empty.
The scalar action is the retraction
\[
 R\longrightarrow B
    \xrightarrow{\ m\ }\End(\mathcal G)
    =\Gamma(T,H^0K)=R,\qquad m(f)(x)=fx .
\]
It defines $b':T\to U$, and $(q')^*\ker(m)\cdot\mathcal G=0$.
For the closed slice $j_{b'}:Y'_{b'}\hookrightarrow Y'_T$ this yields
a unique $T$-flat sheaf $\mathcal V$ with
$\mathcal G=(j_{b'})_*\mathcal V$. Conversely its annihilator on $U$
is precisely $\ker(m)$, by faithfulness of the scalar action.
This gives uniqueness, descent of $b'$ on overlaps, and arbitrary
base change by the identity for $H^0(\LL f^*K)$. The identical
construction recovers $b:T\to C$ from an ordinary stable family.

\emph{Inverse functors and arrows.} Near $b_0$, the contracted
coordinate $t$ and the interior support mark $z$ give the chart
\[
 [\AA^1_a\times\GG/\GG]\ \simeq\ \AA^1_t,
 \qquad t=az,\qquad \lambda(a,z)=(\lambda a,\lambda^{-1}z).
\]
Indeed $(a,z)\mapsto(az,z)$ trivializes the torsor and $t\mapsto(t,1)$
is its section. Pull this chart back in an \'etale coordinate on $C$;
away from $b_0$ it is the identity. It identifies a section $b$
with an expanded curve marked by $b'$, uniquely up to unique
isomorphism over any $T$. The contraction induces the commutative diagram
\[
\begin{tikzcd}[column sep=large]
 Y_{b'}'\arrow[r,hook,"j_{b'}"]\arrow[d,"\wr"']&
 Y\times_C C'_T\arrow[d,"c"]\arrow[r,"q'"]&
 C'_T\arrow[d,"c_C"]\\
 Y_b\arrow[r,hook,"j_b"]&Y\times T\arrow[r,"q_T"]& C\times T .
\end{tikzcd}
\]
The right square is Cartesian, and restriction to the marked slice gives
\[
 Y'_{b'}=(Y\times_C C'_T)\times_{C'_T,b'}T=Y\times_{C,b}T=Y_b .
\]
The inverse functors send $(j_b)_*\mathcal V\longleftrightarrow(j_{b'})_*\mathcal V$ on these identified slices. The whole inverse image of the contracted slice can include a bubble. On the marked slice, pushforward along either closed immersion is fully faithful. The unique isomorphism of marked targets thus identifies arrows, including scalars, and the two constructions commute with base change.

Finally let $G$ be semistable and $G'\subsetneq G$ have the same
reduced Hilbert polynomial. The preceding effective-cycle argument
applied on the ordinary target gives $[m^2]P_G=2r$ and
$[m^2]P_{G'}=2r'$, $0<r'\le r$. For $r'=r$ one has
$P_G-P_{G'}=P_{G/G'}\ne0$, contradicting equality of reduced
polynomials. Otherwise
\[
 r\chi(G')=r'c,\quad \gcd(r,c)=1\quad\Longrightarrow\quad
 r\mid r',\qquad 0<r'<r.
\]
This proves the numerical criterion. If semistability equals stability,
the projective Simpson moduli of \cite[Theorem~1.21]{Simpson}
is the scalar-rigidified stable moduli; the groupoid equivalence
transfers projectivity and both valuative properties.
\end{proof}
The equivalence concerns the classical moduli at the outer divisor. Virtual comparison also requires an identification of obstruction morphisms, as in the rank-three pencil calculation below. At the internal cut $Q\cup_S W$ a root surface can meet $S$ in a curve, whereas the induction above starts with zero boundary restriction.

\subsection{The outer-relative structure-sheaf comparison}\label{sec:K3verifiedcase}
Since $K_Y=-F_{\mathrm{fib}}$, we compare the global obstruction morphisms of the outer-relative and absolute pencil theories for the structure-sheaf charge.

\begin{theorem}\label{thm:K3structurecomparison}
Let $\kappa:Y\to\PP^1$ be the quartic pencil of Section~\ref{sec:seconddegfull}, and choose a smooth fiber $D_0$.  In the character of $\OO_{Y_b}$, the ordinary stable-sheaf moduli and the scalar-rigidified relative structure-sheaf moduli of $(Y,D_0)$ with zero outer boundary are both $\PP^1$.  The latter inserts one bubble at $b=\kappa(D_0)$ and is proper and separated.

Its relative determinant obstruction theory agrees with the obstruction theory of vertical sheaves on the K3-fibered threefold in \cite[Theorem 2.3]{GS}.  Their virtual class and degree are
\begin{equation}\label{eq:K3verifieddegree}
 [\PP^1]^{\vir}=c_1(\OO_{\PP^1}(-1))\cap[\PP^1],
 \qquad \deg[\PP^1]^{\vir}=-1.
\end{equation}
\end{theorem}
\begin{proof}
Write the sheaf as $i_*V$ on an integral K3 fiber and set $R=V^{**}$. Here $\deg R=0$ and $\chi(V)=2$. If $R\ne\OO$, the divisorial section argument of Lemma~\companionref{lem:vacuumgenericfull} gives $H^0(R)=H^0(R^*)=0$. Since each fiber is normal Gorenstein with trivial dualizing sheaf, $H^2(R)=0$ and
\[
       2=\chi(V)=\chi(R)-\length(R/V)\le0,
\]
a contradiction. For $R=\OO$ the same identity gives $R/V=0$. This includes the nodal members. The cyclic Cartier deformation argument and $H^1(\OO_{Y_b})=0$ identify the ordinary moduli with the pencil $\PP^1$.

The normal bundle of $D_0$ in $Y$ is trivial.  In a local coordinate $b$ on the pencil centered at $D_0$, the relative expansion uses a parameter $a$ and a finite bubble coordinate $z$, with the support graph
\[
 b=az,\qquad z\in\GG.
\]
The rescaling $(a,z)\mapsto(\lambda a,\lambda^{-1}z)$ is free. Its quotient has coordinate $az=b$, including at $a=0$, where the support is an interior copy of $D_0$. It misses both divisors and admits no empty bubble. On the support $(b,a)=(az,a)=(a)$, so the blowup leaves the abstract support family unchanged.

The identity chart away from $D_0$ glues this quotient to $\PP^1$, so the relative functor is proper and separated. Its normal tangent sequence is
\begin{equation}\label{eq:K3normalchart}
 0\longrightarrow\OO
 \xrightarrow{(a,-z)}\OO^{\oplus2}
 \xrightarrow{(z,a)}\OO\longrightarrow0.
\end{equation}
The arrows are the infinitesimal rescaling and $d(az)$. For $(u,v)$ in the kernel, $u=-av/z$, hence $(u,v)=(a,-z)(-v/z)$. Since $z$ is invertible the sequence is exact, including at $a=0$.

On the graph of the pencil the Cartier resolution gives
\[
 \mathcal Ext_p^1(\OO_{\mathcal Y},\OO_{\mathcal Y})=T_{\PP^1},\qquad
 \mathcal Ext_p^2(\OO_{\mathcal Y},\OO_{\mathcal Y})
   =R^2\kappa_*\OO_Y=(\kappa_*\omega_{Y/\PP^1})^\vee=\OO(-1).
\]
The supporting family on the expanded chart is this same graph. Equation~\eqref{eq:K3normalchart} identifies its base direction after rescaling, and Part~I, Theorem~\companionref{thm:zeroboundaryPOT}, removes the scalar and its top dual. The resulting two middle bundles agree with \cite[Theorem~2.3]{GS}.

To identify the obstruction \emph{morphisms}, not merely these two bundles, let $E\to L_{\PP^1}$ be either of the two perfect obstruction theories.  Smoothness gives
\[
 h^0(E)=\Omega^1_{\PP^1}=\OO(-2),\qquad
 h^{-1}(E)=\OO(1).
\]
Its truncation triangle is
\[
 \OO(1)[1]\longrightarrow E\longrightarrow\OO(-2)
       \xrightarrow{\epsilon}\OO(1)[2].
\]
The extension class lies in
\[
 \Ext^2(\OO(-2),\OO(1))=H^2(\PP^1,\OO(3))=0.
\]
The ambiguity of the splitting lies in
$\Ext^1(\OO(-2),\OO(1))=H^1(\PP^1,\OO(3))=0$.
Thus the identified cohomology bundles determine the splitting uniquely.
Moreover $\Hom(\OO(1)[1],\Omega^1_{\PP^1})=0$, so the morphism to the cotangent complex is the projection to $\OO(-2)$, after its required isomorphism on $h^0$.  Both theories therefore give the same complex and the same obstruction morphism.

The virtual class on the smooth curve is the top Chern class of its obstruction line bundle, yielding \eqref{eq:K3verifieddegree}.
\end{proof}

\begin{corollary}\label{prop:quarticleadingfull}
For the pencil of Theorem~\ref{thm:K3structurecomparison}, the raw neutral series begins $K^{\rm raw}_0=-q^{-1}+O(1)$.
\end{corollary}
\begin{proof}
At $d=0$, the Fourier power is $1-\chi$. The degree-zero hull argument of the preceding proof excludes $\chi>2$ and identifies the entire $\chi=2$ term with the structure-sheaf pencil. Equation~\eqref{eq:K3verifieddegree} gives its degree.
\end{proof}

\subsection{A complete rank-three fiber class and its internal limit}
\label{sec:nonabelianfiber}
Let $h$ be the primitive quartic polarization on the fibers of
$\kappa:Y\to\PP^1$. Write $A$ for the pullback of the quartic
threefold's hyperplane class. The fibers of a general Lefschetz
pencil are smooth or have one ordinary node.

We compute the full
outer virtual class for one character of generic rank three, then
describe the sheaf on both branches of the internal cut.
\begin{theorem}\label{thm:rankthreepencilspherical}
On a fiber $T_b\subset\Pi_b\simeq\PP^3$, let
$U_b=T_{\Pi_b}(-1)|_{T_b}$. Fix the character
$v=\ch((i_b)_*U_b)$ on $Y$.

Its divisorial component is
$3F_{\mathrm{fib}}$ and its fiber Hilbert polynomial for $h$ is
$6m^2+4m+4$. All semistable sheaves in this character are stable.
Their ordinary moduli and their scalar-rigidified outer relative moduli
from Theorem~\ref{thm:outerbitrarymultiplicity} are both $\PP^1$.

The obstruction theory of sheaves on $Y$ and the obstruction theory relative to the outer K3 fiber have
the same virtual class
\[
       [\PP^1]^{\vir}
       =c_1(\OO_{\PP^1}(-1))\cap[\PP^1]
       \quad\hbox{in }A_0(\PP^1),\qquad
       \deg[\PP^1]^{\vir}=-1 .
\]
\end{theorem}
\begin{proof}
The ambient hyperplanes supply a rank-four bundle
$\mathcal V=\kappa_*\OO_Y(A)$ on the base. Evaluation and its dual
give a bundle $U$ on $Y$ by the exact sequence
\[
 0\longrightarrow\OO_Y(-A)\longrightarrow
   \kappa^*\mathcal V^\vee\longrightarrow U\longrightarrow0 .
\]
On $T=T_b$ this is the restricted Euler sequence. It gives
$\operatorname{rk}U_b=3$, $c_1(U_b)=h$, $\chi(U_b)=4$ and
$H^0(T,U_b^\vee)=0$. These hold on a nodal member as well, by the
hypersurface sequence for the quartic in its hyperplane.

We check slope stability of $U_b$. A proper torsion-free quotient
$G$ of rank $s=1,2$ and slope at most $4/3$ would have
$\deg_h\det G\le\lfloor4s/3\rfloor\le2$.
The reflexive determinant $\det G=(\bigwedge^sG)^{**}$ is generated by the maximal minors in codimension one. Positive degree would therefore give a moving linear system of degree one or two. Every irreducible line or conic has a smooth rational strict transform $C$ on the minimal resolution; blowing up a point of a smooth curve leaves it unchanged. Adjunction gives $C^2=-2$, so two distinct irreducible curves cannot represent its class. The K3 Picard group is countable, hence there are only countably many such cycles, also allowing pairs of lines and double lines. They cannot contain a moving complex linear system.

In degree zero a nonzero determinant section gives $\det G=\OO_T$. Wedging with $s-1$ independent generating sections yields $G\to\OO_T$ nonzero, contradicting $H^0(T,U_b^\vee)=0$.

A negative determinant degree is excluded by
global generation. This proves slope stability.

The support argument of Theorem~\ref{thm:outerbitrarymultiplicity} applies to each stable Jordan--H\"older factor $G_i$: it lies on one fiber and is killed by its defining parameter. If its fiber rank is $s_i$, its leading coefficient is $2s_i$ and equality of reduced Hilbert polynomials gives
\[
 \frac{[m]P_{G_i}}{s_i}=\frac43,\qquad
    3\mid s_i,\qquad \sum_i s_i=3 .
\]
There is exactly one factor. Hence every semistable sheaf of this character is stable.

Let $G$ be a stable sheaf of the stated polynomial on $T$. It is
torsion-free of rank three. The Euler sequence gives, on the smooth
or the nodal Gorenstein surface,
\[
 \chi_T(U_b,G)=4\chi(G)-\chi(G(h))=16-14=2,\qquad
 \Ext^2_T(U_b,G)^\vee=\operatorname{Hom}_T(G,U_b).
\]
Only degrees zero through two occur because $U_b$ is locally free.
If $G\not\simeq U_b$, stability gives vanishing of both Hom groups.
It would follow that $2=-\dim\Ext^1_T(U_b,G)$. Setting $G=U_b$
now gives $h^0(\operatorname{End}U_b)=h^2(\operatorname{End}U_b)=1$
and $h^1(\operatorname{End}U_b)=0$. The trace splits the scalar unit
after division by three, so
\[
 R\kappa_*\operatorname{End}^{0}U=0,\qquad
 R\kappa_*\OO_Y\xrightarrow{\ \cdot1_U\ }
 R\kappa_*\operatorname{End}U\quad\text{is an isomorphism}.
\]
The first assertion follows by derived base change on the proper flat
pencil and Nakayama. It is preserved by base change.

The graph family defines a proper bijective morphism from the pencil
to its ordinary projective stable moduli. The Cartier resolution of
the graph has two terms. In the spectral sequence for its self-Ext groups on $Y$
the two rows are the cohomology of $\operatorname{End}U_b$ and of
$\operatorname{End}U_b\otimes N_{T/Y}$, with the latter row shifted
by one. The normal line bundle is the constant fiber of $T_{\PP^1}$.

Multiplication by $1_U$ identifies both rows, hence the entire
self-Ext complex after pushforward, with the corresponding rows for
$\OO_T$. The global graph and its motion give the resulting tangent
line bundle $T_{\PP^1}$. At each point the local moduli has dimension one
and tangent dimension one, hence is regular. The proper quasi-finite
map is finite and birational in characteristic zero and the target
is normal. It is an isomorphism.

Theorem~\ref{thm:outerbitrarymultiplicity}
gives the relative moduli; its numerical coprimality condition
also holds since $\gcd(3,4)=1$.

On the outer chart $b=az$ the slice and $U$ pull back under contraction. The scalar inclusion induces a quasi-isomorphism on the two Cartier self-Ext rows, since $R\kappa_*\operatorname{End}^0 U=0$. In the smoothing direction its differential is $(z,a)$, as in \eqref{eq:K3normalchart}. On scalars, $\operatorname{tr}(c\,1_U)=3c$. Thus the scalar and dual top reductions of Part~I, Theorem~\companionref{thm:zeroboundaryPOT}, agree up to multiplication by the unit $3$.

Both theories therefore have tangent $T_{\PP^1}$ and obstruction $R^2\kappa_*\OO_Y=\OO(-1)$. Since $\Ext^2(\Omega^1_{\PP^1},\OO(1))=0$ and $\Hom(\OO(1)[1],\Omega^1_{\PP^1})=0$, their cotangent presentations are
\[
  \OO(1)[1]\oplus\Omega^1_{\PP^1}\xrightarrow{\operatorname{pr}_2}\Omega^1_{\PP^1}.
\]

The obstruction-bundle formula gives the displayed virtual class.
\end{proof}

\begin{proposition}\label{prop:sphericalinternalcut}
Fix a smooth member $T=T_b$ of the pencil different from the distinguished
member $S$. In the internal degeneration of $Y$, the sheaf $(i_b)_*U_b$
has a flat perfect extension. On the central support $T\cup_C R_*$
its restrictions are $U_b$ and $p^*(\OO_C\oplus D)$, where $C=C_*$
is the plane quartic base curve, $p:R_*\to C$ is the ruling and
$D=T_{\PP^2}(-1)|_C$. The glued sheaf is stable for suitable
compatible ample polarizations on $Q\cup_{S_\infty}W$.

Its $W$
factor is decomposable, and its boundary sheaf on $S_\infty$
is the pushforward of $\OO_C\oplus D$.
\end{proposition}
\begin{proof}
In the chart of the quartic threefold write $x,u$ for the two
generators of the pencil and $x=bu$ on $T$, with $b\ne0$ finite
after choosing pencil coordinates. Let $s$ be the degeneration
parameter. The first blowup, along $(x,s)$, restricts on the strict
transform of $T\times\mathbb A^1$ to the blowup along $(u,s)$.
In its chart $x=sz$ the transform is $sz=bu$.

The second center has
ideal $(u,z)$, whose restriction is $(z)$; its strict transform
does not meet the complementary part of the exceptional divisor in
the other chart. Thus the final transform is
$\mathcal T=\operatorname{Bl}_{C\times0}(T\times\mathbb A^1)$.
Its central components are $T$ and $R_*$ joined along $C$.

Pull $U_b$ back by contraction $\mathcal T\to T$ and push it forward
from this relative Cartier divisor. Its local resolution has length
one, and it is flat over the parameter because $\mathcal T$ is flat
and the pullback is locally free. Restriction of the Euler
presentation to the plane of $C$ gives
$U_b|_C=\OO_C\oplus T_{\PP^2}(-1)|_C$ upon choosing a complement
to the three plane coordinates. These are precisely the branch and
boundary restrictions stated above.

Write $s_-=R_*\cap S_\infty$ for the section of square $-4$;
the other distinguished section is $s_+=s_-+p^*h|_C$. The normal
of the first section is $\OO_C(-h)$, and
$K_{R_*}=-2s_-$ by the canonical divisor formula for a ruled
surface and $K_C=h|_C$. For rational $0<\epsilon\ll1$ take
\[
 A_T=(1-\epsilon)h,\qquad A_R=p^*(h|_C)+\epsilon s_- .
\]
These are induced by an ample rational class on the two ambient
components. On $Q$ use $(1-\epsilon)H_Q$. On $W$, with blowdown
$\pi:W\to P$, projection $p_P:P\to S$ and exceptional divisor
$E_P$, use
\[
 \pi^*((1+\delta)p_P^*h+(\epsilon+\delta)S_\infty)-\delta E_P,
 \qquad 0<\delta\ll\epsilon .
\]
At $\delta=0$ the class before pullback is ample for small
$\epsilon>0$ (the section is relatively ample). For small $\delta$
the displayed rational divisor class is ample on the blowup, since its perturbation
in the $-E_P$ direction is relatively ample. It restricts to
$(1-\epsilon)h$ on $S_\infty$. Since $E_P|_{R_*}=s_+$ it restricts
to $A_R$ as claimed.

Denote the glued bundle on its supporting surface by $B$. Write
$P_J(m)=a_\epsilon(J)m^2+b_\epsilon(J)m+O(1)$ for a subsheaf
$J\subset B$, interpreting the coefficients for a rational
polarization by Riemann--Roch. Replacing it by an integral multiple
preserves their slope inequalities. The normalization sequence for
$B$ gives
\begin{align*}
 a_\epsilon(B)&=6(1-\epsilon)^2+3(4\epsilon-2\epsilon^2)=6,\\
 b_\epsilon(B)&=4(1-\epsilon)+(12-8\epsilon)-12(1-\epsilon)=4 .
\end{align*}
We prove $b_\epsilon(J)/a_\epsilon(J)<2/3$ for every subsheaf
whose rank is smaller than three on at least one occupied component.
Let $q,s'$ be its ranks on $T,R_*$, respectively. For the contraction
$c:T\cup R_*\to T$, the normalization sequence and
$Rp_*\OO_{R_*}=\OO_C$ give $Rc_*B=U_b$. In particular
$c_*J\subset U_b$. At $\epsilon=0$, projection formula with $h$
on $T$ gives
\[
 a_0(J)=2q,\qquad
 b_0(J)=h\cdot c_1(c_*J)-4\,\operatorname{rk}_C R^1c_*J .
\]
If $0<q<3$, slope stability of $U_b$ gives
$b_0(J)<4q/3$. If $q=3$ and $s'<3$, the values of $c_*J$ along
the generic point of $C$ lie in a subspace of $U_b|_C$ of dimension
at most $s'$. Indeed on the generic ruling the image on $R_*$
is a rank-$s'$ subsheaf of the trivial rank-three bundle; it is
locally free there and its evaluation at the attaching section has
rank at most $s'$. Matching the two evaluations proves the assertion.

It follows over the DVR of the generic point of $C$ in $T$ that
$U_b/c_*J$ has positive length there. Thus
$h\cdot c_1(c_*J)\le4-h\cdot C=0$.

It remains to bound $J$ with $q=0$. It lies in
$p^*(U_b|_C)(-s_-)$ on $R_*$, the kernel of restriction to the
other branch. Put $J'=J(s_-)$ and let $t$ be its degree on the
generic ruling. Since $J'$ embeds in a trivial bundle there,
$t\le0$. Riemann--Roch on $R_*$ gives $b_0(J)=4t$. If $t=0$
the generic subbundle of rank $s'$ is a constant subspace:
on $\PP^1_{\mathbb C(C)}$ each splitting summand has nonpositive
degree, and their degrees sum to zero.

Let $V\subset U_b|_C$ be
the saturated subbundle extending that generic subspace. Then
$J'\subset p^*V$, since the composite to the locally free
$p^*((U_b|_C)/V)$ vanishes generically. The determinant of $J'$
is $p^*M$ for a line bundle $M$ on $C$ because its fiber degree
is zero. Taking determinants gives $\deg M\le\deg V$.

We have $\deg V\le 2s'$. In fact $D$ is slope-stable of slope
two. A destabilizing subline would give a globally generated
quotient line bundle of degree at most two, since $D$ is globally generated.
A smooth plane quartic is a nonhyperelliptic curve of genus three.
It has no moving linear system of degree one or two, and the
degree-zero possibility $\OO_C$ is excluded by the Euler sequence
$H^0(D^\vee)=0$.

The splitting $U_b|_C=\OO_C\oplus D$ thus
gives the bound on $V$. Using $K_{R_*}=-2s_-$ yields
\[
 b_\epsilon(J)=\epsilon\deg M,\qquad
 a_\epsilon(J)=s'(4\epsilon-2\epsilon^2),\qquad
 \frac{b_\epsilon(J)}{a_\epsilon(J)}
       \le\frac{2}{4-2\epsilon}<\frac23
       \quad(0<\epsilon<1/2).
\]

For the other rank profiles put $u_\epsilon(J)=b_\epsilon(J)-2a_\epsilon(J)/3$. Integrality gives $u_0(J)\le-1/3$. Fix an ample rational parameter $\epsilon_0>0$. The maximal Harder--Narasimhan slope at $\epsilon_0$ bounds $b_{\epsilon_0}(J)/a_{\epsilon_0}(J)$. Since $a_\epsilon(J)$ depends only on $(q,s')$, a finite set, there is a uniform $B$ with $u_{\epsilon_0}(J)\le B$. Write $\theta=\epsilon/\epsilon_0$. Affine linearity of $b_\epsilon$ gives
\[
 u_\epsilon(J)\le -(1-\theta)/3+\theta B+
 \tfrac23\big((1-\theta)a_0(J)+\theta a_{\epsilon_0}(J)-a_\epsilon(J)\big).
\]
The parenthesis tends uniformly to zero by finiteness of the rank profiles. Hence $u_\epsilon(J)<0$ for all sufficiently small $\epsilon>0$. For $q=s'=3$ a proper subsheaf has a nonzero quotient of lower dimension, so its reduced Hilbert polynomial is strictly smaller. The glued bundle is pure since its supporting surface is reduced and Cohen--Macaulay.
For the two nonzero summands $B_1,B_2$ of $B_R=p^*(\OO_C\oplus D)$, write $a_i$ for the leading coefficients of $P_{B_i}$. Then
\[
 p_{B_R}=\frac{a_1p_{B_1}+a_2p_{B_2}}{a_1+a_2},\qquad
 \max_i p_{B_i}\ge p_{B_R}.
\]
Thus the stable glued sheaf has a ruled factor which fails stability for every ample polarization. In the normalization triangle for $\operatorname{End}B$, the restriction from $R_*$ to $C$ is a quasi-isomorphism, since $Rp_*\OO_{R_*}=\OO_C$. Its homotopy fiber gives
\[
 R\Gamma(T\cup_C R_*,\operatorname{End}B)\simeq R\Gamma(T,\operatorname{End}U_b).
\]
Gluing along this boundary must therefore retain the matching of automorphisms of the two summands.
\end{proof}

\begin{corollary}\label{cor:sphericalinternalext}
Put $X'_0=Q\cup_{S_\infty} W$ and let $i:Z=T\cup_C R_*\hookrightarrow X'_0$ be the Cartier support of the internal limit of Proposition~\ref{prop:sphericalinternalcut}. Then
\[
 \dim\Ext^j_{X_0'}(i_*B,i_*B)=1\qquad(0\le j\le3).
\]
The moving support gives the generator in degree one. On its ruled
branch, with $j:R_*\hookrightarrow W$ and $V=U_b|_C$, the natural
map
\[
 R\operatorname{Hom}_C(V,V)\longrightarrow
 R\operatorname{Hom}_W(j_*p^*V,j_*p^*V)
\]
is a quasi-isomorphism.
\end{corollary}
\begin{proof}
The branch normal line bundle is $N_{R_*/W}=\OO_{R_*}(-s_-)$ by the
divisor calculation in Section~\ref{sec:seconddegfull}. Its
restriction to a ruling is $\OO_{\PP^1}(-1)$, and therefore
$Rp_*N_{R_*/W}=0$. The degree-zero and degree-one sheaf Ext of the
Cartier pushforward on $W$ are $\operatorname{End}(p^*V)$ and
$\operatorname{End}(p^*V)\otimes N_{R_*/W}$, supported on $R_*$;
the higher sheaf Ext vanish. Projection formula on these two rows
gives the asserted quasi-isomorphism, including the scalar maps.

Write $N=N_{Z/X_0'}$. This line bundle restricts to $\OO_T(h)$ on $T$
and $\OO_{R_*}(-s_-)$ on $R_*$. On the cut both restrictions
equal $\OO_C(h)$. Tensor the normalization sequence with
$\operatorname{End}B\otimes N$. The direct summand of its middle
term from $R_*$ has zero derived global sections, so its homotopy
fiber after applying $R\Gamma$ is
\[
 \begin{aligned}
 R\Gamma(Z,\operatorname{End}B\otimes N)
 &\simeq\operatorname{fib}\big(
 R\Gamma(T,\operatorname{End}U_b(h))\longrightarrow
 R\Gamma(C,\operatorname{End}V(h))\big)\\
 &\simeq R\Gamma(T,\operatorname{End}U_b(h-C))
 \simeq R\Gamma(T,\operatorname{End}U_b).
 \end{aligned}
\]
The local-to-global spectral sequence for the Cartier pushforward $i_*B$ has two rows,
\[
 E_2^{p,q}=
 \begin{cases}
 H^p(Z,\operatorname{End}B),&q=0,\\
 H^p(Z,\operatorname{End}B\otimes N),&q=1,\\
 0,&q\ge2.
 \end{cases}
 \qquad E_2^{p,q}\Longrightarrow \Ext^{p+q}_{X'_0}(i_*B,i_*B).
\]
Since $Rp_*\OO_{R_*}=\OO_C$, restriction to the attaching section is a quasi-isomorphism on ruled-component endomorphisms. Normalization and Proposition~\ref{prop:sphericalinternalcut} give $E_2^{p,0}=H^p(T,\operatorname{End}U_b)$. The normal-line calculation gives $E_2^{p,1}=E_2^{p,0}$. Sphericity makes each row $(\CC,0,\CC)$.

The only possible differential is $d_2^{0,1}:\CC\to\CC$. Moving $b$ in the family of Proposition~\ref{prop:sphericalinternalcut} gives a first-order deformation of $i_*B$. Its normal component on $T$ is the derivative of $x-bu$, namely $-u\,\mathrm{id}_{U_b}$ in a local trivialization. It is nonzero and vanishes along $C$. Hence the edge map
\[
 \Ext^1_{X'_0}(i_*B,i_*B)\longrightarrow
 H^0(Z,\operatorname{End}B\otimes N)
\]
is surjective. Thus $d_2^{0,1}=0$, and each total degree $0,1,2,3$ contains one surviving copy of $\CC$.
\end{proof}

\subsection{The relative Hilbert-scheme factors}
For a smooth projective surface $T$, a smooth divisor $C\subset T$, and a line bundle $M$, the relative Carlsson--Okounkov class of Gholampour and Sheshmani \cite{GSHilb} is the perfect $K$-class on the relative Hilbert scheme of points obtained from
\[
 R\Gamma(T,M)-R\operatorname{Hom}_T(I_Z,I_Z\otimes M).
\]
Its relative construction is defined on the expanded surface and is compatible with the admissibility of the finite quotient.  The relative surface formula of \cite[Theorem~2]{GSHilb} states
\begin{equation}\label{eq:relativeCOfull}
\begin{aligned}
 \sum_{n\ge0}q^n\int_{(T/C)^{[n]}}c_{2n}(\mathbb E^{[n]}(M))
 &=\prod_{m\ge1}(1-q^m)^{-\delta(T,C;M)},\\
 \delta(T,C;M)&=\int_Tc_2(T_T(-\log C)\otimes M).
\end{aligned}
\end{equation}

\begin{proposition}\label{prop:punctualfactorsfull}
Let $T$ be a smooth Cartier support in a threefold satisfying $K_Y=-S$ near $T$, with $C=T\cap S$ smooth.  For $N=N_{T/Y}$ one has
\[
 N=K_T(C),\qquad \delta(T,C;N)=e(T)-e(C).
\]
For the plane transform in $B$ and the ruled support $R_*\subset W$, relative to their gluing curves, the factors in \eqref{eq:relativeCOfull} are respectively
\begin{equation}\label{eq:twofactorfull}
 \Xi(q)^{-27},\qquad\Xi(q)^4,
 \qquad\Xi(q)=\prod_{m\ge1}(1-q^m).
\end{equation}
Their ratio is $\Xi^{-31}$. For a smooth quartic support with hyperplane boundary the same calculation gives $\Xi^{-28}$. Consequently the indicated punctual exponents satisfy $\Xi^{-28}\Xi^{4}=\Xi^{-24}$ and $\Xi^{-28}\Xi^{-27}=\Xi^{-55}$. These exponents enter the comparison of the two relative evaluations.
\end{proposition}
\begin{proof}
Adjunction gives $K_T=(K_Y+T)|_T=-C+N$.  The logarithmic tangent sequence gives
\[
 c_1(T_T(-\log C))=-K_T-C=-N,
 \qquad
 \int_Tc_2(T_T(-\log C))=e(T)-e(C).
\]
For any rank-two bundle $E$ and line bundle of first Chern class $n$,
\[
 c_2(E\otimes M)=c_2(E)+c_1(E)n+n^2.
\]
Substituting $n=N$ cancels the last two terms, proving the exponent formula.  The same calculation applies on $W$ because the additional divisor $B_0$ is disjoint from $R_*$.

For the plane transform and the ruled surface the exponents are
\[
 \begin{aligned}
 e(\operatorname{Bl}_{20}\PP^2)-e(C)&=(3+20)-(-4)=27,\\
 e(R_*)-e(C_*)&=2e(C_*)-e(C_*)=-4.
 \end{aligned}
\]
This proves \eqref{eq:twofactorfull} and gives $\Xi^{-27}/\Xi^{4}=\Xi^{-31}$.
\end{proof}

Equation~\eqref{eq:relativeCOfull} computes the top Chern numbers on fixed smooth supports. Their contribution on the threefold also requires its virtual insertion, the sum over line-sheaf charges and the terms at singular supports, with their K3 boundary values.

In eta notation $\Xi^{-31}=q^{31/24}\eta^{-31}$; the prefactor is absorbed only after the specified charge normalization is checked.

\section{General divisor charge and completed modular transformations}\label{sec:generalmodularity}
In class $rH$, stable factors can have multiplicities $r_i>0$ with $\sum_i r_i=r$. Equality of reduced Hilbert polynomials restricts the Fourier exponents at which such filtrations occur. The Hall logarithm supplies their weights. For $r=1$ there is one stable factor, and Theorem~\ref{thm:twodegmodularityfull} gives a criterion for the proposed $S$-transformation.
\subsection{Charge coordinates for general divisor classes}
For $v=(0,rH,\beta,n)$ set
\begin{equation}\label{eq:generalEMcharges}
\begin{gathered}
 d=H\beta,\qquad q_1=-d,\qquad q_0=n+\frac{25r}{12},\\
 \widehat q_0=q_0-\frac{q_1^2}{10r},\qquad
 \mu=q_1-\frac{5r^2}{2}\pmod{5r}.
\end{gathered}
\end{equation}
The neutral residue is $q_1-5r^2/2\equiv0\pmod{5r}$. On the quartic pencil the residue is the fiber degree modulo four.

The coordinates arise by multiplying the Chern character by $\sqrt{\operatorname{td}(X)}=1+c_2(X)/24$:
\[
 \Gamma(F)=rH+d\frac{H^2}{5}
                    +\left(n+\frac{25r}{12}\right)[\mathrm{pt}].
\]
The sign in $q_1=-d$ is the same convention as in \cite{AFKPS,AFKPmock}.

\begin{proposition}\label{prop:generalspectralflow}
The class $\mu$ is an integral residue class.  Tensoring by $\OO_X(aH)$ preserves Gieseker (semi)stability and the ordinary or generalized DT invariant, and gives
\begin{equation}\label{eq:generalspectralflow}
 q_1'=q_1-5ra,\qquad
 q_0'=q_0-aq_1+\frac{5r}{2}a^2.
\end{equation}
It preserves $\widehat q_0$ and $\mu$.  Put
\begin{equation}\label{eq:generalcentralcharges}
 c_r=5r^3+50r,\qquad a_r=\frac{5r^3+25r}{6}=\chi(\OO_D).
\end{equation}
For a representative $q_1=\mu+5r^2/2$,
\begin{equation}\label{eq:generalTphase}
 -\widehat q_0
 =-\frac{c_r}{24}+\frac{\mu^2}{10r}+\frac{r\mu}{2}
             +a_r-\chi(F).
\end{equation}
Thus the formal $T$-phase is determined without a modularity assumption.
\end{proposition}
\begin{proof}
By \textup{(\companionref{eq:generalpolynomial})}, $d+5r/2\in\mathbb Z$.  Hence
\[
 \mu=-(d+5r/2)-\frac{5r(r-1)}2\in\mathbb Z.
\]
Expanding $e^{aH}\ch(F)$ gives
\[
 d'=d+5ra,\qquad n'=n+ad+\frac{5r}{2}a^2.
\]
These are \eqref{eq:generalspectralflow}.  Direct subtraction gives
\begin{align*}
 q_0'-\frac{(q_1')^2}{10r}
 &=q_0-aq_1+\frac{5r}{2}a^2
    -\frac{q_1^2-10raq_1+25r^2a^2}{10r}\\
 &=q_0-\frac{q_1^2}{10r}.
\end{align*}
Also $q'_1-q_1=-5ra$, so $\mu'=\mu$. Twisting translates the Hilbert polynomial, identifies the endomorphism complexes with units and traces, and preserves short exact sequences. It therefore preserves both the stable obstruction theories and the semistable Hall logarithms with their Behrend functions.

For a smooth member $D\in|rH|$, adjunction gives $K_D=D|_D$ and
\[
 c_2(T_D)=(c_2(X)+D^2)|_D,
 \qquad e(D)=D^3+D c_2(X)=5r^3+50r.
\]
The sequence $0\to\OO_X(-D)\to\OO_X\to\OO_D\to0$ and threefold Riemann--Roch give
\[
 \chi(\OO_D)=\frac{D^3}{6}+\frac{D c_2(X)}{12}
              =\frac{5r^3+25r}{6}.
\]
This Euler characteristic is the same for every Cartier member, including singular and nonreduced members, because it depends on its exact sequence.  Finally substitute $n=\chi(F)-25r/6$ and $q_1=\mu+5r^2/2$ into \eqref{eq:generalEMcharges}.  The result is
\[
 -\widehat q_0=-\chi(F)+\frac{25r}{12}
         +\frac{5r^3}{8}+\frac{\mu^2}{10r}+\frac{r\mu}{2},
\]
which is \eqref{eq:generalTphase}.  In particular $a_r$ is an integer.
\end{proof}

\begin{proposition}\label{prop:generalvacuumbound}
For every Gieseker-semistable pure two-dimensional sheaf of class $(0,rH,\beta,n)$ on a smooth quintic,
\begin{equation}\label{eq:generalvacuumbound}
 n\le\frac{5r^3}{24}+\frac{d^2}{10r},
 \qquad \widehat q_0\le\frac{5r^3+50r}{24}.
\end{equation}
The assertion includes all support types and does not use the existence of an invertible hull.
\end{proposition}
\begin{proof}
We use the rank-zero inequality of \cite[Theorem 1.1(i) and Lemma 3.7]{FeyzOSV}.  Its stated quintic case applies to slope-semistable surface sheaves.  In its notation
\[
 \frac12\left(\frac{DH^2}{H^3}\right)^2
 +6\left(\frac{\beta H}{DH^2}\right)^2
 -\frac{12n}{DH^2}\ge0.
\]
Gieseker semistability implies this slope semistability.  Substituting $D=rH$, $DH^2=5r$ and $d=\beta H$ gives
\[
 \frac{r^2}{2}+\frac{6d^2}{25r^2}-\frac{12n}{5r}\ge0.
\]
Multiply by $5r/12$ and rearrange to obtain the first inequality in \eqref{eq:generalvacuumbound}.  Adding $25r/12-d^2/(10r)$ gives the second.  
\end{proof}

\subsection{A polar range free of strictly semistable surface sheaves}
The first Chern class is $rH$ while the rank on the threefold is zero.
For a semistable surface sheaf put
\[
 c_r=5r^3+50r,\qquad
 \varepsilon(F)=c_r/24-\widehat q_0(F)\ge0 .
\]
The inequality is Proposition~\ref{prop:generalvacuumbound}. It gives a
stronger restriction on a wall where a pure surface sheaf is strictly
semistable.
\begin{proposition}\label{prop:step22-jhgap}
Let $F_i$ be the Jordan--H\"older factors of $F$ in its Gieseker
chamber, with divisorial multiplicities $r_i$. Then
\[
 \varepsilon(F)=\sum_i\varepsilon(F_i)
       +\frac5{24}\biggl(r^3-\sum_i r_i^3\biggr).
\]
In particular, every strictly semistable $F$ with $r\ge2$ satisfies
$\varepsilon(F)\ge 5r(r-1)/8$. In multiplicity two the bound
is attained, in a charge which also contains stable sheaves.
\end{proposition}
\begin{proof}
All $F_i$ are pure of dimension two, $r_i\ge1$, and
$\sum_i r_i=r$. Equality of the linear terms of the reduced
Hilbert polynomials gives a common $u=q_1(F_i)/r_i=q_1(F)/r$.
Additivity of $q_0$ gives
\[
 \widehat q_0(F)=\sum_i q_0(F_i)-ru^2/10
                 =\sum_i\widehat q_0(F_i).
\]
Subtracting from $c_r/24$ gives the first formula. For positive $a,b$, $(a+b)^3-a^3-b^3=3ab(a+b)>0$. For $1\le a\le r-1$,
\[
 r^3-a^3-(r-a)^3=3ra(r-a)\ge 3r(r-1).
\]
Merging parts until two remain proves the gap. For $F=\OO_T^{\oplus2}$ one has $\varepsilon(\OO_T)=0$ and $\varepsilon(F)=5(8-2)/24=5/4$.

Proposition~\companionref{prop:determultiple} provides a DVR family of stable
sheaves specializing to this very object.
\end{proof}

\begin{theorem}[Two simple limits on the smooth quintic]\label{thm:step22-nonsep}
For every $m\ge1$ the scalar-rigidified moduli functor which retains
all simple Gieseker-semistable sheaves of character
$2(1-e^{-mH})$ on $X$ fails the uniqueness criterion for a separated
moduli space, even for a family with geometrically stable generic fiber.
\end{theorem}
\begin{proof}
Put $L=\OO_X(mH)$. Take four general sections $s_1,s_2,a,b$ of $L$
and let $D_i=V(s_i)$, $E_i=\OO_{D_i}$. Put $R=\CC[[t]]$ and
$K=\CC((t))$. On $X_R$ take the two cokernels of
$(L^{-1})^{\oplus2}\to\OO^{\oplus2}$ given by
\[
 M_+=\begin{pmatrix}s_1&a\\t^2 b&s_2\end{pmatrix},
 \qquad M_-=\begin{pmatrix}s_1&t^2 a\\ b&s_2\end{pmatrix}.
\]
Their common determinant $s_1s_2-t^2ab$ is regular on each fiber,
so each matrix is injective fiberwise. The two cokernels are
$R$-flat by the local criterion for flatness. Over $K$ the matrices
are conjugate by $\operatorname{diag}(t^2,1)$.
The generic determinant surface is smooth and geometrically integral.
Indeed for $t=1$ evaluation of the four entries onto the fiber of
$L^{\oplus 4}$ is surjective at each point of $X$.

The rank-zero
stratum has codimension four and is avoided on the threefold; the
rank-one stratum of the determinant is smooth, so transversality of
a general section in characteristic zero gives a smooth divisor.
Smoothness is open in the pencil parameter and a smooth ample divisor
is geometrically connected by the divisor sequence and Kodaira
vanishing.

Over the generic surface one matrix entry is a unit
locally; row and column operations identify the cokernel with a line
bundle there. Its pushforward is stable by purity and generic rank one.

The two special fibers fit into
\[
 0\to E_1\to F_+\to E_2\to0,qquad
 0\to E_2\to F_-\to E_1\to0 .
\]
Resolving $E_i$ by $[L^{-1}\xrightarrow{s_i}\OO_X]$ identifies
their extension classes, up to sign, with $a|_{D_1\cap D_2}$ and
$b|_{D_1\cap D_2}$ in $H^0(D_1\cap D_2,L)$, both nonzero.
The stable factors have distinct supports and equal reduced
polynomials. Each extension is strictly semistable.

For an extension $0\to E_i\to F\to E_j\to0$ with class $\xi\ne0$, $\Hom(E_i,E_j)=\Hom(E_j,E_i)=0$ gives
\[
 0\to\End F\to\End E_i\oplus\End E_j=\mathbb C^2
   \xrightarrow{(u,v)\mapsto(u-v)\xi}\Ext^1(E_j,E_i).
\]
Both limits are simple. Applying $\Hom(E_1,-)$ to the two sequences gives
\[
 \Hom(E_1,F_+)=\mathbb C,\qquad
 \Hom(E_1,F_-)=\ker\big(\mathbb C\xrightarrow{1\mapsto\xi_-}
                                    \Ext^1(E_1,E_2)\big)=0 .
\]

The special fibers remain distinct after
scalar rigidification or a finite base change of traits. Their
common associated graded is $E_1\oplus E_2$.
\end{proof}
For the ordinary stable functor, Proposition~\companionref{prop:noordinaryproperextension} exhibits a missing limit. Enlarging it to retain every simple semistable sheaf as a distinct point introduces the two limits $F_+$ and $F_-$ above. They give the same point in Simpson's projective moduli space \cite[Theorem~1.21]{Simpson} because both have associated graded $E_1\oplus E_2$.

Thus retaining every simple semistable sheaf violates separatedness in each even divisor multiplicity. In the Hall logarithm the extensions and their common associated graded belong to the same semistable stack; the integral of that logarithm defines the coefficient of \eqref{eq:generalseries}.
\subsection{The generalized DT series and its predicted modular transformations}
Fix $r$ and a representative $\mu\in\{0,\ldots,5r-1\}$. Take $q_1=\mu+5r^2/2$. For each allowed $\widehat q_0$, choose the character $v$ with these coordinates and let $\overline\Omega_{r,\mu}(\widehat q_0)$ be its generalized Donaldson--Thomas invariant for $H$-Gieseker semistability. Then
\begin{equation}\label{eq:generalseries}
 h^{\mathrm{DT}}_{r,\mu}(\tau)=
 \sum_{\widehat q_0}\overline\Omega_{r,\mu}(\widehat q_0)
          q^{-\widehat q_0},\qquad q=e^{2\pi i\tau}.
\end{equation}
For example $F$ and $F(H)$ have curve charges $q_1$ and $q_1-5r$. Tensoring identifies their moduli and preserves $\widehat q_0$ by Proposition~\ref{prop:generalspectralflow}. The representative just chosen assigns their common invariant to a single coefficient of \eqref{eq:generalseries}.

For $r=1$ these coefficients are ordinary stable-sheaf virtual degrees. For $r>1$ they also include the Hall weights of all strictly semistable objects of the same character. Proposition~\ref{prop:generalvacuumbound} makes the part with negative exponents finite in each case.

For the large-volume D4--D2--D0 index, the proposed completion adds terms depending on invariants in smaller positive divisor multiplicities. In the normalization of \cite[\S\S2.4--2.5, (2.34)]{AFKPmock} its two transformations are
\begin{align}
 \widehat h_{r,\mu}(\tau+1,\bar\tau+1)
  &=e^{2\pi i\alpha_{r,\mu}}\widehat h_{r,\mu}(\tau,\bar\tau),
   \qquad
   \alpha_{r,\mu}=-\frac{c_r}{24}
           +\frac{\mu^2}{10r}+\frac{r\mu}{2},
       \label{eq:generalcompletedT}\\
 \widehat h_{r,\mu}(-1/\tau,-1/\bar\tau)
  &=\frac{(-1)^{a_r+1}(-i\tau)^{-3/2}}{\sqrt{5r}}
       \sum_{\nu\bmod5r}e^{-2\pi i\mu\nu/(5r)}
                         \widehat h_{r,\nu}(\tau,\bar\tau).
       \label{eq:generalcompletedS}
\end{align}
All roots use the continuous branch on the upper half-plane.  The $T$-factor is well defined on the residue class because
\[
 \alpha_{r,\mu+5r}-\alpha_{r,\mu}
      =\mu+\frac{5r(r+1)}2\in\mathbb Z.
\]
For $r=1$, $a_1=5$ and $c_1=55$; the two equations reduce to \eqref{eq:Tmodnew}--\eqref{eq:Smodnew}.

For completeness, the phase conversion from \cite[Equation (2.34)]{AFKPmock} is explicit.  Its $S$-matrix multiplies $\tau^{-3/2}$ by $(-1)^{a_r}/\sqrt{i5r}$.  Since
\[
 \frac{\tau^{-3/2}}{\sqrt i}
      =-(-i\tau)^{-3/2},
\]
this is exactly \eqref{eq:generalcompletedS}.  Its $T$-exponent is
\[
 \frac{(\mu+5r^2/2)^2}{10r}+\frac{50r}{24},
\]
whose difference from $\alpha_{r,\mu}$ is $a_r\in\mathbb Z$.  The two multipliers therefore agree.

\begin{conjecture}\label{conj:generalSduality}
The rational series \eqref{eq:generalseries} admit the completed transformation laws \eqref{eq:generalcompletedT}--\eqref{eq:generalcompletedS}.  The primitive series is weakly holomorphic.  For $r>1$ the holomorphic series is the holomorphic part of a mock modular object of depth at most $r-1$, with the completion prescribed by the smaller divisor-charge series in the S-duality formulation of \cite{AFKPmock}.
\end{conjecture}

For instance, the charge-two formulation in \cite[Equations (2.36)--(2.40)]{AFKPmock}, specialized to the quintic degree $\kappa=5$, predicts for the geometric series
\begin{equation}\label{eq:charge2anomaly}
 h^{\mathrm{DT}}_{2,\mu}=h^{(0)}_{2,\mu}
       +(-1)^{\mu+5}\sum_{a\bmod5}
          G^{(5)}_{\mu-2a+5}\,h^{\mathrm{DT}}_{1,a}h^{\mathrm{DT}}_{1,\mu-a}.
\end{equation}
Here $h^{(0)}_2$ is a weakly holomorphic vector of weight $-3/2$ for the charge-two representation, and $G^{(5)}$ has weight $3/2$ and the specified mock shadow.  Its index is modulo $10$; the $h^{\mathrm{DT}}_1$ indices are modulo $5$.  One concrete normalization is $G^{(5)}_\lambda=\tfrac12(V_5\mathcal H)_\lambda$, as in the cited source.  Here
\[
 \mathcal H_\epsilon(\tau)
 =\sum_{\substack{D\ge0\\D\equiv-\epsilon^2\!\!\!\pmod4}}
       H(D)q^{D/4},\qquad \epsilon\in\mathbb Z/2,
\]
where $H(0)=-1/12$, and for $D>0$ the Hurwitz class number $H(D)$ is the sum over $\operatorname{SL}_2(\mathbb Z)$ classes of integral positive definite binary quadratic forms of discriminant $-D$, including imprimitive forms, with weight $2/|\operatorname{Aut}_{\operatorname{SL}_2}(Q)|$.  Values with no such forms are zero.  For a two-component vector $\phi$, define
\begin{equation}\label{eq:HurwitzHecke5}
 (V_5\phi)_\lambda(\tau)
 =5\sum_{\substack{a,d>0\\ad=5}}d^{-2}
       \sum_{b=0}^{d-1}{\bf1}_{a\mid\lambda}
       e^{\pi i b\lambda^2/(10a)}
       \phi_{\lambda/a\bmod2}\!\left(\frac{a\tau+b}{d}\right).
\end{equation}
The index $\lambda$ is taken modulo ten. Formula~\eqref{eq:charge2anomaly} describes the proposed modular anomaly. Each quadratic term has total weight $3/2-3/2-3/2=-3/2$. Once $h^{\mathrm{DT}}_1$ is known, this formula prescribes the quadratic term. The weakly holomorphic summand $h^{(0)}_2$ still has polar coefficients to be determined at charge $2H$.

For $r=1$ the criterion in Theorem~\ref{thm:twodegmodularityfull} uses the multiplicity-one, four-sector K3 series of Section~\ref{sec:modularity}. For $r>1$ the fiber-charge normalization, support profiles and boundary stability enter separately (Proposition~\companionref{prop:generalprofiles} and Theorem~\companionref{thm:generalpolynomialstability}). A derivation of \eqref{eq:charge2anomaly} and its higher-depth analogues requires relative evaluations for all boundary classes, including semistable flags.

\subsection{Charge conjugation for the sheaf series}
\label{sec:hall-charge-conjugation}
The primitive reconstruction uses equality of opposite D2 residues. On the smooth quintic we prove this symmetry for generalized Donaldson--Thomas invariants of $H$-Gieseker-semistable sheaves in every positive divisor class, including strictly semistable sheaves. For the large-volume index see \cite[discussion following (2.56)]{AFKPS}.

\begin{theorem}\label{thm:rank-zero-charge-conjugation}
For every $r>0$ and every numerical Chern character $v=(0,rH,\beta,n)$, put $v^\natural=(0,rH,-\beta,n)$. With the generalized Donaldson--Thomas invariants of $H$-Gieseker-semistable sheaves of \textup{(\companionref{eq:Hallloggeneral})},
\[
 \overline\Omega_X(v)=\overline\Omega_X(v^\natural),\qquad
 h^{\mathrm{DT}}_{r,\mu}=h^{\mathrm{DT}}_{r,-\mu}.
\]
These are identities of individual invariants and of formal series. For $r=1$ the invariants on both sides are degrees of virtual zero-cycles.
\end{theorem}

\begin{proof}
For a pure surface sheaf put $r(F)H=\ch_1(F)$, $d(F)=H\ch_2(F)$ and $s(F)=d(F)/r(F)$. Slope semistability is equivalent to
\[
 s(A)\le s(F)\quad\text{whenever}\quad
 0<r(A)<r(F),\qquad F/A\text{ is pure of dimension two}.
\]
Indeed saturation preserves $r(A)$ and adds the degree of an effective curve cycle to $d(A)$. For $r(A)=r(F)$ the inequality is automatic.

\smallskip\noindent\emph{Duality of slopes.}
Set $D(K)=\RR\mathcal Hom_X(K,\OO_X)[1]$. By Theorem~\companionref{thm:generalCMhull},
\[
 \mathcal H^0(DF)=F^{\mathrm d},\quad
 \mathcal H^1(DF)=T_F,
\]\[
 F\subset (F^{\mathrm d})^{\mathrm d},\quad
 \dim\Supp(T_F),\ \dim\Supp((F^{\mathrm d})^{\mathrm d}/F)\le0 .
\]
Hence $r(F^{\mathrm d})=r(F)$ and $d(F^{\mathrm d})=-d(F)$. If $F$ is slope semistable and $F^{\mathrm d}$ is destabilized by a saturated subsheaf, the complementary pure quotient $C$ has
$s(C)<s(F^{\mathrm d})$. Left exactness of $(-)^{\mathrm d}$ on such quotient sequences follows from $\mathcal Hom_X(A,\OO_X)=0$ for a torsion sheaf $A$. Consequently
\[
 C^{\mathrm d}\lhook\joinrel\longrightarrow (F^{\mathrm d})^{\mathrm d},
 \qquad
 s(C^{\mathrm d}\cap F)=-s(C)>s(F).
\]
The intersection has the same $r,d$ as $C^{\mathrm d}$ since their quotient is finite. This contradicts semistability of $F$. Any finite-colength subsheaf of $F^{\mathrm d}$ is again pure and slope semistable: it has the same $r,d$ and its subsheaves are subsheaves of $F^{\mathrm d}$.

\smallskip\noindent\emph{An open heart for the Hall operations.}
Let $\mathcal T=\operatorname{Coh}_0(X)$ and $\mathcal U=\{A\in\operatorname{Coh}(X):\Hom(\mathcal T,A)=0\}$. Quotient by the maximal finite-length subsheaf gives the torsion pair $(\mathcal T,\mathcal U)$. Its right tilt
\[
 \mathcal A=\langle\mathcal U,\mathcal T[-1]\rangle
 =\{E:\mathcal H^0(E)\in\mathcal U,\ \mathcal H^1(E)\in\mathcal T,\ \mathcal H^i(E)=0\ (i\ne0,1)\}
\]
is a bounded heart. We verify the openness required in \cite[\S2.5, (2.6)]{TodaHall22}. For $K$ perfect put $K^*=\RR\mathcal Hom(K,\OO_X)$. The following conditions are equivalent to $K\in\mathcal A$:
\begin{align*}
 &\mathcal H^{>1}(K)=0,\quad \dim\Supp\mathcal H^1(K)\le0,\\
 &\mathcal H^{\ge3}(K^*)=0,\quad
   \operatorname{codim}\Supp\mathcal H^i(K^*)\ge i\quad(i=1,2).
\end{align*}
For $K\in\mathcal A$, dualize the ordinary truncation triangle. The dual of $\mathcal H^1(K)$ is concentrated in degree three, while $\mathcal Ext^3(\mathcal H^0(K),\OO_X)=0$ by depth. Regularity gives the other codimension bounds. Conversely, in
\[
 E_2^{p,q}=\mathcal Ext^p(\mathcal H^{-q}(K^*),\OO_X)
          \ \Longrightarrow\ \mathcal H^{p+q}(K)
\]
the grade bound $p\ge\operatorname{codim}\Supp\mathcal H^{-q}(K^*)$ for a nonzero term excludes $p+q<0$. The truncation triangle now gives $\mathcal Ext^3(\mathcal H^0(K),\OO_X)=\mathcal H^3(K^*)=0$. By Auslander--Buchsbaum, $\mathcal H^0(K)$ has depth at least one wherever it is nonzero at a closed point, as required.

For a perfect family $P$ on $X\times B$, let $L_i(P)$ be the locus where $H^j(P\otimes^{\LL} k(x,b))\ne0$ for some $j\ge i$. It is closed by the rank conditions on a finite free representative. On each fiber, Nakayama in the largest nonvanishing degree gives
\[
 L_i(P)\cap X_b=\bigcup_{j\ge i}\Supp\mathcal H^j(P_b).
\]
The conditions above are therefore upper bounds on the fiber dimensions of the proper morphisms $L_i(P),L_i(P^*)\to B$. Upper semicontinuity gives openness. For $X$ smooth, a relatively perfect family is locally perfect on $X\times B$. This proves the openness used for the Behrend integration of \cite[Theorem~2.9]{TodaHall22}.

\smallskip\noindent\emph{The two filtrations.}
For a fixed rational $s$ write $U_s$ for the Hall element of slope-semistable pure surface sheaves of slope $s$, including the zero object, and $Z$ for that of $\mathcal T[-1]$, also including zero. The product records the subobject first. The middle objects in $U_s*Z$ are exactly those with a unique ordinary truncation
\[
 F=\mathcal H^0(E)\longrightarrow E\longrightarrow Q[-1],
 \qquad Q=\mathcal H^1(E).
\]
Let $Q^\dagger=\mathcal Ext^3(Q,\OO_X)$. As $D(Q[-1])=Q^\dagger[-1]$, the dual triangle gives the sequence
\[
 0\to \mathcal H^0(DE)\to F^{\mathrm d}\to Q^\dagger
   \to\mathcal H^1(DE)\to T_F\to0 .
\]
Thus $\mathcal H^0(DE)$ has finite colength in $F^{\rm d}$ and $\mathcal H^1(DE)$ has finite length. The first is zero or slope semistable of slope $-s$. Biduality gives the converse. Dualizing the ordinary truncation of $DE$ gives the unique filtration of $E$ with factors in $Z,D_*U_{-s}$, in that order. Hence
\begin{equation}\label{eq:dual-charge-hall}
 U_s*Z=Z*D_*U_{-s}.
\end{equation}
Here $D$ acts on the stack of complexes by sending a family to its relative dual shifted by one, and an isomorphism to the inverse of its dual. This commutes with base change and preserves the intrinsic Behrend function. For the Hall equality, on each finite-type family stratify until the cohomology sheaves of the family and its dual are flat, using generic flatness and Noetherian induction on the reduced base.

The cohomology and the two truncation triangles then commute with base change on each stratum. Uniqueness identifies the automorphism groups of the filtered and middle objects. Equivalently the forgetful maps from the filtration stacks are geometric bijections onto the same constructible strata, which suffice by \cite[\S2.5, relation~(2)]{TodaHall22}.

For fixed output $(r,d,n)$ at slope $s$, the factors satisfy $\sum_i r_i=r$, $d_i=sr_i$, $n=\sum_i n_i-\ell$, $\ell\ge0$. Proposition~\ref{prop:generalvacuumbound} applies to each slope-semistable factor. With $B_i=5r_i^3/24+d_i^2/(10r_i)$ it gives
\[
 n-\sum_{j\ne i}B_j\le n_i\le B_i,\qquad
 0\le\ell\le\sum_i B_i-n.
\]
The ranks have finitely many partitions and the numerical charges lie in a discrete lattice. For each character the Gieseker-semistable stack is of finite type, as is the stack of points of fixed length. The same bounds applied to the Gieseker Harder--Narasimhan filtration give finitely many types and finite-type extension stacks for the slope locus. This permits coefficientwise Hall logarithms and integration.

\smallskip\noindent\emph{Integration and strictly semistable sheaves.}
Write $\mathbb L=[\mathbb A^1]$, $u_s=(\mathbb L-1)\log_*U_s$ and $z=(\mathbb L-1)\log_*Z$. Factor $U_s$ by the ordinary Gieseker Harder--Narasimhan filtration:
\[
 U_s=\prod_{p}^{\searrow}G_p ,
\]
Here $G_p$ comprises zero and the Gieseker-semistable sheaves of reduced polynomial $p$ at slope $s$. For the HN factors of a slope-semistable $F$,
\[
          s_1\ge\cdots\ge s_k,\qquad
          s_1\le s(F)=\frac{\sum r_i s_i}{\sum r_i}.
\]
Hence every $s_i=s(F)$. Conversely, extensions of pure semistable sheaves of this common slope remain pure and slope semistable.

The Joyce absence-of-poles theorem for semistable loci \cite[Theorem~8.7]{JoyceNoPoles22} gives regular elements $g_p=(\mathbb L-1)\log_*G_p$ and $z$ (the latter in the Hall algebra of finite-length sheaves shifted by $-1$), as in \cite[Theorem~3.14, (3.19)--(3.20)]{TodaHall22} and the definition \textup{(\companionref{eq:Hallloggeneral})}. The extension stacks used for a fixed $p$ consist of sheaves in both $\operatorname{Coh}(X)$ and $\mathcal A$, and agree in the two hearts. Thus the same regularity and integration apply to those stack functions here.

Baker--Campbell--Hausdorff expresses $u_s$ in terms of the $g_p$ and iterated brackets $[a,b]_*/(\mathbb L-1)$. In the numerical Poisson algebra of \cite[Theorem~2.9]{TodaHall22}, the bracket of two basis elements contains the numerical Euler form as a factor. For two surface charges $v_i=(0,r_iH,\beta_i,n_i)$ of slope $s$ and any shifted point charge $w=(0,0,0,-\ell)$, Riemann--Roch gives
\[
 \chi(v_1,v_2)=r_2d_1-r_1d_2=0,
 \qquad \chi(w,v_i)=0 .
\]
All the nested brackets in the integrated Harder--Narasimhan product therefore vanish, as in the slope-to-Gieseker comparison of \cite[Lemma~5.2]{FeyzOSV}. With the stack Behrend sign used in \cite[(3.21)--(3.22)]{TodaHall22}, the integration $I$ satisfies
\[
 I(u_s)=\sum_p I(g_p)
       =-\sum_{d/r=s}\overline\Omega_X(0,rH,\beta,n)e_{(0,rH,\beta,n)} .
\]
The equality includes every term of the Gieseker Hall logarithm in a strictly semistable character. Taking logarithms in \eqref{eq:dual-charge-hall}, multiplying by $\mathbb L-1$ and passing to the semiclassical Hall algebra gives
\[
 u_s=\exp\bigl(\{z,-\}\bigr) D_*u_{-s}.
\]
Duality reverses products but fixes the order of the identical factors in each power of $U_{-s}$. It therefore commutes with its logarithm and preserves regularity. By $\chi(w,v_i)=0$, application of $I$ gives $I(u_s)=I(D_*u_{-s})$. Duality preserves the Behrend function and sends $v$ to $v^\natural$, proving the coefficient identity.

Finally $q_1\mapsto-q_1$ and $q_0\mapsto q_0$. Thus $\widehat q_0$ is fixed and $\mu\mapsto-\mu-5r^2\equiv-\mu\pmod{5r}$. Summing in a fixed residue proves the asserted formal-series symmetry.
\end{proof}

\subsection{Two sources of higher divisor-charge coefficients}\label{sec:hall22}
Let $X$ be the smooth quintic. We first determine the part of its rational invariants forced by charge divisibility. We then evaluate full invariants on a range of numerical classes by a wall-crossing theorem on $X$. Both calculations use the generalized Donaldson--Thomas series in \eqref{eq:generalseries}.

\begin{proposition}[The contribution forced by divisibility]\label{prop:hall22multi}
For $v$ in the numerical Grothendieck group, set
\[
 b(v)=\sum_{m\mid v}\frac{\mu_{\rm ar}(m)}{m^2}\,\overline\Omega(v/m),
 \qquad \overline\Omega(v)=\overline{\dt}_H(v).
\]
Here $m\mid v$ means $v/m\in K_{\rm num}(X)$ and $\mu_{\rm ar}$ is the arithmetic M\"obius function. Write $h^b_{r,\mu}$ for the series \eqref{eq:generalseries} with coefficients $b$. In divisor charge $2H$ one has
\begin{equation}\label{eq:hall22double}
 h^{\mathrm{DT}}_{2,\mu}(\tau)-h^b_{2,\mu}(\tau)=
 \begin{cases}
 \tfrac14 h^{\mathrm{DT}}_{1,(\mu+5)/2}(2\tau), & \mu\ \text{odd},\\
 0,&\mu\ \text{even}.
 \end{cases}
\end{equation}
The index on the right belongs to $\mathbb Z/5$. Each equality is coefficientwise, including at exponents which do not arise by doubling.
\end{proposition}
\begin{proof}
The definition is \cite[Definition 6.10]{JS}. The convolution identity
$\sum_{d\mid n}\mu_{\rm ar}(d)=\mathbf{1}_{n=1}$ gives
$\overline\Omega(v)=\sum_{m\mid v}m^{-2}b(v/m)$. For a class $w$ of divisor multiplicity $s$, the coordinates \eqref{eq:generalEMcharges} give
\[
 \widehat q_0(mw)=m\widehat q_0(w),\qquad
 \mu(mw)\equiv m\mu(w)-\tfrac52m(m-1)s^2\pmod{5ms}.
\]
In charge $2H$, divisibility greater than one forces $m=2$ and $s=1$. Primitive divisor charge implies $b(w)=\overline\Omega(w)$. The map $a\mapsto 2a-5$ identifies $\mathbb Z/5$ with the odd residue classes modulo ten. It is injective on charge representatives modulo twists as well: $2w'=2(e^{tH}w)$ implies $w'=e^{tH}w$ in the torsion-free numerical charge lattice. Summing the correction $\overline\Omega(w)q^{-2\widehat q_0(w)}/4$ proves the assertion.
\end{proof}
For $w=1-e^{-H}$ the quintic invariant is $5$. Thus \eqref{eq:hall22double} gives $b(2w)+5/4$ at $(\mu,x)=(5,-55/12)$. The term $5/4$ is the divisibility correction, and $b(2w)$ requires independent input. The contribution of the full strictly semistable locus in this charge is computed in Corollary~\ref{cor:strict22number} by integrating its Hall logarithm.

\begin{theorem}[A range of exact polar coefficients by wall crossing]\label{thm:hall22strip}
Let $\ell\ge 0$, $k\ge 6\ell+1$, and write
\[
 v_{k,\ell}=1-e^{-kH}-\ell[\mathrm{pt}],\qquad
 \sum_{j\ge0}I_{j,0}z^j=\prod_{m\ge1}\bigl(1-(-z)^m\bigr)^{200m}.
\]
Every Gieseker-semistable sheaf of class $v_{k,\ell}$ is stable. Its full sheaf invariant is
\begin{equation}\label{eq:hall22strip}
 [q^{\ell-c_k/24}]\,h^{\mathrm{DT}}_{k,0}
 =\overline\Omega(v_{k,\ell})
 =(-1)^{a_k-\ell-1}(a_k-\ell) I_{\ell,0}.
\end{equation}
This identity ranges over all sheaves of the specified charge. It imposes no assumption on the support or on restriction to a K3 surface.
\end{theorem}
\begin{proof}
Throughout, $H^3=5$ and $\operatorname{Pic}(X)=\mathbb ZH$. We have $d=-5k^2/2$, $n=5k^3/6-\ell$, $\mu=0$ and $-\widehat q_0=\ell-c_k/24<0$. The discriminant in \cite[Theorem 1.1]{FeyzOSV} equals
\[
 Q=\frac{k^2}{2}+6\left(\frac{d}{5k}\right)^2-\frac{12n}{5k}=\frac{12\ell}{5k}.
\]
The hypothesis of \cite[Theorem 1.1(ii), Lemma 3.7]{FeyzOSV} is
\[
 25Q < 5k+\frac{2}{5k}-\frac52-\frac{2}{25k^2}.
\]
For $\ell=0$ its right-hand side is positive. For $\ell\ge1$, $k\ge7$, its left-hand side is $60\ell/k<10$ and its right-hand side exceeds $32$. The theorem gives stability and a signed sum of products $P_{-m_1,\beta_1}I_{m_2,\beta_2}$ for rank $-1$ and $1$ factors
\[
 v_1=-e^{D_1}(1,0,-\beta_1,-m_1),\quad
 v_2=e^{D_2}(1,0,-\beta_2,-m_2),\qquad v_1+v_2=v_{k,\ell}.
\]
The determinant torsion factor has order one. Each nonzero summand supplies actual objects of the indicated stable-pair and ideal-sheaf types. By \cite[Propositions~3.5--3.6 and (31)]{FeyzOSV}, the effective curve classes of such a summand satisfy
\[
 0\le H\beta_i\le\frac52\Bigl(k^2-k\sqrt{k^2-24\ell/(5k)}\Bigr).
\]
Indeed the two terms involving $D_i$ on the left of (31) cancel for $D_i\in\mathbb Z H$. By $k\ge6\ell+1$ we have $\ell<k/6-1/(30k)$, equivalently $k^2-24\ell/(5k)>(k-2/(5k))^2$. The upper bound is strictly less than one, so the integral degrees vanish and ampleness gives $\beta_1=\beta_2=0$. Write $D_1=tH$, $D_2=(t+k)H$. The degree-four component gives
\[
 5kt+\tfrac52 k^2=d=-\tfrac52 k^2,\qquad t=-k.
\]
A stable pair with zero curve class has zero target, hence $P_{-m_1,0}=\mathbf 1_{m_1=0}$. The degree-six component gives $m_2=\ell$. For the two remaining triples the inequalities of \cite[Definition~3.4]{FeyzOSV} reduce, after canceling determinant terms, to
\[
 0\le k/10-1/25,\qquad \max(0,-\ell)\le k(k+1)/6 .
\]
Both right sides are positive.

Riemann--Roch now gives, for the ideal $I_Z$ of a length-$\ell$ subscheme,
\[
 \chi(v_2,v_1)=-\chi(I_Z,\mathcal O_X(-kH))
 =-\int_X(1+\ell[\mathrm{pt}])e^{-kH}\operatorname{td}(X)=a_k-\ell.
\]
The signed wall factor is $(-1)^{\chi(v_2,v_1)-1}\chi(v_2,v_1)$. Finally $\chi_{\rm top}(X)=-200$, and the degree-zero DT formula \cite[Theorem 4.12]{BFpoints22} gives the displayed product for $I_{\ell,0}$.
\end{proof}
In particular $[q^{-35/6}]h^{\mathrm{DT}}_{2,0}=15$. For $k=7$, $a_7=315$, $c_7=2065$, and $I_{1,0}=200$, whence $[q^{-2041/24}]h^{\mathrm{DT}}_{7,0}=-62\,800$. These identities use wall crossing on the smooth quintic. The degeneration comparison table records coefficients obtained from the relative moduli spaces.

The recursion of \cite[Theorem 1.2]{FeyzOSV} also applies beyond this range. Its inputs include threefold stable-pair invariants and invariants of smaller divisor multiplicity. In \cite[Theorems 1 and 2]{FTone22}, rank refers to $\operatorname{ch}_0$; all the classes considered here have $\operatorname{ch}_0=0$. Thus a comparison of primitive divisor-charge series supplies only part of the input to a recursion for the higher divisor charges.

\subsection{A range with nonzero curve class}\label{sec:lines22}
We next determine polar Fourier coefficients in curve-charge residue $1$. At the stability wall used in Theorem~\ref{thm:hall22strip}, the rank-one factor now carries a curve. Write $[L]$ for the class of a line in the quintic, normalized by $H[L]=1$.

\begin{theorem}[Line and point charges on the quintic]\label{thm:hall22lines}
Suppose $j\ge1$, $k\ge\max(18,10j)$, and set $a=a_k-k-j$. For
\[
 w_{k,j}=1-e^{-kH}-[L]-j[\mathrm{pt}],\qquad
 x_{k,j}=-c_k/24+k/2+1/(10k)+j,
\]
all Gieseker-semistable sheaves in class $w_{k,j}$ are stable, and
\begin{align}
 [q^{x_{k,j}}]h^{\mathrm{DT}}_{k,1}
 &=(-1)^{a-1}a\cdot 2875\,[z^{j-1}]
 (1+z)^{-2}\prod_{m\ge1}\bigl(1-(-z)^m\bigr)^{200m}.
 \label{eq:hall22lines}
\end{align}
\end{theorem}
\begin{proof}
The charge map gives $\mu=1$ and $-\widehat q_0=x_{k,j}$. Put $d=-5k^2/2-1$ and $n=5k^3/6-j$. Then
\[
 Q=\tfrac65+\tfrac6{25k^2}+\tfrac{12j}{5k},\qquad
 25Q=30+6/k^2+60j/k<5k+2/(5k)-5/2-2/(25k^2).
\]
Indeed $Q<1.45$ while $5k-5/2\ge 87.5$. By \cite[Theorem 1.1(ii), Lemma 3.7]{FeyzOSV} this implies stability and the two-factor formula used in Theorem~\ref{thm:hall22strip}. Set $b_i=H\beta_i$ for the curve degrees on its two sides. Equations (31)--(32) of \cite[Proposition 3.5]{FeyzOSV} give
\begin{align*}
 0\le b_i&\le\frac{5Q}{1+\sqrt{1-2Q/k^2}}<4,\\
 m_1&\le \tfrac23 b_1(b_1+5/2),\qquad
 -m_2\le \tfrac23 b_2(b_2+5/2).
\end{align*}
The square root is greater than $0.99$ for $k\ge18$ and $Q<1.45$. Write $D_2=tH$, $D_1=(t-k)H$. The degree-four identity is
\[
 -5k^2/2-1=5kt-5k^2/2+b_1-b_2.
\]
Since $0\le b_i\le3$, $5kt=b_2-b_1-1$ forces $t=0$ and $b_2=b_1+1$. The degree-six identity now reads $m_1-m_2=kb_1-j$. If $b_1=1$, the two bounds on $m_i$ give $k-j\le25/3$, contrary to $k-j\ge0.9k\ge16.2$. If $b_1=2$, they give $2k-j\le17$, contrary to $2k-j\ge1.9k$. Thus $(b_1,b_2)=(0,1)$, the zero-class stable pair forces $m_1=0$, and $m_2=j$. Conversely these data satisfy \cite[Definition 3.4]{FeyzOSV}, since $1/5\le k/10-1/25$ and $-j\le0$. The Euler form is
\[
 \chi(v_2,v_1)=-\int_X(1-[L]+j[\mathrm{pt}])e^{-kH}\operatorname{td}(X)
 =a_k-k-j.
\]
It remains to compute $I_{j,[L]}$. By deformation invariance of the threefold ideal-sheaf invariants, take a general quintic. Its scheme of lines is the zero scheme of a general section of the globally generated vector bundle $\operatorname{Sym}^5 U^\vee$ on $\operatorname{Gr}(2,5)$; $U$ denotes the tautological rank-two bundle. This bundle has rank six on a smooth sixfold, and the restriction $H^0(\mathbb P^4,\mathcal O(5))\to H^0(L,\mathcal O_L(5))$ is onto for every line.

Generic smoothness for the incidence of pairs $(L,X)$ therefore makes the line scheme finite and reduced. Its length, with $x,y$ the Chern roots of $U^\vee$, is
\[
 \int_{\operatorname{Gr}(2,5)}c_6(\operatorname{Sym}^5U^\vee)
 =[x^4y^3](x-y)\prod_{i=0}^{5}(ix+(5-i)y)=2875.
\]
Every effective one-cycle of degree one is a single line. Each contributes the series $z(1+z)^{-2}\sum_{s\ge0}I_{s,0}z^s$ by the smooth-curve DT/PT formula \cite[Theorem 1.1 and Section 1.1]{Ric22}: the curve has genus zero, its point in the Hilbert scheme of lines is smooth of dimension zero, and the corresponding Behrend weight is $1$. Summing over lines gives \eqref{eq:hall22lines}.
\end{proof}
For $k=20$, the coefficients of $h^{\mathrm{DT}}_{20,1}$ at $x_{20,1}$ and $x_{20,2}=x_{20,1}+1$ are respectively $19\,345\,875$ and $-3\,829\,914\,000$.

The first counts the class with a single line in the rank-one factor. The second also permits one zero-dimensional defect, including a defect supported on the line. The factor $[z](1+z)^{-2}\sum I_{s,0}z^s=200-2=198$ includes that contribution. These are higher-charge wall-crossing evaluations.

\subsection{The complete strictly semistable locus in a doubled charge}\label{sec:kronecker22}
Hecht proposes the five-arrow Kronecker quiver for the polar class $(\mu,\nu)=(5,5/4)$ in \cite[Section~6, Table~5]{Hecht22}, obtaining the rational quiver index $-155/4$. Here $\nu=c_2/24-\widehat q_0$ with $c_2=140$ denoting the divisor-charge central charge. We identify exactly which part of the sheaf stack this quiver parametrizes and determine its contribution to the sheaf invariant.

Write $w=1-e^{-H}$, and let $K_5$ be the quiver with five arrows from one vertex to another. A dimension-$(n,n)$ representation is an $n$-by-$n$ matrix $M$ of linear forms in $H^0(X,\mathcal O_X(H))$ after bases have been chosen. A subrepresentation is a pair of subspaces $A,B$ such that every coefficient matrix carries $A$ into $B$.

Use the stability condition $\dim A\le\dim B$ on subrepresentations, with a strict inequality for every nonzero proper subrepresentation in the stable case.

\begin{theorem}[The Hall contribution in the class $2(1-e^{-H})$]\label{thm:kronecker22}
The semistable stack of sheaves on $X$ of charge $2w$ is a disjoint union of an open and closed substack $\mathcal K$ and a stable locus $\mathcal U$. Taking $\operatorname{coker}(M:\mathcal O_X(-H)^2\to\mathcal O_X^2)$ identifies $\mathcal K$ with the semistable $(2,2)$-representation stack of $K_5$. The contribution of $\mathcal K$ to the generalized sheaf invariant is $-155/4$. Consequently, writing $U$ for the rigidified stable moduli of $\mathcal U$,
\begin{equation}\label{eq:Kronfull22}
 [q^{-55/12}]h^{\mathrm{DT}}_{2,5}=-\tfrac{155}{4}+\chi(U,\nu_U).
\end{equation}
Here $U$ is a proper union of components of the Simpson moduli space. The residual weighted Euler characteristic is a stable-sheaf invariant.
\end{theorem}
\begin{proof}
\smallskip\noindent\emph{1. Stability of a matrix and of its cokernel.}
Write $M=\sum_{j=0}^4 x_jM_j$, where $(x_0,\ldots,x_4)$ is a basis of $H^0(X,\OO_X(H))$. For any subspaces $A,B\subset\CC^2$ the condition on a subrepresentation is $M_j A\subset B$ for all $j$. The dimension pairs violating semistability are
\[
              (\dim A,\dim B)=(1,0),\ (2,1),\ (2,0).
\]
Each such pair forces $\det M\equiv0$. Conversely, suppose the span of the $M_j$ lies in the rank-one cone. For two nonzero members $u\otimes v$ and $u'\otimes v'$, the rank of their sum is two unless $u\parallel u'$ or $v\parallel v'$. A pair sharing only its first factor forces every other member to share that factor too; likewise for the second. Thus all $M_j$ have a common image line or a common kernel line.

Hence a semistable $M$ has nonzero determinant and gives the resolution
\[
 0\longrightarrow \OO_X(-H)^2\xrightarrow{M}\OO_X^2
       \longrightarrow F_M\longrightarrow0 .
\]
The sheaf $F_M$ is Cohen--Macaulay of dimension two.

For stable $M$ set $W=\langle M_0,\ldots,M_4\rangle\subset\operatorname{Mat}_2(\CC)$. If $\dim W=1$, every source line maps into a line. If $\dim W=2$, choose a singular $A\in\PP W$ and $\ell\subset\ker A$. For a second generator $B$ one has $W\ell=B\ell$, of dimension at most one. Both cases contradict stability, so $\dim W\ge3$.

If $\dim W=3$, the restricted determinant can degenerate only when $\PP W$ is tangent to the smooth determinant quadric in $\PP^3$. After changes of the two bases this means that the same matrix entry vanishes on all of $W$, giving a $(1,1)$ subrepresentation. Thus in the stable case $\det|_W$ has rank three or four, so $\det M$ is irreducible.

On $X$ a decomposition of its divisor into effective summands would give $2H=H+H$ since $\Pic(X)=\mathbb ZH$. Both summands would be hyperplane sections, including when they coincide. The restriction $H^0(\PP^4,\OO(2))\to H^0(X,\OO_X(2H))$ is an isomorphism, so the quadratic polynomial would factor. It follows that $F_M$ has rank one on its integral support, and is stable.

For a strictly semistable matrix, choose bases adapted to a $(1,1)$ subrepresentation. The two diagonal entries $s_1,s_2$ are nonzero linear forms, since $\det M\ne0$. Writing $D_i=(s_i=0)$, the snake lemma on the two triangular blocks gives
\[
 0\longrightarrow\OO_{D_1}\longrightarrow F_M\longrightarrow\OO_{D_2}
 \longrightarrow0 .
\]
The factors have the same reduced Hilbert polynomial. This proves strict semistability of $F_M$.

\smallskip\noindent\emph{2. Recovering the matrix from a family.}
Let $p:X\times T\to T$ and let $F_T$ be flat, with fiber $F_M$ at $t$. The resolution and Kodaira vanishing give $h^0(F_M)=2$, $h^1(F_M)=0$. Near $t$, base change and Nakayama give
\[
 V=p_*F_T\ \text{locally free of rank two},\qquad p^*V\twoheadrightarrow F_T .
\]
Its kernel $K$ is $T$-flat and $K_t=\OO_X(-H)^2$.

The same argument for $K(H)$ gives a rank-two bundle $U=p_*K(H)$ and an isomorphism $p^*U\xrightarrow{\sim}K(H)$. This yields
\[
 0\longrightarrow p^*U\otimes\OO_X(-H)
   \xrightarrow{M_T}p^*V\longrightarrow F_T\longrightarrow0 .
\]
Both bundles and the arrow commute with base change. The frames of $U$ and $V$ give the two changes of basis of the quiver representation. Thus the cokernel construction identifies the quiver stack with an open substack $\mathcal K$ of the sheaf stack.

\smallskip\noindent\emph{3. Which sheaves lie outside this substack?}
A strictly semistable sheaf of class $2w$ has two stable factors of multiplicity one. Equality of their reduced Hilbert polynomials with the reduced polynomial of $2w$ gives character $w$ for each. Lemma~\companionref{lem:vacuumgenericfull} identifies them with $\OO_{D_1},\OO_{D_2}$. For their extension $F$, the section $1\in H^0(\OO_{D_2})$ lifts since $H^1(\OO_{D_1})=0$. The two resulting generators give the commutative diagram with exact rows and columns
\[
\begin{tikzcd}[column sep=small,row sep=small]
 &0\arrow[d]&0\arrow[d]&0\arrow[d]&\\
0\arrow[r]&\OO_X(-H)\arrow[r]\arrow[d]&K\arrow[r]\arrow[d]&\OO_X(-H)\arrow[r]\arrow[d]&0\\
0\arrow[r]&\OO_X\arrow[r]\arrow[d]&\OO_X^2\arrow[r]\arrow[d]&\OO_X\arrow[r]\arrow[d]&0\\
0\arrow[r]&\OO_{D_1}\arrow[r]\arrow[d]&F\arrow[r]\arrow[d]&\OO_{D_2}\arrow[r]\arrow[d]&0\\
 &0&0&0&
\end{tikzcd}
\]
The top row splits because $\Ext^1_X(\OO_X(-H),\OO_X(-H))=H^1(\OO_X)=0$. This puts $F$ in $\mathcal K$.

It follows that $\mathcal K$ is saturated for $S$-equivalence: it contains every strictly semistable sheaf, and a stable sheaf forms its own equivalence class. Its image in the Simpson quotient is open. It is also the image of the projective quiver GIT quotient in the separated Simpson scheme. Indeed the affine invariant quotient is a point, since scaling one vertex contracts the representation space to zero. Thus the complementary space $U$ is a proper stable locus.

\smallskip\noindent\emph{4. Integrating the Hall logarithm.}
If $\delta_{nw}$ denotes the semistable Hall element of class $nw$, the term in class $2w$ in its logarithm is
\[
                \epsilon_{2w}=\delta_{2w}-\tfrac12\delta_w*\delta_w .
\]
On $\mathcal K$ a flag is exactly a quiver subrepresentation of dimensions $(1,1)$ with semistable factors. In one direction this is the snake lemma for the injective matrices. In the other direction $H^1(\OO_{D_1})=0$ makes sections exact in the extension, and taking the kernels of the evaluations recovers the source flag. The same two sequences of bundles in families identify the flag stacks over $\mathcal K$.

Put $u=\mathbb L$ and $g(u)=[\mathrm{GL}_2]^2=((u^2-1)(u^2-u))^2$. For a nonzero tuple the common kernel line, if present, is unique, as is the common image line. With one line fixed, the parameter space has dimension $5\cdot2$; with both fixed it has dimension $5\cdot1$. Removing the zero tuple in each case and applying inclusion--exclusion gives
\[
 V(u)=u^{20}-1-2(u+1)(u^{10}-1)+(u+1)^2(u^5-1).
\]
For tuples with a flag, choose the source and target lines, the two nonzero diagonal five-tuples and the arbitrary off-diagonal five-tuple. The bundles are Zariski locally trivial on the product of the flag varieties. Their total class is
\[
 W(u)=(u+1)^2(u^5-1)^2u^5.
\]
Dividing by the two frame groups and regularizing gives
\[
 \begin{aligned}
 B(u)&=\frac{(u-1)(V(u)-W(u)/2)}{g(u)}\\
 &=\frac{(u^4-u^3+u^2-u+1)(u^4+u^3+u^2+u+1)}{2(u+1)}\\[-2pt]
 &\hspace{19mm}\cdot(2u^6+4u^5+6u^4+7u^3+6u^2+4u+2).
 \end{aligned}
\]
The stack $\mathcal K$ is smooth of dimension $20-8=12$, so its intrinsic Behrend function is $+1$. The no-poles statement used in Theorem~\ref{thm:rank-zero-charge-conjugation} permits specialization of the regularized Hall logarithm at $u=1$. Its weighted stack integral is $B(1)=155/4$. With the sign of \textup{(\companionref{eq:Hallloggeneral})} this gives $\overline\Omega_{\mathcal K}(2w)=-155/4$. Subtracting the divisibility term $5/4$ gives $-40$, in agreement with the quiver formula of \cite[Theorem~5.2]{Reineke22}.

Every object in $\mathcal U$ is stable, so the same-polynomial two-step flag correspondence is empty there. Integration on that component gives $\chi(U,\nu_U)$. Adding the two contributions proves \eqref{eq:Kronfull22}.
\end{proof}

\begin{corollary}\label{cor:strict22number}
Restrict the charge-$2w$ Hall logarithm to the locally closed locus of strictly semistable middle sheaves, using the Behrend function of the stack of sheaves on $X$. Its integrated contribution is $-155/4$.
\end{corollary}
\begin{proof}
The unframed stable matrix moduli is smooth of dimension
$5\cdot4-(2\cdot4-1)=13$, with Behrend function $-1$. Its ordinary compactly supported Euler characteristic is zero by \cite[Corollary~6.14]{Weist22}.

Hence its integrated contribution to the charge-$2w$ logarithm is zero. All strictly semistable objects lie in $\mathcal K$ by Theorem~\ref{thm:kronecker22}. Additivity of the Hall integration after stratifying by the middle object now gives the number.
\end{proof}

\subsection{Principal parts after fixing a higher-charge shadow}
For each $r$, the transformation problem has two steps. The first is to construct the correction to the holomorphic sheaf series and prove its modular transformation. The second is to identify the resulting function from the finite list of negative Fourier coefficients.

\begin{proposition}\label{prop:step22-shadowunique}
Write two completions in the form
\[
 \widehat f(\tau,\bar\tau)=f(\tau)+R(\tau,\bar\tau),\qquad
 \widehat g(\tau,\bar\tau)=g(\tau)+R(\tau,\bar\tau),
\]
with the same correction $R$. Assume $f$ and $g$ are holomorphic on the upper half-plane, their Fourier series have finite principal parts at infinity, and both completions satisfy \eqref{eq:generalcompletedT} and \eqref{eq:generalcompletedS}. Equality of all coefficients of $q^x$ for $x<0$ in each component implies $\widehat f=\widehat g$.

For $r=2$, the bound of Proposition~\ref{prop:generalvacuumbound} allows six negative exponents in each component. Up to the charge symmetry of Theorem~\ref{thm:rank-zero-charge-conjugation}, there are six components $\mu=0,1,2,3,4,5$, hence $36$ possible polar positions.
\end{proposition}
\begin{proof}
Put $u=f-g=\widehat f-\widehat g$. Cancellation of principal parts makes $u$ bounded at the cusp. The unitary multiplier and weight $-3/2$ make $y^{-3/4}\|u(\tau)\|$ invariant, with limit zero at infinity. If $u\ne0$, this function has a positive interior maximum. On its nonzero locus,
\[
 \Delta\log\big(y^{-3/4}\|u\|\big)
 =\frac{3}{4y^2}+\Delta\log\|u\|>0,
\]
since $u$ is holomorphic. This excludes such a maximum and gives $f=g$.
For $r=2$, $c_2/24=35/6$. The spectral-flow phase gives all candidate
exponents as
\[
 x=-35/6+\{\mu^2/20\}+j,\qquad j\in\mathbb Z_{\ge0}.
\]
Here $\{\cdot\}$ denotes fractional part; the six fractions are
$0,1/20,1/5,9/20,4/5,1/4$. All admit $j=0,\ldots,5$.
\end{proof}
At charge two we can determine which polar coefficients in \eqref{eq:charge2anomaly} require additional sheaf data. Write $A_{2,\mu}$ for its quadratic expression in $h^{\mathrm{DT}}_1$. The primitive lower bound is $-55/24$ and the Hurwitz series in \eqref{eq:HurwitzHecke5} have nonnegative exponents. Therefore
\[
            [q^x]A_{2,\mu}=0\qquad (x<-55/12).
\]
Among the $36$ allowed polar positions up to charge symmetry this applies to nine: $j=0$ for all six representatives, and $j=1$ for $\mu=0,1,2$. Indeed the test is $j+\{\mu^2/20\}<5/4$. At each such position the proposed formula assigns the full coefficient to $h^{(0)}_2$.

For $\mu=0$ and $x=-35/6$, Theorem~\ref{thm:hall22strip} gives $15$. The sheaf calculation varies tilt stability on the derived category of the smooth quintic and applies the wall-crossing formula of \cite[Theorem~1.1(ii)]{FeyzOSV} to $v_{2,0}$. This agrees with the most polar index in \cite[\S4.4, (4.4)]{Hecht22}.

At the endpoint $(\mu,x)=(5,-55/12)$, equality in the bound requires both primitive factors at $q^{-55/24}$. Proposition~\ref{prop:polar-from-specializations} gives this coefficient only in component zero, with value $5$. Since $G^{(5)}_0=\tfrac12(V_5\mathcal H)_0$, the two divisors of five in \eqref{eq:HurwitzHecke5} give
\[
 [q^0]G^{(5)}_0=\tfrac12\bigl(5\cdot5^{-2}\cdot5+5\bigr)H(0)
                =-\tfrac14 .
\]
The sign $(-1)^{5+5}$ is positive, so $[q^{-55/12}]A_{2,5}=-25/4$. Combining this with the geometric Hall calculation \eqref{eq:Kronfull22}, the charge-two prediction requires
\[
 [q^{-55/12}]h^{(0)}_{2,5}
       =\left(-\frac{155}{4}+\chi(U,\nu_U)\right)
          +\frac{25}{4}=\chi(U,\nu_U)-\frac{65}{2}.
\]
This computation distinguishes the three inputs to an induction in $r$: the lower-multiplicity series entering $A_2$, the stable sheaves on the proper space $U$, and the identity identifying the completed geometric series with a modular form. The uniqueness proposition applies once the last identity and the polar data are available.

\section{Charges and the normalized K3-fibration theorem}\label{sec:modularity}
Both degenerations involve the quartic relative class. For comparison with its pairing against the ruled threefold, we compute the absolute vertical series of the pencil in fiber degrees $1,2,3\bmod4$. In degree $0\bmod4$ the correction is expressed by point ideals on the fibers. Section~\ref{sec:twodegoperators} specifies the relative identity needed on the quintic.
\subsection{The quintic charge and its five sectors}
For $\ch(F)=(0,H,\beta,n)$ put
\[
 d=H\beta,\qquad q_1=-d,\qquad q_0=n+\frac{25}{12},\qquad
 \widehat q_0=q_0-\frac{q_1^2}{10},\qquad
 \mu=q_1-\frac52\pmod5.
\]
The degree-six Chern character $n$ and the Euler characteristic satisfy $\chi(F)=n+25/6$.  Integrality of the first difference of the Hilbert polynomial gives $d+5/2\in\mathbb Z$, so $\mu$ is an integral residue class.

\begin{proposition}\label{prop:spectralflownew}
Tensoring by $\OO_X(aH)$ acts by
\begin{equation}\label{eq:spectralflowfull}
 (q_1,q_0)\longmapsto
 \left(q_1-5a,\ q_0-aq_1+\frac52a^2\right).
\end{equation}
It preserves the sheaf obstruction theory, $\widehat q_0$, and $\mu$.  On a representative with $q_1=\mu+5/2$,
\begin{equation}\label{eq:exponentcongruence}
 -\widehat q_0\equiv-\frac{55}{24}+\frac{\mu^2}{10}+\frac\mu2
 \equiv-\frac{55}{24}+\frac{3\mu^2}{5}\pmod{\mathbb Z}.
\end{equation}
\end{proposition}
\begin{proof}
Apply Proposition~\ref{prop:generalspectralflow} with $r=1$. Its tensor equivalence identifies the obstruction theories and gives \eqref{eq:spectralflowfull}. Since $a_1=5$, \eqref{eq:generalTphase} becomes
\[
 -\widehat q_0=-\frac{55}{24}+\frac{\mu^2}{10}+\frac{\mu}{2}
                      +5-\chi(F),\qquad
 \frac{3\mu^2}{5}-\left(\frac{\mu^2}{10}+\frac{\mu}{2}\right)
                   =\frac{\mu(\mu-1)}{2}\in\mathbb Z.
\]
This gives \eqref{eq:exponentcongruence}.
\end{proof}

Choose one representative of each spectral-flow orbit and define the formal series
\begin{equation}\label{eq:hseriesnew}
 h^{\mathrm{DT}}_{1,\mu}(\tau)
 =\sum_{\widehat q_0}\Omega_\mu(\widehat q_0)q^{-\widehat q_0},
 \qquad q=e^{2\pi i\tau}.
\end{equation}
The desired S-duality statement is convergence to weakly holomorphic functions satisfying
\begin{align}
 h^{\mathrm{DT}}_{1,\mu}(\tau+1)&=e^{2\pi i(-55/24+3\mu^2/5)}h^{\mathrm{DT}}_{1,\mu}(\tau),\label{eq:Tmodnew}\\
 h^{\mathrm{DT}}_{1,\mu}(-1/\tau)&=\frac{(-i\tau)^{-3/2}}{\sqrt5}
 \sum_{\nu\bmod5}e^{-2\pi i\mu\nu/5}h^{\mathrm{DT}}_{1,\nu}(\tau).\label{eq:Smodnew}
\end{align}
The first equation is already a formal consequence of the exponent congruence.  The second and the analytic statement are not consequences of spectral flow.

Let $D_5=(\mathbb Z/5,Q_5)$ with $Q_5(\mu)=3\mu^2/5$.  Its bilinear form is $\mu\nu/5$.  The relevant numerator representation has matrices
\begin{equation}\label{eq:Weil5new}
 \rho_5(T)e_\mu=e(3\mu^2/5)e_\mu,\qquad
 \rho_5(S)e_\mu=-\frac1{\sqrt5}\sum_\nu e(-\mu\nu/5)e_\nu,
 \quad e(x)=e^{2\pi ix}.
\end{equation}

\begin{proposition}\label{prop:numeratorsuffices}
Suppose an independently constructed $f\in M_{26}(\rho_5)$ has coefficientwise expansion $f=\eta^{55}h^{\mathrm{DT}}_1$.  Then \eqref{eq:hseriesnew} converges and satisfies \eqref{eq:Tmodnew}--\eqref{eq:Smodnew} with a finite principal part.
\end{proposition}
\begin{proof}
The eta function is nonvanishing on the upper half-plane and has expansion
$\eta=q^{1/24}\prod_{m\ge1}(1-q^m)$.  Thus $f/\eta^{55}$ is holomorphic there, has the prescribed convergent Fourier expansion, and has only a finite pole at the cusp.  The identities
\[
 \eta(\tau+1)=e(1/24)\eta(\tau),\qquad
 \eta(-1/\tau)=(-i\tau)^{1/2}\eta(\tau)
\]
and \eqref{eq:Weil5new} give the $T$ equation.  For $S$, the numerator has scalar factor $-\tau^{26}$ times the positive Fourier matrix.  Since $(-i)^{26}=-1$,
\[
 \frac{-\tau^{26}}{(-i\tau)^{55/2}}=(-i\tau)^{-3/2}.
\]
This proves the displayed $S$ equation and the weight $26-55/2=-3/2$.
\end{proof}

\subsection{The four-sector input from the quartic pencil}
Use the polarization $L=A+F_{\mathrm{fib}}$ of Proposition~\ref{prop:pencilcomplete}.  Its square on a fiber is $4$.  For a fiber sheaf let
\[
\begin{aligned}
 P_{d,c}(m)&=2m^2+dm+c,\\
 K_a^{\mathrm{raw}}(\tau)&=q^{1+d^2/8}\sum_c\dt_Y(d,c)q^{-c},
 \qquad a=d\pmod4.
\end{aligned}
\]
Here $\dt_Y(d,c)$ includes all fiber-sheaf characters with the indicated Hilbert polynomial.  Tensoring by $L$ sends
\[
 (d,c)\longmapsto(d+4,c+d+2),\qquad
 1+\frac{(d+4)^2}{8}-(c+d+2)=1+\frac{d^2}{8}-c.
\]
We may choose one $d$ per residue $a\bmod4$, in particular $d=0$ for $a=0$. The correspondence of Section~\ref{sec:twodegoperators} must carry these four components to the five quintic components $\mu\bmod5$.

The nodal K3-fibration comparison in \cite[Theorem~2 and the definition of $Z(X,q)$ in \S1.4]{GS} involves both a period series and a convention for the zero class. To compare it with the sheaf theory of our second degeneration, put $\Delta=\eta^{24}$. Let $p:\widetilde C\to\PP^1$ be the degree-two cover used for simultaneous resolution of the nodal fibers, and let $\widetilde\kappa:\widetilde Y\to\widetilde C$ be the resulting smooth quasi-polarized K3 family. Define
\begin{equation}\label{eq:Phihatfull}
 \widehat\Phi=\tfrac12\Phi^{\widetilde\kappa}_{\rm NL}
 \in M_{21/2}(\rho_4),\qquad D_4=(\mathbb Z/4,a^2/8),
\end{equation}
using the period normalization of \cite[\S4.4, Corollary~3]{MPNL}. Write
\begin{equation}\label{eq:Kcorrectednew}
 K^{\rm per}=\Delta^{-1}\widehat\Phi .
\end{equation}
The raw series above is defined by the obstruction theory on vertical sheaves. The coefficient at $q^{-1}$ fixes the normalization of the summand with trivial hull on a fiber.

\begin{proposition}[Vacuum normalization of the pencil]\label{prop:pencil-hodge-audit}
For the hyperplane pencil on the quartic threefold, the two independently defined series satisfy
\[
 [q^{-1}]K^{\rm raw}_0=[q^{-1}]K^{\rm per}_0=-1.
\]
In particular $[q^{-1}](K^{\rm raw}+e_0/\Delta)_0=0$.
\end{proposition}
\begin{proof}
At $d=0$ the exponent is $1-c$. For $c\ge2$ a degree-zero hull $R$ satisfies $\chi(R)\ge c$, so Theorem~\ref{thm:K3structurecomparison} forces $R=\OO$. Then $c=2-\length(R/V)$ gives $c=2$, $V=R$. The entire coefficient at $q^{-1}$ is consequently $\deg R^2\kappa_*\OO_Y=\deg\OO_{\PP^1}(-1)=-1$.

For the period form let $\lambda_{\widetilde C}=\widetilde\kappa_*\omega_{\widetilde Y/\widetilde C}$ be its Hodge line bundle. The relative dualizing sheaf of $Y/\PP^1$ is $\kappa^*\OO(1)$; the morphism is flat Gorenstein, so this formula pulls back to the branched cover. The local threefold singularities on the cover are ordinary double points, and the simultaneous small resolution is crepant: it is an isomorphism in codimension one between Gorenstein spaces. As the fibers are connected, the projection formula yields
\[
 \lambda_{\widetilde C}=p^*\OO_{\PP^1}(1),\qquad \deg\lambda_{\widetilde C}=2 .
\]
In the discriminant-zero convention of \cite[\S4.4, immediately before Corollary~3]{MPNL}, the neutral constant coefficient of the Noether--Lefschetz form is the degree of the dual Hodge line bundle. Hence $[q^0]\widehat\Phi_0=-2/2=-1$. Since $\Delta=q+O(q^2)$, this proves the second equality. The last assertion follows by adding the principal coefficient of $1/\Delta$.
\end{proof}

\begin{theorem}[Absolute vertical series in nonzero degree residues]
\label{thm:pencil-all-charge-normalization}
For the very general quartic hyperplane pencil, for $a=1,2,3$,
\[
 K^{\rm raw}_a=K^{\rm per}_a.
\]
In particular if $d=a$, and $b_N=e(\Hilb^N(\mathrm{K3}))$ for $N\ge0$, $b_N=0$ otherwise, then
\begin{equation}\label{eq:pencil-all-charge-convolution}
 \dt_Y(d,c)=\frac12\sum_{h\le1+d^2/8}
    \mathrm{NL}^{\widetilde\kappa}_{h,d}\,b_{h+1-c}.
\end{equation}
\end{theorem}
\begin{proof}
The general hyperplane fiber has Picard group generated by the quartic polarization $l$. At each nodal fiber $S_b$, its minimal resolution $j:\widetilde S_b\to S_b$ has
\[
 \Pic(\widetilde S_b)=\mathbb Z l\oplus\mathbb Z e,
 \quad l^2=4,\quad l\cdot e=0,\quad e^2=-2,
\]
where $e$ is the exceptional curve and $l$ is the pullback of the hyperplane class. This is the generic one-node quartic lattice of \cite[\S5.6, (27)]{MPNL}. Its genericity also applies to a hyperplane pencil on a very general quartic threefold. For a fixed hyperplane, restriction of quartic equations from $\PP^4$ onto it is surjective, and a general extension of a one-node quartic defines a smooth threefold: the normal derivative can be chosen nonzero at the node and Bertini applies elsewhere.

The universal incidence of such hyperplanes and threefolds is irreducible, by marking the unique node and imposing vanishing of the first jet of the restricted equation. Every additional Picard condition has positive codimension on that incidence.

For a very general threefold, these countably many loci have codimension at least two in its space of hyperplanes (the one-node locus itself is a divisor); a very general line avoids all of them. The identical restriction argument on the smooth locus shows that the generic fiber has Picard rank one.

A pure multiplicity-one vertical sheaf is the pushforward of a rank-one torsion-free sheaf $V$ on its integral fiber. Its Hilbert coefficient $d$ is the degree of its divisorial class. On a nodal fiber take the Weil divisor of $V^{**}$ and its strict transform $B$ on $\widetilde S_b$. Restriction to a smooth hyperplane curve avoiding the node and the finite defect of $V$ gives
\[
 d=\deg(V|_{l})=l\cdot B\in4\mathbb Z .
\]
Indeed the hyperplane section is a plane quartic of genus three, and $P_V(m)-P_V(m-1)=4m+d-2$, which yields the degree in this equality. The same divisibility holds on the geometric generic fiber. Thus for $4\nmid d$ the sheaves in a fixed polynomial occur only on smooth fibers in isolated Noether--Lefschetz loci. The resolved period family likewise has no classes of this degree above any singular member.

On a smooth such fiber $S_t$ write $V=M\otimes I_Z$, $\xi=c_1(M)$, $\xi^2=2h-2$, $l\cdot\xi=d$. Riemann--Roch and the Hodge index theorem give
\[
 c=\chi(V)=2+\xi^2/2-\length Z=h+1-\length Z,\qquad h\le1+d^2/8.
\]
Fix the class $\xi$ and put $N=h+1-c$. Since $\Pic^0(S_t)=0$, its moduli on $S_t$ is $S_t^{[N]}$, of dimension $2N$. At a transverse period intersection the type-II formula gives $\int_{S_t^{[N]}}c_{2N}(T_{S_t^{[N]}})=b_N$. Local deformation of an isolated intersection of multiplicity $m$ gives $m b_N$ \cite[Propositions~2.13--2.14 and proof of Theorem~1, (11)]{GS}.

The ordinary moduli is proper because every pure multiplicity-one fiber sheaf is Gieseker-stable and the semistable moduli is projective. There are no contributions at the omitted nodal fibers.

Over the smooth support locus the simultaneous cover has two unramified sheets, on which the isolated class and its period multiplicity are unchanged. The sum of these multiplicities on $Y$ is therefore $\mathrm{NL}^{\widetilde\kappa}_{h,d}/2$, with its quasi-polarized degree label $d$ as in \cite[\S1.6, (8), and \S5.6]{MPNL}. This proves \eqref{eq:pencil-all-charge-convolution}. The factor $\Delta^{-1}=\sum_{N\ge0}b_Nq^{N-1}$ converts it to the series identity since $1+d^2/8-h+N-1=1+d^2/8-c$. Tensoring translates to all representatives of $a$.
\end{proof}

\begin{proposition}[One point in fiber degree zero on the absolute pencil]
\label{prop:pencil-neutral-onepoint}
For the pencil of Proposition~\ref{prop:pencilcomplete},
\[
 [q^0]K^{\rm raw}_0=\dt_Y(0,1)=84=[q^0]K^{\rm per}_0 .
\]
\end{proposition}
\begin{proof}
Write the sheaf as $i_*V$ on $S_b$, with $P_V(m)=2m^2+1$. A nonzero section would inject $\OO_{S_b}$ into $V$ with quotient of polynomial $-1$. Hence $h^0(V)=0$. Serre duality on the Gorenstein quartic gives
\[
 \dim\Hom(V,\OO_{S_b})=h^2(V)=1+h^1(V)>0.
\]
A nonzero map is injective by torsion-freeness. Its cokernel has polynomial $P_{\OO_{S_b}}-P_V=1$, so $V=I_{y/S_b}$. The ideal sequence gives $h^1(V)=0$ and the inclusion is unique up to scalar.

This identifies the moduli in families, including at a node. A flat family of simple multiplicity-one vertical sheaves has a relative supporting fiber: locally on the base of the pencil its parameter acts as a scalar on the family, by simplicity and cohomology and base change for the relative endomorphism complex. This cuts out the graph of the corresponding map of parameter schemes, as in Theorem~\ref{thm:outerbitrarymultiplicity}.

On its family of quartics $f:S_T\to T$, the sheaf $V_T$ has $Rf_*V_T=A[-2]$ for a line bundle $A$ on $T$, by the fiber cohomology calculation and base change. Adjunction gives $\omega_{S_T/T}=f^*B$ for a line bundle $B$ on $T$. Relative duality yields
\[
 Rf_*\mathbf R\mathcal Hom(V_T,\OO_{S_T})
   \simeq (Rf_* V_T)^\vee\otimes B^\vee[-2]
   =A^\vee\otimes B^\vee .
\]
A local generator gives a fiberwise injection $V_T\hookrightarrow\OO_{S_T}$ with $T$-flat cokernel. Its image is independent of the generator. The quotient defines a section of $S_T/T$, whose ideal recovers $V_T$ up to a line from $T$. Thus the moduli is $\Hilb^1(Y/\PP^1)=Y$ in families. Its universal sheaf $\mathcal U$ is the relative diagonal ideal on $Y\times_{\PP^1}Y$, pushed into $Y\times Y$.

Write $p_1,p_2$ for the space and moduli projections of $Y\times Y$. On the support of $\mathcal U$, one has $p_1^*K_Y=p_2^*K_Y$ since $K_Y=\kappa^*\OO(-1)$. Threefold Serre duality, the sheaf-moduli tangent identification, and \cite[Theorem~2.3]{GS} give
\[
 \mathcal Ext^2_{p_2}(\mathcal U,\mathcal U)
   =\big(\mathcal Ext^1_{p_2}(\mathcal U,\mathcal U\otimes p_1^*K_Y)\big)^\vee
   = (T_Y\otimes K_Y)^\vee .
\]
The obstruction theory has virtual dimension zero on the smooth threefold $Y$, so its virtual degree is the third Chern number of this obstruction bundle. Put $f=c_1(\kappa^*\OO(1))$. Then $f^2=0$, $\int_Y c_3(T_Y)=-60$ by the blowup computation of Section~\ref{sec:seconddegfull}, and $\int_Yf\,c_2(T_Y)=24$ by restriction to a smooth K3 fiber with trivial normal bundle. Hence
\[
 \int_Y c_3(\Omega_Y\otimes K_Y^\vee)
       =\int_Y\{-c_3(T_Y)+f\,c_2(T_Y)\}=84 .
\]

For the period coefficient write $\widehat\Phi_0=-1+a_1q+O(q^2)$. A square-$(-2)$ class of degree zero on a smooth quartic contradicts ampleness, since K3 Riemann--Roch makes one sign effective. At each node the exceptional class $e$ and $-e$ both have period multiplicity one on the double cover. Indeed $x_i=\sqrt t\,u_i$ on $\sum_i x_i^2=t$ gives
\[
  \int_{\delta_t}\operatorname{Res}\frac{dx_1\,dx_2\,dx_3}{\sum x_i^2-t}
    =\sqrt{t}\int_{\delta_1}\operatorname{Res}\frac{du_1\,du_2\,du_3}{\sum u_i^2-1},
\]
with nonzero last integral, as in \cite[\S5.6, (27) and the proof of Corollary~2]{MPNL}. Therefore $a_1=\tfrac12(2\cdot108)=108$. Multiplying by $\Delta^{-1}=q^{-1}+24+O(q)$ gives $[q^0]K^{\rm per}_0=108-24=84$.
\end{proof}

Thus $c_Y=K^{\rm per}-K^{\rm raw}$ vanishes off $a=0$. In residue zero, $(c_Y)_0=O(q)$ by Propositions~\ref{prop:pencil-hodge-audit} and~\ref{prop:pencil-neutral-onepoint}. The exceptional classes $\pm je$ enter at the squares $j^2$.

\begin{figure}[htbp]\centering
\begin{tikzpicture}[font=\footnotesize,x=1cm,y=1cm,line join=round,line cap=round]
\node[font=\small\bfseries] at (2.55,3.07) {Resolved quartic};
\node[font=\small\bfseries] at (9.55,3.07) {Nodal quartic};
\path[fill=qIIteal!9] (.60,.61) .. controls (1.22,1.53) and (2.00,1.20)..(2.58,1.20)
 .. controls (3.22,1.20) and (3.93,1.53)..(4.55,.61)
 -- (4.55,2.56) .. controls (3.90,1.74) and (3.20,1.94)..(2.58,1.94)
 .. controls (1.93,1.94) and (1.20,1.74)..(.60,2.56)--cycle;
\draw[qIIteal] (.60,.61)..controls (1.22,1.53) and (1.96,1.20)..(2.58,1.20)..controls (3.20,1.20) and (3.93,1.53)..(4.55,.61);
\draw[qIIteal] (.60,2.56)..controls (1.22,1.74) and (1.96,1.94)..(2.58,1.94)..controls (3.20,1.94) and (3.93,1.74)..(4.55,2.56);
\draw[qIIclay,line width=1.2pt] (2.58,1.57) ellipse (.26 and .37);
\node[text=qIIclay] at (2.60,.70) {$e\simeq\PP^1$};
\draw[-{Stealth[length=6pt]},qIIblue,line width=.9pt] (5.10,1.58)--(7.27,1.58);
\node[text=qIIblue,above] at (6.18,1.68) {contraction};
\path[draw=qIIteal,fill=qIIteal!9] (7.69,.66)..controls (8.48,1.34) and (9.04,1.46)..(9.57,1.57)
 ..controls (9.05,1.70) and (8.48,1.83)..(7.69,2.53);
\path[draw=qIIteal,fill=qIIteal!9] (11.45,.66)..controls (10.66,1.34) and (10.10,1.46)..(9.57,1.57)
 ..controls (10.08,1.70) and (10.66,1.83)..(11.45,2.53);
\fill[qIIclay] (9.57,1.57) circle (2.6pt);
\node[text=qIIclay,below] at (9.57,1.37) {$p$};
\node at (2.58,.05) {$h\cdot e=0,\quad e^2=-2$};
\node at (9.57,.05) {$xy=z^2$ at $p$};
\end{tikzpicture}
\caption{The exceptional class has degree zero against the pulled-back quartic polarization. Its signed multiples $\pm je$ have square $-2j^2$, accounting for the exponents $j^2$ in \eqref{eq:pencil-neutral-reduction}. The drawing denotes the contraction on one fiber.}\label{fig:II-node-resolution}
\end{figure}
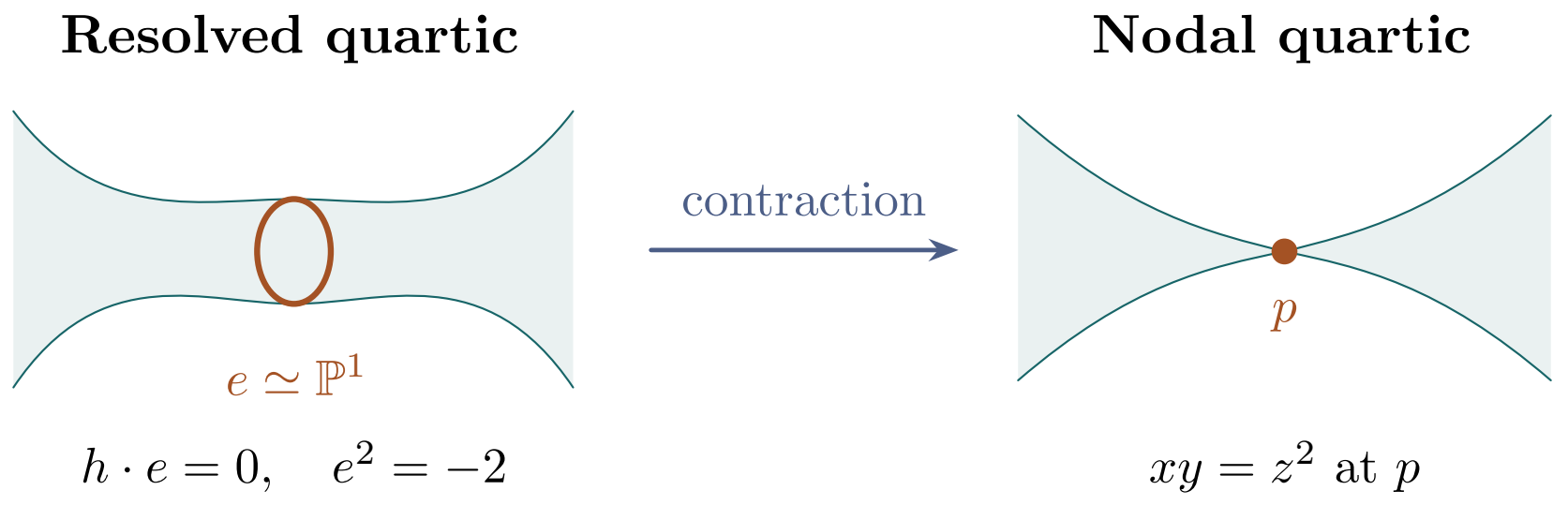

Figure~\ref{fig:II-node-resolution} records the exceptional degree entering the neutral period series.

\begin{corollary}[The series of point ideals in fiber degree zero]\label{cor:pencil-ideal-reduction}
For $N\ge0$ let $M^{\rm id}_N$ be the open and closed part of the absolute vertical sheaf moduli of $P_{0,2-N}$ whose geometric sheaves have trivial reflexive hull on their supporting fiber. Endow it with the full induced scheme structure and threefold obstruction theory. These loci exist, and with
\[
 \mathcal I(q)=\sum_{N\ge0}q^{N-1}\deg[M^{\rm id}_N]^{\vir},
 \qquad
 \Theta_+(q)=\sum_{j\ge1}q^{j^2},
\]
the exact remaining correction is
\begin{equation}\label{eq:pencil-neutral-reduction}
 (K^{\rm per}-K^{\rm raw})_0
       =\frac{-1+108\Theta_+(q)}{\Delta(q)}-\mathcal I(q).
\end{equation}
\end{corollary}
\begin{proof}
Fix $c=2-N$. On a nodal fiber with resolution $\widetilde S\to S$, strict transform of Weil divisors gives
\[
          \mathbb Z e\ \longrightarrow\ \Pic(\widetilde S)
                  \ \longrightarrow\ \operatorname{Cl}(S)\ \longrightarrow\ 0 .
\]
Indeed the resolution is an isomorphism off $e$, and principal divisors come from the common function field. The lattice computed above gives $\operatorname{Cl}(S)=\mathbb Z l$, so a degree-zero rank-one sheaf $V$ has $V^{**}=\OO_S$. On a smooth fiber put $M=V^{**}$ and $\xi=c_1(M)$. For $\xi\ne0$ the Hodge index theorem and Riemann--Roch give
\[
 l\cdot\xi=0,\quad \xi^2=2h-2<0,\quad
 V=M\otimes I_Z,\qquad \length Z=h+1-c\ge0.
\]
These inequalities give $c-1\le h\le0$. For fixed $h$ the complete resolved period curve meets the corresponding Noether--Lefschetz divisor in finitely many points, since the geometric generic Picard group is $\mathbb Z l$. At each fiber negative definiteness of $l^\perp\cap\Pic(S)$ gives finitely many classes of the fixed norm. For $\xi=0$ one has $V=I_Z$, $\length Z=N$.

Write $M=M(P_{0,c})$. For $\xi\ne0$ let $H_{t,\xi}=\Hilb^{h+1-c}(S_t)$ and $M_0$ be the image of the relative Hilbert scheme of $N$ points. Properness gives a disjoint closed decomposition on underlying spaces:
\[
 |M|=|M_0|\ \sqcup\!\bigsqcup_{(t,\xi)}\operatorname{im}(|H_{t,\xi}|\to|M|).
\]
Here the second map sends $Z$ to $M_\xi\otimes I_Z$ on $S_t$, with $c_1(M_\xi)=\xi$. The union is finite by the preceding bounds. Restricting the structure sheaf and the obstruction theory of $M$ to the open and closed subset $|M_0|$ defines $M_N^{\rm id}$ with its full scheme structure.

Theorem~\ref{thm:pencil-all-charge-normalization} gives the same multiple of $b_{h+1-c}$ for each nonzero class on a smooth fiber. The remaining period terms have constant coefficient $-1$ and nodal classes $\pm je$ with $1-h=j^2$, $j\ge1$. Multiplication by $\pm j$ preserves the vanishing order of the period. Dividing by the degree of the cover gives
\[
 -1+\tfrac12\sum_{t\ {\rm nodal}}\sum_{j\ge1}2q^{j^2}
       =-1+108\Theta_+(q)
\]
as in \cite[\S5.6, proof of Corollary~2]{MPNL}. Dividing by $\Delta$ and subtracting the remaining sheaf series $\mathcal I(q)$ proves \eqref{eq:pencil-neutral-reduction}.
\end{proof}

The local construction in \cite[Remark~2.15 and \S2.3]{GS} requires a polarization near each resolved fiber. If $L$ were such a line bundle on a simultaneous resolution over a disk, flatness of $c_1(L)$ and $\Pic(S_t)=\mathbb Z l$ at a very general $t$ would give
\[
                c_1(L_0)\cdot e=ml\cdot e=0 .
\]
Relative ampleness requires $c_1(L_0)\cdot e>0$. In \eqref{eq:pencil-neutral-reduction} only the periods are computed on the resolution. The obstruction theory on $M_N^{\rm id}$ comes from the original projective pencil.

\begin{proposition}[Two points on the degree-zero component of the pencil]
\label{prop:pencil-ideal-two}
For the very general Lefschetz pencil, the morphism
$H=\Hilb^2(Y/\PP^1)\to M^{\rm id}_2$ defined by the universal ideal is an isomorphism. This component has virtual degree $2268$. In particular $(K^{\rm per}-K^{\rm raw})_0=O(q^2)$.
\end{proposition}
\begin{proof}
First $H$ is smooth of dimension five. At two distinct points the equation in an ordered local chart is $f(y_1)=f(y_2)$, with $f$ a local pencil parameter. It is smooth since at most one of the points can be critical.

At a length-two jet in a smooth threefold use affine coordinates $x\in\CC^3$, a direction $u=(1,u_2,u_3)$ in an affine chart of $\PP^2$ and a parameter $t$. The universal subscheme is given by $x+su$ with $s^2=t$. These are coordinates on the Hilbert square: on the open where the first coordinate generates its rank-two coordinate algebra, subtract half its trace to obtain $s$, then express the other two coordinates in the basis $(1,s)$. The requirement that $f$ be scalar on this algebra is
\[
 \frac{f(x+su)-f(x-su)}{2s}=0 ,
\]
The quotient is even in $s$. At a critical point take $f$ quadratic by the Morse lemma. The derivative of the defining equation in $x$ is $\operatorname{Hess}(f)(-,u)\ne0$. At a regular point take $f$ itself as coordinate. This proves smoothness. At a double point on the node the map to the base and its differential are
\[
   (x,u,t)\longmapsto f(x)+t f(u),\qquad d\pi_H|_{x=t=0}=f(u)\,dt.
\]
It is submersive there precisely when $f(u)\ne0$.

At two distinct points including the node the first derivative of the base value is zero. In the remaining cases the support lies in the smooth locus of $f$, and the differential surjects.

We compare tangents to the fixed-fiber moduli. Write $S=Y_b$, $J=I_{Z/S}$ and $Q=\OO_Z$. Since $S$ is normal and $J=\OO_S$ off finitely many points, $\Hom_S(J,\OO_S)=\CC$. Applying $\Hom_S(J,-)$ to the ideal sequence thus gives
\[
 0\to \Hom_S(J,Q)\to\Ext^1_S(J,J)\to
 \Ext^1_S(J,\OO_S)\to\Ext^1_S(J,Q).
\]
The last arrow is injective. Indeed applying $\Hom_S(-,\OO_S)$ to the ideal sequence and using $H^1(\OO_S)=0$ identifies $\Ext^1_S(J,\OO_S)$ with the kernel $K$ of $\Ext^2_S(Q,\OO_S)\to H^2(\OO_S)$. By Serre duality this kernel is the functionals in $H^0(Q)^*$ vanishing on $1$. Naturality of connecting homomorphisms for the maps from $\OO_S$ to $Q$ shows that it suffices to check injectivity of
\[
 K\lhook\joinrel\longrightarrow\Ext^2_S(Q,\OO_S)\longrightarrow\Ext^2_S(Q,Q).
\]
Both global Ext groups in this display are direct sums of the local Ext groups at the finite support of $Q$, by the local-to-global spectral sequence. At a complete intersection of codimension two the whole latter map is injective by the Koszul resolution: the induced map on top Ext is the identity of $Q$.

Here are the local calculations in the two remaining cases. At a reduced point of $Z$ at the node write $R=\CC[[x,y,z]]/(xy-z^2)$ and $Q=k$. At a length-two jet at the node for which $f(u)\ne0$, the ideal in $R$ is generated by two independent linear forms vanishing on the direction, since restriction of $f$ to that axis is a nonzero multiple of its coordinate squared.

It is a complete intersection. For an isotropic direction choose linear coordinates with that direction the $x$-axis, so $Q=R/(x^2,y,z)$. The reduced and isotropic cases have the uniform description
\[
 Q=R/(x^a,y,z),\quad a=1,2.
\]
Over $A=\CC[[x,y,z]]$ resolve $Q$ by the Koszul complex on $g=(x^a,y,z)$. The identity $xy-z^2=(0,x,-z)\cdot g$ gives the degree-one homotopy of multiplication by $xy-z^2$ via wedge with $xe_2-ze_3$. By \cite[Theorem~7.2]{EisenbudHyp22}, over $R$ its resolution is the Koszul exterior algebra on $e_1,e_2,e_3$ with an additional divided-power generator $t$ of degree two and $\partial t=xe_2-ze_3$. In degree two define a cochain $\varphi$ over $R$ by
\[
 \varphi(e_1e_2)=-z,\qquad \varphi(e_1e_3)=-x,\qquad
 \varphi(e_2e_3)=\varphi(t)=0 .
\]
It is closed: the four values of $d^*\varphi$ on $e_1e_2e_3,te_1,te_2,te_3$ are $xy-z^2,xz-xz,0,0$ (up to the choice of Koszul signs). For $a=2$ its value on $e_1e_3$ modulo $Q$ is $-x$, and every boundary has zero value on that pair modulo $Q$ since $g=0$ there.

By local duality $\Ext^2_R(Q,R)=Q^*$ is cyclic over $Q$ of length $a$. For $a=2$ the one-dimensional subspace in the global kernel $\Ext^2_S(Q,\OO_S)\to H^2(\OO_S)$ is the functionals vanishing on $1$. The functional with value one on $x$ and zero on $1$ spans this line and generates $Q^*$ as a $Q$-module.

We show that the cocycle $\varphi$ represents a generator by proving that it does not belong to $x\Ext^2_R(Q,R)$. Any cocycle $(A,B,C,D)$, with coordinates its values on $(e_1e_2,e_1e_3,e_2e_3,t)$, satisfies $x^2 C-yB+zA=0$ in $R$. Comparing the terms of degree one forces $A,B$ to have vanishing constant terms. Its reduction evaluated on $e_1e_3$ is thus divisible by $x$. As boundaries reduce to zero on this pair, the image there of $x$ times \emph{any} class is zero.

Our displayed cocycle has nonzero image there, so is a generator. The trace-zero generator differs from a nonzero constant multiple by a socle element, whose image on the same pair is zero. Its image is therefore nonzero as required. For a reduced node paired with a distinct regular point, local duality identifies the global trace kernel with the anti-diagonal in the duals of the two residue fields. Its projection to the regular point is nonzero, where the Koszul calculation applies.

In particular the asserted fixed-fiber tangent comparison holds in all cases; local injectivity on the entire dual of $Q$ is only needed in the complete-intersection cases.

A first-order sheaf deformation induces motion of its supporting fiber by Theorem~\ref{thm:outerbitrarymultiplicity}. Where $H\to\PP^1$ is submersive, the fixed-fiber tangent isomorphism and the common base direction identify the total tangents.

In the remaining cases $Z$ contains a reduced point at the node or an isotropic double there. The displayed resolution of $Q$ is minimal and infinite, hence $J$ has infinite projective dimension over $R=P/(f)$, $P=\CC[[x,y,z]]$, $f=xy-z^2$. A lift in a nonzero base direction would be flat over $\mathbb C[\epsilon]/(\epsilon^2)$ with $\epsilon=f$ in $P/(f^2)$. If $(D_\bullet,d)$ is the minimal $R$-resolution, the obstruction operator on the lifted differentials is
\[
 t_i=\frac{\widetilde d_{i-1}\widetilde d_i}{f}\bmod f:
                   D_i\longrightarrow D_{i-2}.
\]
A flat lift supplies a free resolution over $P/(f^2)$ reducing to a resolution over $R$, so $[t]=0$ up to chain homotopy. But \cite[\S1 and Theorem~6.1(i),(iii)]{EisenbudHyp22} identifies $t$ on the periodic tail with an isomorphism. After tensoring the minimal resolution with $k$, the differential vanishes and that isomorphism survives. Thus the transverse tangent is zero. Both moduli have the same fixed-fiber tangent in this case as well.

The Hilbert morphism is proper and bijective on geometric points, hence finite. At $m\in M_2^{\rm id}$ the dimension and tangent calculations give $\dim\OO_{M_2^{\rm id},m}=5=\dim T_m M_2^{\rm id}$ for each closed $m$. Thus the target is regular. The finite map is birational on each component in characteristic zero, hence an isomorphism by normality.

By Proposition~\ref{prop:pencil-neutral-onepoint}, the obstruction bundle is $\Omega_H\otimes\pi_H^*\OO(1)$ for $\pi_H:H\to\PP^1$. Let $S$ be a one-nodal quartic. Off the diagonal in $S^{[2]}$ we count unordered pairs; over the diagonal the punctual fiber is $\PP(T_pS)$. Thus
\[
 e(S^{[2]})=\binom{23}{2}+2(23-1)+3=300,\qquad
 e(H)=2\cdot324+108(300-324)=-1944 .
\]
Here $324$ is the Euler number on a smooth fiber.

With $f=\pi_H^*c_1(\OO(1))$, $f^2=0$ and trivial normal line of a smooth fiber, we obtain
\[
 \int_H c_5(\Omega_H\otimes\pi_H^*\OO(1))
    =-e(H)+\int_H f c_4(T_H)=1944+324=2268 .
\]
This is also $[q^1](-1+108\Theta_+)/\Delta=-324+108\cdot24$. Equation~\eqref{eq:pencil-neutral-reduction} finishes the proof.
\end{proof}

Write $\mathcal I_N=\deg[M_N^{\rm id}]^{\vir}$. The expected values in \eqref{eq:pencil-neutral-reduction} through $N=3$ are
\[
\begin{array}{c|rrrr}
N & 0&1&2&3\\ \hline
[q^{N-1}]\,(-1+108\Theta_+)/\Delta & -1 & 84 &2268&31792
\end{array}
\]
because $\Delta^{-1}=q^{-1}+24+324q+3200q^2+O(q^3)$.
For $N=3$ this is $108\cdot324-3200=31792$. We compute the virtual degree on the relative Hilbert scheme, including its nodal fibers.

\begin{proposition}[The length-three virtual degree on the quartic pencil]\label{prop:pencil-ideal-three-m27}
Let $H_3=\Hilb^3(Y/\PP^1)$ and let $\pi_3:H_3\to\PP^1$ send a subscheme to its supporting fiber. The universal ideal gives an isomorphism $H_3\simeq M^{\rm id}_3$. This is a smooth projective sevenfold and
\[
             \mathcal I_3=\int_{H_3}c_7\big(\Omega_{H_3}\otimes\pi_3^*\OO(1)\big)
                      =31792 .
\]
\end{proposition}
\begin{proof}
Work near an affine pencil coordinate $f$. On $\Hilb^3(Y)$ the tautological function defines
\[
 H_3=Z(\bar f),\qquad
 \bar f\in\Gamma\big(\OO_Y^{[3]}/\OO\big),\qquad
 \operatorname{rk}(\OO_Y^{[3]}/\OO)=2 .
\]
Away from the node the relative Hilbert scheme is smooth over the pencil. At disjoint support the local Hilbert functor is the product of the punctual deformation functors over the base. There is at most one node on each fiber, so only the factor supported there requires a calculation. Proposition~\ref{prop:pencil-ideal-two} covers its lengths one and two.

For a curvilinear triple choose $P=\CC[[t,y,z]]$, $I=(y,z,t^3)$ and $A=P/I$. The tangent space $\Hom_P(I,A)=A^3$ has dimension nine. Smoothing three points along the axis gives local dimension at least nine, hence smoothness of the ambient Hilbert scheme. At the node $df=0$ and $f(t,0,0)=O(t^3)$ imply $\operatorname{Hess}_f(\partial_t,\partial_t)=0$. Choose $y$ so that $\operatorname{Hess}_f(\partial_t,\partial_y)\ne0$. Then
\[
  (\partial_y f)|_A=t\cdot u,\qquad u\in A^\times,\qquad
  (\partial_y f)A=(t)\xrightarrow{\sim}A/\mathbb C .
\]
Moving the generator $y$ realizes this surjection as the differential of $\bar f$.

For a noncurvilinear triple choose Morse coordinates $x,y,z$ with its tangent plane $W$ given by $z=0$. Its ideal and algebra are
\[
 I=(z,x^2,xy,y^2),\qquad A=\CC\oplus\mathfrak m,qquad
 \mathfrak m=\langle x,y\rangle,\quad\mathfrak m^2=0.
\]
The two syzygies $y\cdot x^2=x\cdot xy$ and $y\cdot xy=x\cdot y^2$ force the images of the three quadrics in $A$ to lie in $\mathfrak m$. There are no further constraints once $z$ has arbitrary image in $A$. Thus
\[
 \Hom_P(I,A)=A\oplus\mathfrak m^{\oplus3},\qquad
 \dim\Hom_P(I,A)=3+6=9 .
\]
The family $(z,x^2-sx,xy,y^2-sy)$ smooths $Z$ to three distinct points. The nondegenerate quadratic $f$ cannot vanish on a two-plane, so $f|_{W}$ has a nonzero coefficient of $x^2,xy$ or $y^2$. Varying that generator in $\mathfrak m$ makes $d\bar f$ onto $A/\mathbb C$, proving smoothness of dimension seven.

\smallskip\noindent\emph{Euler characteristics of the fibers.}
For a nodal K3 fiber $S_0$ with node $o$, $e(S_0\setminus\{o\})=22$. Write $H^j_o$ for the punctual locus of length $j$. A double point has a tangent direction in $\PP^2$. A curvilinear triple has an isotropic tangent on a smooth conic.

Its remaining transverse second derivative gives an affine plane locally over that conic, since a change of parameter on $\mathbb C[t]/(t^3)$ removes the tangential second derivative. The equation $f=0$ imposes only isotropy. The complementary square-zero triple is determined by its plane. Hence
\[
 e(H_o^0)=e(H_o^1)=1,\quad e(H_o^2)=e(\PP^2)=3,\quad
 e(H_o^3)=e(\PP^1)e(\mathbb A^2)+e(\operatorname{Gr}(2,3))=5.
\]

Stratify the full Hilbert scheme by the length supported at the node. For the smooth locus the Euler series to order three is $\prod_{m\ge1}(1-z^m)^{-22}$, by stratifying by support (equivalently the Euler product of G\"ottsche \cite{GottscheBetti}). The first four coefficients are $1,22,275,2530$. Therefore
\[
 e(\Hilb^3 S_0)=2530+275+3\cdot22+5=2876,\qquad
 e(\Hilb^3\mathrm{K3})=3200.
\]
Summing over the pencil with its $108$ nodal fibers gives $e(H_3)=2\cdot3200+108(2876-3200)=-28592$. Put $u=\pi_3^*c_1(\OO(1))$. Since $u^2=0$ and a smooth fiber has trivial normal bundle,
\[
 \int_{H_3}c_7(\Omega_{H_3}\otimes\pi_3^*\OO(1))
  =-e(H_3)+\int_{H_3}u\,c_6(T_{H_3})
  =28592+3200 .
\]

\smallskip\noindent\emph{Infinitesimal deformations on a fixed fiber.}
Write $S=Y_b$, $A=\OO_Z$ and $J=I_{Z/S}$. The universal ideal gives a proper morphism $H_3\to M_3^{\rm id}$, bijective on geometric points because $J^{**}=\OO_S$ and an embedding $J\to\OO_S$ is unique up to scalar. We compare the deformation functors at first order. On a fixed fiber the map of tangent spaces is the connecting homomorphism
\[
 \Hom_S(J,A)\lhook\joinrel\longrightarrow \Ext^1_S(J,J)
\]
of the ideal sequence. Its injectivity follows from $\Hom_S(J,\OO_S)=\CC$. As in Proposition~\ref{prop:pencil-ideal-two}, this is an isomorphism if the map
\[
 \operatorname{Ext}^2_S(A,\OO_S)\longrightarrow\operatorname{Ext}^2_S(A,A)
\]
is injective on $\ker(H^0(A)^*\xrightarrow{\lambda\mapsto\lambda(1)}\CC)$. Indeed $\Ext^1_S(J,\OO_S)$ injects into $\Ext^2_S(A,\OO_S)$ with precisely this image by $H^1(\OO_S)=0$ and Serre duality. The two connecting maps for $\OO_S\to A$ commute. Both Ext groups in the display decompose as sums of local groups over the support.

At a regular point, duality identifies the dual local map with $\operatorname{End}(A)\to A$, $\phi\mapsto\phi(1)$. Its inverse is $a\mapsto(\cdot a)$ since $A$ is a quotient algebra. This version also covers the noncurvilinear triple at a regular point. At a complete intersection on a nodal fiber the Koszul resolution identifies the displayed map with the identity of $A$. We treat the other punctual algebras at a node explicitly.

At the node put $R=\CC[[u,v]]^{\{\pm1\}}$, with $x=u^2$, $y=v^2$, $z=uv$, and let $M=Ru+Rv$ be the odd summand. For $\phi\in\Hom_R(M,R)$, the elements $p=\phi(u)$, $q=\phi(v)$ satisfy $vp=uq$. Thus $(p,q)=(u\alpha,v\alpha)$ for a unique odd series $\alpha$, identifying $\Hom_R(M,R)=M$ by multiplication.

A curvilinear quotient at the node of length $a\le3$ with isotropic tangent may be written, after a linear change preserving the quadric, with $x$ its parameter and $y,z$ vanishing to order at least two. For $a=3$, write $y=\alpha x^2$, $z=\beta x^2$ modulo $x^3$. Replacing $z$ by $z-\beta x^2$ and $y$ by $y-2\beta xz+\beta^2x^3$ preserves $xy-z^2$.

The ideal for $\alpha\ne0$ is the complete intersection $(z,y-\alpha x^2)$ in the new coordinates. For $a=2$ a nonisotropic tangent likewise gives a complete intersection, since the two linear equations of the axis cut out its double point.

For $a=1$, for an isotropic double, and for $a=3$, $\alpha=0$, the ideal is $(x^a,z,y)$. For a noncurvilinear triple its tangent plane meets the quadric in two distinct lines or in a double line. Up to linear isometry the plane is respectively $z=0$ or $y=0$.

The remaining ideals have the following presentations $R^3\twoheadrightarrow J$. Each ordered triple in the center gives its generators. In the last column, $b_1,b_2$ generate the module of relations in the sense that it equals $M b_1\oplus M b_2\subset R^3$.
\[
\begin{array}{c|c|c}
 &\text{ordered generators of }J&(b_1,b_2)\\ \hline
 A=\CC[x]/(x^a)&(u^{2a},uv,v^2)&((v,-u^{2a-1},0),(0,v,-u))\\
 W=(z=0)&(u^4,uv,v^4)&((v,-u^3,0),(0,v^3,-u))\\
 W=(y=0)&(u^4,u^3v,v^2)&((v,-u,0),(0,v,-u^3))
\end{array}
\]
To check the relations, fix a monomial in the image ideal. Its coefficient relations are generated by differences between two generators dividing that monomial. If the first and third divide it, so does the middle generator. It therefore suffices to compare consecutive generators. In each pair the common multiples give exactly the odd multiples of the listed column. The two submodules have zero intersection by their first and third coordinates. Grouping coefficients proves the assertion also for formal power series.

Dimension shifting gives $\Ext_R^2(A,R)=\Ext_R^1(J,R)$. An element of $\Hom_R(M b_1\oplus M b_2,R)$ is a pair of odd series $(p,q)$. Its four values on $u b_1,v b_1,u b_2,v b_2$ are $(up,vp,uq,vq)$. The image of the class in $\Ext_R^1(J,A)$ is obtained by reducing those values modulo $J$ and quotienting by the restrictions of maps $R^3\to A$. The following evaluations detect its rank.

\begin{itemize}
\item For $(x^a,z,y)$, the entries of $u b_1=(z,-x^a,0)$ vanish in $A$. Every map from $R^3$ vanishes on it modulo $J$. The pairs $(u x^k,0)$ for $0\le k\le a-2$ have values $x^{k+1}$ there. These give $a-1$ independent images.
\item For $W=(z=0)$, both $u b_1=(z,-x^2,0)$ and $v b_2=(0,y^2,-z)$ vanish modulo $J$. On these two vectors the pairs $(u,0),(0,v)$ have values $(x,0),(0,y)$.

The images are independent.
\item For $W=(y=0)$, use $v b_2$ and $v b_1+u b_2$. Both vanish modulo $J$, having respective entries $(0,y,-xz)$ and $(y,0,-x^2)$. Evaluating the pairs $(u,0),(0,u)$ gives $(0,z),(z,x)$, again independent.
\end{itemize}

By local duality $\Ext_R^2(A,R)=A^*$. The calculations bound the local kernel by one dimension. If nonzero it equals $\operatorname{Soc}(A^*)=(A/\mathfrak m_A)^*$ and pairs nontrivially with $1$. All regular-point maps are injective and there is at most one node. Hence the total kernel meets $\{\lambda:\lambda(1)=0\}$ trivially, proving the fixed-fiber tangent isomorphism.

\smallskip\noindent\emph{Moving the supporting fiber.}
A first-order deformation of a simple fiber sheaf induces a first-order motion of its support in the pencil, as in Proposition~\ref{prop:pencil-neutral-onepoint}. If $\pi_3$ is submersive at $Z$, the fixed-fiber tangent isomorphism and the common one-dimensional support direction give an isomorphism on total tangent spaces. This covers subschemes supported at regular points and the complete intersections at the node. For the nonisotropic double the axis restriction of $f$ has nonzero quadratic coefficient.

For the complete-intersection triple above it has nonzero cubic coefficient in adapted curvilinear coordinates. In either case, varying the constant term of the last curvilinear generator gives the required base direction.

For each other algebra in the table the first syzygy of $J$ is $M\oplus M$. The odd summand $M\subset\CC[[u,v]]$ has the regular sequence $u^2,v^2$, hence depth two. It has rank one and two minimal generators. Auslander--Buchsbaum would force a finite-projective-dimension $M$ to be free. Thus $J$ has infinite projective dimension over $R$.

At such a subscheme a nonzero base direction in the sheaf functor would give a flat lift of $J$ across the first-order smoothing $f=\epsilon$. The periodicity operator on its minimal resolution obstructs this lift by the argument in Proposition~\ref{prop:pencil-ideal-two} using \cite[\S1 and Theorem~6.1(i),(iii)]{EisenbudHyp22}. Consequently both total tangent spaces are the fixed-fiber tangent spaces. This exhausts all supports of length three.

Let $j:H_3\to M^{\rm id}_3$ be the universal-ideal morphism. It is finite by properness and quasi-finiteness, and bijective on geometric points. At a closed $m=j(Z)$,
\[
 \dim\OO_{M^{\rm id}_3,m}=\dim\OO_{H_3,Z}=7
                  =\dim T_ZH_3=\dim T_m M^{\rm id}_3 .
\]
Thus the target is regular. In characteristic zero $j$ is birational on each component, so normality gives an isomorphism. By relative duality as in Proposition~\ref{prop:pencil-neutral-onepoint}, its obstruction bundle is $\Omega_{H_3}\otimes\pi_3^*\OO(1)$.
The Chern number above therefore proves
\[
                \mathcal I_3=31792,\qquad (K^{\rm per}-K^{\rm raw})_0=O(q^3).
\]
\end{proof}
Equation~\eqref{eq:pencil-neutral-reduction} concerns the absolute pencil. The relative transfer is the additional hypothesis of Theorem~\ref{thm:twodegmodularityfull}.

\section{The primitive sheaf moduli for point and line defects}\label{sec:initial-checks}
For $i:T=X\cap\Pi\hookrightarrow X$ integral, Proposition~\companionref{prop:primitive} writes $F=i_*G$ with $G$ rank one and torsion free. The classification below includes singular hyperplane sections and all hulls in the charges used for specialization.

\subsection{Rank-one sheaves on the hyperplane sections}
Write $h=\mathcal O_T(1)$. The degree-five equation of $T$ gives
$\omega_T=h$, $\chi(\mathcal O_T)=5$, and $h^2=5$. For a rank-one
torsion-free sheaf $G$, define its degree by
\begin{equation}\label{eq:polar-surface-polynomial}
 \chi(G(th))=\frac52t^2+
 \left(e(G)-\frac52\right)t+\chi(G).
\end{equation}
For a general smooth plane section $C$ avoiding the nonfree locus, take first differences:
\[
 \chi(G(th))-\chi(G((t-1)h))=5t+e(G)-5=\chi(C,G|_C(th)).
\]
Since $g(C)=6$, Riemann--Roch gives $e(G)=\deg(G|_C)$.

\begin{lemma}\label{lem:planequinticseries}
Let $C\subset\mathbb P^2$ be a smooth plane quintic.  Its gonality is four.  A degree-five line bundle with three independent sections is $\OO_C(1)$.  Line bundles of degrees one and two have at most one section, and line bundles of degrees five, six and seven have at most three sections.
\end{lemma}
\begin{proof}
Adjunction gives $g(C)=6$, $K_C=\OO_C(2)$. The canonical embedding factors through the quadratic Veronese surface $V\subset\PP^5$. Since multiplication of conics surjects onto quartics and $I_{C/\PP^2}(4)=\OO_{\PP^2}(-1)$, restriction identifies
\[
 H^0(I_{C/\PP^5}(2))=H^0(I_{V/\PP^5}(2)),
 \qquad h^0(I_{V/\PP^5}(2))=21-15=6 .
\]
The symmetric $2\times2$ minors generate $I_V$. For a general fiber $D$ of a trigonal pencil, Clifford and Riemann--Roch give
\[
 h^0(D)=2,\qquad h^0(K_C-D)=4,\qquad \dim\langle D\rangle_{K_C}=1 .
\]
A quadric through $C$ vanishes at three distinct points of this line, hence contains it. The line would lie on $V$, where every curve has even degree. The canonical embedding excludes degree two, and projection from a point of $C$ has degree four.

Let $\deg M=5$ and $h^0(M)\ge3$. Clifford gives $h^0(M)=3$. After removing the base divisor, the plane map has degree $k$ onto a nondegenerate curve of degree $e$, with $ke\le5$. If $ke\le4$, then either $k=1$ and $g(C)\le (e-1)(e-2)/2\le3$, or $(k,e)=(2,2)$ and $C$ is hyperelliptic. Thus $ke=5$, $(k,e)=(1,5)$. The plane image has arithmetic genus $6=g(C)$, so it is smooth, and adjunction gives $M^2=K_C$.

This plane model gives the same canonical $C\subset\PP^5$. Its canonical quadrics recover the same $V$. On $V$ the ample line bundle whose square is $\OO_{\PP^5}(1)|_V$ is unique, since $\Pic(V)=\mathbb Z$. Its restriction gives $M=\OO_C(1)$.

Gonality gives $h^0\le1$ in degrees $1,2,3$. For degree $6$, four sections would give equality in Clifford for a noncanonical special bundle on a nonhyperelliptic curve. For degree $7$ use $h^0(M)=2+h^0(K_C\otimes M^{-1})\le3$. We use Clifford and its equality condition as in \cite{ACGH}.
\end{proof}

Put $R=G^{\vee\vee}$ and $k=\length(R/G)$.  Normality implies that
$R$ is rank-one reflexive and that $k$ is finite; in particular
\[
 e(R)=e(G),\qquad \chi(R)=\chi(G)+k.
\]
A rank-one reflexive sheaf on the normal Gorenstein surface $T$ is
maximal Cohen--Macaulay.  Duality therefore gives
\begin{equation}\label{eq:polar-reflexive-duality}
 e(R^\vee)=-e(R),\qquad
 \chi(R^\vee)=\chi(R)+e(R).
\end{equation}
Indeed, maximal Cohen--Macaulay duality identifies the derived dual
with $R^\vee$, and Serre duality gives
$\chi(R^\vee(h))=\chi(R)$.  Applying
\eqref{eq:polar-surface-polynomial} to $R^\vee$ at $t=1$ yields the
second equality.  Also
\[
 H^2(T,G)^\vee=\operatorname{Hom}(G,\omega_T)
             =H^0(T,R^\vee(h)).
\]

\begin{lemma}\label{lem:polar-degree-zero}
If $R$ is rank-one reflexive, $e(R)=0$, and
$R\not\cong\mathcal O_T$, then $\chi(R)\leq2$.
Consequently, if $G$ is rank-one torsion-free of degree zero and
$\chi(G)\geq3$, then
\[
 G\cong I_{Z/T},\qquad \length Z=5-\chi(G)\leq2.
\]
\end{lemma}

\begin{proof}
A nonzero section of a rank-one reflexive sheaf of degree zero has
an effective zero divisor of degree zero, hence no divisorial zero.
It is an isomorphism outside finitely many points, and reflexivity
on the normal surface makes it an isomorphism everywhere.  Thus
\[
 H^0(T,R)=H^0(T,R^\vee)=0
\]
when $R\not\cong\mathcal O_T$.

Suppose for a contradiction that $\chi(R)\geq3$, and put
$M=R^\vee(h)$.  Serre duality gives $h^0(T,M)\geq3$.
Choose a general smooth plane section $C$ avoiding the singular points. The sequence
\[
 0\longrightarrow R^\vee\longrightarrow M
   \longrightarrow M|_C\longrightarrow0
\]
and Lemma~\ref{lem:planequinticseries} give
\[
 3\le h^0(T,M)\le h^0(C,M|_C)\le3,\qquad M|_C\simeq\OO_C(1).
\]
Thus $H^0(T,M)\to H^0(C,M|_C)$ is an isomorphism. A fixed curve would intersect $C$ and give a base point of $|\OO_C(1)|$. The base locus is therefore finite, and the rational map $\varphi:T\dashrightarrow\PP^2$ restricts to the plane embedding on $C$.

If $\dim\varphi(T)=2$, a general section of $M$ gives $I_A=R(-h)\subset\OO_T$. The depth lemma gives $\depth\OO_{A,x}=1$ at closed points of $A$. Bertini, applied off the finite base locus, gives an irreducible and generically reduced divisor. Thus $A$ is integral of degree five, and
\[
 0\longrightarrow R\longrightarrow\mathcal O_T(h)
 \longrightarrow\mathcal O_A(h)\longrightarrow0
\]
gives
\[
 p_a(A)=1+\chi(R)\geq4.
\]
An integral nonplanar curve of degree five in $\mathbb P^3$ has
arithmetic genus at most two by the Castelnuovo bound.
Thus $A$ lies in a plane.  Since it has degree five and lies on $T$,
it is the full plane section of $T$.  It follows that
$I_A=\mathcal O_T(-h)$ and $R\cong\mathcal O_T$, a contradiction.

If $\dim\varphi(T)=1$, the image is $\varphi(C)\simeq C$. For the closure $F_y$ of a general fiber,
\[
 h\cdot [F_y]=\deg(C\longrightarrow\varphi(C))=1 .
\]
Indeed $C$ misses the base locus. Each component has positive degree, so the general fiber closure is a line with multiplicity one. These lines form an irreducible family covering $T$.

A member of such a moving family cannot avoid $\operatorname{Sing}T$.
If a line $\ell$ lies entirely in the smooth locus of $T$, adjunction
and $K_T=h$ give
\[
 N_{\ell/T}\cong\mathcal O_{\mathbb P^1}(-3),
\]
so $H^0(N_{\ell/T})=0$ and $\ell$ is isolated in the Hilbert scheme
of curves on $T$.  Since the singular locus is finite, irreducibility
of the parameter curve implies that every general ruling line passes
through one fixed singular point $p$.  Hence $T$ is the cone over
its general plane section $C$, with vertex $p$.

Take $C$ in a plane missing $p$. A general target line cuts five collinear points $x_1,\ldots,x_5$ on $C$, since $M|_C\simeq\OO_C(1)$. For the corresponding divisor $A$,
\[
 A_{\mathrm{red}}=\bigcup_{i=1}^5\langle p,x_i\rangle
   \subset\langle p,\langle x_1,x_2\rangle\rangle .
\]
The same depth argument excludes an embedded point at the vertex. The five rulings have multiplicity one, so $A$ is the full plane section. Again $R\simeq\OO_T$, proving $\chi(R)\le2$ for nontrivial $R$.

Finally, if $\chi(G)\geq3$, then
$\chi(R)=\chi(G)+k\geq3$, so $R\cong\mathcal O_T$.
The inclusion $G\subset R$ is an ideal of a finite subscheme $Z$,
and its Euler characteristic gives
$\length Z=5-\chi(G)\leq2$.
\end{proof}

\begin{lemma}\label{lem:polar-small-negative-degree}
Let $G$ be rank-one torsion-free on $T$.
\begin{enumerate}
 \item If $e(G)=-1$ or $e(G)=-2$, then $\chi(G)\leq4$.
 \item If $e(G)=-1$ and $\chi(G)=4$, then
 $G=I_{\ell/T}$ for a line $\ell\subset T$.
 \item If $e(G)=1$, then $\chi(G)\leq3$; equality holds precisely
 for $G\cong(I_{\ell/T})^\vee$ with $\ell\subset T$ a line.
 \item If $e(G)=2$, then $\chi(G)\leq2$.
\end{enumerate}
\end{lemma}

\begin{proof}
Write $R=G^{\vee\vee}$ and $k=\length(R/G)$.  Suppose first that
$e(G)=-j$ with $j=1$ or $2$.  A sheaf of negative degree has no
nonzero section, since its image would have an effective divisorial
quotient of negative degree.  Thus $h^0(G)=0$.
Choose a general smooth plane section $C$ avoiding the singularities.
Since $e(R^\vee(-h))=j-5<0$, restriction to $C$ gives
\[
 h^0(R^\vee)\leq h^0(R^\vee|_C)\leq1.
\]
The sequence
\[
 0\to R^\vee\to R^\vee(h)\to R^\vee(h)|_C\to0
\]
and the degree-six and degree-seven bounds on the plane quintic
then give $h^0(R^\vee(h))\leq1+3=4$.  Serre duality yields
\[
 \chi(G)=-h^1(G)+h^0(R^\vee(h))\leq4.
\]
This proves (i).

For $j=1$, $\chi(G)=4$ forces $h^0(R^\vee)=1$. The divisorial zero locus of its section is a degree-one curve $\ell$. By reflexivity and the depth lemma for $0\to I_{\ell/T}\to\OO_T\to\OO_\ell\to0$,
\[
 R^\vee=\OO_T(\ell),\qquad R=\OO_T(-\ell)=I_{\ell/T},\qquad
 k=\chi(I_{\ell/T})-\chi(G)=5-1-4=0 .
\]
This also applies at a singular point of $T$ where $I_{\ell/T}$ is noninvertible.

For positive degree $j=1,2$, the dual reflexive sheaf $R^\vee$
has degree $-j$.  By (i) and \eqref{eq:polar-reflexive-duality},
\[
 \chi(G)+k+j=\chi(R^\vee)\leq4.
\]
This proves the bounds in (iii)--(iv).  In degree one, equality
$\chi(G)=3$ forces $k=0$ and $\chi(R^\vee)=4$.
Part (ii) then gives $R^\vee=I_{\ell/T}$, as claimed.
Conversely its reflexive dual has degree one and Euler characteristic
three by \eqref{eq:polar-reflexive-duality}.
\end{proof}

For a sheaf $F=i_*G$ on a hyperplane support $T\subset X$, put
$e=e(G)$.  Its charges satisfy
\begin{equation}\label{eq:polar-charge-euler}
\begin{gathered}
 H\beta=e-\frac52,\qquad n=\chi(G)-\frac{25}{6},\\
 \widehat q_0=\chi(G)-\frac{25}{12}
                   -\frac{(e-5/2)^2}{10},\qquad
 \mu=-e\pmod5.
\end{gathered}
\end{equation}
Tensoring by $\mathcal O_X(H)$ permits one to take
$e\in\{-2,-1,0,1,2\}$.  Proposition~\ref{prop:generalvacuumbound} and
\eqref{eq:polar-charge-euler} leave exactly the following possibilities
for a polar term.  The table gives the full sheaf moduli in each case:
\[
\begin{array}{c|c|c|l}
 e&\mu&\chi(G)&G\\\hline
 0&0&5,4,3&I_{Z/T},\quad\length Z=0,1,2\\
 -1&1&5&\text{empty}\\
 -1&1&4&I_{\ell/T},\quad \ell\text{ a line}\\
 1&4&4&\text{empty}\\
 1&4&3&(I_{\ell/T})^\vee,\quad \ell\text{ a line}\\
 -2&2&6,5&\text{empty}\\
 2&3&4,3&\text{empty}.
\end{array}
\]
The assertions are Lemmas~\ref{lem:polar-degree-zero} and
\ref{lem:polar-small-negative-degree}; in particular, no charge-
conjugation conjecture or restriction to smooth supports is used.
\begin{proposition}\label{prop:polar-incidence-functors}
For $0\leq k\leq2$, the stable-sheaf moduli of character
$\operatorname{ch}(\mathcal O_T)-k[\mathrm{pt}]$ is isomorphic to
$\operatorname{Hilb}^k(\mathcal U_H/|H|)$ as a scheme.
The moduli of the degree-minus-one polar sheaves is the incidence
scheme of pairs $(\ell,T)$ with $\ell\subset X$ a line and
$T\in|H|$ containing it.  Its opposite polar sector is isomorphic
to the same incidence scheme under relative duality.
\end{proposition}

\begin{proof}
The preceding classification proves that every geometric point of the
first moduli problem is $F=i_*I_{Z/T}$.  We explain why the description
also holds on nonreduced bases.  For such an $F$, adjunction gives
\[
 \operatorname{Ext}^j_X(F,\mathcal O_X(-H))
 =\operatorname{Ext}^{j-1}_T(I_{Z/T},\mathcal O_T).
\]
Applying $R\operatorname{Hom}_T(-,\mathcal O_T)$ to
$0\to I_{Z/T}\to\mathcal O_T\to\mathcal O_Z\to0$ gives
\[
 0\longrightarrow\operatorname{Ext}^1_T(I_{Z/T},\mathcal O_T)
 \longrightarrow\mathbb C^k\longrightarrow H^2(T,\mathcal O_T)
 \longrightarrow\operatorname{Ext}^2_T(I_{Z/T},\mathcal O_T)
 \longrightarrow0.
\]
The middle map is dual to
$H^0(T,\omega_T)\to H^0(Z,\omega_T|_Z)$.
Since $\omega_T=\mathcal O_T(H)$ and $H$ separates every subscheme
of length at most two, that restriction is surjective.  As
$h^2(\mathcal O_T)=4$, we obtain
\begin{equation}\label{eq:polar-extension-groups}
\begin{gathered}
 \operatorname{Ext}^1_X(F,\mathcal O_X(-H))=\mathbb C,
 \qquad \operatorname{Ext}^2_X(F,\mathcal O_X(-H))=0,\\
 \dim\operatorname{Ext}^3_X(F,\mathcal O_X(-H))=4-k.
\end{gathered}
\end{equation}
Also $\operatorname{Hom}_X(F,\mathcal O_X(-H))=0$.
The unique nonzero extension class up to scalar is represented by
\[
 0\longrightarrow\mathcal O_X(-H)\longrightarrow I_{Z/X}
 \longrightarrow i_*I_{Z/T}\longrightarrow0.
\]
It is nonzero because its middle sheaf is torsion-free.

Now let $F$ be a flat family over a base $B$ and let
$p:X\times B\to B$.  The vanishing in
\eqref{eq:polar-extension-groups} and cohomology and base change
show that
\[
 E=\mathcal Ext^1_p(F,p_X^*\mathcal O_X(-H))
\]
is an invertible sheaf, compatible with arbitrary base change.
The identity of $E$ gives a universal extension
\[
 0\longrightarrow p_X^*\mathcal O_X(-H)\otimes p^*E^\vee
 \longrightarrow\mathcal I\longrightarrow F\longrightarrow0.
\]
Since $\mathcal Hom_p(F,p_X^*\OO_X(-H))=0$, the local extensions glue uniquely. Their middle term $\mathcal I$ is flat, with fiber $I_{Z/X}$. Put $L_B=\mathcal Hom_p(\mathcal I,\OO_{X\times B})$. By depth, $\mathcal Ext^j_X(\OO_Z,\OO_X)=0$ for $j<3$. The ideal sequence and $H^1(X,\OO_X)=0$ give
\[
 \Hom_X(I_{Z/X},\OO_X)=\mathbb C,\qquad \Ext^1_X(I_{Z/X},\OO_X)=0 .
\]
The relative Hom complex is perfect, with no negative fiber cohomology. Base change in degree zero gives a line bundle $L_B$ compatible with arbitrary pullback. The canonical evaluation
\[
 \mathcal I\otimes p^*L_B\longrightarrow\OO_{X\times B}
\]
is injective on every fiber. The local criterion of flatness gives an injection with flat cokernel, whose support is proper and quasi-finite over $B$. Thus $\mathcal I=I_{\mathcal Z}\otimes p^*B_0$ for a finite flat $\mathcal Z/B$ of length $k$ and $B_0=L_B^\vee$. This also gives the ideal reconstruction used in \cite[Proposition~2.1]{GST}.

The injection defining the extension corresponds to
\[
 s\in H^0\!\left(X\times B,I_{\mathcal Z}(H)\otimes p^*(B_0\otimes E)\right),
 \qquad \mathcal Z\subset\mathcal T=V(s).
\]
It is fiberwise nonzero, so $\mathcal T$ is a relative hyperplane and the cokernel recovers $F$. This gives the inverse incidence construction after arbitrary base change. Twisting the universal sheaf changes only the base line factors, so the inverse descends to the moduli scheme.

For the line sector let $R=I_{\ell/T}$ and $F=i_*R$.
The depth lemma and Gorenstein duality give
$\mathcal Ext_T^{>0}(R,\mathcal O_T)=0$. The restriction bound in
Lemma~\ref{lem:polar-small-negative-degree} and the inclusion
$R\subset\OO_T$ give $h^0(R^\vee)=1$. The sequence
\[
 0\to I_{\ell/T}(h)\to\OO_T(h)\to\OO_\ell(1)\to0
\]
gives $h^0(R(h))=4-2=2$. Hence $h^2(R^\vee)=2$ by Serre
duality, and $\chi(R^\vee)=3$ gives $h^1(R^\vee)=0$.
Adjunction on $X$ now reads
\[
 \dim\Ext_X^j(F,\OO_X(-H))
 =\dim H^{j-1}(T,R^\vee)
 =\begin{cases}1&j=1,\\2&j=3,\\0&j=0,2.\end{cases}
\]
The same universal extension therefore applies. On each fiber its
middle term is $I_{\ell/X}$ by
$0\to\OO_X(-H)\to I_{\ell/X}\to i_*I_{\ell/T}\to0$.
For a line on the smooth quintic, local duality for its regular embedding gives $\mathcal Ext_X^2(\OO_\ell,\OO_X)=\omega_\ell=\OO_\ell(-2)$ and $\mathcal Ext_X^{<2}(\OO_\ell,\OO_X)=0$. Consequently $\Hom_X(I_{\ell/X},\OO_X)=\mathbb C$ and $\Ext_X^1(I_{\ell/X},\OO_X)=H^0(\ell,\omega_\ell)=0$. The same evaluation gives a flat family of ideals of lines and the extension arrow cuts out the containing hyperplane. Both constructions commute with arbitrary base change.

Finally $R$ is maximal Cohen--Macaulay on $T$, hence $F$ has local
projective dimension one on the smooth threefold $X$ along its support.
For a flat family of these sheaves the derived dual into $\OO_X(-H)$
has fibers concentrated in degree one. By the fiberwise criterion for
a relatively perfect complex, its first cohomology is flat and commutes
with base change. Adjunction gives
\[
 \mathcal Ext_X^1(i_*R,\OO_X(-H))=i_*R^\vee.
\]
Twice applying this operation is the identity by biduality.
It therefore identifies both line sectors as schemes.
\end{proof}

\section{Specialization for punctual defects}
\subsection{A proper expanded sheaf space in the one-point charge}
\label{sec:one-defect-global}
For $v_1=\ch(\OO_T)-[\mathrm{pt}]$, $T\in|H|$, the generic sheaf is $I_{p/T}$. If $p$ specializes to a joining K3, its residue field has infinite projective dimension on the singular target. We retain the point on an expansion of the limiting hyperplane.

Let $P=\mathscr P$ be the blowup in
\textup{(\companionref{eq:vacuumblowup})}. Write $Z_P\to P$ for its universal minimal
expanded target and $U\hookrightarrow Z_P$ for its universal Cartier
surface. A refinement $c:Z'\to Z$ over an object $(Z,T)$ of $P$
inserts additional ruled components at joining divisors. A
\emph{one-point refinement} consists of such a refinement and a
section $p'$ of the smooth locus of $Z'$ with
$p'\in c^*T$.

Every new ruled component on a geometric fiber must
contain $p'$ in its interior. An arrow preserves the contraction,
the embedded Cartier divisor and the section, and is an
isomorphism of ordered expansions. This is the one-mark stability
condition of \cite[Definitions~1.2.1--1.2.2]{AF} relative to the
given minimal support expansion.

\begin{lemma}\label{lem:one-defect-universal}
Contraction of the section identifies the groupoid of one-point
refinements with $U$, over $P$. In particular it is proper and
separated over the degeneration curve. A fiber has at most two
inserted components, including the component required by the
minimal support.
\end{lemma}
\begin{proof}
Away from the joining divisors the inverse uses the original target
and the given section. At a join use an ordered expansion chart
$xy=d$ and a section with $x(p)=a$, $y(p)=b$. Then $ab=d$ in
the base ring, including over a nonreduced base. Replace that factor
of the target by the standard expansion with parameters $(a,b)$.
In two charts it is
\[
 x=az,\quad zy=b,\qquad xw=a,\quad y=wb,\qquad w=z^{-1}
 \text{ on the overlap}.
\]
Lift $p$ by $z(p')=1$. Then $c(p')=p$ and $p'\in c^*T$. The equation of $c^*T$ restricts to the old Cartier equations on the old components and to that of $\rho^{-1}(T\cap S)$ on a new ruled component. It is fiberwise regular, hence cuts out a flat relative Cartier divisor.

For the inverse, let $(u,v)$ be the parameters of a marked level and $z_0=z(p')\in\OO^\times$. The normalization
\[
 (u,v,z)\longmapsto (uz_0,v/z_0,z/z_0)
\]
gives $(a,b,1)$ and contracts the level if $a$ or $b$ is invertible. Under a change $x\mapsto hx$ of branch equation the normalized coordinates transform by
\[
 a\mapsto h(p)a,\qquad z\mapsto h z/h(p).
\]
The transition law is multiplicative. It descends the construction for the branch line bundles, whose tensor product has smoothing section $d$. Distinct joining divisors are disjoint. The resulting contraction and lift are inverse on each chart and respect its arrows. A dilation fixing $z(p')=1$ is the identity, and the minimal support has trivial stabilizer by Theorem~\companionref{thm:verifiedcompactification}.

Finally $U$ is closed in the proper family $Z_P$, and $P$ is
projective over the curve. The bound follows from one mark and at
most one old bubble.
\end{proof}

\begin{figure}[htbp]
 \centering\begingroup
 \definecolor{mpS}{HTML}{176568}
 \definecolor{mpT}{HTML}{B15C25}
 \definecolor{mpI}{HTML}{53618A}
 \begin{tikzpicture}[x=1cm,y=1cm,line join=round,line cap=round,
     branch/.style={draw=mpI!85!black,line width=1.1pt},
     guide/.style={draw=mpT!68,line width=.45pt},
     leader/.style={draw=black!57,line width=.4pt},
     joinpt/.style={circle,draw=mpS!75!black,line width=.7pt,fill=white,inner sep=0pt,minimum size=4pt},
     txt/.style={font=\small,inner sep=2pt}]
   \node[txt,anchor=west] at (-.35,2.29) {\textbf{(a)}\ \ The ruling over $p$};
   \node[txt,anchor=west] at (7.65,2.29) {\textbf{(b)}\ \ The contracted join};
   \filldraw[fill=mpT!11,draw=mpT,line width=.95pt]
      (2.07,.00) circle (1.20);
   \draw[guide,dash pattern=on 2pt off 2pt]
      (.87,0) arc[start angle=180,end angle=0,x radius=1.2,y radius=.29];
   \draw[guide] (.87,0) arc[start angle=180,end angle=360,x radius=1.2,y radius=.29];
   \draw[guide] (2.07,-1.20) .. controls (1.18,-.99) and (1.18,1.00) .. (2.07,1.20);
   \draw[guide,dash pattern=on 2pt off 2pt]
      (2.07,-1.20) .. controls (2.92,-1.00) and (2.92,.99) .. (2.07,1.20);
   \draw[branch] (2.07,1.20)..controls (1.51,1.31) and (1.09,1.67)..(.45,1.76);
   \draw[branch] (2.07,-1.20)..controls (1.42,-1.32) and (1.16,-1.66)..(.45,-1.76);
   \node[joinpt] at (2.07,1.20) {};
   \node[joinpt] at (2.07,-1.20) {};
   \node[txt,anchor=south west] at (2.11,1.18) {$w=0$};
   \node[txt,anchor=north west] at (2.11,-1.18) {$z=0$};
   \node[txt,anchor=east,text=mpI!75!black] at (.54,1.48) {$y=0$};
   \node[txt,anchor=east,text=mpI!75!black] at (.54,-1.48) {$x=0$};
   \draw[leader] (2.64,-.25)--(3.42,-.58)--(3.73,-.58);
   \node[txt,anchor=west] at (3.73,-.58) {$z(p')=1$};
   \filldraw[draw=mpT!55!black,fill=white,line width=.7pt]
               (2.64,-.25) circle (3.2pt);
   \fill[mpT!60!black] (2.64,-.25) circle (1.55pt);
   \node[txt,anchor=north] at (2.16,-1.94) {$c^{-1}(p)\simeq\PP^1$};
   \draw[-{Stealth[length=2.0mm,width=1.5mm]},draw=mpI!85!black,line width=.8pt]
      (5.40,0)--(6.90,0) node[midway,above,txt] {$c$};
   \begin{scope}[shift={(9.61,0)}]
     \draw[branch] (-1.63,0)--(1.58,0);
     \draw[branch] (-.60,-1.27)--(.60,1.27);
     \node[txt,anchor=west,text=mpI!75!black] at (.63,1.23) {$x=0$};
     \node[txt,anchor=south east,text=mpI!75!black] at (1.67,.03) {$y=0$};
     \draw[leader] (-.06,-.06)--(-.63,-.58)--(-1.04,-.58);
     \node[txt,anchor=east] at (-1.04,-.58) {$p$};
     \filldraw[fill=white,draw=mpS!75!black,line width=.75pt] (0,0) circle (3.2pt);
     \fill[mpS!75!black] (0,0) circle (1.55pt);
     \node[txt,anchor=north] at (0,-1.94) {$xy=0$};
   \end{scope}
 \end{tikzpicture}\endgroup
 \caption{A point of $T\cap S$ on the special fiber of the chart
 $xy=ab$, with $a=b=0$. The refinement inserts a ruling over $p$
 in the pulled-back surface $T'=c^*T$. The points $z=0$ and
 $w=0$ attach to the branches shown in the normal slice (b). In the coordinate fixed by the
 lift, $p'$ lies at $z=1$. The length-one quotient
 $\OO_{T'}\to k(p')$ is supported there.}
 \label{fig:one-point-refinement}
\end{figure}
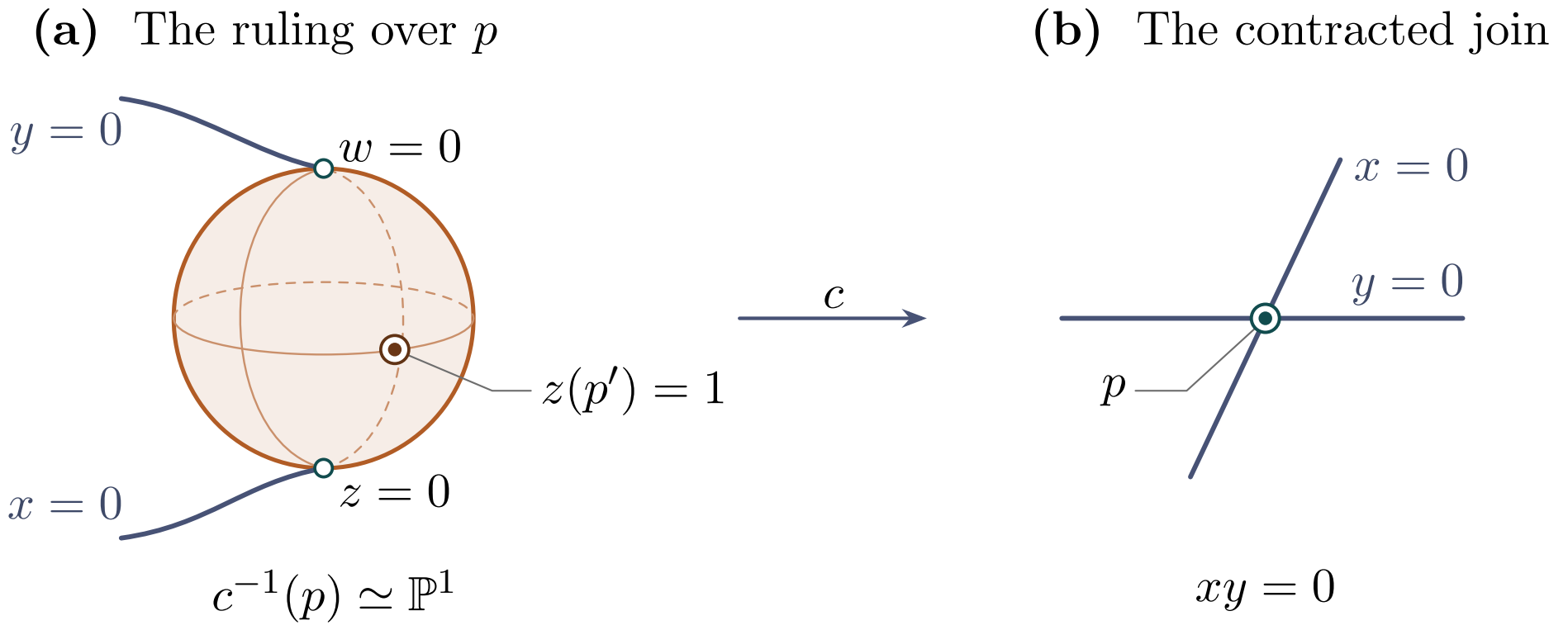

Figure~\ref{fig:one-point-refinement} locates the quotient on the new ruling in the normalized chart.

\begin{theorem}\label{thm:one-defect-global}
Under the hypotheses of Theorem~\companionref{thm:verifiedcompactification},
$U$ represents a proper, separated scalar-rigidified expanded
sheaf moduli space $M_1$ with universal sheaf, up to a base
invertible sheaf,
\[
                 F=I_{p'/T'},\qquad T'=c^*T.
\]
It contains the complete generic moduli of character $v_1$.
For each represented sheaf, every infinitesimal deformation
as a sheaf on a deforming expanded target comes uniquely from
the indicated incidence construction, up to the stipulated
isomorphisms. Its relative obstruction theory over $\EE$ has
tangent and obstruction groups $\Ext^1_{Z'}(F,F)$ and
$\Ext^2_{Z'}(F,F)$, respectively.
\end{theorem}

\begin{proof}
On each fiber the sequences
\[
 0\to\OO_{Z'}(-T')\to\OO_{Z'}\to\OO_{T'}\to0,\qquad
 0\to F\to\OO_{T'}\to\OO_{p'}\to0
\]
remain exact under base change. The first is a Cartier resolution and $p'$ is a relative regular immersion of codimension three in the smooth locus. Thus $F$ is flat over the base and perfect on the target. Its inclusion in the Cohen--Macaulay surface sheaf gives purity. Near a join $F=\OO_{T'}$, with boundary value $0$ or $\OO_A$ for the integral Cartier curve $A\in|h|$.

For the nef polarization a new ruled support over $A$ adds $P_{\OO_A}$ to each sum in
\[
 P_{\OO_{T'}}=\sum_i P_{\OO_{T'_i}}-\sum_e P_{\OO_{A_e}} .
\]
Thus the strict linear cut inequalities in
\textup{(\companionref{eq:vacuumordinarypolys})} and \textup{(\companionref{eq:vacuumgraphpolys})} persist. For the new component alone the residual polynomial is
$(D_{\rm right}-B_{\rm right})-D_{\rm left}=-B_e$ by \textup{(\companionref{eq:interiorres})}.

For every proper interior interval its
linear residual is
$[z](D_{\rm right}-B_{\rm right})-[z]D_{\rm left}<0$.
All these strict linear bounds persist under a sufficiently small
ample perturbation. Removing the point changes only the
constant terms, by Proposition~\companionref{prop:hull-punctual-DVR}.
On a single integral support stability follows from rank one.

This gives stability on expansion charts; a polarization furnishing
strict inequalities at a fiber remains ample and furnishes the
same strict inequalities near that fiber in the corresponding
branches of the chart.

To construct the inverse on infinitesimal families, fix a complete geometric
fiber and write $W=Z'$, $D=T'$, $p=p'$, and $N=\OO_W(D)$. The
normalization sequence and Lemma~\companionref{lem:vacuumlinearsystems} give
\begin{equation}\label{eq:one-defect-hull-sections}
 H^{>0}(W,N)=0,\qquad h^0(W,N)=5,\qquad
 H^0(W,N)\otimes\OO_W\twoheadrightarrow N .
\end{equation}
For the last assertion, in the ordinary profile the sections of
$\OO_Q(H_Q)$ and $\OO_B(H_B)$ restrict surjectively to
$H^0(S,\OO_S(h))$. In the empty-left-root profile the generators
on $B$ come from $f_5$ and $f_4H^0(\PP^3,\OO(1))$ of
$I_\Gamma(5)$, with $f_5$ nonzero on the strict transform of $S$.
In the graph profile the two boundary restrictions of a bubble
section are its two entries in
$H^0(\OO_S)\oplus H^0(\OO_S(h))$.

They extend independently
to the two roots, and $\OO_\Delta(D_+)$ is generated by those
sections. On an additional ruled component,
$N|_\Delta=\rho^*(N|_S)$ and
$\RR\rho_*\OO_\Delta=\OO_S$. Both joining restrictions are the
identity on $\RG(S,N|_S)$, so normalization cancels the new
component against one copy of this complex. It also gives the
asserted global generation on that component.

The expanded threefold is Gorenstein and $\omega_W\simeq\OO_W$.
From $0\to N^{-1}\to\OO_W\to\OO_D\to0$, Serre duality,
\eqref{eq:one-defect-hull-sections} and
Proposition~\companionref{prop:structurecohom} one obtains
$H^0(D,\OO_D)=\CC$, $H^1(D,\OO_D)=0$ and
$\omega_D=N|_D$, generated at $p$.
Apply $\RHom_D(-,\OO_D)$ to the point-ideal sequence.
The local ring of $D$ at $p$ is Gorenstein of dimension two,
even at a singularity of $D$. Its local residue-field Ext
groups are concentrated in degree two. The exact segment is
\[
 0\to\Ext^1_D(I_{p/D},\OO_D)\to
 \Ext^2_D(k(p),\OO_D)\longrightarrow H^2(D,\OO_D).
\]
The last arrow is dual to evaluation of $H^0(D,\omega_D)$ at
$p$, hence injective. Cartier adjunction now gives
\begin{equation}\label{eq:one-defect-extension}
 \Ext^j_W(F,N^{-1})=\Ext^{j-1}_D(I_{p/D},\OO_D),\qquad
 (\ext^0,\ext^1,\ext^2)_W(F,N^{-1})=(0,1,0).
\end{equation}
A generator is
$0\to N^{-1}\to I_{p/W}\to F\to0$. Splitting would inject
the nonzero divisor-torsion sheaf $F$ into the ideal of $p$ in
$\OO_W$, whereas a local equation of $D$ is a non-zero-divisor
on $\OO_W$.

Let $F_A$ be a flat deformation on $W_A$ over a local Artin algebra. At a small extension with kernel $J$, the obstruction and ambiguity for lifting $N$ lie in $H^2(\OO_W)\otimes J$ and $H^1(\OO_W)\otimes J$; both vanish. Write $N_A$ for its lift. By \eqref{eq:one-defect-extension} and perfect base change,
\[
 \Ext^1_{W_A}(F_A,N_A^{-1})\simeq A,\qquad
 0\to N_A^{-1}\to J_A\to F_A\to0
\]
for a generator. Units of $A$ act transitively on the generators. To lift the inclusion of the middle term, apply $\RHom_W(-,\OO_W)$ to the ideal sequence of $p$. Smoothness at $p$ gives
\[
 \Ext_W^{<3}(k(p),\OO_W)=0,\quad
 \Hom_W(I_{p/W},\OO_W)=\CC,\quad
 \Ext_W^1(I_{p/W},\OO_W)=0.
\]
Thus $J_A\hookrightarrow\OO_{W_A}$ lifts uniquely up to a unit. The fiber criterion for flatness gives a flat quotient. Its support is proper and quasi-finite, hence finite. Its coordinate algebra is locally free of rank one, recovering a section $p_A$. The composite $N_A^{-1}\to J_A\to\OO_{W_A}$ recovers a relative Cartier divisor $D_A$ containing $p_A$.

To recover the minimal target, work on the expansion
chart centered at $W$. If the mark stabilized an additional
bubble of the closed fiber, multiply the two adjacent smoothing
parameters and let $c_A:W_A\to W_A^0$ be that elementary
contraction. Otherwise take the identity. The old closed support
line bundle on the contracted target lifts to $N_A^0$ by the
same $H^1,H^2$ vanishing. Its pullback is $N_A$ up to a base
factor, by uniqueness. The map
\[
 H^0(W_A^0,N_A^0)\longrightarrow H^0(W_A,c_A^*N_A^0)
\]
is a map between free $A$-modules of rank five, by
\eqref{eq:one-defect-hull-sections} and base change. Its
reduction is the isomorphism from the normalization cancellation
above. It is consequently an isomorphism. The equation of $D_A$
descends and its descended divisor is relatively Cartier by the
fiber criterion. This is an object of $P$ by
Theorem~\companionref{thm:verifiedcompactification}. The point on it and
Lemma~\ref{lem:one-defect-universal} recover the original
refinement, section, and ideal.

The geometric inverse reads
\[
 I_D=\operatorname{Ann}_{Z'}F,\qquad (F|_D)^{\vee\vee}/F=\OO_p ,
\]
where the double dual is taken on $D$; indeed $\mathcal Hom_D(I_{p/D},\OO_D)=\OO_D$ by depth two. Contraction to the minimal support distinguishes the ruled level carrying only $p$. Thus the incidence map to the scalar-rigidified simple-sheaf stack is injective on geometric points and induces isomorphisms on stabilizers.

The map is representable and locally of finite presentation on the sheaf charts. The Artin inverse proves formal \'etaleness; universal injectivity then gives an open immersion. Its image inherits the obstruction theory of Theorem~\companionref{thm:ambientPOT}. Proposition~\ref{prop:polar-incidence-functors} identifies the generic fiber, and Lemma~\ref{lem:one-defect-universal} gives properness and separatedness.
\end{proof}

\begin{corollary}\label{cor:one-defect-specialization}
Let $M_1^\dagger=M_1\times_{\EE}\EE_0^\dagger$, where a point of
$\EE_0^\dagger$ is an expansion together with a chosen joining
divisor. If $g:M_1^\dagger\to(M_1)_0$ forgets the choice, then
\[
 \spc[M_X(v_1)]^\vir=g_*[M_1^\dagger]^\vir
       \quad\text{in } A_0((M_1)_0)_{\QQ}.
\]
Here $X$ is the geometric generic fiber with the corresponding
specialization convention. The degree of this equality is $-800$.
\end{corollary}
\begin{proof}
Writing $b:M_1\to\EE$ and using the virtual pullback of its sheaf
obstruction theory, the specialization and the normalization of
the expansion divisor give
\[
 0^!b^![\EE]=b_0^![\EE_0]
   =b_0^!\nu_*[\EE_0^\dagger]
   =g_*(b^\dagger)^![\EE_0^\dagger].
\]
The proper maps and obstruction morphisms used here are those
of the theorem and Lemma~\companionref{lem:virtualnormalization}. At a
chosen join, with restrictions $F_L,F_R$ and common perfect
restriction $E$, the marked sheaf obstruction theory uses
\[
 K_W\longrightarrow K_L\oplus K_R
 \xrightarrow{\,\mathrm{res}_L-\mathrm{res}_R\,}K_S
 \longrightarrow K_W[1].
\]
Here $K_W=\RHom_W(F,F)$, $K_i=\RHom_{W_i}(F_i,F_i)$ and
$K_S=\RHom_S(E,E)$. Apply $\RHom_W(F,-)$ to the normalization
triangle tensored with $F$ and use perfect adjunction. Scalar
rigidification and the global determinant give the reduction of
Proposition~\companionref{prop:POTbasechange}. The equality therefore uses
the obstruction theory of the surface sheaf even on a newly
inserted bubble. Finally the smooth-fiber incidence is a
$\PP^3$-bundle over $X$, of dimension six. Its symmetric
obstruction theory gives degree
$(-1)^6 e(\PP^3)e(X)=4(-200)=-800$.
\end{proof}

At a positive join the marked obstruction theory is determined by the difference of the two restriction maps. The factorization question is whether this homotopy fiber is realized on the product of two independently proper relative sheaf spaces. The distinction can be seen on the plane component in Proposition~\ref{prop:one-defect-plane-extra}.

The marked pushforward in Corollary~\ref{cor:one-defect-specialization}
defines a Chow class on the entire proper special fiber. To express
each summand as a diagonal pullback from independently compactified
relative sheaf spaces one also needs the corresponding deformation
identifications on each half. The evaluation argument in
\eqref{eq:one-defect-extension} uses the dualizing sheaf of the
\emph{complete} support.

For example, on a ruled surface with one
joining section $s$ at the end of a relative chain, the dualizing
sheaf of the union restricts to $\omega_R(s)$, of degree $-1$ on a
ruling. Its global sections vanish along that ruled component.
This specifies the additional extension calculation required at an
outer-boundary defect.

\begin{corollary}\label{cor:one-defect-zero-relative}
The summand of the preceding Chow identity having zero restriction
at the marked K3 is the pushforward of an independently defined
proper relative sheaf class. Write $D_B=5H_B-E_B$ and
$V_B=H^0(B,\OO_B(D_B))$. Its relative moduli scheme is
\[
 R_{B,1}^0 =
 \{(b,[s])\in B\times\PP(V_B) : s(b)=0\}
 =
 \PP_B\!\left(\ker(V_B\otimes\OO_B\longrightarrow\OO_B(D_B))\right),
\]
with one-dimensional-subspace convention for projectivization.
If $g_0:R_{B,1}^0\to(M_1)_0$ attaches the empty quartic side,
this summand is $(g_0)_*[R_{B,1}^0]^\vir$, with the
boundary-trivialized determinant of Theorem~\companionref{thm:zeroboundaryPOT}.
\end{corollary}
\begin{proof}
On the $B$ side use the relative minimal Cartier supports of
class $D_B$ and one marked point, with an extra expansion whenever
the point hits an internal joining divisor. A zero outer
restriction means that the support misses the outer divisor by
Nakayama's lemma. There is consequently no refinement based at
the outer divisor.

By Theorem~\companionref{thm:verifiedcompactification}, the minimal support is the universal divisor on $B\times\PP(V_B)$; its center ideal stays invertible over the hyperplane of divisors containing $S$. Lemma~\ref{lem:one-defect-universal} identifies the marked incidence at each internal join. On $B$, \eqref{eq:one-defect-hull-sections} gives
\[
 V_B\otimes\OO_B\twoheadrightarrow\OO_B(D_B),\qquad
 \operatorname{rk}\ker(\mathrm{ev})=4,
\]
by the generators $f_5$ and $f_4H^0(\PP^3,\OO(1))$.

Let $M_{1,0}^{\dagger}$ be the marked piece with zero restriction. The maps on families satisfy
\[
 R_{B,1}^0\ \underset{\mathrm{restriction}}{\overset{\mathrm{attachment}}{\rightleftarrows}}\ M_{1,0}^{\dagger},
 \qquad F_B\longmapsto 0\cup F_B,\qquad (0\cup F_B)|_{W_B}=F_B .
\]
The marked join has an empty quartic root and no stable empty bubble on that side. Since the relative support misses the marked K3, restriction and extension by zero are inverse. Nakayama preserves this disjointness in infinitesimal families, to which Theorem~\ref{thm:one-defect-global} applies after attachment. On the marked K3 the matching object is $0$, with $\mathbf R\!\operatorname{Hom}(0,0)=0$. Theorem~\companionref{thm:zeroboundaryPOT} identifies the two scalar and determinant reductions, giving the asserted virtual class.
\end{proof}

\begin{proposition}\label{prop:one-defect-zero-count}
The relative degree of the zero-boundary summand is $-120$.
The sum of the degrees of the positive-boundary marked summands is
therefore $-680$.
\end{proposition}
\begin{proof}
Write $R=R_{B,1}^0$ and $[c:\ell]$ for the parameters of
\textup{(\companionref{eq:vacuumBzeroeq})}. The complement of the unexpanded
locus of $R$ is $D_1\cup D_2$, with
\[
 D_1=S\times\PP^3,\qquad
 D_2=\{(b,[\ell]):\ell(b)=0\}\subset B\times\PP^3,
 \qquad D_{12}=D_1\cap D_2.
\]
The last two schemes are $\PP^2$-bundles over $B$ and $S$,
respectively, since four linear forms generate $\OO_B(H_B)$.
Locally along $D_{12}$ write the incidence as $c=-u\ell(b)$
as in \textup{(\companionref{eq:vacuumBzeroeq})}. The equations $u=0$ and
$\ell(b)=0$ have independent differentials: $u$ is a parameter
normal to $S$, and varying the plane varies its value at $b$.
Thus $D=D_1+D_2$ has simple normal crossings on the smooth
sixfold $R$.

Let $\mathcal E_{B,S}$ be the stack of expansions of the pair.
At a point lying on exactly $k$ irreducible branches of $D$,
the chart of $R\to\mathcal E_{B,S}$ is the map to
$[\AA^k/\GG^k]$ given by the corresponding local equations
$d_1,\dots,d_k$. Indeed at a general boundary point the parameter
is $c$; at $D_{12}$ the graph-to-root join has coordinates
$u,w$ with $uw=c$, and the contracted mark has
$(u(b),w(b))=(u(b),-\ell(b))$.

The one-mark expansion factors
$c$ by precisely these two entries. Since $d_1,\ldots,d_k$ are
part of a regular system of parameters, the map to the expansion
stack is smooth. Its ordinary relative tangent bundle $K$ is
the kernel in the locally split exact sequence
\[
 0\to K\to T_R\oplus\OO_R^k
   \xrightarrow{(v,\alpha)\mapsto(v(d_j)-\alpha_jd_j)_j}
   \OO_R^k\to0 .
\]
Projection to $T_R$ identifies $K=T_R(-\log D)$: each
$\alpha_j$ is uniquely determined by $v(d_j)\in(d_j)$.
This description is unchanged on replacing $d_j$ by a unit
multiple and therefore is global.

The support misses $D_{\rm out}$. Write $p$ for the universal target
over $R$. Relative Serre duality gives
\[
 F(-D_{\rm out})\simeq F,\qquad
 \mathcal Ext^2_p(F,F)\simeq
       \bigl(\mathcal Ext^1_p(F,F)\bigr)^\vee .
\] Indeed
$\omega_{W/\mathcal E_{B,S}}\simeq\OO_W(-D_{\rm out})$ by
relative log adjunction, after choosing the trivialization of
$K_B+S$. On a bubble the invariant logarithmic differential
in the ruling direction identifies these trivializations
under rescaling. The boundary determinant construction removes
only the scalar and the top Ext.

Smoothness over
$\mathcal E_{B,S}$ identifies its degree-zero tangent sheaf
with $K$, so its obstruction bundle is $K^\vee$.
Smoothness also makes the ordinary tangent morphism to the expansion
stack surjective in degree zero. The transitivity triangle therefore
identifies the absolute obstruction sheaf with the same bundle.
Consequently
\[
 [R]^\vir=c_6(\Omega_R^1(\log D))\cap[R].
\]
For each smooth stratum $D_I=\bigcap_{i\in I}D_i$, adjunction
gives $c(T_{D_I})=c(T_R)|_{D_I}/\prod_{i\in I}(1+D_i)|_{D_I}$.
The exact logarithmic tangent sequence gives
$c(T_R(-\log D))=c(T_R)/((1+D_1)(1+D_2))$.
Expanding $\prod_j(1-D_j/(1+D_j))$, pushing forward from the
strata and applying Gauss--Bonnet yields
\[
 \deg[R]^\vir=e(R)-e(D_1)-e(D_2)+e(D_{12})
  =4e(B)-4e(S)-3e(B)+3e(S)=e(B)-24 .
\]
Here the dualization sign is $(-1)^6=1$. The center
$\Gamma\subset\PP^3$ is a complete intersection of degrees
$4,5$, so $2g(\Gamma)-2=20(4+5-4)=100$. The blowup formula
gives $e(B)=4+(2-1)e(\Gamma)=4-100=-96$.
This proves the first value. Subtract it from the degree in
Corollary~\ref{cor:one-defect-specialization} to obtain
$-800-(-120)=-680$ for the sum of the other marked degrees.
\end{proof}

\begin{proposition}\label{prop:one-defect-positive-marked}
The two positive-boundary marked pieces, distinguished by the
side carrying the point defect, have intrinsic marked virtual
degrees $-320$ (quartic side) and $-360$ (blown-up projective
side). Together with Proposition~\ref{prop:one-defect-zero-count}
this computes the degree of every summand of the one-point
specialization directly on its marked space. More precisely,
write $V_0=R_{B,1}^0$, $D_0=D_1+D_2$ for its boundary in
Proposition~\ref{prop:one-defect-zero-count}, and let
$g_i:V_i\to(M_1)_0$ forget the chosen join. Then
\begin{equation}\label{eq:one-defect-three-Chern}
 \spc[M_\eta(v_1)]^\vir=
 \sum_{i=L,R,0}(g_i)_*
     \left(c_6(\Omega^1_{V_i}(\log D_i))\cap[V_i]\right)
      \quad\text{in }A_0((M_1)_0)_{\QQ},
\end{equation}
where the two positive boundaries are specified below.
\end{proposition}
\begin{proof}
Put $A=\operatorname{Bl}_{[S]}|H_Q|$, $a:A\to|h|=\PP^3$.
The morphism $a$ is a $\PP^1$-bundle, and its exceptional
divisor $E_A$ is a section. On the positive matching an
underlying minimal relative support is a point $d\in A$
together with the unique plane of parameter $a(d)$ on $B$, by
\textup{(\companionref{eq:vacuumpositiveproduct})}. On the side carrying the
defect the mark may reach the outer divisor of this support
pair.

A section meeting the outer divisor with local normal
coordinate $x$ has the unique relative expansion with parameter
$a_0=x(p)$ and interior lifted coordinate $z(p')=1$ in
$x=a_0z$. This is the one-ended version of
Lemma~\ref{lem:one-defect-universal}: replacing the initial
interior coordinate $z_0$ by $z/z_0$ changes the parameter by
the unit $z_0$. The other chart is $a_0=xw$. The construction
commutes with arbitrary base change.

Its contraction leaves
the underlying incidence of the point with the abstract support
unchanged. The ideal of the mark recovers the point by
Theorem~\ref{thm:one-defect-global}. At the marked join the
boundary sheaf is $\OO_C$ for the same plane curve $C$ on both
sides. Its scalar identifications are absorbed by the scalar
automorphisms on the two simple branches. Consequently as
\emph{marked} spaces with their classical scheme structures the
two sectors are
\begin{align*}
 V_L&=\{(d,q)\in A\times Q : q\in T_d\},\qquad
        T_d\in |H_Q|\text{ the contracted hyperplane},\\
 V_R&=A\times_{\PP^3}
       \{(b,[\ell])\in B\times\PP^3 : \ell(b)=0\}.
\end{align*}
For $V_L$ the isomorphism of the abstract expanded Cartier support
with its contracted surface follows from the invertible ideal on
the support in the relative chart proof of Theorem~\companionref{thm:verifiedcompactification}. To verify the
identifications over the marked stack one may equivalently cut
the one-point refinement at the distinguished join and contract
each half independently to its minimal Cartier support, retaining
the contracted mark on its half.

The inverse glues the reconstructed
one-ended refinements. If the mark of the minimal support was on
the old join, its added bubble has two ends: marking one end puts
the defect on the left and marking the other puts it on the right.
These two choices belong to different marked spaces.

We compute the smoothing parameters on the marked expansion stack. In $V_L$ near $E_A$ write the hyperplane equation
$u+s\ell=0$ where $u$ cuts out $S\subset Q$, $s$ cuts out
$E_A$ and $[\ell]=a(d)$. At $s=0$ the equation cuts out the
smooth $S$ with nonvanishing normal differential $du$, so $V_L$
is smooth there. Away from $E_A$ it is an open of the smooth
ordinary incidence, a $\PP^3$-bundle over $Q$.

At $s=0$ the
graph has normal coordinate to the outer divisor proportional
to $\ell(q)$, by \textup{(\companionref{eq:vacuumstrict})}. Thus in addition to
the fixed marked node the possible unsmoothed parameters are
$s$ (the graph-to-empty-root node) and $\ell(q)$ (the added
node when the mark reaches the outer divisor). On $s\ne0$ the
latter is a unit multiple of the usual outer coordinate.
Their zero loci meet transversely.

At $s=\ell(q)=0$, the
coordinate $s$ is free and $\ell(q)$ varies with the plane
parameter; the values of linear forms on $S$ generate
$\OO_S(h)$. Away from $s=0$ the only boundary locus is the
smooth incidence of a plane through a point of $S$ with the
remaining hyperplane parameter. In $V_R$ the analogous two
divisors are the pullback of $E_A$ and $b\in S$.

Indeed the
ordinary plane incidence is a $\PP^2$-bundle over $B$, also
after restricting to $S$; $V_R$ is a $\PP^1$-bundle over that
incidence and $E_A$ pulls back to its section. The outer marking
coordinate is a local equation of $S$ at $b$. Thus $V_R$ is
smooth and the two boundary divisors meet transversely as well.

Let $D_i$ be the reduced union of these divisors on $V_i$, $i=L,R$. The face marking the $j$th join in $t=t_0\cdots t_k$ is
\[
 t_j=0,\qquad \chi_\ell=\lambda_\ell^{-1}\lambda_{\ell+1}\quad
  (\ell\ne j),\qquad \lambda_0=\lambda_{k+1}=1.
\]
The inverse character change is
\[
 \lambda_m=\prod_{\ell=0}^{m-1}\chi_\ell\quad(m\le j),
 \qquad \lambda_m=\prod_{\ell=m}^k\chi_\ell^{-1}\quad(m>j).
\]
Hence this face is $[\AA^k/\GG^k]$ with coordinatewise scaling. The remaining smoothing parameters cut out the transverse divisors of $D_i$; the kernel computation of Proposition~\ref{prop:one-defect-zero-count} gives relative tangent $T_{V_i}(-\log D_i)$ and smoothness over the face.

The marked obstruction theory is the restriction of the complete-chain theory of Corollary~\ref{cor:one-defect-specialization}. The relative Calabi--Yau volume form extends over every ruling by the invariant logarithmic differential. Relative Serre duality and smoothness over the expansion stack give the same absolute and relative obstruction bundle. Since $\dim V_i=6$,
\[
 \mathrm{Ob}_{V_i}=\Omega_{V_i}^1(\log D_i),\qquad
 \deg[V_i]^\vir=\int_{V_i}c_6(\Omega_{V_i}^1(\log D_i))
                =e(V_i\setminus D_i).
\]

For $V_L\setminus D_L$ one picks $q\in Q\setminus S$ and a
hyperplane through $q$. The excluded hyperplane $S$ does not
contain such a point. The fiber is $\PP^3$ and the family
is a projective bundle. For $V_R\setminus D_R$ one picks
$b\in B\setminus S$, a plane through its image in $\PP^3$,
and a point of $a^{-1}([\ell])\setminus E_A$. The latter
two choices give a $\PP^2$-bundle followed by a Zariski locally
trivial affine-line bundle. Finally
$c(T_Q)=(1+H)^5/(1+4H)$ gives
$e(Q)=-14\int_Q H^3=-56$. Therefore
\[
 \deg [V_L]^\vir=4(-56-24)=-320,\qquad
 \deg [V_R]^\vir=3(-96-24)=-360 .
\]
Their sum is $-680$, and including the zero-boundary class gives
$-320-360-120=-800$.
\end{proof}

\begin{proposition}\label{prop:one-defect-plane-extra}
On the positive plane side $(B,S)$ there are singular plane
transforms $T$ and points $p\in T\setminus S$ for which
\[
 \dim\Ext^1_T(I_{p/T},I_{p/T})=4,\qquad
 \dim T_p\Hilb^1(T)=3 .
\]
In particular even with the plane support fixed, its one-point
incidence does not exhaust the deformations as a relative surface
sheaf.
\end{proposition}
\begin{proof}
Choose a point $q$ of $\Gamma$ at which its curvature is nonzero
and a plane containing its tangent line with contact exactly two
and distinct from the tangent plane of $S$. Such $q$ exist:
if the derivative of the tangent-line map were identically zero
in characteristic zero, the tangent line would be constant and
$\Gamma$ would be a line.

At such $q$ at most one plane in the
pencil contains the second osculating direction. Write $s$ for
a parameter on $\Gamma$ and $x,y$ for transverse coordinates so
that $I_\Gamma=(x,y)$ and the plane equation is $x+s^2$. Indeed
its restriction has a double zero, whose unit factor can be
absorbed in the parameter, and its differential is a nonzero
conormal at $q$. In the chart $x=yv$ of the blowup, the plane
transform has ordinary double point
\[
 R=\widehat{\OO}_{T,p}=\CC[[y,v,s]]/(yv+s^2),\qquad
 \mathfrak m=(y,v,s).
\]
The point $p$ on the exceptional ruling represents the normal
direction annihilated by the plane conormal. The strict transform
of $S$ meets the ruling in the direction annihilated by the
conormal of $S$. The distinct tangent planes imply $p\notin S$.

On $T$ the sheaf $I=I_{p/T}$ has $\mathcal Hom(I,\OO_T)=\OO_T$
by depth two and $H^{>0}(T,\OO_T)=0$ by the divisor sequence on
$B$. Gorenstein duality and the point-ideal sequence yield
\begin{align*}
 \Hom_T(I,I)&=\Hom_T(I,\OO_T)=\CC,\\
 \Hom_T(I,k(p))&=(\mathfrak m/\mathfrak m^2)^*,\\
 \Ext_T^1(I,\OO_T)&=\Ext_R^1(\mathfrak m,R)=\CC .
\end{align*}
The map $\Ext_R^1(\mathfrak m,R)\to\Ext_R^1(\mathfrak m,k)$
is zero, as follows directly from a presentation. Let $e_y,e_v,e_s$
map to the generators of $\mathfrak m$, and put
$N=\ker(R^3\to\mathfrak m)$. It is generated by the three
Koszul relations $K_{ij}=x_j e_i-x_i e_j$ and
$g=v e_y+s e_s$. To see generation, lift a relation to
$\CC[[y,v,s]]$, subtract a multiple of the lift of $g$ to remove
its multiple of $yv+s^2$, and use Koszul exactness on $y,v,s$.

For any $\phi:N\to R$ the Koszul relation among the $K_{ij}$
implies that each $\phi(K_{ij})$ lies in $\mathfrak m$:
reducing the relation modulo $\mathfrak m^2$ gives a linear
relation among the independent classes $y,v,s$ with those
constant terms as coefficients. Further
$x_jg=vK_{yj}+sK_{sj}$ in $N$.

Applying $\phi$ gives
$x_j\phi(g)\in\mathfrak m^2$ for every $j$, whence
$\phi(g)\in\mathfrak m$. Thus each map representing an
$\Ext_R^1(\mathfrak m,R)$ class becomes zero after postcomposition
with $R\to k$. Globally the target $\Ext_T^1(I,k(p))$ is local
at $p$, so the same vanishing holds. Applying $\RHom_T(I,-)$ to
the ideal sequence now gives
\[
 0\to(\mathfrak m/\mathfrak m^2)^*
 \to\Ext^1_T(I,I)\to\CC\to0 .
\]
The first arrow moves the length-one quotient on the fixed
surface. Since $\Hilb^1(T)=T$ and its tangent dimension at $p$
is three, this proves the assertion. All these deformations remain
admissible at the outer divisor over a square-zero parameter:
admissibility is open in flat families and holds on its unique
geometric fiber. Pushing a deformation on the fixed closed
subscheme $T$ to $B$ retains its extension class, since any
$\OO_B$-linear splitting between sheaves annihilated by the ideal
of $T$ is $\OO_T$-linear.
\end{proof}

\begin{proposition}[Deformation with prescribed structure-sheaf boundary]\label{prop:plane-boundary-rigidity}
Let $W$ be one of the positive relative chains just constructed,
let $S_o$ be its outer K3 divisor, and let $D\subset W$ be its
Cartier support on that expansion, including the pullback to a
one-point refinement. Put $C=D\cap S_o$, $N=\OO_W(D)$, and
$F=i_*I_{p/D}$ for a point $p\in D\setminus S_o$ smooth on
$W$. The point may be singular on $D$.

Given a local Artin
$\CC$-algebra $A$, a deformation $(W_A,S_{o,A})$ as an
expanded pair, and a relative Cartier deformation
$C_A\subset S_{o,A}$, any flat sheaf deformation $F_A$ equipped
with $F_A|_{S_{o,A}}\simeq\OO_{C_A}$ comes from a relative
Cartier divisor $D_A$ meeting $S_{o,A}$ in $C_A$ and a section
$p_A\in D_A\setminus S_{o,A}$:
\[
                 F_A\simeq I_{p_A/D_A}.
\]
The divisor and section are unique. The isomorphisms are taken
with the scalar rigidification of Section~\companionref{sec:reduction}.
In particular, for the singular plane transform of
Proposition~\ref{prop:one-defect-plane-extra} the restriction on
fixed-support first-order deformations has rank exactly one,
\[
 \ker\left[\Ext_T^1(I_{p/T},I_{p/T})
   \longrightarrow H^1(C,\OO_C)\right]=T_pT .
\]
\end{proposition}
\begin{proof}
On the minimal positive branch the abstract support is a
hyperplane of $Q$, a graph isomorphic to $S$, or a plane
transform on $B$. The divisor sequences used in
Propositions~\ref{prop:Bplanefamilyfull} and
\ref{prop:Qstructurefamilyfull} give
$H^0(\OO_D)=\CC$ and $H^1(\OO_D)=0$; the graph gives
the same identities since it is a K3. A point refinement adds
a ruled component $P_C\to C$ of the support, attached along a
section, at the end or at an existing join.

Its direct image
of $\OO_{P_C}$ is $\OO_C$ with vanishing higher images. The
normalization sequence therefore cancels the new component
against one copy of $R\Gamma(\OO_C)$, also for a singular
integral $C$. The relative support line bundle on the minimal
branch is generated at $p$: it is the hyperplane line bundle on
the roots, while on the graph it restricts to $\OO_S(h)$.
Contraction to the minimal branch pulls this bundle back on a
refinement. Logarithmic adjunction gives
\[
 H^0(\OO_D)\xrightarrow{\sim}H^0(\OO_C),\qquad
 H^0(\OO_D(-C))=H^1(\OO_D(-C))=0,\qquad
 \omega_D(C)=N|_D .
\]
Here $C$ is an integral Cartier curve and $D$ is Gorenstein.
Apply $R\operatorname{Hom}_D(-,\OO_D(-C))$ to the ideal
sequence of $p$. Serre duality identifies the dual of
\[
 \Ext^2_D(k(p),\OO_D(-C))\longrightarrow H^2(D,\OO_D(-C))
\]
with the evaluation $H^0(\omega_D(C))\to\omega_D(C)|_p$,
which is onto. The preceding groups and Gorenstein local duality
give $\operatorname{Hom}_D(I_{p/D},\OO_D(-C))=
\Ext^1_D(I_{p/D},\OO_D(-C))=0$. Cartier adjunction yields
\[
  \Ext_W^j(F,N^{-1}(-S_o))
      =\Ext_D^{j-1}(I_{p/D},\OO_D(-C))=0
        \qquad (j=1,2).
\]

By $H^1(\OO_W)=H^2(\OO_W)=0$ the line bundle $N$ lifts uniquely
up to isomorphism to $N_A$; its restriction to the K3 equals
$\OO_{S_{o,A}}(C_A)$ up to isomorphism, since $H^1(\OO_{S_o})=0$.
For $j=1,2$ the above vanishings persist with target
$N_A^{-1}(-S_{o,A})$. Indeed, filter this line bundle by the
powers of the maximal ideal of $A$.

Each graded quotient is a
sum of copies of $N^{-1}(-S_o)$ on the closed fiber, and derived
adjunction with that fiber computes Ext from $F$ because $F_A$
is $A$-flat. The long exact sequences give the assertion by
induction. The relative Cartier sequence of $S_{o,A}$ and
derived adjunction now yield an isomorphism
\[
 \Ext^1_{W_A}(F_A,N_A^{-1})\xrightarrow{\sim}
 \Ext^1_{S_{o,A}}(\OO_{C_A},N_A^{-1}|_{S_{o,A}}).
\]
There are no extra Tor terms on the right: multiplication by
the outer equation is injective on the closed-fiber $F$, so
the relative flatness criterion gives injection on $F_A$ and
flat quotient. Lift the extension
$0\to\OO_{S_{o,A}}(-C_A)\to\OO_{S_{o,A}}\to\OO_{C_A}\to0$.
The middle term $J_A$ of its lift is $A$-flat and reduces to
$I_{p/W}$, because the extension
$0\to N^{-1}\to I_{p/W}\to F\to0$ has precisely this boundary.

Changing either boundary identification by a unit gives an
isomorphic extension after rescaling. Indeed
$H^0(\OO_{S_{o,A}})=A$ and $H^0(\OO_{C_A})=A$, by the flatness
and geometric connectedness of the family of integral curves
in $|h|$ (or the Cartier sequence on the constant K3).

Finally $\operatorname{Hom}_W(I_{p/W},\OO_W)=\CC$ and
$\Ext^1_W(I_{p/W},\OO_W)=0$ by the point-ideal sequence,
$H^1(\OO_W)=0$, and the depth-three smooth local ring at $p$.
Successive square-zero liftings therefore give
$J_A\hookrightarrow\OO_{W_A}$ uniquely up to a unit of $A$:
a lift exists at each step and the differences between two
lifts are scalar multiples of the inclusion.

The fiberwise
injectivity criterion gives a flat cokernel, which is the algebra
of a unique section $p_A$ of length one. The composite
$N_A^{-1}\to J_A\to\OO_{W_A}$ cuts out $D_A$ since its
reduction is a regular equation. Near $S_{o,A}$ the ideal
$J_A$ is the unit ideal, and the restricted extension shows
$D_A\cap S_{o,A}=C_A$ as relative Cartier divisors. The quotient
is $I_{p_A/D_A}$, as asserted.

To identify the kernel in the last display of the proposition,
work on the fixed surface $T$ and put $I=I_{p/T}$. Tensor its
ideal sequence by $\OO_T(-C)$ and use the two vanishings above:
\[
 \Ext^1_T(I,I(-C))
   =\operatorname{Hom}_T(I,k(p))\otimes\OO_T(-C)|_p .
\]
The exact sequence $0\to I(-C)\to I\to\OO_C\to0$, adjunction
and $H^0(\OO_T)\simeq H^0(\OO_C)$ identify the image of this
group injectively with the kernel of restriction to
$\Ext^1_C(\OO_C,\OO_C)=H^1(\OO_C)$. Since $p\notin C$, the
outer equation trivializes $\OO_T(-C)|_p$. The ideal-quotient
tangent is $\operatorname{Hom}_T(I,k(p))=T_pT$, of dimension
three at the ordinary double point. Proposition~\ref{prop:one-defect-plane-extra}
gives dimension four in the source. Thus the additional
direction changes the boundary sheaf already to first order.
\end{proof}
\begin{corollary}[Finite defects with prescribed boundary]\label{cor:finite-defect-boundary}
Let $(W,D,C,N)$ satisfy the hypotheses and cohomology identities of
Proposition~\ref{prop:plane-boundary-rigidity}. Let $Z\subset D$ be a
finite subscheme of length $m\ge1$ in the smooth locus of $W$,
disjoint from the outer divisor and all joining divisors. Assume that
\begin{equation}\label{eq:finite-defect-evaluation}
 H^0(D,\omega_D(C))\longrightarrow H^0(Z,\omega_D(C)|_Z)
       \quad\hbox{is surjective}.
\end{equation}
At $i_*I_{Z/D}$, every deformation of the sheaf on an expanded pair
over a local Artin base with prescribed structure-sheaf boundary
comes from a unique flat Cartier support containing a finite flat
subscheme of degree $m$, up to the scalar convention above. On an
unexpanded quartic root this holds for all $m=2$, including a
nonreduced length-two subscheme at a singular point of its
hyperplane surface, provided the subscheme avoids the outer divisor.
\end{corollary}
\begin{proof}
Let $I=I_{Z/D}$. Serre duality on the proper Gorenstein surface, or
local Gorenstein duality followed by the local-to-global sequence,
gives $\Ext_D^{j}(\OO_Z,\OO_D(-C))=0$ for $j<2$ and identifies
the dual of
\[
 \Ext^2_D(\OO_Z,\OO_D(-C))\longrightarrow H^2(D,\OO_D(-C))
\]
with \eqref{eq:finite-defect-evaluation}. Its injectivity gives
$\Ext^{j}_D(I,\OO_D(-C))=0$ for $j=0,1$, by the ideal sequence
and $H^0(\OO_D(-C))=H^1(\OO_D(-C))=0$. Cartier adjunction,
devissage over an Artin base and the outer Cartier sequence give
exactly the relative extension isomorphism in
Proposition~\ref{prop:plane-boundary-rigidity}. The canonical
lifted boundary extension has middle term reducing to $I_{Z/W}$.
Since $\OO_Z$ has finite length in the smooth threefold locus,
$\Ext_W^j(\OO_Z,\OO_W)=0$ for $j<3$. The source ideal therefore
still satisfies
\[
 \operatorname{Hom}_W(I_{Z/W},\OO_W)=\CC,\qquad
 \Ext_W^1(I_{Z/W},\OO_W)=0 .
\]
Its flat lift embeds in $\OO_{W_A}$ uniquely up to a scalar, with
flat quotient supported on finitely many points. This quotient is
the algebra of a finite flat subscheme $Z_A$ of degree $m$, and
the composite from $N_A^{-1}$ determines the unique $D_A$
containing $Z_A$. The inverse sends a flag to its ideal sheaf
with the canonical boundary identification. In the unexpanded
quartic case $\omega_D(C)=\OO_D(H_Q)$ is very ample:
the hyperplane linear forms separate every length-two subscheme
on the embedded surface. This verifies
\eqref{eq:finite-defect-evaluation} there.
\end{proof}

\begin{theorem}[The positive one-point summands as refined diagonal pullbacks]\label{thm:point-open-diagonal}
Use the notation of \eqref{eq:one-defect-three-Chern} and write
$J=J_{\chi(\OO_C)}$. On either side $a=Q,B$, prescribe the
minimal positive hyperplane or plane support, its one-point
refinements, and the boundary condition $F|_{S_o}\simeq\OO_C$,
where $C$ varies in $|h|$. The resulting one-point incidence
spaces are proper. Writing $A=\operatorname{Bl}_{[S]}|H_Q|$, they are
\begin{align*}
 \mathcal P_Q^\sigma&=\{(d,q)\in A\times Q:q\in T_d\},\\
 \mathcal P_B^\sigma&=\{([\ell],b)\in |h|\times B:\ell(b)=0\}.
\end{align*}
The first space includes the graph over $[S]$. On a fixed
expansion their formal deformations with this specified boundary
condition are the deformations of the sheaf just described.

For each $i=L,R$ there exist finite-type open neighborhoods
$U_{i,Q},U_{i,B}$ in the scalar-rigidified relative simple
perfect-sheaf stacks, with the usual two-term obstruction
theories and boundary maps to $J$. The boundary on these opens
is allowed to vary in $J$. For the defect-free branch we may
use its proper structure-sheaf component. There is an open
immersion
\[
 j_i:V_i\longrightarrow U_{i,Q}\times_J U_{i,B}
\]
and an equality in $A_0(V_i)_{\QQ}$:
\begin{equation}\label{eq:one-defect-open-diagonal}
 [V_i]^\vir=j_i^*\Delta_J^!
       \left([U_{i,Q}]^\vir\times[U_{i,B}]^\vir\right).
\end{equation}
In particular $(g_i)_*$ of the right hand side is defined
regardless of properness of $U_{i,Q}$ and $U_{i,B}$. The
left hand side is the class computed in
\eqref{eq:one-defect-three-Chern}.
\end{theorem}
\begin{proof}
The one-ended lifting and contraction of the mark in the proof of
Proposition~\ref{prop:one-defect-positive-marked} apply on each
pair independently. A mark at an internal join uses
Lemma~\ref{lem:one-defect-universal}. Forgetting the
refinement gives a point of the abstract minimal support, and
the normalized interior coordinate of the lifted point recovers
the refinement. The minimal relative support is the family over
$A$ or $|h|$ of Theorem~\companionref{thm:verifiedcompactification}.

Both displayed incidence spaces are closed in projective
products; the one-ended normalization is unique over arbitrary
base rings. The independence of the target coordinate follows
from the formula $z'=hz/h(p)$ for a change of branch equation
$x'=hx$, as in Lemma~\ref{lem:one-defect-universal}.
Proposition~\ref{prop:plane-boundary-rigidity} recovers the
divisor and point from a deformation of the sheaf with the
prescribed boundary.

Let $c:W'\to W$ contract the extra ruled component. Since $H^1(\OO_W)=H^2(\OO_W)=0$, the support line bundle on $W'$ is $c^*L$ also over an Artin base. On the closed fiber the extra ruled term cancels in the normalization sequence, giving
\[
 H^j(W,L)\xrightarrow{\sim}H^j(W',c^*L),\qquad j=0,1.
\]
For a small extension with kernel $J$ the obstruction to lifting a section lies in $H^1(L)\otimes J$ and the lifts form a torsor under $H^0(L)\otimes J$. Induction therefore gives the same isomorphism on sections over the base. Their descended zero divisor is Cartier by fiberwise regularity.

To obtain virtual classes with unrestricted boundary, take the
relative stacks of simple perfect sheaves which are normal at the
outer divisor, have stable boundary in $J$, and have trivial
infinitesimal target stabilizer after scalar rigidification.
These are open conditions near the restricted one-point sheaves.
Indeed each nonzero branch is a Cartier sheaf or its point
ideal.

For either branch sheaf $G$ on its connected surface $D$,
\[
 \mathcal End_D(G)=\OO_D,\qquad \End_D(G)=\CC,\qquad G|_C=\OO_C
\]
by depth at the point and integrality of the plane quartic $C$. A graph fixes its dilation, and on a point refinement the recovered mark fixes the dilation. Perfect base change for self-Ext, the fiberwise flatness criterion, openness of boundary stability and upper semicontinuity of stabilizer dimension give the required open neighborhoods. The bounded charts of Section~\companionref{sec:compactification} are locally of finite type with finite-type stabilizers and quasi-compact diagonal. Finitely many cover the quasi-compact images of the $V_i$.

Their union gives $U_{i,Q},U_{i,B}$.
On the branch without a point,
Propositions~\ref{prop:Bplanefamilyfull}
and~\ref{prop:Qstructurefamilyfull} give the proper open
component: flat deformations of $\OO_D$ still have their
unique generating section up to a unit, since $H^1(\OO_D)=0$.
Lifting that section over small extensions and applying Nakayama
recovers the relative Cartier quotient. Proposition~\companionref{prop:relativePOT}
gives the branch obstruction theories, with absolute virtual
dimension three.

On classical truncations Proposition~\ref{prop:one-defect-positive-marked} gives
\[
 V_L=\mathcal P_Q^\sigma\times_{|h|}|h|,\qquad
 V_R=A\times_{|h|}\mathcal P_B^\sigma .
\]
Pull the open immersion of Theorem~\ref{thm:one-defect-global} back to the marked expansion face. Branch charges are locally constant on this face; moving the point across a new level changes the marked node. Hence the two positive sectors are disjoint opens. By derived descent and Proposition~\companionref{prop:redglue}, the scalar- and determinant-reduced product over $J$ is the corresponding marked matching space. Restriction to these opens gives $j_i$.

Finally the derived diagram is the open restriction of
\[
 U_{i,Q}^{\rm der}\times^R_J U_{i,B}^{\rm der}
       \longrightarrow U_{i,Q}^{\rm der}\times U_{i,B}^{\rm der}
       \longrightarrow J\times J
\]
formed by pullback of the diagonal of the smooth sixfold $J$.
The cotangent triangle \textup{(\companionref{eq:POTdiagram})}, its compatibility
with the ordinary cotangent triangle, and smooth-base
transitivity apply on these opens. Refined pullback along
$\Delta_J$ gives the virtual class of the entire product, and
flat pullback along the open immersion gives
\eqref{eq:one-defect-open-diagonal}.

Restriction of virtual
classes commutes with a further shrinking of either neighborhood.
Properness of $V_i$ and $g_i$ is what permits the final
pushforward. Pushing each branch class separately to $J$ would
add a properness requirement on its boundary evaluation.
\end{proof}

\begin{corollary}[Proper branch classes with structure-sheaf boundary]\label{cor:one-point-independent-sigma}
Let $P=|h|$, let $\sigma:P\hookrightarrow J$ be the section
$[C]\mapsto[\OO_C]$, and set $L=b_P^*\OO_{\PP^4}(-1)$ on $A$.
On the two proper relative incidence spaces of
Theorem~\ref{thm:point-open-diagonal}, the boundary constraint
defines independent virtual classes of dimension zero
\[
 \xi_a=\iota_a^*\sigma^![U_a]^\vir
         \in A_0(\mathcal P_a^\sigma)_{\QQ},
       \qquad a=Q,B .
\]
Here $U_a$ is an unrestricted relative neighborhood of the
one-point branch, and
$\iota_a:\mathcal P_a^\sigma\hookrightarrow U_a\times_JP$
is the open inclusion. The classes are independent of the
neighborhood chosen around these incidences. If
$r':V_R=A\times_P\mathcal P_B^\sigma\to
\mathcal P_B^\sigma$ is projection and
$\mathrm{pr}_A:V_R\to A$ is the other projection, then
\[
 [V_L]^\vir=\xi_Q,\qquad
 [V_R]^\vir=c_1(\mathrm{pr}_A^*L)\cap(r')^*\xi_B,
 \qquad
 \deg\xi_Q=-320,\qquad \deg\xi_B=360.
\]
\end{corollary}
\begin{proof}
The map $\sigma$ is a regular immersion of codimension three
since $P$ and $J$ are smooth of dimensions three and six. The
proper point-incidence morphism to the constrained sheaf stack
is injective on geometric points and on rigidified stabilizers:
$D$ is recovered as the annihilator of the ideal and
$\OO_D/I_{p/D}$ recovers the point, with the expansion determined
by the support and that point. Boundary automorphisms of the
simple sheaf $\OO_C$ are scalars.

Proposition~\ref{prop:plane-boundary-rigidity} supplies the Artin inverse. Thus $\iota_a$ is representable, \'etale and universally injective, hence an open immersion between the indicated locally finitely presented moduli spaces. On the quartic side it is $j_L$ of Theorem~\ref{thm:point-open-diagonal}. For two choices of neighborhood, flat base change to their intersection identifies the refined pullbacks.

The defect-free plane component is $P$, with virtual class $[P]$
and evaluation $\sigma$. Thus the diagonal pullback in
\eqref{eq:one-defect-open-diagonal} for $V_L$ is precisely
$\xi_Q$. The defect-free quartic factor is the smooth
$\PP^1$-bundle $a:A\to P$, with evaluation $\sigma a$ and
virtual class $c_1(L)\cap[A]$ by
Proposition~\ref{prop:Qstructurefamilyfull}. Smooth pullback of
the diagram defining $\sigma^!$ by $a$, followed by the
projection formula for $c_1(L)$, gives the asserted class on
$V_R$. Finally $L$ has degree $-1$ on each projective-line
fiber, so the projective bundle formula gives
\[
 (r')_*[V_R]^\vir=-\xi_B .
\]
The branch classes themselves have intrinsic Chern formulas.
For $\xi_Q$ use $V_L=\mathcal P_Q^\sigma$ in
\eqref{eq:one-defect-three-Chern}. For $\xi_B$ write $D_B^p$ for
the inverse image of $S\subset B$ in the incidence
$\mathcal P_B^\sigma$. This is a smooth divisor, a
$\PP^2$-bundle over $S$. On the projective-line bundle
$r':V_R\to\mathcal P_B^\sigma$ the boundary $D_R$ is the
union of $(r')^{-1}D_B^p$ and one section $E$. The exact
sequence of logarithmic differentials gives
\[
 c_6(\Omega^1_{V_R}(\log D_R))
  =(r')^*c_5(\Omega^1_{\mathcal P_B^\sigma}(\log D_B^p))
     \,c_1(\Omega^1_{V_R/\mathcal P_B^\sigma}(E)).
\]
The last factor has degree $-1$ on each fiber. Applying $(r')_*$
to \eqref{eq:one-defect-three-Chern} on this marked component
therefore gives
\[
 \xi_B=c_5(\Omega^1_{\mathcal P_B^\sigma}(\log D_B^p))
        \cap[\mathcal P_B^\sigma],\qquad
 \deg\xi_B=-3\big(e(B)-e(S)\big)=360,
\]
where logarithmic Gauss--Bonnet was proved in
Proposition~\ref{prop:one-defect-zero-count}. The same
calculation on $\mathcal P_Q^\sigma$ has degree
$4(e(Q)-e(S))=-320$ by
Proposition~\ref{prop:one-defect-positive-marked}.
\end{proof}

\begin{lemma}[Properness of a projected matching space]\label{lem:proper-projected-matching}
Let $B$ be a noetherian base and let $N\to B$ be a separated algebraic space of finite type. Suppose a proposed compactification $N$ of a relative sheaf locus admits a $B$-morphism $\pi:V\to N$ from a proper marked space $V\to B$. Then $\pi$ is proper and its image is closed after any base change. If $\pi$ is surjective on geometric points, $N\to B$ is proper.
\end{lemma}
\begin{proof}
The graph and projection factor $\pi$ as
\[
 V\ \xrightarrow{\ \Gamma_\pi\ }\ V\times_B N
      \ \xrightarrow{\ \mathrm{pr}_2\ }\ N .
\]
The graph is a closed immersion by separatedness and the projection is a base change of $V\to B$. For every $B'\to B$ and closed $Z\subset |N_{B'}|$, surjectivity on geometric points gives
\[
 \operatorname{im}(Z\to |B'|)
 =\operatorname{im}\big(\pi_{B'}^{-1}(Z)\to |B'|\big).
\]
The right side is closed since $V_{B'}$ is proper. Thus $N\to B$ is universally closed; the other two conditions for properness are hypotheses.
\end{proof}
Applying Lemma~\ref{lem:proper-projected-matching} to a matching over a separated boundary moduli space $J$ requires that every geometric branch sheaf $F_1$ admit a branch sheaf $F_2$ on the opposite side with an isomorphism $\LL i_1^*F_1\simeq\LL i_2^*F_2$ and stable glued object in $V$.

For the positive one-point classes of Corollary~\ref{cor:one-point-independent-sigma} the independently constructed spaces impose $\sigma(|h|)\subset J$. Their properness supplies proper pushforwards for that boundary condition. Extending this argument to an unrestricted branch requires a separated finite-type relative moduli space and a covering by proper matching spaces at all its geometric points.

\subsection{Specialization of the first point-defect virtual class}
For the charge $v_\bullet=1-e^{-H}-[\mathrm{pt}]$, the incidence
description gives a proper model on which the specialization can be
computed as a zero-cycle. Write $C$ for a smooth pointed base curve of
the first degeneration and $g\colon\XX\to\PP^4$ for its contraction
followed by the projection. Let $V=H^0(\PP^4,\OO(1))$ and put
\[
 0\longrightarrow K\longrightarrow V\otimes\OO_{\XX}
     \xrightarrow{\mathrm{ev}}g^*\OO(1)\longrightarrow0,
 \qquad \mathcal I=\PP_{\XX}(K).
\]
Here $\PP(K)$ parametrizes lines in the fibers of $K$. For $Y=Q,B,S$
write $I_Y=\PP_Y(K|_Y)$, $p_Y\colon I_Y\to Y$, and
$i_Y\colon I_Y\hookrightarrow\mathcal I_0$.

\begin{theorem}\label{thm:first-point-log-specialization}
The morphism $\mathcal I\to C$ is projective and separated. Over a
smooth quintic fiber its fiber is the entire sheaf moduli of charge
$v_\bullet$. On the central fiber the specialization of its intrinsic
Donaldson--Thomas virtual class satisfies
\begin{equation}\label{eq:first-point-log-specialization}
 \spc_{\mathcal I}[M_\eta(v_\bullet)]^{\vir}
 =\sum_{Y=Q,B} i_{Y*}\!
  \left(c_6\big(\Omega^1_{I_Y}(\log I_S)\big)\cap[I_Y]\right)
 \quad\text{in }A_0(\mathcal I_0).
\end{equation}
The degrees of the two summands are $-320$ on the quartic component
and $-480$ on the blown-up projective-space component.
\end{theorem}
\begin{proof}
Evaluation is surjective since $\OO_{\PP^4}(1)$ is generated by $V$.
Thus $K$ has rank four. The projective bundle is smooth over $\XX$
of relative dimension three and
\[
 \mathcal I_0=I_Q\cup_{I_S}I_B,\qquad
 \mathcal I_\eta=\{(x,[s]):x\in X_\eta,\ s(x)=0\}.
\]
The second equality identifies the generic fiber with $M_\eta(v_\bullet)$
by Proposition~\ref{prop:polar-incidence-functors}. It includes
singular hyperplanes through $x$. On this smooth six-dimensional
scheme the symmetric obstruction theory has obstruction bundle
$\Omega^1_{\mathcal I_\eta}$. The obstruction bundle formula \cite[Proposition~5.6]{BF} gives
\[
 [M_\eta(v_\bullet)]^{\vir}
   =c_6(\Omega^1_{\mathcal I_\eta})\cap[\mathcal I_\eta].
\]
Let $E$ be the cokernel of the logarithmic differential
$\pi^*\Omega_C^1(\log 0)\to\Omega^1_{\mathcal I}(\log\mathcal I_0)$,
where $\pi\colon\mathcal I\to C$. At a point of $I_S$, take local
coordinates $xy=t,z_1,\ldots,z_5$. Then
\[
 \mathrm d\log t\longmapsto\mathrm d\log x+\mathrm d\log y,
 \qquad E=\langle\mathrm d\log x,\mathrm dz_1,\ldots,\mathrm dz_5\rangle,
 \qquad \mathrm d\log y=-\mathrm d\log x .
\]
Away from $I_S$ this quotient equals the relative cotangent bundle.
The branch restrictions identify canonically with
$E|_{I_Y}=\Omega^1_{I_Y}(\log I_S)$ and
$E|_{\mathcal I_\eta}=\Omega^1_{\mathcal I_\eta}$.
Both branches have multiplicity one. Specializing the preceding
top Chern class and using the projection formula gives
\[
 c_6(E|_{\mathcal I_0})\cap
      (i_{Q*}[I_Q]+i_{B*}[I_B])
  =\sum_{Y=Q,B}i_{Y*}\big(c_6(E|_{I_Y})\cap[I_Y]\big),
\]
as asserted, by \cite[\S20.3]{Fulton}.

For each $Y$, the sequence of differentials is
\[
 0\to p_Y^*\Omega^1_Y(\log S)
    \to\Omega^1_{I_Y}(\log I_S)\to\Omega^1_{I_Y/Y}\to0 .
\]
Both end terms have rank three. On a projective three-space fiber
$\int c_3(\Omega^1_{\PP^3})=-4$, hence the pushforward along $p_Y$
is $-4c_3(\Omega^1_Y(\log S))\cap[Y]$. Adjunction on $S$ and the
identity $c(T_Y(-\log S))=c(T_Y)/(1+[S])$ give its degree
$4(e(Y)-e(S))$. For $Q$, $c_3(T_Q)=-14H^3$ and $\int_Q H^3=4$.
For $B$, the blowup replaces $\Gamma$ by a projective line bundle.
Since $2g(\Gamma)-2=4\cdot5(4+5-4)=100$,
\[
 e(Q)=-56,\qquad e(B)=4+e(\Gamma)=-96,\qquad e(S)=24 .
\]
The required degrees are $4(-56-24)$ and $4(-96-24)$.
\end{proof}
The proper maps and their targets fit into the diagram
\[
 \begin{tikzcd}[column sep=3em,row sep=1.35em]
 M_1^\dagger \ar[r,"g"] \ar[d] & (M_1)_0 \ar[r,"u_0"] \ar[d] &
      \mathcal I_0 \ar[d] \\
 \EE_0^\dagger \ar[r,"\nu"'] & \EE_0 \ar[r] & \{0\}.
 \end{tikzcd}
\]
Here $g$ forgets the chosen joining K3, $\nu$ is the normalization of
the expansion divisor, and $u_0$ contracts the expansion and records
the hyperplane and the image of the point. The left square is
Cartesian. The next corollary computes the composite pushforward.

The log Chern classes in \eqref{eq:first-point-log-specialization} lie in $A_0(\mathcal I_0)$. Comparison with a relative sheaf virtual class takes place after pushforward along $u_0g$.

\begin{corollary}\label{cor:point-defect-contraction}
Let $M_1$ and $g\colon M_1^\dagger\to(M_1)_0$ be as in
Corollary~\ref{cor:one-defect-specialization}. There is a proper
contraction $u\colon M_1\to\mathcal I$, an isomorphism on the
generic fiber, such that
\[
 (u_0 g)_*[M_1^\dagger]^{\vir}
   =\sum_{Y=Q,B}i_{Y*}
      \big(c_6(\Omega^1_{I_Y}(\log I_S))\cap[I_Y]\big).
\]
\end{corollary}
\begin{proof}
Contract the marked point to $x\in\XX$ and take the image $[s]$ of
the support parameter in $|H|$ under the blowdown
$\mathscr P\to |H|\times C$. On the unmodified chart the support
is $s=0$. In the charts \textup{(\companionref{eq:vacuumtchart})} and
\textup{(\companionref{eq:vacuumsection})} the contracted hyperplane equation is
the support section multiplied by the exceptional factor, in the
chosen trivialization.

Thus $s(x)=0$ as a section on the base
scheme, including on nonreduced bases. This gives a morphism to
the incidence scheme. It is proper since $M_1$ is proper over $C$
and $\mathcal I$ is separated over $C$. On the generic fiber the
two operations are inverse to $(x,[s])\mapsto I_{x/(s=0)}$.
Compatibility of specialization with proper pushforward gives
\[
 u_{0*}\spc_{M_1}[M_\eta]^{\vir}
 =\spc_{\mathcal I}(u_\eta)_*[M_\eta]^{\vir}
 =\spc_{\mathcal I}[M_\eta]^{\vir}.
\]
Substitute the two specialization formulas.
\end{proof}

\subsection{The double-point contribution under specialization}\label{sec:double-point-direct-d}
Write $X$ for the generic quintic after base change to a strictly henselian trait, and let $V=X^{[2]}$. The surface moduli of class $\ch(\OO_T)-2[\mathrm{pt}]$ is the plane bundle $p:M\to V$ whose fiber over $Z$ consists of hyperplanes containing $Z$. The description includes length-two schemes supported at one point, by Proposition~\ref{prop:polar-incidence-functors}. We construct a representative of its virtual cycle which specializes to actual stable surface sheaves on the two roots.
\begin{theorem}\label{thm:two-points-direct}
The special-fiber contributions to the degree in this charge, indexed by the root configuration of the limiting length-two scheme, are
\[
\begin{array}{c|r}
\text{two distinct points on }Q&9600\\
\text{one point on each root}&28800\\
\text{two distinct points on }B&21600\\
\text{length two at one point on }Q&-600\\
\text{length two at one point on }B&-900
\end{array}
\qquad
 [q^{-7/24}]h^{\mathrm{DT}}_{1,0}=58\,500 .
\]
The sector labels refer to the specified moving representative in the proper Simpson closure of the generic surface-sheaf space.
\end{theorem}
\begin{proof}
Let $Y=X\times X$, $b:U=\operatorname{Bl}_\Delta Y\to Y$, $i:\Delta\hookrightarrow Y$, $j:E=\PP(T_X)\hookrightarrow U$, $p_E:E\to X$, and $q:U\to V$ the degree-two quotient. Locally a length-two scheme supported at a point is specified by its center, a tangent direction $d=(1,a,b)$ up to normalization on the chosen chart, and one discriminant coordinate $v$.

The ordered points in that chart are $m+sd,m-sd$ and $v=s^2$, where $m$ is the center. This is the double-cover model $U\to V$ in \'etale coordinates about the diagonal. In particular $V$ is smooth, $q$ is branched along $E$, with $q|_E$ an isomorphism onto its image. Put $\tau(Z)=c_{\dim Z}(T_Z)\cap[Z]$ for smooth $Z$.
We use the following identities of zero-cycles:
\begin{align}
 b_*\tau(U)&=\tau(Y)+2i_*\tau(X),\label{eq:two-ordered-blowup}\\
 2\tau(V)&=q_*\tau(U)+(qj)_*\tau(E),\qquad
 (p_E)_*\tau(E)=3\tau(X).\label{eq:two-branch}
\end{align}
For a smooth variety $W$ with smooth divisor $D$ the logarithmic tangent sequence gives
$\tau(W)=c_{\dim W}(T_W(-\log D))\cap[W]+i_{D*}\tau(D)$.
Use the deformation to the normal cone $\operatorname{Bl}_{\Delta\times0}(Y\times\mathbb A^1)$.
Its central fiber is $U\cup_E P$ with $P=\PP(T_X\oplus\OO)$. Let $\tau_{\log}(W,D)=c_{\dim W}(T_W(-\log D))\cap[W]$. On a fiber of $P\to X$ one has $c(T_{\PP^3}(-\log\PP^2))=(1+h)^3$, hence
\[
 (P\to X)_*\tau_{\log}(P,E)=\tau(X),\qquad (p_E)_*\tau(E)=3\tau(X).
\]
Specializing the relative log tangent bundle and pushing to $Y$ gives
\[
 \tau(Y)=b_*(\tau(U)-j_*\tau(E))+i_*\tau(X).
\]
This proves \eqref{eq:two-ordered-blowup}. For the double cover, $v=s^2$ gives $\mathrm d\log v=2\,\mathrm d\log s$, hence $q^*\Omega^1_V(\log qE)=\Omega^1_U(\log E)$. Push forward $\tau_{\log}(U,E)=\tau(U)-j_*\tau(E)$ and use $q|_E:E\simeq qE$ to obtain \eqref{eq:two-branch}.

The blowup formula gives $b_*:A_0(U)\to A_0(Y)$ an isomorphism, since its additional summands would be Chow groups of negative dimension on $X$. Similarly $(p_E)_*:A_0(E)\to A_0(X)$ is an isomorphism. Choose three rationally equivalent representatives $\alpha,\beta,\delta$ of $\tau(X)$, with pairwise disjoint finite support. Then $\alpha\times\beta$ is off the diagonal and has a strict lift $\widetilde{\alpha\times\beta}$ to $U$. Choose any lift $\widetilde\delta$ of $\delta$ to $E$ by choosing a tangent direction at each support point after a finite splitting extension. Equations \eqref{eq:two-ordered-blowup}--\eqref{eq:two-branch} give
\[
 2\tau(V)=q_*\widetilde{\alpha\times\beta}
                     +5(qj)_*\widetilde\delta .
\]
On $M$ put $a=c_1(\OO_{M}(1))$ from the containing hyperplane.
Smoothness and the symmetric obstruction theory identify the generic virtual class with $c_8(\Omega^1_M)\cap[M]$.
In the cotangent sequence of the plane bundle only the top vertical degree survives after multiplication by $p^*c_6(\Omega^1_V)$. The Euler sequence therefore gives
\[
 [M]^\vir=3a^2\cap p^*\tau(V).
\]

It remains to choose the representatives over the degeneration. Let $\mathcal E$ be the rank-three relative logarithmic tangent bundle of the smooth semistable total threefold family. On each root $Y_i$ it restricts to $T_{Y_i}(-\log S)$, and over the generic point it restricts to $T_X$. For a relatively very ample bundle of class $z$ and sufficiently high degree,
\[
 c_3(\mathcal E)=c_3(\mathcal E(z))-z\,c_2(\mathcal E(z))
                      +z^2c_1(\mathcal E(z))-z^3 .
\]
Each term on the right is represented by degeneracy loci of generally chosen sections of the generated rank-three bundle and independent hyperplane divisors. Choose the three tuples generally so that the intersections have dimension zero on each special component, are smooth there, and avoid the joining divisor, exceptional divisors, and the finite supports from the preceding tuples.

This is possible since evaluation of the sections is surjective, rank-drop loci have their determinantal codimensions, and any imposed bad locus on a component has dimension at most two. Thus they extend as signed finite horizontal cycles and give $\alpha,\beta,\delta$ of the required kind.

If pairs of their special points are required to avoid an additional proper closed subset of a product of components, choose the tuples independently: over that subset the two sets of rank-drop and hyperplane conditions have total codimension six, exceeding its dimension. Their degrees on each root are independent of the tuple:
\[
 d_Q=\int_Q c_3(T_Q(-\log S))=e(Q)-e(S)=-80,\qquad
 d_B=e(B)-e(S)=-120 .
\]
Over each section of $\delta$ take a general tangent direction whose image in $\PP^4$ is nonzero. In smooth relative coordinates its double point is given by
\[
 x_i\longmapsto x_i(p)+v_i\epsilon,\qquad \epsilon^2=0 .
\]
The ordered-pair term uses two disjoint sections. In either case the central length-two scheme $Z$ is contained in the smooth locus and spans a line, so the hyperplanes through $Z$ form a plane. Two further general point conditions representing $a^2$ select one member.

Bertini gives smooth branch sections away from the secant or tangent line. At the finite base points, singularity requires the chosen hyperplane to contain the tangent space there. If the secant lies in the plane root, its section is a plane transverse to $\Gamma$ for a general choice. Thus the support is a smooth quartic surface and a plane transform meeting in a smooth quartic $C$. The ideal sequence gives a flat family, with
\[
 I_{Z/T}|_C=\OO_C,\qquad P_{I_{Z/T}}=P_{\OO_T}-2 .
\]
A punctual modification preserves the strict linear stability inequalities. These families therefore specialize the chosen virtual zero-cycle in the relative Simpson closure.

For this representative the ordered-pair and tangent-direction terms contribute
\[
 \tfrac32(d_Q+d_B)^2+\tfrac{15}2(d_Q+d_B)
 =9\,600+28\,800+21\,600-600-900=58\,500 .
\]
The five summands describe this moving representative; their sum is the specialized degree.
\end{proof}

\subsection{Two different points on a single exceptional ruling}
\label{sec:double-point-rank-drop}
The plane-bundle proof of the smooth-fiber two-point calculation has two distinct boundary difficulties: projected points can cease to impose independent hyperplane conditions, and the vanishing that recovers the point presentation from arbitrary infinitesimal sheaf deformations can fail. Both occur at two reduced, distinct points.

\begin{proposition}\label{prop:exceptional-pair-hidden-tangent}
Let $\gamma\in\Gamma$ be general and choose two distinct points $z_1,z_2$ on its exceptional ruling $f_\gamma\subset B$, off $S$. There exists an ordinary positive support $T$ containing them with smooth component surfaces, smooth joining quartic, and stable admissible sheaf $F=i_*I_{Z/T}$ on $W=Q\cup_S B$, where $Z=z_1+z_2$. On the fixed target $W$ the incidence of hyperplanes and reduced unordered pairs is smooth of dimension eight at $(T,Z)$. The obstruction space for the extension presentation used in the one-point charge has dimension one:
\[
 \dim\Ext_W^2(F,\OO_W(-T))=1 .
\]
\end{proposition}
\begin{proof}
Choose a general plane $\Pi\subset\PP^3$ through $\gamma$ transverse to $\Gamma$, its plane section of $S$ smooth. One can choose such $\Pi$: at $\gamma$ only the planes containing $T_\gamma\Gamma$ are excluded by the first condition, and only the tangent plane to $S$ is excluded by the second. Away from its single base point use Bertini. Extend $\Pi$ to an ordinary hyperplane in $\PP^4$ varying the lift generally, so its quartic surface section is smooth.

Bertini on the whole system through $\gamma$, and avoidance of $T_\gamma Q$, give this last assertion simultaneously. Its transform in $B$ is the blowup of $\Pi$ at its twenty transverse intersections with $\Gamma$, containing $f_\gamma$. Thus $Z$ lies in the smooth surface and threefold loci, off their joins. The ideal restricts to $\OO_C$ at the joining curve $C$.

The strict linear cut inequalities for the structure sheaf persist after subtracting two points. Near each $z_i$ the surface is smooth; the complete fixed-target flag space is thus smooth over the four smooth hyperplane parameters there, with fiber dimension four. This proves its dimension and smoothness.

The evaluation of hyperplane sections on $Z$ has rank one. In local blowup coordinates the map to $\PP^3$ reads $(a,b,c)\mapsto(a,ab,c)$ with distinct ruling coordinates $b_1,b_2$ at the two points. In the ordered-pair chart the ideal where their projections coincide is
\[
 (a_1-a_2,\ a_1b_1-a_2b_2,\ c_1-c_2)
       =(a_1,a_2,c_1-c_2);
\]
the equality uses that $b_1-b_2$ is a unit. Thus even at distinct points there is a smooth codimension-three locus on which the two projected point conditions coalesce.

Adjunction and the hyperplane sequence give $\omega_T=\OO_T(T)$ and $h^0(\omega_T)=4$; its sections are ambient forms modulo the support equation. Normalization along the connected quartic gives $H^1(\OO_T)=0$. Since $Z$ is in the smooth locus, dualizing its ideal sequence gives
\[
 0\longrightarrow\Ext_T^1(I_{Z/T},\OO_T)\longrightarrow
 \Ext_T^2(\OO_Z,\OO_T)\xrightarrow{d}H^2(\OO_T).
\]
The dual map $d^\vee:H^0(\omega_T)\to H^0(\omega_T|_Z)$ evaluates sections at two points projecting to $\gamma$. Hence $\operatorname{rk}d^\vee=1$. By Cartier adjunction,
\[
 \ext^2_W(F,\OO_W(-T))
   =\ext^1_T(I_{Z/T},\OO_T)
   =2-\operatorname{rk}d^\vee=1 .
\]
The obstruction associated to a given deformation is the image of its tangent class in this one-dimensional space.

\end{proof}

\section{Specialization of a line defect}
The line ideal changes in two different ways under the first degeneration. The clean excess loci of the classical line degeneration were studied by Katz \cite{KatzLines}. We compute their degrees here and keep track of the containing divisor, its sheaf and the class in the proper sheaf moduli.
\subsection{A line defect and the specialization of its surface-sheaf virtual class}
\label{sec:direct-line}
Take the general one-parameter degeneration $f_1f_4=tf_5$ with
semistable model $\mathcal X$, quartic component $Q$, and component
$B=\operatorname{Bl}_{\Gamma}\mathbb P^3$ where
$\Gamma=V(f_1,f_4,f_5)$. Assume $S=V(f_1,f_4)$ contains no line.
Work over the strictly henselian trait at $t=0$, with residue field
$\mathbb C$. All numerical identities also hold after passage to its
completion.

Write $H$ for the pullback of the ambient hyperplane
class, and set $A_m=mH+Q$. This is a relative polarization for all
sufficiently large $m$: on $Q$ its restriction is $(m-1)H_Q$,
and on $B$ it is $(m+4)H_B-E_B$. The latter is ample for $m\gg0$
since $-E_B$ is ample for $B\to\mathbb P^3$. Ampleness on every
irreducible component implies ampleness on the union.

On the smooth generic fiber $X_\eta$, let $M^-_\eta$ be the
moduli of the line-defect surface sheaves $i_*I_{L/T}$ and
$M^+_\eta$ its dual sector $i_*(I_{L/T})^\vee$, with $T$
a hyperplane section and $L\subset T$ a line. These descriptions
exhaust the sheaves in their respective characters, as schemes, by
Proposition~\ref{prop:polar-incidence-functors}. Denote by
$\overline M^\pm$ the schematic closures of these moduli in the
corresponding projective relative Simpson spaces for $A_m$.

On the generic fiber stability for $A_m=mH$ equals stability for
$H$. Pure two-dimensional sheaves of the indicated character have
fixed Hilbert polynomial; specifying the full generic numerical
character picks an open and closed union inside that generic moduli.
The notation $\spc$ below means Chow specialization to
$A_0(\overline M^\pm_0)_{\mathbb Q}$. We describe its value on
the complete generic virtual class by choosing a moving
representative of that class.

\begin{theorem}[Specialization for the two primitive line-defect charges]
\label{thm:direct-line}
Choose a general line $\Lambda\subset\mathbb P^4$. There are
two reduced finite sets $D_Q,D_B$ on the special fiber. They and
their associated surface sheaves are as follows.
\begin{enumerate}
\item $D_Q$ parametrizes lines $L\subset Q$ such that their
point $p=L\cap S$ lies in $\Gamma$. Let
$U_L=\langle\Lambda,L\rangle$, a hyperplane. In the central
surface $T_L=(Q\cap U_L)\cup(B\cap U_L)$ take
$Z_L=L\cup E_p$, with $E_p$ the exceptional fiber of $B$ over $p$.
\item $D_B$ parametrizes lines $L\subset\mathbb P^3$ for which
$f_5|_L$ is divisible by $f_4|_L$. In the same surface $T_L$
take $Z_L$ to be the strict transform of $L$ on $B$.
\end{enumerate}
With $i_L:T_L\hookrightarrow\mathcal X_0$ the inclusion, one has
\begin{align}
 \spc[M^\pm_\eta]^\vir
 &=3\sum_{L\in D_Q}[i_{L*}\OO_{T_L}(\pm Z_L)]
   +3\sum_{L\in D_B}[i_{L*}\OO_{T_L}(\pm Z_L)]
 \quad\text{in }A_0(\overline M^\pm_0)_{\mathbb Q}.
 \label{eq:direct-line-chow}
\end{align}
The point notation denotes the class of the corresponding closed
point of the special moduli. The degrees of the two sums with
their coefficients are $4800$ and $3825$. In particular the
quintic DT coefficient at $q^{-83/120}$ is $8625$ in each of
residues $1$ and $4$.
\end{theorem}
\begin{proof}
\emph{The moving virtual cycle.}
Let $G=\operatorname{Gr}(2,5)$, $U$ its tautological rank-two
bundle and $E=\operatorname{Sym}^5U^\vee$. A generic quintic section
is transverse on $G$, so its vanishing scheme, which describes the
possible defects, is finite and \'etale. After splitting that scheme,
each connected piece of the surface-sheaf moduli is precisely the
plane of hyperplanes containing its line.

The symmetric obstruction
theory of stable sheaves on the smooth Calabi--Yau fiber restricts
on a smooth such piece to obstruction bundle $\Omega^1_{\mathbb P^2}$.
For a smooth moduli space the virtual class is the top Chern class
of its obstruction bundle capped with its fundamental class. Thus,
if $a$ is the pullback of $c_1(\OO(1))$ under the support map
to $(\mathbb P^4)^\vee$, the Euler sequence gives
\[
 [M^\pm_\eta]^\vir=c_2(\Omega^1_{M^\pm_\eta})\cap[M^\pm_\eta]
 =3a^2\cap[M^\pm_\eta].
\]
This equality can equally be obtained from the relative Euler
sequence over the finite \'etale scheme of defects. Choose two points
spanning $\Lambda$. The two conditions that the support hyperplane
contain these points cut each plane of supports in one reduced point,
as soon as $\Lambda$ and the defect line are disjoint. This choice
of representative retains the virtual multiplicity $3$. The number
of components of $M_\eta^\pm$ has not entered the computation.

\emph{All horizontal limits of the selected supports.}
On $G$ the equations for the defect are the six coefficients of
$s=f_1|_U f_4|_U-t f_5|_U\in E$. At $t=0$ its zero scheme has two
disjoint smooth clean components
\[
 Z_Q=Z(f_4|_U),\qquad Z_B=\operatorname{Gr}(2,\ker f_1).
\]
They are disjoint since $S$ contains no line. The quartic is
general and $\operatorname{Sym}^4U^\vee$ is generated by ambient
quartic forms; hence $Z_Q$ is smooth of the expected dimension one.
Along $Z_Q$ the linearization has image
$f_1|_U\operatorname{Sym}^4U^\vee$ of rank five; along $Z_B$
it has image $f_4|_U U^\vee$ of rank two. In each case the
normal bundle of the respective zero locus maps isomorphically
onto this image.

Write $V$ for the cokernel. The restriction of
$f_5$ induces a section of $V$, generated by the ambient quintic
parameters. Its zeros in each smooth component are therefore
reduced and zero-dimensional for general $f_5$.
On $Z_Q$ the cokernel fiber is $\OO_p(5)$ at $p=L\cap S$,
and on $Z_B$ vanishing is the stated divisibility condition.

For the scheme structures and the multiplicities of specialization,
choose local coordinates $(y,z)$, where $z$ are normal coordinates
to either component. Resolving the equations in the image of the
linearization by the formal implicit function theorem gives
$z=t b(y,t)$. After substitution the residual equations are
\[
 t\bigl(-\overline f_5(y)+t h(y,t)\bigr)=0
 \quad\text{in the cokernel coordinates}.
\]
On the horizontal closure, obtained by saturation by $t$, the
residual ideal is the one generated by
$-\overline f_5(y)+t h(y,t)$ near its zeros. Indeed the Jacobian
with respect to $y$ there is invertible, so these equations give a
smooth section over the trait. Away from the zeros one equation in
parentheses is a unit, so there is no horizontal limit.

Properness
of $G$ exhausts the limits. The horizontal defect scheme is
therefore finite \'etale over the trait after shrinking, with one
section through each point of $D_Q\amalg D_B$. In particular all
multiplicities at this stage are one. The resulting finite list
allows $\Lambda$ to be chosen simultaneously so that $U_L$ is
defined for every limit and $U_L\ne\{f_1=0\}$.

\emph{Sizes of the two special loci.}
On $Z_Q$ the map $e:Z_Q\to S$ sending a line to its intersection
with $S$ satisfies $e^*\OO_S(1)=\det U^\vee$: in
$0\to\ker(f_1|_U)\to U\to\OO\to0$ the kernel represents that point.
For $x,y$ the formal Chern roots of $U^\vee$ the Grassmannian
integration formula gives
\begin{align*}
 |D_Q|
 &=5\int_G c_5(\operatorname{Sym}^4U^\vee)(x+y)\\
 &=5[x^4y^3](x-y)32xy(x+y)^2(3x+y)(x+3y)=1600.
\end{align*}
On $G_B=\operatorname{Gr}(2,4)$ there is a vector bundle sequence
\[
 0\longrightarrow U^\vee\xrightarrow{\cdot f_4|_U}
       \operatorname{Sym}^5U^\vee\longrightarrow V\longrightarrow0 .
\]
Multiplying the Chern roots $1+ix+(5-i)y$ for $0\le i\le5$ and
dividing by $(1+x)(1+y)$ gives in degree four
\[
 c_4(V)=120(x^4+y^4)+1290(x^3y+xy^3)+2565x^2y^2 .
\]
Thus $|D_B|=[x^3y^2](x-y)c_4(V)=1275$.
All lines in $D_B$ meet $S$ in four distinct points for general data:
the tangent-line locus is a proper closed subset of $G_B$, and the
incidence where the globally generated rank-four section vanishes
over this subset has dimension strictly less than the parameter
space of quintics. All four points lie on $\Gamma$.

\emph{Extension as surface sheaves.}
For each central line $L$ the containing hyperplanes form a $\PP^2$. On a quartic-root line the tangent hyperplanes $T_xQ$, $x\in L$, occupy at most a curve in this plane. Bertini away from $L$ and the tangent test at the finite base points on $S$ give smooth sections $K\subset Q$ and $C\subset S$. For a plane-root line the same test applies at its transverse intersections with $S$. The plane section on the plane root meets $\Gamma$ transversely: at a required point the line is transverse to $S$, or the plane varies through that single point. Thus
\[
 P=\operatorname{Bl}_{\Pi\cap\Gamma}\Pi,\qquad
 \#(\Pi\cap\Gamma)=20,\qquad K\cap P=C .
\]
The general $\langle\Lambda,L\rangle$ range over a dense open of the containing hyperplanes, so these transverse choices hold simultaneously for the finitely many lines.

At $C$, smoothness of the section gives local coordinates with target $xy=t$ and hyperplane equation a coordinate along the double locus. Away from $C$ the component surfaces are smooth. Thus the closure $\mathcal Z$ of the moving line is a prime Cartier divisor in the regular $\mathcal T$.

For a quartic-root line, the center ideal $(f_1,t)$ restricts on its moving line surface to the parameter ideal at $(p,0)$ and to the unit ideal elsewhere. The strict transform is the blowup there, with $\mathcal Z_0=L\cup E_p$, both of multiplicity one. For a plane-root line,
\[
 f_1|_{L_t}=t u_t,\qquad u_0=f_5|_L/f_4|_L\in H^0(L,\OO_L(1)).
\]
The center ideal is $(t)$, so the lift is the graph of $u_0$ on the central fiber, the strict transform of the four-secant disjoint from $S$.

Since $\mathcal Z$ contains no component of $\mathcal T_0$, its restriction $Z_L$ is Cartier on $K\cup_C P$, also at the node over $p$. The sheaves $i_*\OO_{\mathcal T}(\pm\mathcal Z)$ are flat; on either branch their common restrictions are $\OO_C(\pm p)$ for a quartic-root line and $\OO_C$ for a plane-root line.

\emph{Stability and Chow specialization.}
The sheaves just obtained are pure because $K\cup_C P$ is a Cartier
surface and their restrictions to it are invertible. For their
two possible generic multiranks write $Z_K,Z_P$ for the two pieces
of $Z_L$, allowing $0$, and $d_i=H\cdot Z_i$. Then $d_K+d_P=1$.
Since both components are integral, a subsheaf supported on one
component is contained in $\OO_i(\pm Z_i-C)$; this follows by
restricting the glued line bundle to the complementary component.

A proper subsheaf of full multirank has quotient of dimension at
most one and hence has strictly smaller reduced Hilbert polynomial.
For the one-component possibilities it suffices to compare the
linear terms after dividing the quadratic leading coefficient by
a common convention. Explicitly, put
$\lambda_H(F)=a_1/a_2$ for
$P_F(k)=a_2 k^2/2+a_1 k+O(1)$. The formal coefficients at the
nef pullback $H$ give
\begin{align*}
 \lambda_H(\OO_{T_L}(\pm Z_L))&=(\pm1-5/2)/5,\\
 \lambda_H(\OO_K(\pm Z_K-C))&=(\pm d_K-4)/4,\\
 \lambda_H(\OO_P(\pm Z_P-C))&=\pm d_P-5/2 .
\end{align*}
Indeed $H_K^2=4$, $H_P^2=1$, $H_K K_K=0$, $H_P K_P=-3$,
and $H\cdot C=4$; subtract $\chi(\OO_C(\pm Z_L)(kH))$
once in the normalization sequence. Both branch slopes are
strictly smaller than the global slope for either sign and either
$(d_K,d_P)=(1,0),(0,1)$. Continuity gives these inequalities for
$A_m/m=H+Q/m$ when $m\gg0$.

Relative ampleness gives Gieseker stability. All the moving representatives therefore extend to the stable locus of $\overline M^\pm$. Intersecting their closures with $t=0$ gives \eqref{eq:direct-line-chow}, of degree $3(1600+1275)=8625$.
\end{proof}
Equation~\eqref{eq:direct-line-chow} is a Chow specialization in the projective Simpson closure. A factorization into unrestricted relative virtual classes further requires identification with the marked diagonal pullbacks.

\begin{figure}[htbp]\centering
\begin{tikzpicture}[x=1cm,y=1cm,font=\footnotesize,line join=round,line cap=round]
\node[font=\small\bfseries] at (2.79,3.37) {A quartic line};
\node[font=\small\bfseries] at (9.45,3.37) {A four-secant};
\path[draw=qIIteal!48,fill=qIIteal!9] (2.19,.98)--(3.09,1.41)--(3.44,2.89)--(2.54,2.46)--cycle;
\draw[qIIblue,line width=1.3pt] (.34,1.16)--(2.78,1.92);
\draw[qIIclay,line width=1.35pt] (2.78,1.92)..controls (3.25,2.11) and (4.25,2.18)..(5.03,2.53);
\fill[qIIteal] (2.78,1.92) circle(2.8pt);
\node[below,text=qIIblue] at (1.22,1.41) {$L\subset Q$};
\node[below,text=qIIclay] at (4.06,2.08) {$E_p\subset B$};
\node[above left,text=qIIteal] at (2.77,1.94) {$p$};
\node[text=qIIteal] at (3.27,2.97) {$S$};
\draw[qIIteal,line width=1.05pt] plot[domain=7.06:11.92,samples=100]
 (\x,{1.85+.45*sin((\x-7.46)*140)});
\draw[qIIblue,line width=1.2pt] (7.02,1.85)--(11.93,1.85);
\foreach \x in {7.46,8.745714,10.031429,11.317143}{
 \fill[qIIclay] (\x,1.85) circle(2.25pt);}
\node[below,text=qIIblue] at (7.25,1.68) {$L$};
\node[text=qIIteal,anchor=south] at (9.43,2.52) {$L\cap S\subset\Gamma$};
\node at (2.79,.53) {$p=L\cap S\in\Gamma$};
\node at (9.47,.53) {$\widetilde L\cap S=\varnothing$ on $B$};
\end{tikzpicture}
\caption{The two central line defects of Theorem~\ref{thm:direct-line}. The first acquires an exceptional ruling. In the second panel the four marked intersection points lie on the blowup center. Their exceptional lifts separate the strict transform of the line from the strict transform of $S$.}\label{fig:II-special-line-defects}
\end{figure}
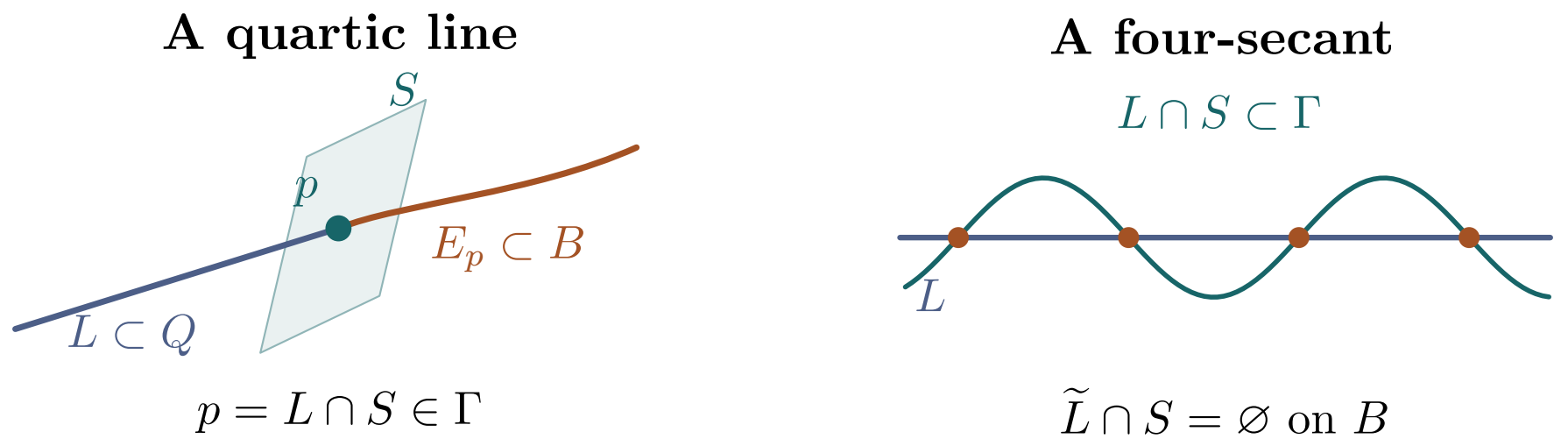

\subsection{Deformations of the line ideal at the boundary}
In this subsection $W$ is a fiber of the expanded target, $T$ an
active hyperplane support on it and $D\subset T$ the defect curve.
The ordinary split curve is $D=L\cup f_p$, with $L\subset Q$ a
line, $p=L\cap S\in\Gamma$, and $f_p$ the exceptional fiber of $B$.
The other curves are strict transforms $D=L'$ of four-secant lines
in $B$ and are disjoint from $S$. We work with general defining
forms and $\operatorname{Pic}(S)=\mathbb Zh$.

\begin{lemma}\label{lem:linedefect-target-rigidity}
Each of these curves is a regular immersion of codimension two in $W$.
It deforms uniquely with a given infinitesimal deformation of the
target over the expansion charts.
\end{lemma}
\begin{proof}
For a split curve the two branches meet $S$ with contact one.
In the completed local threefold
$\mathbb C[[x,y,u,v]]/(xy)$ their union can be straightened to
$u=v=0$: subtract the $x$- and $y$-dependent graph coordinates
of its two branches from each of $u,v$. This change of coordinates
has invertible linear part; on $x=y=0$ choose $u,v$ as parameters
on $S$. The equations $u,v$ form a regular sequence. Restricting
the normal sheaf to the two branches and using normalization gives
\[
 0\to N_{D/W}\to N_{L/Q}\oplus N_{f_p/B}\to T_pS\to0
\]
with the branch bundles understood as their pushforwards. The last map is $\operatorname{ev}_L-\operatorname{ev}_{f_p}$. Transversality of the section of $\operatorname{Sym}^4 U^\vee$ on $\operatorname{Gr}(2,5)$ gives $h^0(N_{L/Q})=1$. Riemann--Roch gives $\chi(N_{L/Q})=-K_Q L=1$. For the exceptional ruling the normal sequence has the form
\[
 0\to\OO_{f_p}\to N_{f_p/B}\to\OO_{f_p}(-1)\to0.
\]
Thus both branch bundles have $(h^0,h^1)=(1,0)$. The map on sections is the derivative of the intersection defining $F_1(Q)\times_S\Gamma$. A general divisor in $|5h|$ cuts the nonconstant evaluation branches transversely and misses each constant branch. Therefore
\[
 0\to H^0(N_{D/W})\to\CC^2\xrightarrow{\sim}T_pS
       \to H^1(N_{D/W})\to0.
\]

On $B\setminus S$ the lift of a projected line is the graph of $f_5/f_4$. In an infinitesimal deformation the projection is still a line embedding, by base change for its complete pencil $\OO(1)$. Since $S$ contains no line, multiplication by $f_4$ is a subbundle injection on $\operatorname{Gr}(2,4)$. The lifts are given by
\[
 [f_5|_L]=0\ \text{in }(\operatorname{Sym}^5U^\vee/f_4U^\vee)|_L,
 \qquad f_5|_L=f_4|_L\,u_L,\quad u_L\in H^0(L,\OO_L(1)).
\]
This general section has transverse zeros. Hence $H^0(N_{L'/B})=0$ and $\chi(N_{L'/B})=(-K_B)\cdot L'=S\cdot L'=0$, giving $H^1(N_{L'/B})=0$. The expansions are isomorphisms near $L'$. For both regular immersions the obstruction and the torsor of lifts across a small extension with kernel $J$ are governed by $H^1(N)\otimes J$ and $H^0(N)\otimes J$, respectively; both vanish.
\end{proof}

\begin{lemma}\label{lem:lineideal-recovery-on-expanded-support}
At each hyperplane support containing one of these curves, every
infinitesimal flat deformation of the sheaf $I_{D/T}$ on a deforming
target recovers a unique deformation $D'\subset T'$ of the flag.
The only ambiguity in its inclusion is multiplication by a unit on
the parameter base.
\end{lemma}
\begin{proof}
The regular immersion in Lemma~\ref{lem:linedefect-target-rigidity}
and the Cartier divisor sequence imply that $F=I_{D/T}$ has projective
dimension one on $W$. In local coordinates, write
$I_{D/W}=(u,v)$ and write a local equation of $T$ as $au+bv$.
With $t=au+bv$ the presentation and its determinant are
\[
 \OO_W^2\xrightarrow{M}\OO_W^2\twoheadrightarrow F,\qquad
 M=\begin{pmatrix}-v&a\\ u&b\end{pmatrix},\qquad
 M\,\operatorname{adj}(M)=-(au+bv)\,\mathrm{id}.
\]
Since $t$ is regular, $M$ is injective. Flatness and Nakayama lift this resolution on a deformed target. Its determinant lifts $-t$ and defines a relative Cartier divisor $T'$ annihilating the deformed sheaf. On the fiber, $F$ is faithful on $T$: it is invertible away from $D$, and $\OO_T$ is Cohen--Macaulay with no associated point in $D$. Flatness gives faithfulness on successive graded pieces of the Artin base; induction on its length then gives $\operatorname{Ann}(F')=I_{T'}$.

Put $A=\omega_T|_D$. In the ordinary split profile it has multidegree
$(1,0)$ on $L\cup f_p$, in the other ordinary profile it is
$\OO_{L'}(1)$. In the exceptional one-bubble chart the extra twisting
divisor of the support lies on the unused quartic root and misses
$L'$. In the $t$-chart with support of class $5H_B-E_B$,
adjunction gives $\omega_T=\OO_T(H_B)$. Thus in all cases
$h^0(A)=2$ and $H^1(A)=0$. Further,
$H^0(T,\omega_T)\to H^0(D,A)$ is onto.

On the ordinary
profile linear forms restrict to the complete pencil on the
degree-one component and to its value on an exceptional tail.
On the exceptional graph and $t$-charts the four plane-coordinate
sections restrict to the complete pencil on $L'$; equivalently
restrict the basis of the hyperplane system in
Lemma~\companionref{lem:vacuumlinearsystems}. The divisor sequence on the
ambient Gorenstein complete target with trivial dualizing sheaf
and $H^1(\OO_W)=0$ identifies its sections modulo the defining
section with $H^0(T,\omega_T)$.

The curve is Cohen--Macaulay, so duality on $T$ gives
$R\mathcal Hom_T(\OO_D,\OO_T)=(\omega_D\otimes A^{-1})[-1]$.
Here $H^0(\omega_D\otimes A^{-1})=H^1(A)^\vee=0$, and the dual of
$H^1(\omega_D\otimes A^{-1})\to H^2(\OO_T)$ is the restriction
surjection just proved. Applying $R\operatorname{Hom}_T(-,\OO_T)$
to the ideal sequence, with $H^1(\OO_T)=0$, gives
\[
 \operatorname{Hom}_T(I_{D/T},\OO_T)=\mathbb C,\qquad
 \operatorname{Ext}^1_T(I_{D/T},\OO_T)=0 .
\]
At each square-zero step these two equalities lift the inclusion to
the lifted Cartier support and give uniqueness up to scalar.
The criterion for flat cokernels makes it injective with flat
cokernel; its image is an ideal. This recovers the flag over all
Artin test rings and also proves its recoverability on geometric
points. No invertibility of the ideal on a singular hyperplane was
assumed.
\end{proof}

For a four-secant line branch, trivialize the rank-two bundle of the
moving lines over the local base. Let $u_t=f_1|_{L_t}/t$, which
is regular since $L_0\subset\{f_1=0\}$. The condition on the
unexpanded hyperplane $f_1+\sum b_i\ell_i$ to contain the line is
\[
 t u_t+\sum_i b_i\ell_i|_{L_t}=0
 \quad\text{in }H^0(L_t,\OO(1)).
\]
The map from the four plane coordinates to this rank-two module is
surjective. Hence one can eliminate two variables, and the
restriction of $(t,b_1,\ldots,b_4)$ to this incidence is generated
by $t$ and the two remaining variables. The proper transform is
the blowup of the smooth $\mathbb P^2$-bundle at its central point.
On the $s,v$ chart the two equations after division by $s$ become
$v u_t+\ell_c|_{L_t}=0$.

In particular at $s=v=0$ they are
$\ell_c|_{L_0}=0$, so that the bubble graph meets the plane along
its quartic curve while the line defect remains on $B$ away from
that join. In the $t$-chart the equations are $u_t+\ell_c|_{L_t}=0$.
This verifies that the strict transform represents the actual
incidence on both exceptional charts, including nilpotent bases.
For the split curve the image line lies on $Q$ and cannot lie in
the hyperplane of $S$. The incidence avoids the center; its parameter
is a smooth $\mathbb P^2$-bundle over the base.

Write $q/M$ for the normalized polarization on the
quartic or active graph branch. With the notation of
\textup{(\companionref{eq:vacuumgraphpolys})}, the degree of the defect on the whole
curve is $M$. Its degree $d_1$ on the left active branch equals
$q$ for the split curve and $0$ for a four-secant defect.
The new linear cut coefficient for the ideal equals
\[
 [m]D_F=\frac{2q^2}{M}-d_1+\frac{4q^2}{5M^2}M
         =\frac{14q^2}{5M}-d_1.
\]
For instance choose $5/14<q_0/M<q_1/M<1$. Then
$0<[m]D_F<4q$ in both profiles, so the cut inequalities are
strict. A zero-boundary support is integral and its rank-one ideal
is stable. In the dual case $[m]D_{F^\vee}=6q^2/(5M)+d_1$, with
the same strict bounds. A punctual quotient does not change the
linear cut coefficient.

\begin{theorem}\label{thm:expanded-line-compact}
Let $R$ be the strictly henselian trait of
Theorem~\ref{thm:direct-line}. There is a proper separated
scalar-rigidified moduli space $M_L\to\operatorname{Spec}R$ on the
stack of expansions, with full generic stable-sheaf moduli in the
line-ideal character.

Its connected pieces over $R$ are the
$\PP^2$-bundles in the split case and the blowups, at one closed
point, of such bundles in the four-secant case.
Every represented flat deformation on an expanded target is the
deformation of the surface sheaf itself. The marked boundary
therefore carries the restriction of its determinant-reduced
threefold obstruction theory.
\end{theorem}
\begin{proof}
By Theorem~\ref{thm:direct-line} the horizontal defect lines give distinct \'etale sections. Write $p:C_R\to\Spec R$ for the embedded family over one section. On a split fiber $H$ has multidegree $(1,0)$ and its complete pencil contracts the exceptional tail. On a four-secant fiber it embeds the projected line. In both cases $h^0(H|_C)=2$ and $h^1(H|_C)=0$. The containing hyperplanes therefore form $\PP_R(E)$, where
\[
 0\to E\to H^0(\PP^4,\OO(1))\otimes R\to p_*(H|_{C_R})\to0,
 \qquad \operatorname{rk}E=3 .
\]
For a split family this incidence misses the expansion center. For a four-secant, the two preceding charts identify its strict transform with the blowup of $\PP_R(E)$ at one reduced point of the closed fiber.

The defect family is disjoint from the expansion centers in
this latter case (the strict transform is disjoint from $S$ on the
closed fiber); its inverse image is still the defect family.
The scheme obtained over all sections is proper and separated.
Its relative Cartier support and regular defect family define a
flat perfect sheaf by the ideal sequence.

On a split join the
coordinates of Lemma~\ref{lem:linedefect-target-rigidity} make
the curve sequence admissible: the restricted one-branch
sequence is cut out by the regular parameter normal to $S$ on
that branch. The relative boundary of the ideal is consequently
the ideal of the contact point in the plane quartic. At other
joins it is the structure sheaf of the boundary quartic, or it
vanishes at an empty outer end.

The strict cut inequalities above
apply even at singular points on those quartics. A graph bubble
has trivial dilation stabilizer since its support is a nonzero
graph. The connected rank-one ideal has only scalar endomorphisms
by stability. Thus the relative sheaf defines a morphism of
this scheme to the rigidified expanded simple-sheaf moduli.

Every closed point recovers its support and defect by
Lemma~\ref{lem:lineideal-recovery-on-expanded-support}. On an
Artin lifting of the expanded target the recovered curve has a
unique lift by Lemma~\ref{lem:linedefect-target-rigidity}.
On the split chart this is the pullback of $C_R$; on the other
charts it is the same pullback on the unmodified neighborhood
of $L'$. The recovered support parameters and their changes
of gauge are those of the minimal-support chart: the determinant
line bundle lifts uniquely since $H^1(\OO_W)=H^2(\OO_W)=0$,
and its sections have the bases of
Lemma~\companionref{lem:vacuumlinearsystems}.

Imposing containment of
the unique lifted curve gives exactly the incidence equations
already constructed, including the divided equation on the graph
chart. Conversely these equations produce the ideal under
arbitrary base change. The morphism is representable, injective
on geometric points and formally \'etale by this inverse. Hence it
is an open immersion into the stack of perfect sheaves over
expansions after scalar rigidification.

Its proper source, with
the induced sheaf groupoid, is the asserted compactification.
This argument also specifies the boundary objects of that
groupoid. Surjectivity onto the entire stack of stable sheaves
on arbitrary expansions is unnecessary. The complete generic
identification comes from Proposition~\ref{prop:polar-incidence-functors}.
Each chart has a universal flat perfect sheaf locally up to tensoring by a base line bundle.

The complex of self-extensions and its trace descend under that ambiguity. The recoverability and lifting arguments give the requisite openness; the plane-bundle and blowup descriptions give boundedness, existence and uniqueness over traits for this chosen groupoid. The cohomology of the expanded targets is given by Proposition~\companionref{prop:structurecohom}. Thus Theorem~\companionref{thm:main} supplies the relative perfect obstruction theory on this space. Applying its virtual pullback to the normalization
of $t=0$ gives
\[
 \spc[M_{L,\eta}]^\vir=g_*[M_L^\dagger]^\vir
   \quad\text{in }A_0(M_{L,0})_{\QQ}.
\]
Here $M_L^\dagger$ marks a join and $g$ forgets it.
\end{proof}

\begin{corollary}\label{cor:line-opposite-expanded}
The opposite polar line class has the same compactification and marked
zero-cycle degrees. 
\end{corollary}
\begin{proof}
Let $j:T\hookrightarrow W$ be the Cartier support and $G=I_{C/T}$. The ideal sequence and the depth lemma make $G$ maximal Cohen--Macaulay. Gorenstein duality on $T$ gives
\[
 R\mathcal Hom_T(G,\OO_T)=G^\vee[0],\qquad G\simeq G^{\vee\vee}.
\]
For the square presentation $0\to P_1\xrightarrow{M}P_0\to j_*G\to0$, Cartier adjunction gives
\[
 \operatorname{coker}(M^\vee:P_0^\vee\to P_1^\vee)
 =\mathcal Ext^1_W(j_*G,\OO_W)=j_*(G^\vee\otimes N_{T/W}).
\]
On every fiber $\det M$ is a nonzerodivisor. Hence $M^\vee$ is injective on fibers and its cokernel is base-flat with formation commuting with base change. Applying this construction twice recovers $G$ on all nearby flat deformations. On the positive split join it restricts to curve duality for $I_{p/C}$; Serre duality and $\deg\omega_C=4$ give
\[
 \chi(I_{p/C}^\vee)=-\chi(I_{p/C})-4=-1 .
\]
These rank-one torsion-free boundary sheaves are stable. The strict cut inequalities for the dual were verified on the expanded support. Duality identifies the two smooth proper split-branch loci and their complete marking divisors, preserving the logarithmic tangent bundle. Relative Serre duality on the Calabi--Yau expanded targets then identifies the obstruction bundles with the corresponding dual logarithmic tangent bundles.
\end{proof}

\subsection{Proper relative factors for the split line}
Write $P=|h|$. Let $P_Q^L$ parametrize $(L,T_Q)$ where
$L\subset Q$ is a line and $T_Q$ is a hyperplane section containing it.
Let $P_B^f$ parametrize $(p,\Pi)$ with $p\in\Gamma$ and $\Pi$
a plane in $\PP^3$ containing $p$. Put $T_B=\rho^{-1}\Pi$ and denote
the exceptional ruling by $f_p$. Both incidences are projective-plane
bundles over smooth projective curves and carry universal surface
ideals $I_{L/T_Q}$ and $I_{f_p/T_B}$.

\begin{proposition}\label{prop:proper-split-line-factors}
The two incidences are independently proper relative stable sheaf
spaces. They represent the indicated incidences as proper open loci of the respective relative sheaf stacks. Their virtual classes are their fundamental classes. With
the scalar and determinant conventions of the matching theorem,
their product over $J_{-3}$ is a projective-plane bundle over
$F_1(Q)\times_S\Gamma$. Its diagonal pullback has degree $4800$.
\end{proposition}
\begin{proof}
The ideals are perfect on the separate roots by the square presentation
in Lemma~\ref{lem:lineideal-recovery-on-expanded-support}.
On $T_Q$, adjunction and the divisor sequence give
$\omega_{T_Q}=\OO_{T_Q}$ and $H^1(\OO_{T_Q})=0$.
On $T_B$ the divisor sequence and
$R\rho_*\OO_B=\OO_{\PP^3}$ give $H^1(\OO_{T_B})=H^2(\OO_{T_B})=0$,
also for singular plane transforms. Writing $E$ for the exceptional
divisor gives $\omega_{T_B}=(E-3H_B)|_{T_B}$, hence
$\omega_{T_B}|_{f_p}=\OO_{f_p}(-1)$. Thus for either curve $D=L,f_p$
one has
\[
 H^1(\omega_T|_D)=0,\qquad H^0(\omega_T)\twoheadrightarrow H^0(\omega_T|_D).
\]
The support and inclusion recovery in
Lemma~\ref{lem:lineideal-recovery-on-expanded-support} applies on
each root separately. The root has $H^1(\OO)=H^2(\OO)=0$ and its
support line bundle therefore has no deformation parameter. The
Hilbert scheme of the line is the smooth Fano scheme on $Q$.
On $B$ its ruling family over $\Gamma$ is the full embedded deformation
space since $h^0(N_{f_p/B})=1$, $h^1(N_{f_p/B})=0$ and motion over
$\Gamma$ spans the tangent.

Containment in $T_B$ is exactly the
plane-incidence equation $p\in\Pi$ over every base. Hence the maps
of the two incidence spaces to the respective scalar-rigidified
sheaf stacks are injective on geometric points and formally \'etale.
They preserve the trivial stabilizers and are representable open
immersions.

The supports on the roots are integral: the class on $Q$
is irreducible, and the total transform of a plane has no exceptional
component because $\Gamma$ is not contained in a plane. Their
rank-one torsion-free ideals are stable.

The line or ruling meets $S$ simply, so its ideal sequence restricts exactly to $C$. Both restrictions are $I_{p/C}$, stable with $\chi(I_{p/C})=-3$, also at a singular $p$. Each smooth branch has dimension and virtual dimension three, hence zero obstruction sheaf and virtual class equal to its fundamental class.

For the Abel map $(C,p)\mapsto I_{p/C}$ the ideal sequence gives
\[
 0\to\CC\to\Hom_C(I_p,\OO_C)\to\Ext^1_C(k_p,\OO_C)
          \longrightarrow H^1(C,\OO_C).
\]
The last map is injective: its dual is evaluation of $\omega_C=\OO_C(h)$ at $p$. Thus every map $I_p\to\OO_C$ is a multiple of the inclusion. This persists over Artin bases for maps of an ideal with its given inclusion: across a small extension with kernel $J$, subtract a lift of the old multiple; the difference belongs to $\Hom_C(I_p,\OO_C)\otimes J$. The annihilator recovers $C$, and the cokernel recovers $p$. This gives universal injectivity and injectivity on tangent spaces, including at singular $p$. Properness and representability then give a closed immersion.
Matching the sheaves therefore matches the flags scheme-theoretically. The boundary plane for $T_Q$ determines $\Pi$ and requires $L\cap S=p\in\Gamma$. This gives the stated bundle over a transverse set of $1600$ points.

For the excess calculation let $U=\{(p,C):p\in C\}\subset J_{-3}$. Over $S$ this is the projective-plane bundle of hyperplanes through the given point. Inside the smooth six-dimensional holomorphic symplectic variety $J_{-3}$ it has codimension two.

The restriction of the symplectic form $\omega$ to $U$ is pulled back from $S$: in the filtration of $\Omega_U^2$ by relative degree, all terms of positive relative degree have zero direct image, by $H^0(\PP^2,\Omega^j)=0$ for $j=1,2$. Hence $T_{U/S}\subset (T_U)^\perp$. Both bundles have rank two, and contraction with $\omega$ identifies
\[
 T_{U/S}=(T_U)^\perp,\qquad N_{U/J_{-3}}\simeq\Omega^1_{U/S}.
\]
The two factors map to $U$ by base change of $F_1(Q)\to S$ and $\Gamma\to S$. Indeed a plane through $p$ in $\PP^3$ spans a unique hyperplane with a given quartic-root line through $p$. The two maps intersect transversely in $U$, since their tangent maps on $S$ span $T_pS$ by Lemma~\ref{lem:linedefect-target-rigidity}. For the diagonal in $J_{-3}$ the excess bundle is therefore $N_{U/J_{-3}}$ on their product. The excess-intersection formula gives
\[
 \Delta^!_{J_{-3}}([P_Q^L]\times[P_B^f])
   =\sum_{L\in D_Q} c_2(\Omega^1_{\PP^2})\cap[\PP^2_L].
\]
The verified sheaf and boundary identifications above identify this pullback with the marked virtual class by Theorem~\companionref{thm:diagonal}. Its degree is $3\cdot1600=4800$.
\end{proof}

\subsection{The two complementary evaluations on a family of Abel cycles}
\label{sec:route19-abel}
We compare the two complementary components on boundary sheaves
$I_{p/C}$, where $C$ varies in $|h|$ and $p$ varies on $C$.
This comparison works before choosing the relative quartic cycle
or taking the degree of its virtual class.

Let $J_{-3}=M_S(0,h,-3)$ and denote its support map by
$\rho:J_{-3}\to |h|$. With the assumptions on $(S,h)$ of
Proposition~\ref{prop:proper-split-line-factors}, consider the maps
\[
 S\ \xleftarrow{\ \pi\ }\ U=\{(p,C):p\in C\}
 \ \xrightarrow{\ i\ } J_{-3},\qquad i(p,C)=[I_{p/C}].
\]
The projective-plane bundle $\pi$ is smooth, and $i$ is a regular
closed immersion of codimension two. For $D\in A^1(S)_{\QQ}$ set
$A_D=i_*(\pi^*D\cap[U])$. Write
$f_{-1}=[\Pic^{-1}(C_*)]$ for the class of the smooth fiber
$\rho^{-1}([C_*])$. Here $[C_*]\in|h|$ denotes the point specifying
the curve; its divisor class on $S$ is $h$.

\begin{proposition}\label{prop:route19-abel}
The split-ruling family on the component $B$ and the degree-minus-one
line-bundle family on the ruled support in $W$ have boundary virtual
pushforwards $A_{\Gamma}$ and $f_{-1}$, respectively. For every
$D\in A^1(S)_{\QQ}$ their pairings are
\begin{equation}\label{eq:route19-abel}
 \int_{J_{-3}} A_D\,A_{\Gamma}
  =3\int_S D[\Gamma]
  =15\int_S Dh
  =15\int_{J_{-3}}A_D\,f_{-1}.
\end{equation}
For the quartic-root split-line family the two degrees are $4800$ and
$320$, respectively.
\end{proposition}

\begin{proof}
The family of surfaces in $B$ carrying the ideal of the exceptional
ruling over $p\in\Gamma$ is $\pi^{-1}(\Gamma)$. Its boundary
restriction is $I_{p/C}$. Proposition~\ref{prop:proper-split-line-factors}
identifies its virtual class with its fundamental class, giving
the pushforward $A_{\Gamma}$.

On the rigid ruled support $j:R_*\hookrightarrow W$, write
$b:R_*\to C_*$ for the ruling. For $L\in\Pic^{-1}(C_*)$, let
$F_L=j_*b^*L$. Its restriction at $S_\infty$ is the line bundle $L$
supported on $C_*$, and $R_*$ misses the other relative divisor.
The two sheaf Ext terms of the Cartier immersion are
$\OO_{R_*}$ and $N_{R_*/W}$. The second has zero cohomology in
all degrees by Proposition~\ref{prop:rigidruledfull}, so
\[
 \operatorname{Ext}^k_W(F_L,F_L)=H^k(C_*,\OO_{C_*}).
\]
Varying $L$ realizes the entire tangent $H^1(C_*,\OO_{C_*})$, and
the obstruction group vanishes. The proper degree-minus-one Picard
torsor thus gives a smooth sheaf component as in
Proposition~\ref{prop:Wlinefamilyfull}, with virtual boundary
pushforward $f_{-1}$.

For any $D,E\in A^1(S)_{\QQ}$ the self-intersection formula gives
\[
 \int_{J_{-3}}A_D A_E
 =\int_U \pi^*(DE)c_2(N_{U/J_{-3}})
 =3\int_S DE .
\]
Indeed $N_{U/J_{-3}}\simeq\Omega^1_{U/S}$ by
Proposition~\ref{prop:proper-split-line-factors}, and
$\int_{\PP^2} c_2(\Omega^1_{\PP^2})=3$. To compute the
other pairing, set $g=\rho i$. The universal curve $g:U\to |h|$
is flat (its defining relative Cartier equation is regular on
every surface fiber). Over the smooth member $C_*$ the map
$\rho$ is a relative Picard morphism and is smooth. Consequently
$f_{-1}=\rho^*[\mathrm{pt}]$ for this choice of point, and
\[
 \int_{J_{-3}}A_D f_{-1}
 =\int_U \pi^*D\,g^*[\mathrm{pt}]
 =\int_{C_*}D|_{C_*}=\int_S Dh .
\]
Now $[\Gamma]=5h$. Finally, for the quartic Fano curve take
$D=(L\mapsto L\cap S)_*[F_1(Q)]$. The relative incidence is the
base change of $\pi$ along this map, so its boundary cycle is $A_D$.
The intersection in Theorem~\ref{thm:direct-line} gives
$\int_S D[\Gamma]=1600$, and hence $\int_S Dh=320$.
Substitution proves the last assertion.
\end{proof}
Equation~\eqref{eq:route19-abel} compares the two stated complementary
families on the entire span of the $A_D$. A full Fourier coefficient
also includes the contributions of the other marked families of that
charge.

\begin{lemma}[Length-two contact with the two curves]
\label{lem:route19-two-contact}
For a smooth divisor $D\subset S$, put $D^{[2]}=\Hilb^2(D)\subset
S^{[2]}=\Hilb^2(S)$. Let $r:\widetilde{S^2}=\operatorname{Bl}_{\Delta}(S^2)
\to S^{[2]}$ be the map forgetting the ordering, $E$ its exceptional
divisor, and $h_j$ the pullback of $h$ from the $j$th factor of $S^2$.
Then in $A^2(\widetilde{S^2})$,
\[
 r^*[\Gamma^{[2]}]=25r^*[C_*^{[2]}]+20[E]h_2.
\]
The last term has nonzero numerical class.
\end{lemma}
\begin{proof}
For a line bundle $L$ write $L^{[2]}$ for its rank-two tautological
bundle on $S^{[2]}$, with fiber $H^0(Z,L|_Z)$ at $Z$. A defining
section for $D\in|L|$ gives a section of $L^{[2]}$ with zero
scheme $D^{[2]}$. The latter is the smooth surface
$\operatorname{Sym}^2D$ in the smooth fourfold $S^{[2]}$, so
$[D^{[2]}]=c_2(L^{[2]})$. On $\widetilde{S^2}$ write $L_j$ for the
pullbacks from the two factors. The fundamental exact sequence on the ordered Hilbert square is described in \cite[\S1, (1)]{SchlickeweiHilb}. Evaluation at the second point gives
\[
 0\longrightarrow L_1(-E)\longrightarrow r^*L^{[2]}
 \longrightarrow L_2\longrightarrow0 .
\]
Indeed the ordered universal subscheme is the union of its two graphs,
meeting along $E$. Locally take the first point as origin, choose
blowup coordinates in which the second is $(u,au)$, and call the
relative surface coordinates $(x,y)$. The union has equations
$y=ax$, $x(x-u)=0$, including its flat double point at $u=0$.
The kernel of evaluation on the second graph is the sheaf of sections
on the first graph vanishing at $u=0$. Tensoring by $L$ gives the
displayed sequence. Hence for $L=\OO_S(kh)$ we have
\[
 c_2(r^*L^{[2]})=k^2 h_1h_2-k[E]h_2 .
\]
Substitution of $k=5,1$ gives the identity. For the blowdown
$d:\widetilde{S^2}\to S^2$ the projective-bundle formula gives
$d_*[E]^2=-[\Delta]$. Thus
$\int_{\widetilde{S^2}}([E]h_2)^2=-\int_S h^2=-4$.
\end{proof}
This identifies the additional cycle supported on the exceptional
divisor that a comparison of length-two contact conditions detects.
Applying it to threefold virtual classes requires identifying
the boundary morphisms with these incidence conditions and
computing the virtual multiplicities along the collision locus.

\begin{proposition}\label{prop:four-secant-sheaf-marked}
For each four-secant defect $L'$ the two marked special-fiber spaces
are $V_+=\operatorname{Bl}_{o}\PP^2$ and $V_0=\PP^2$, with common
divisor $D_+=D_0=\PP^1$. Their marked classes are
$c_2(\Omega^1_{V_i}(\log D_i))\cap[V_i]$, of degrees $2$ and $1$.
The positive term is the restriction to its proper marked space
of the refined diagonal on relative sheaf stacks.
The zero-boundary term is a proper one-sided surface-sheaf class.
\end{proposition}
\begin{proof}
The incidence blowup gives the central fiber $V_+\cup_D V_0$.
Its expansion parameters $s,v$ satisfy $sv=t$. On a normalized
marked face the unmarked parameter $z$ is a local equation for
the smooth divisor $D_i$. Thus its morphism to the normalized
expansion stack is smooth, with relative tangent
$T_{V_i}(-\log D_i)$: in the chart
$[\mathbb A^1/\mathbb G_m]$ a vector field satisfies $\xi(z)=az$.

The intrinsic flag inverse identifies this with the relative
sheaf deformations on the complete target. Its relative Serre
duality is compatible with dilation on the bubble since the
logarithmic ruling differential is invariant under rescaling.
Thus the obstruction bundle is the dual logarithmic bundle,
by the cotangent triangle over the smooth marked base, precisely
as in Proposition~\ref{prop:one-defect-positive-marked}.

Logarithmic Gauss--Bonnet gives
$e(\operatorname{Bl}_o\PP^2)-e(\PP^1)=4-2=2$ and
$e(\PP^2)-e(\PP^1)=3-2=1$. At the zero marked join one
relative sheaf is empty and the other is the entire expanded
flag on $B$. It misses its outer end, including along the graph
bubble. One-ended contraction and support/ideal recovery identify
the latter relative space with $V_0$.

For the positive marked join set $K_{L'}=\PP^1$, the pencil of
planes containing the image line. Put
$A=\operatorname{Bl}_{[S]}|H_Q|$. Here the common boundary is
the structure sheaf of the plane quartic, so its class is
$\sigma(C)=[\OO_C]\in J_{-2}$. Around the one-ended plane sheaves
of this pencil take the relative unrestricted locus $U_B$ of simple
perfect sheaves normal at $S$, with stable boundary in $J_{-2}$
and no infinitesimal target stabilizer.

There is such a finite-type
open union of neighborhoods by the openness and bounded-neighborhood
argument in Theorem~\ref{thm:point-open-diagonal}. Cutting the
complete open immersion of the expanded line-ideal space and applying
derived matching gives an open immersion
\[
 V_+=A\times_P K_{L'}\xrightarrow{j}A\times_{J_{-2}}U_B,\qquad
 [V_+]^\vir=j^*\Delta_{J_{-2}}^!([A]^\vir\times[U_B]^\vir).
\]
This uses the earlier proper quartic structure-sheaf component.
After one prescribes $F|_C\cong\OO_C$, the plane pencil itself
has the full constrained sheaf deformations. Indeed set $I=I_{L'/T_B}$.
The curve $C$ misses $L'$, and
\[
 H^1(\OO_{T_B}(-C))=0,\qquad
 \operatorname{Ext}^1_{T_B}(I,\OO_{T_B}(-C))=0.
\]
The first assertion follows from $H^0(\OO_{T_B})\cong H^0(\OO_C)$
and $H^1(\OO_{T_B})=0$. For the second apply the ideal sequence
and duality: the possible kernel is dual to the cokernel of
$H^0(\omega_{T_B}(C))\to H^0(\omega_{T_B}(C)|_{L'})$;
this is onto since $\omega_{T_B}(C)=H_B|_{T_B}$. Also
$H^1(L',\OO_{L'}(1))=0$ gives
$\operatorname{Hom}_{T_B}(I,\OO_{T_B}(-C))=0$.
These groups give existence and uniqueness for lifting the inclusion
with its prescribed boundary map by the restriction triangle, as in
Proposition~\ref{prop:plane-boundary-rigidity}.

This identifies
$K_{L'}$ as an open subspace of $U_B\times_{J_{-2}}P$.
The induced zero-cycle on this proper pencil has degree $-2$, since
$V_+\to K_{L'}$ is a projective-line bundle and the quartic factor
has virtual class $c_1(b_P^*\OO(-1))\cap[A]$ with fiber degree $-1$.
These identifications give the positive marked comparison without
assuming properness of the unrestricted $U_B$.
\end{proof}

\section{A conic on a moving hyperplane surface}\label{sec:marked-conics}
For the classical degeneration of conics see \cite{KatzConicLimit}. We compute the surface-sheaf contributions of conics on either root and of a split conic whose two lines have a moving attachment point on $S$. As before, $f_1,f_4,f_5$ are general and $\Pic(S)=\mathbb Zh$, so $S$ contains neither a line nor a plane conic.

\begin{lemma}\label{lem:wholeconic-18}
The stable sheaves in the two primitive charges $(e,\chi)=(-2,4),(2,2)$ on a general smooth quintic are respectively the ideals $i_*I_{D/T}$ and their surface duals, where $D$ is a smooth plane conic, $T$ a hyperplane section containing it, and $i:T\hookrightarrow X$. Each of the two sheaf moduli schemes is the pencil of containing hyperplanes over the finite \'etale scheme of conics.
\end{lemma}
\begin{proof}
Write $R=G^{\vee\vee}$ and $k=\length(R/G)$. In degree $-2$, Lemma~\ref{lem:polar-small-negative-degree} gives
\[
 4+k=\chi(R)\le4,\qquad
 4\le h^0(R^\vee(h))\le h^0(R^\vee)+3\le4.
\]
For the middle inequality, restrict to a general plane quintic $C$: the degree-seven restriction has $h^0\le3$ by Riemann--Roch and gonality four. Thus $k=0$ and the section of $R^\vee$ identifies $R=I_{D/T}$ for a Cohen--Macaulay curve with $(\deg D,\chi(\OO_D))=(2,1)$, by the depth lemma.

An integral such curve is a plane conic. Disjoint lines have Euler characteristic two; meeting lines depend on $4+3+3=10$ parameters and impose eleven independent quintic conditions. For a doubled line $L$, purity extends the generic equality $I_{L/X}^2\OO_D=0$ across $L$. Its nilradical is a torsion-free quotient $N^\vee_{L/X}=\OO_L(1)^{\oplus2}\twoheadrightarrow\OO_L(a)$, since $L$ is isolated and reduced with normal degree $-2$. Hence $a\ge1$ and $\chi(\OO_D)=1+(a+1)\ge3$. Only smooth conics remain. In degree two, $\chi(R^\vee)=\chi(R)+2=4+k$ reduces to the previous case. Relative Gorenstein duality identifies the functors as in Corollary~\ref{cor:line-opposite-expanded}.

For the scheme structure, the Hilbert scheme of plane conics is smooth of dimension eleven, a projective bundle over $\operatorname{Gr}(3,5)$. A quintic equation induces a general section of the generated rank-eleven bundle with fibers $H^0(D,\OO_D(5))$, so its zeros are transverse. For $F=i_*I_{D/T}$ a square presentation on $X$ (the conic is lci of codimension two and $T$ is Cartier containing it) recovers its Cartier support under every infinitesimal deformation as in Lemma~\ref{lem:lineideal-recovery-on-expanded-support}. On $T$, $H^1(\OO_T)=0$ and
\[
 H^1(D,\omega_T|_D)=0,\qquad
 H^0(T,\omega_T)\twoheadrightarrow H^0(D,\omega_T|_D)
                         =H^0(\PP^1,\OO(2)).
\]
The ideal sequence and Gorenstein duality give the two groups used in lifting the inclusion:
\[
 \Hom_T(I_{D/T},\OO_T)=\CC,\qquad \Ext_T^1(I_{D/T},\OO_T)=0 .
\] The embedded conic itself is rigid; deformations of the ideal thus give exactly containing hyperplanes, also on Artin bases. Properness, injectivity on geometric points and this inverse identify the rigidified sheaf scheme.
\end{proof}

\begin{lemma}[Three special loci of the conic defect]\label{lem:threeconics-18}
Over the strictly henselian trait of the first degeneration the smooth-fiber conic scheme extends in the projective plane-conic scheme to a finite \'etale horizontal closure. Its special points have the following cardinalities, computed on the special loci:
\[
\begin{array}{c|r}
 A:\text{ conic on }Q\text{ meeting }\Gamma\text{ at its two points on }S&258200\\
 B':\text{ conic in }\PP^3\text{ meeting }\Gamma\text{ at its eight points on }S&187850\\
 M:\text{ a line on each root, with one matching point on }S&163200 .
\end{array}
\]
Here the plane-root line in $M$ meets $\Gamma$ at the other three points of its intersection with $S$.
\end{lemma}
\begin{proof}
On $\mathcal H=\PP(\operatorname{Sym}^2 U^\vee)$ over $\operatorname{Gr}(3,5)$ use the lines convention and tautological plane bundle $U$. For $d\ge2$ write
\[
0\longrightarrow \OO_{\mathcal H}(-1)\otimes\operatorname{Sym}^{d-2}U^\vee
\xrightarrow{\cdot q}\operatorname{Sym}^dU^\vee\longrightarrow E_d\longrightarrow0
\]
and put $E_1=U^\vee$. Its sections induced by ambient forms generate $E_d$, of rank $2d+1$. The zero scheme of $f_1f_4$ in $E_5$ has three disjoint smooth components. On conics contained in $Q$, multiplication by $f_1:E_4\to E_5$ is injective on a neighborhood, since no component of such a conic can be in $S$. The zero scheme there is $Z(f_4\in E_4)$, smooth of dimension two. On conics in the plane root it is the bundle over $\operatorname{Gr}(3,4)$ cut by $f_1\in E_1$, since $f_4:E_1\to E_5$ is injective there. The dimension is eight.

At a remaining zero the conic splits in its plane as $q=f_1y$, with $M_0=(f_1=0)$ on the plane root and $L=(y=0)$ on $Q$. These lines are distinct. On an Artin deformation write $f_1f_4=qh$. Modulo $f_1$, the closed-fiber value $\bar h=f_4|_{M_0}/y$ is nonzero since $S$ contains no line. Flatness of the moving line makes $\bar h$ a nonzerodivisor. Hence
\[
 \bar q\bar h=0\ \Longrightarrow\ q=f_1 y\ \Longrightarrow\ f_4=y h .
\]
Thus the remaining zero scheme is the $\PP^2$-bundle over the Fano curve of $Q$ parametrizing the lines $M_0$ through $p=L\cap S$. This proves cleanness along the three loci.

On each such locus $Y$, the derivative identifies $N_{Y/\mathcal H}$ with a subbundle of $E_5|_Y$. Its quotient has rank $\dim Y$ and is generated by quintics. Hence a general $f_5$ has transverse, isolated residual zeros. As in Theorem~\ref{thm:direct-line}, solve the normal equations formally and divide the residual equations by $t$. Transversality gives the \'etale horizontal closure.

On $A,B'$ the residual section is $f_5$ restricted to the contact scheme with $S$. On $M$, put $g=f_4|_{M_0}/y$. The normal image restricts to zero modulo $g$, by differentiating $f_1f_4-qh$ on $M_0$, where $f_4=yg$ and $h=g$. Restriction identifies the rank-three residual quotient and its zero condition:
\[
 H^0(M_0,\OO(5))/g H^0(M_0,\OO(2)),\qquad
 f_5|_{M_0}\in g H^0(M_0,\OO(2)).
\]
Indeed the restriction onto this quotient is surjective and both sides have rank three.

A general residual zero avoids any fixed proper closed subset of these smooth loci since the respective excess ranks equal their dimensions. On $A$ and $B'$ singular conics and conics meeting $S$ nontransversely form proper subloci: for $A$, use general $f_4$ and then general $f_1$ (the singular conics are a divisor in $\mathcal H$); for $B'$ use general planes and conics on them.

In the mixed locus a general line through $p$ meets the smooth $S$ simply in four points; the tangent or repeated-residual conditions form proper subloci over each $p$. The additional value $f_5(p)$ is independent of the three imposed residual values since $p$ is distinct from them. Thus $f_5(p)\ne0$ at every mixed solution.

For the intersection numbers one has
\[
 N_A=\int_{\mathcal H}c_9(E_4)c_2(E_5-E_4),\qquad
 N_{B'}=\int_{\mathcal H_3}c_8(E_5-E_1),
\]
where $\mathcal H_3$ is the same conic bundle for $\PP^3$.
For explicit evaluation take distinct weights $w_i=11^i$ ($0\le i\le4$ or $0\le i\le3$). For a three-element subset $I$ use $S_d(I)$ for the multiset of its degree-$d$ monomial weights (sums with repetition). For each $a\in S_2(I)$ set
\[
 R_d=S_d(I)\setminus(a+S_{d-2}(I)),\quad R_1=I,\qquad
 T(I,a)=\prod_{i\in I,\,j\notin I}(i-j)\prod_{b\in S_2(I),\ b\ne a}(b-a).
\]
Here $I$ and its complement are sets of weights, and $S_0=\{0\}$. Equivariant integration on the Grassmannian and projective fibers gives the finite sums
\[
\begin{aligned}
 N_A&=\sum_{I,a}\frac{\prod R_4}{T(I,a)}[z^2]
       \frac{\prod_{b\in R_5}(1+bz)}{\prod_{b\in R_4}(1+bz)}=258200,\\
 N_{B'}&=\sum_{I,a}\frac{1}{T(I,a)}[z^8]
       \frac{\prod_{b\in R_5}(1+bz)}{\prod_{b\in I}(1+bz)}=187850 .
\end{aligned}
\]
The two sums run over $60$ and $24$ fixed points, respectively. In $R_4$ the notation $\prod R_4$ means the product of all weights with multiplicity.

On the mixed locus let $P$ be the tautological point line of $p$ on the Fano curve, let $V$ be the rank-two tautological bundle of $M_0$, and let $J=V/P$. Write $x=c_1(P^\vee)$, $y=c_1(J^\vee)$. The residual cubic $g$ has its scale prescribed by $f_4$: the canonical linear equation of $p$ on $M_0$ is $V\to J$, so $g$ takes values in $\operatorname{Sym}^3V^\vee\otimes J^\vee$. Its image in $\operatorname{Sym}^5 V^\vee$ is the subbundle $J\otimes\operatorname{Sym}^2 V^\vee$. The degree-three Chern class of the quotient, obtained by dividing the two Chern polynomials, is
\[
 60x^3+495x^2y+750xy^2+240y^3 .
\]
The mixed locus is $\PP(\CC^4/P)$ over the quartic Fano curve; $J$ is the tautological line there. Integration on the plane fibers gives $p_*y^2=1$, $p_*y^3=-c_1(\CC^4/P)=-x$. Finally $\int x=320$ by the quartic Fano intersection in Theorem~\ref{thm:direct-line}. Hence $N_M=(750-240)320=163200$. These are the three counts used in the marked specialization.
\end{proof}

\begin{theorem}[Marked virtual specialization in the conic charges]\label{thm:marked-conics-18}
In each of the two characters of Lemma~\ref{lem:wholeconic-18} there is a proper scalar-rigidified expanded sheaf space over the strictly henselian trait with complete generic fiber. Its marked virtual zero-cycles have the following degrees:
\[
\begin{array}{c|r}
 A&-2\cdot258200\\
 M&-2\cdot163200\\
 B',\ \text{positive marked boundary}&-187850\\
 B',\ \text{zero marked boundary}&-187850 .
\end{array}
\]
Consequently $[q^{13/120}]h^{\mathrm{DT}}_{1,2}=[q^{13/120}]h^{\mathrm{DT}}_{1,3}=-1218500$.
\end{theorem}
\begin{proof}
\emph{The embedded special curves.}
Lift the sheets of Lemma~\ref{lem:threeconics-18} to the blowup of $(f_1,t)$. For type $A$ its restriction to the moving conic surface has base scheme $\{p_1,p_2\}$, the two transverse contacts on the central conic. The blowup adds the exceptional rulings over both points. For type $B'$ one has $(f_1,t)=(t)$ on the moving surface. The strict transform misses $S$, since the regular ratio $f_1/t=f_5/f_4$ is finite at all eight base points by divisibility.

For type $M$ write $q_t h_t=f_1f_4-t f_5$ on the moving plane. At the node $p$ the factors $f_1,f_4,q_0$ vanish, whereas $h_0(p)=g(p)\ne0$. Differentiation at $(p,0)$ gives
\[
 \dot q_0(p)=-f_5(p)/g(p)\ne0,\qquad t=f_1 f_4/f_5
\]
on the target near $p$. Thus the conic surface is smooth there and the blowup is an isomorphism at $p$. Off $p$, $(f_1,t)$ is the unit ideal on the quartic branch and is principal on the plane branch, where $f_1$ vanishes on the central component. The limit is $L\cup\widetilde M_0$: the plane line lifts at its other three roots of $f_4$ and meets $S$ simply at $p$.

At a join the transverse branches have equations $(u,v)$ in $\CC[[x,y,u,v]]/(xy)$ by Lemma~\ref{lem:linedefect-target-rigidity}; elsewhere each branch is a smooth curve in a smooth threefold. Hence each limit $D\subset W$ is a regular immersion of codimension two.

\emph{Rigid defects on the expanded fiber.} For a rational tree with components $D_i$ crossing joins transversely, normalization of its locally free rank-two normal sheaf gives
\[
0\longrightarrow N_{D/W}\longrightarrow
 \bigoplus_i N_{D_i/Y_i}\longrightarrow\bigoplus_{p\,\mathrm{join}} T_pS\longrightarrow0.
\]
Transversality identifies each branch-normal fiber with $T_pS$; the right arrow is the difference of the two attachment displacements. For type $A$ the plane and equation of the conic deform with it, since $H|_{\PP^1}=\OO(2)$ has three ambient generators and $H^1=0$. Each exceptional tail has $0\to\OO\to N\to\OO(-1)\to0$ and can move only over $\Gamma$. After ordering the contacts, the kernel on sections is therefore the tangent to $f_5(p_1)=f_5(p_2)=0$ on the smooth quartic-conic surface. It vanishes by transversality.

For $B'$ the same three sections give the projected embedded conic. A lift to $B\setminus S$ is the graph of $f_5/f_4$, so the eight conditions are
\[
 f_5|_D\in f_4|_D\,H^0(D,\OO_D(1))
 \quad\Longleftrightarrow\quad f_5|_{D\cap S}=0 .
\]
In case $M$ the complete pencil of the plane-root line gives its projection. Write its four simple contacts as $p_1,p_2,p_3,p$, with $p_i\in\Gamma$ and $p\notin\Gamma$. Splitting the contacts \'etale locally, regularity of $f_5/f_4$ at $p_i$ is exactly $f_5(p_i)=0$. Near $p$ the blowup is an isomorphism, and the two branches must have the same moving attachment point there. These are precisely the reduced residual equations of Lemma~\ref{lem:threeconics-18}. Again the tangent space vanishes. For $B'$ all relevant expansions are isomorphisms near the curve.

In every case $H^0(N_{D/W})=0$. Since $D$ is a rational tree and $\det N_{D/W}=\omega_D$ by adjunction,
\[
 \chi(N_{D/W})=2\chi(\OO_D)+\deg\omega_D=2-2=0.
\]
For a square-zero extension with kernel $J$, the local lifts of $(u,v)$ differ by a cocycle in $N_{D/W}\otimes J$. Its obstruction group $H^1(N_{D/W})\otimes J$ and the translation group $H^0(N_{D/W})\otimes J$ both vanish. Thus the embedded curve lifts uniquely with the target.

\emph{The expanded sheaf functor.} On a conic sheet $p:D\to\Spec R$, $H$ has degree two on the conic and zero on its tails, or multidegree $(1,1)$ on the split pair. Normalization gives $R^1p_*H=0$ and $\operatorname{rk}p_*H=3$. Ambient linear forms generate $p_*H$, so the containing pencil is $\PP(\ker(H^0(\PP^4,\OO(1))\otimes R\to p_*H))$.

In cases $A,M$ it misses the distinguished hyperplane $(f_1=0)$; the allowed targets remain unexpanded. In case $B'$ take the strict transform of the pencil in the minimal-support expansion. This is the blowup of the pencil bundle at its unique distinguished closed point. To check the scheme structure, set $u_t=f_1|_{D_t}/t$, a regular section. For plane-coordinate sections $l_i$ the three linear equations
\[
 t u_t+\textstyle\sum b_i l_i|_{D_t}=0
\]
cut out the incidence on the original chart. The four $l_i$ span the three-dimensional section module at $t=0$. Thus the center ideal on the pencil has generators $t$ and one free variable. On the graph chart $(s,v)$ divide the equations by $s$, obtaining $v u_t+l_c|_{D_t}=0$, while on the $t$-chart they read $u_t+l_c|_{D_t}=0$. As in Theorem~\ref{thm:expanded-line-compact}, these are the actual flat incidences on the expanded support charts. The curve is disjoint from the expansion center in case $B'$.

Set $I=I_{D/T}$. Flatness over the incidence base follows from
$0\to I\to\OO_T\to\OO_D\to0$. On a branch put $C=T_i\cap S$,
$Z_i=D_i\cap S$ (possibly empty). Multiplication by the local equation
$s$ of $S$ is injective on both $\OO_{T_i}$ and $\OO_{D_i}$.
Accordingly, at every attachment, including a singular point of $C$,
\[
 \LL i_S^*I_{D_i/T_i}
 =\fib(\OO_C\twoheadrightarrow\OO_{Z_i})
 \simeq I_{Z_i/C}.
\]
Recall the square presentation with columns $(-v,u)$ and $(a,b)$ for
$I$ over $W$, where $I_{D/W}=(u,v)$ and $T=(au+bv=0)$.
Its determinant recovers $T$ over Artin bases as in
Lemma~\ref{lem:lineideal-recovery-on-expanded-support}.
Put $A=\omega_T|_D$. Adjunction gives $A=H|_D$ on ordinary
supports. The additional chart twist misses $D$ on expanded supports,
and $A=H_B|_D$ on the zero-boundary chart.

Thus $H^1(D,A)=0$ and
$h^0(D,A)=3$ by the preceding pencil calculation. The divisor sequence
gives $H^1(T,\OO_T)=0$ on each minimal profile; the basis of
Lemma~\companionref{lem:vacuumlinearsystems} gives a surjection
$H^0(T,\omega_T)\to H^0(D,A)$, already on the plane-coordinate
sections for the plane-root case. Surface duality and the ideal
sequence now give
\[
\begin{aligned}
0&\to H^0(\OO_T)\to\Hom_T(I,\OO_T)\to H^1(D,A)^\vee=0,\\
0&\to\Ext^1_T(I,\OO_T)\to
 H^0(D,A)^\vee\hookrightarrow H^0(T,\omega_T)^\vee.
\end{aligned}
\]
Hence the inclusion of $I$ lifts up to a unit, with flat quotient:
its obstruction group vanishes and $\Hom_T(I,\OO_T)=\CC$.
Together with the embedded-curve rigidity already proved, this
recovers the incident flag $D'\subset T'$ from each Artin sheaf
deformation.

For stability use the polarization of \textup{(\companionref{eq:vacuumlocalpolarization})} with
$5/9<q_0/M<q_1/M<1$. For an ordinary or graph support put $q=q_0$ or $q_1$ as appropriate and $m_1=2,0,1$ for profiles $A,B',M$, respectively, the projected degree on the left active branch. The defect has total polarization degree $2M$ by flatness, and degree $m_1 q$ on that branch.

For example, in case $M$ the right branch meets the exceptional divisor three times, so its polarization degree is $(5M-4q)-3(M-q)=2M-q$. The linear cut coefficients for the ideal and its surface dual follow by subtracting or adding these degrees in \textup{(\companionref{eq:vacuumgraphpolys})}:
\[
[m]D_-=\frac{18q^2}{5M}-m_1q,\qquad
 [m]D_+=\frac{2q^2}{5M}+m_1q .
\]
Since $m_1\in\{0,1,2\}$ and $5/9<q/M<1$, both coefficients lie strictly between $0$ and $4q$. The cut criterion gives stability, including on the punctured chart axes. A rank-one torsion-free sheaf on an integral zero-boundary support is stable, and a graph support fixes its bubble dilation. Transposing the square resolution commutes with base change and gives duality as in Corollary~\ref{cor:line-opposite-expanded}.

The pencil bundles for $A,M$ and their indicated blowups for $B'$ thus represent the scalar-rigidified sheaf groupoid on this locus. They are proper and separated; each has at most one inserted level. Flag recovery gives an open immersion by the infinitesimal criterion used in Theorem~\ref{thm:expanded-line-compact}. Lemma~\ref{lem:wholeconic-18} identifies the generic fiber. Theorem~\companionref{thm:main} therefore applies to these marked spaces.

\emph{Marked virtual degrees.}
Every marked curve is a $\PP^1$. For $A,M$ there is no free smoothing parameter. On either marked branch of a $B'$ sheet the other node has parameter $z$ vanishing at the crossing $d$. The relative vector fields there satisfy $\xi(z)\in(z)$ on the dilation quotient. As in Proposition~\ref{prop:four-secant-sheaf-marked}, flag recovery and relative duality give
\[
\begin{array}{c|c|c|r}
\text{sheet}&T_{\rm rel}&\mathrm{Ob}&\deg[\PP^1]^\vir\\\hline
 A,M&T_{\PP^1}&\Omega^1_{\PP^1}&-2\\
 B',\ \text{either branch}&T_{\PP^1}(-\log d)&\Omega^1_{\PP^1}(\log d)&-1 .
\end{array}
\]
Substitute the three counts of Lemma~\ref{lem:threeconics-18} and the two characters in \eqref{eq:polar-charge-euler} to obtain the stated coefficient. The specialization takes place on the complete marked moduli.

\end{proof}

\subsection{Two boundary insertions for conic ideals}
\label{sec:mark20-contacts}
Set $X_2=S^{[2]}$, $J=J_{-4}=M_S(0,h,-4)$ and write
\[
 S^{[2]}\xleftarrow{\ p\ } U_2=\{(Z,C):Z\subset C\in|h|,\ \length Z=2\}
          \xrightarrow{\ a\ }J,\qquad a(Z,C)=I_{Z/C}.
\]
As throughout the primitive Jacobian calculation, $\Pic(S)=\mathbb Zh$.
For a smooth $C_*\in|h|$ define the following middle-dimensional cycles:
\[
 e=[C_*^{[2]}]\in A^2(X_2),\qquad
 \tau=j_*\varpi^* h\in A^2(X_2),\qquad
 f=[\Pic^{-2}(C_*)]\in A_3(J).
\]
Here $j:\mathcal E\hookrightarrow X_2$ is the exceptional divisor of
Hilbert--Chow, and $\varpi:\mathcal E\to S$ sends a nonreduced
subscheme to its support. Thus $\tau$ imposes a double point over
a hyperplane curve, with its tangent direction free. Write
$A_D=a_*(p^*D\cap[U_2])\in A_3(J)$ for $D\in A^2(X_2)_\QQ$.

\begin{lemma}\label{lem:mark20-abel2}
The map $p$ is a projective-line bundle. The map $a$ is a regular
closed immersion with $N_{U_2/J}\simeq\Omega^1_{U_2/X_2}$. For
$D,H\in A^2(X_2)_\QQ$,
\begin{equation}\label{eq:mark20-abel2}
  \deg_J(A_D A_H)=-2\int_{X_2}DH,\qquad
  \deg_J(A_D f)=\int_{X_2}De .
\end{equation}
\end{lemma}
\begin{proof}
Very ampleness gives a surjection of vector bundles
$H^0(S,h)\otimes\OO_{X_2}\to h^{[2]}$; its kernel has rank two and
its line projectivization equals $U_2$. Each integral plane quartic
$C$ is Gorenstein, with $\omega_C=\OO_C(h)$. In the dual of the
ideal sequence for $Z\subset C$, the arrow
\[
 \Ext^1_C(\OO_Z,\OO_C)\longrightarrow H^1(C,\OO_C)
\]
is dual under Serre duality to
$H^0(C,\omega_C)\twoheadrightarrow H^0(Z,\omega_C|_Z)$.
The surjection holds for every length-two subscheme, including one at a
singular point. Thus $\Hom_C(I_Z,\OO_C)=\CC$, with basis the inclusion.
The annihilator of $I_Z$ on $S$ is the ideal of $C$. This remains valid
over a test family: $I_Z\subset\OO_C$ equals $\OO_C$ away from the
relative finite scheme, and sections of $\OO_C$ supported on that scheme
vanish by relative depth one.

For a first-order deformation mapping to a
constant sheaf the support $C$ is therefore constant. An inclusion on
this constant family reducing to the given inclusion differs from it by
a scalar multiple over the dual numbers, by the Hom calculation.
Its image and hence $Z$ are constant.

This proves injectivity of the
tangent map. The same uniqueness proves injectivity on geometric
points. Since $U_2$ is proper, $a$ is proper, unramified and universally
injective in characteristic zero, hence a closed immersion. The source
and target are smooth of dimensions five and six.

Let $\omega$ denote the symplectic form of $J$. The exterior-power
sequence for the smooth map $p$ and $p_*\Omega^1_{U_2/X_2}=0$ imply
$a^*\omega=p^*\alpha$ for a two-form $\alpha$ on $X_2$. Hence the
line $T_{U_2/X_2}$ equals $(T_{U_2})^\perp$. The symplectic pairing now
gives $N_{U_2/J}=\Omega^1_{U_2/X_2}$. Self-intersection gives the
first identity in \eqref{eq:mark20-abel2} since $p_*c_1(N)=-2$.
For the second, let $\rho:J\to|h|$ be support.

At the smooth curve $C_*$, $f=\rho^*[C_*]$ and
$(\rho a)^*[C_*]=[C_*^{[2]}]$ on $U_2$, where $[C_*]$ on the base
denotes a point. Indeed the Hilbert scheme of two points on the family
of smooth curves over the smooth-member locus is smooth over that
locus. The fiber projects by $p$ isomorphically to $C_*^{[2]}$.
The projection formula proves the identity.
\end{proof}

\begin{proposition}\label{prop:mark20-twoprobe}
Put $g=[\Gamma^{[2]}]\in A^2(X_2)$ for the smooth curve
$\Gamma\in|5h|$. For all $D\in A^2(X_2)_\QQ$,
\begin{equation}\label{eq:mark20-twoprobe}
 \deg_J(A_D A_g)=-50\deg_J(A_D f)+10\deg_J(A_D A_\tau).
\end{equation}
The two functionals on the right are linearly independent on the span
of the $A_D$. The ruled component $W$ realizes $f$ by the sheaves
$i_*b^*L$, $L\in\Pic^{-2}(C_*)$, on its rigid ruled support
$i:R_*\hookrightarrow W$, $b:R_*\to C_*$.
\end{proposition}
\begin{proof}
Use the ordered double cover $r:\widetilde{S^2}\to X_2$ in
Lemma~\ref{lem:route19-two-contact}, with exceptional divisor $E$.
It ramifies along $E$ and maps $E$ isomorphically to $\mathcal E$;
therefore $r^*\tau=2[E]h_2$. For instance, if $H$ denotes the divisor
on $X_2$ pulled back from the symmetric divisor $h_1+h_2$, then
$\tau=[\mathcal E]H/2$ because $H|_{\mathcal E}=2\varpi^* h$.
Writing $u=h_1h_2$, $v=[E]h_2$, we have
\[
 r^*e=u-v,\qquad r^*g=25u-5v,\qquad g=25e+10\tau .
\]
The last equality follows by $r_*r^*=2$ on rational Chow groups.
Now apply \eqref{eq:mark20-abel2}. The ordered intersections satisfy
$\int u^2=16$, $\int uv=0$ and $\int v^2=-4$. Thus
\[
 \begin{pmatrix}
 \int e^2&\int e\tau\\ \int e\tau&\int\tau^2
 \end{pmatrix}
 =\begin{pmatrix}6&4\\4&-8\end{pmatrix}.
\]
Its determinant is $-64$. Taking $D=e,\tau$ in
\eqref{eq:mark20-abel2} proves independence. For $L$ of degree $-2$,
the two local self-Ext sheaves of $i_*b^*L$ are $\OO_{R_*}$ and
$N_{R_*/W}$. The latter has zero cohomology by
Proposition~\ref{prop:rigidruledfull}. As in the proof of
Proposition~\ref{prop:route19-abel}, $\Pic^{-2}(C_*)$ parametrizes a
smooth proper component with its fundamental virtual class and evaluation $f$.
\end{proof}

\begin{corollary}\label{cor:mark20-conictest}
Let $\mathcal F_2(Q)$ be the parameter surface of plane conics on the
quartic threefold $Q$ and $D$ its pushforward to $X_2$ by
$K\mapsto K\cap S$. Then
\begin{align*}
 \deg_J(A_D f)&=5016,&\deg_J(A_D A_\tau)&=-26560,\\
 \deg_J(A_D A_g)&=-516400
   & &=-50(5016)+10(-26560).
\end{align*}
\end{corollary}
\begin{proof}
No irreducible component of a conic on $Q$ lies in $S$. Cutting
the universal conic by $S$ gives a relative Cartier subscheme
finite flat of degree two. Let $U_3$ be the tautological plane
bundle on $\operatorname{Gr}(3,5)$ pulled back to the conic
parameter bundle $\mathcal H$ of Lemma~\ref{lem:threeconics-18}.
Put $S=(s_0=0)$ in $Q$ and $C_*=(s_1|_S=0)$.

For any conic
$K$ the restriction map $U_3^\vee\to H^0(K,\OO_K(1))$ is an
isomorphism, and restriction to $K\cap S$ yields the quotient by the
nowhere-zero section $s_0$. Indeed the Cartier sequence on $K$
has kernel $\OO_K$ and $H^1(K,\OO_K)=0$. Over $\mathcal F_2(Q)$
the pullback of $h^{[2]}$ thus satisfies
\[
 0\to\OO\xrightarrow{s_0}U_3^\vee\to h^{[2]}|_{\mathcal F_2(Q)}\to0 .
\]
The top Chern class of $h^{[2]}$ is $e$. With $E_4$
the quartic restriction bundle on $\mathcal H$,
\[
 \int_{X_2}De=\int_{\mathcal H} c_9(E_4)c_2(U_3^\vee)=5016 .
\]
Explicitly the fixed-point expression, with the weights $I,a,R_d,T(I,a)$
of Lemma~\ref{lem:threeconics-18}, is
\[
 \sum_{I,a}\frac{(\prod R_4)\sum_{i<j,\ i,j\in I}ij}{T(I,a)} .
\]
The condition $K\cap S\subset\Gamma$ is the zero locus of the
induced section of $\OO_S(5h)^{[2]}$. The corresponding intersection
$\int_{X_2}Dg=258200$ is the transverse quartic-root calculation of
that same lemma. We calculate the collision independently on
$\mathcal F_2(Q)$. Set $H=c_1(U_3^\vee)$ and
$\xi=c_1(\OO_{\mathcal H}(1))$ in the lines convention. The conic
equation sequence in its plane gives
\[
 0\to\OO_{\mathcal H}(-1)\to\operatorname{Sym}^2 U_3^\vee
 \to E_2\to0 .
\]
Restriction from $K$ to $K\cap S$ further quotients $E_2$ by
$s_0 U_3^\vee$. Since $H^1(K,\OO_K(1))=0$, the pullbacks of the
tautological first Chern classes are
\[
 c_1(h^{[2]})=H,\qquad c_1((2h)^{[2]})=3H+\xi .
\]
Let $L_{(2)}$ on $X_2$ denote the divisor class with ordered pullback
$h_1+h_2$, and let $\delta=[\mathcal E]/2$. The ordered exact sequence
of Lemma~\ref{lem:route19-two-contact} gives
$c_1((kh)^{[2]})=kL_{(2)}-\delta$. Taking $k=1,2$ yields
$L_{(2)}|_{\mathcal F_2(Q)}=2H+\xi$ and
$\delta|_{\mathcal F_2(Q)}=H+\xi$. As $\tau=\delta L_{(2)}$, it
follows that
\[
 \int_{X_2}D\tau=\int_{\mathcal H}c_9(E_4)(H+\xi)(2H+\xi)=13280.
\]
The fixed-point sum computing $5016$ above now has $\sum_{i<j}ij$
replaced by $(\sum I-a)(2\sum I-a)$, since $\xi|_{(I,a)}=-a$.
Thus $258200=25\cdot5016+10\cdot13280$ holds with its three
numbers evaluated separately. The conclusion follows from \eqref{eq:mark20-abel2}. Here $A_D$
parametrizes the boundary ideals $I_{K\cap S/C}$ on the incidence
of hyperplanes containing $K$: the two points span the intersection
of the conic plane with the hyperplane defining $S$, also for a
nonreduced intersection.

Indeed a plane in the hyperplane defining
$S$ contains $K\cap S$ exactly when it contains this line. It then
spans with the conic plane a unique hyperplane in $\mathbb P^4$
containing $K$. The numerical value $-516400$ equals the
quartic-conic component of Theorem~\ref{thm:marked-conics-18}.
\end{proof}
These calculations concern cycles on $J_{-4}$, including the
collision cycle $A_\tau$. A comparison of complete threefold theories
requires a relative virtual class or insertion realizing its pairing
with the quartic state and a compatible K3-fibration series with that
insertion. Equation~\eqref{eq:mark20-twoprobe} specifies its
coefficient and charge sector. It does not evaluate that series.
There is also a support restriction on its possible realization. For
$\iota:\Pic^{-2}(C_*)\hookrightarrow J$ one has
\[
 \iota_* A_3(\Pic^{-2}(C_*))_\QQ=\QQ f,\qquad
 \big(\deg_J(A_D A_\tau)\big)_D\notin
       \QQ\big(\deg_J(A_D f)\big)_D .
\]
The first equality follows because the fiber is integral of dimension
three; the second follows from the nonsingular matrix in
Proposition~\ref{prop:mark20-twoprobe}. Thus any proper insertion
whose middle-dimensional boundary pushforward still factors through
this single Picard fiber sees only the first of the two measurements.

\begin{proposition}[Moving ruled components for the collision class]
\label{prop:mark20-moving}
Fix a smooth $H\in|h|$. Let $D_H=\varpi^{-1}(H)\subset\mathcal E\subset X_2$
with inclusion $k:D_H\hookrightarrow X_2$. Its points are nonreduced
length-two subschemes $Z$ supported at $x\in H$, with tangent line
$[v]\in\PP(T_xS)$. Put $T=U_2\times_{X_2}D_H$. This smooth projective
threefold carries a family of stable pure surface sheaves on moving ruled
components with boundary morphism $e:T\to J_{-4}$ and
\[
                       e_*[T]=A_\tau .
\]
\end{proposition}
\begin{proof}
Because $D_H=\PP(T_S|_H)$ and $p:U_2\to X_2$ is a projective-line
bundle, $T\to D_H\to H$ is a tower of two projective-line bundles.
Pull back the universal plane curve and the universal length-two subscheme
to obtain $\mathcal Z\subset\mathcal C\subset S\times T$. Both
$\OO_{\mathcal Z}$ and $\OO_{\mathcal C}$ are $T$-flat, whence
$\mathcal I=I_{\mathcal Z/\mathcal C}$ is $T$-flat with fibers $I_{Z/C}$.

Write $P=\PP_S(\OO_S\oplus\OO_S(h))$ with its two disjoint sections
$S_0,S_\infty$. Blow up $P\times T$ along
$\mathcal C\subset S_0\times T$ and denote the result by
$\mathcal W\to T$. Locally at the center the ideal is $(z,g)$,
where $z$ cuts out $S_0$ and $g$ cuts out $\mathcal C$ on $S\times T$.
These form a regular sequence on each fiber.

The blowup is the hypersurface
$zv=gu$ with $[u:v]=[z:g]$; its two affine equations are
$zv=g$ and $z=gu$. Each equation is a nonzerodivisor on every
ambient fiber. This proves flatness over $T$ and compatibility with
base change. The strict transform $j:\mathcal R\hookrightarrow\mathcal W$
of the inverse image of $\mathcal C$ is
\[
 b:\mathcal R=P|_{\mathcal C}\longrightarrow\mathcal C.
\]
Indeed $\mathcal C$ is a Cartier section in this inverse image; in
the chart $u=1$ its strict transform is $v=0$, $g=0$.
It misses the strict transform of $S_0\times T$ and meets
$S_\infty\times T$ along the section $\mathcal C$.

Set $\mathcal F=j_*b^*\mathcal I$. Since $b$ is flat, $b^*\mathcal I$
is $T$-flat; pushforward along the closed immersion preserves
$T$-flatness and commutes with base change. Each $C$ is integral.
Its ideal $I_{Z/C}$ is torsion-free of rank one, as is
$b_C^*I_{Z/C}$ on the integral ruled surface $P|_C$, by the
affine-bundle charts.

Hence $F=j_{C*}b_C^*I_{Z/C}$ is pure and
Gieseker stable for every ample polarization. Near $S_\infty$ the
blowup is an isomorphism. In a ruling chart with coordinate $w$
cutting out this section, the restriction of $\mathcal F$ is computed
by the complex
\[
 [\mathcal I[w]\xrightarrow{w}\mathcal I[w]]\simeq\mathcal I
\]
in degrees $-1,0$, including at singular points of $C$. Thus
$e(Z,C)=[I_{Z/C}]$. The square defining $T$ is a flat base change
of $k$; consequently $[T]\mapsto p^*(k_*[D_H])=p^*(\tau\cap[X_2])$
on $U_2$. Pushforward by $a$ proves the formula. For a smooth $C$ transverse
to $H$, the fiber of $T\to|h|$ consists of the four subschemes
$2x\subset C$ with $x\in C\cap H$.

If $2x\sim2y$ with $x\ne y$, the smooth plane quartic would carry a basepoint-free $g^1_2$. Thus the four values $\OO_C(-2x)$ are distinct, and $\deg_J(A_\tau f)=4$.
\end{proof}
The moving collision gives a Chow class. Comparison with a fixed threefold, whose blowup may be singular above a singular point of $C$, requires the virtual comparison below.

\begin{corollary}[Multiplicity on the smooth moving-curve locus]
\label{cor:mark20-genericslice}
Let $V\subset |h|$ consist of smooth curves transverse to $H$.
In the moving ruled sheaf problem over $V$, the collision incidence
has multiplicity one at each of its four points over every $C\in V$.
If a compactified virtual incidence class in $A_3(T)_\QQ$ restricts to
this class over $V$, then it equals $[T]$.
\end{corollary}
\begin{proof}
Fix $C\in V$. For $L'\in\Pic^{-2}(C)$, set $L_R=b_C^*L'$ and
$G=j_{C*}L_R$. For the smooth Cartier embedding
$j_C:R_C\hookrightarrow W_C$, the sheaf-Ext spectral sequence has rows
\[
 E_2^{p,0}=H^p(R_C,\OO_{R_C}),\qquad
 E_2^{p,1}=H^p(R_C,N_{R_C/W_C})=0 .
\]
The second equality follows from
$N_{R_C/W_C}|_{b_C^{-1}(x)}=\OO_{\PP^1}(-1)$ as in
Proposition~\ref{prop:rigidruledfull}. Also $\mathbf Rb_{C*}\OO_{R_C}=\OO_C$.
Thus $\Ext^1_{W_C}(G,G)=H^1(C,\OO_C)$ and
$\Ext^2_{W_C}(G,G)=0$. In the relative moduli over $V$ the deformation
of $L'$ induces the edge isomorphism on tangent spaces.

The smooth
relative Picard variety gives an \'etale chart of the scalar-rigidified
relative sheaf moduli at these objects. Boundary evaluation on this
chart is the identity of $\Pic^{-2}(\mathcal C_V/V)$, by restriction
to the section. Indeed the maps from line-bundle families may be
constructed \'etale locally on the Picard base and agree after scalar
rigidification.

The intersection $H\cap\mathcal C_V\to V$ is finite \'etale of degree
four, and the length-two subscheme lying on $C$ at such a point $x$
is uniquely $2x\subset C$. This identifies $T_V$ with that cover
scheme-theoretically: the tangent direction of a double point in a
smooth curve is its tangent line, also in families.

On an \'etale splitting of the cover the four values give sections $s_i(C)=\OO_C(-2x_i)$ of the smooth relative Picard variety. Each is a regular immersion of codimension three, hence $s_i^![\operatorname{Pic}]=[V]$ on that chart. The vanishing of $\Ext^2$ above identifies this ordinary calculation with the virtual one.

Finally $T$ is integral of dimension three and $T_V$ is dense in
$T$. The localization sequence
\[
 A_3(T\setminus T_V)_\QQ\longrightarrow A_3(T)_\QQ
 \longrightarrow A_3(T_V)_\QQ\longrightarrow0
\]
has zero left term. The class of any extension on this same $T$ is
therefore uniquely determined. Additional components of a different
virtual compactification require their own specialization calculation.
\end{proof}

\subsection{The virtual insertion on a fixed pencil component}
Write $B=|h|$, $J=J_{-4}$, $\rho:J\to B$. In this subsection
$C\in B$ is the base curve of the pencil on the quartic threefold.
Its complementary component is $W_C$, with the ruled support
$R_C\to C$ of Proposition~\ref{prop:mark20-moving}. Keep the
hyperplane curve $H\in |h|$ which defines the collision condition
fixed. Set $J_C=\rho^{-1}(C)$, $\iota_C:J_C\hookrightarrow J$,
and $i_C:\{C\}\hookrightarrow B$. Write $i_C^!$ for the refined
pullback along this base point, on any space mapping to $B$.

\begin{theorem}[Fixed-pencil collision class]\label{thm:step21-fixedslice}
Let $C$ be any smooth plane section of $S$, and set
$\nu_C:C\to J_C=\Pic^{-2}(C)$, $x\mapsto\OO_C(-2x)$.
Let $\mathcal P_C$ be the ruled sheaf component of the relative
problem on $W_C$ with evaluation $e_C:\mathcal P_C\to J_C$
as in Proposition~\ref{prop:mark20-twoprobe}.
Define an operational class on the smooth $J_C$ by
\[
  \xi_C=(\nu_C)_*\big(c_1(\omega_C)\cap[C]\big)
     \in A_0(J_C)=A^3(J_C).
\]
Then
\begin{equation}\label{eq:step21-fixedslice}
  i_C^! A_\tau
   =(e_C)_*\big(e_C^*\xi_C\cap[\mathcal P_C]^{\vir}\big)
   =\xi_C\qquad\text{in }A_0(J_C).
\end{equation}
If $C\ne H$ this class has the representative
$\sum_{x\in C\cap H}i_x(C,H)[\OO_C(-2x)]$. Its degree is four.
The formula also applies to $C=H$, where the collision fiber has
dimension one.
\end{theorem}
\begin{proof}
Let $B^{\rm sm}\subset B$ be the smooth-curve locus, with universal
curve $\pi:\mathscr C\to B^{\rm sm}$ and projection $u:\mathscr C\to S$.
The first infinitesimal neighborhood of the relative diagonal
in $\mathscr C\times_{B^{\rm sm}}\mathscr C$ is finite flat of degree
two over $\mathscr C$. The ideal sequence of this Cartier double
point defines the relative morphism
\[
 \nu:\mathscr C\longrightarrow J|_{B^{\rm sm}},
 \qquad (x,C')\longmapsto\OO_{C'}(-2x).
\]
In a smooth family of curves the tangent direction of a double point
centered at $x$ is the relative tangent line. Thus, scheme-theoretically
over $B^{\rm sm}$, the incidence $T$ of Proposition~\ref{prop:mark20-moving}
is $u^{-1}(H)$ in $\mathscr C$. Indeed the equation imposing the tangent
condition on $\PP(T_S)$ is the vanishing of the differential of the
curve equation on the tautological line.

On each smooth curve its
kernel is a line subbundle, whose projectivization is the curve itself.
The morphism $u$ is smooth: the full universal plane-section incidence
over $S$ is a projective-plane bundle, and $\mathscr C$ is open in it.
Consequently $T|_{B^{\rm sm}}$ is a Cartier divisor of class $u^*h$.
Over $B^{\rm sm}$ the proper pushforward description reads
\[
 A_\tau=\nu_*(c_1(u^*\OO_S(H))\cap[\mathscr C]).
\]
Proper pushforward and first Chern class commute with refined
pullback \cite[Theorem~6.2 and Proposition~6.3]{Fulton}.
Flatness of $\pi$ gives $i_C^![\mathscr C]=[C]$, and
adjunction on the quartic K3 gives $\OO_S(H)|_C=\omega_C$.
Thus $i_C^!A_\tau=\xi_C$. For $C\ne H$ the section cutting out $H$
restricts to a nonzero section on $C$, whose divisor gives the
indicated multiplicities. For $C=H$, refined intersection gives
$c_1(N_{C/S})\cap[C]$ on the entire excess fiber $C$. This equals
$c_1(\omega_C)\cap[C]$ by the same adjunction calculation.

On the fixed $W_C$ the sheaves are $j_*b_C^*L'$,
$L'\in\Pic^{-2}(C)$. By Corollary~\ref{cor:mark20-genericslice}'s
Cartier self-Ext calculation, which uses smoothness of $C$ alone,
$\Ext^1_{W_C}=H^1(C,\OO_C)$, $\Ext^2_{W_C}=0$ for these
self-Ext groups. The map from the Picard variety is proper and
induces the isomorphism on tangent spaces, and its boundary evaluation
is its inverse on this component. Hence
$[\mathcal P_C]^{\vir}=[\mathcal P_C]$ and
$e_{C*}[\mathcal P_C]^{\vir}=[J_C]$.
The projection formula gives the middle equality of
\eqref{eq:step21-fixedslice}.
\end{proof}

There is also a singular-curve form of the slice calculation
which uses only Chow groups. Write $D_H=\PP(T_S|_H)$, with
projection $q:D_H\to H$. On $D_H$ the rank-two bundle $E$
of evaluations of $\OO_S(h)$ on the universal double point
has an exact sequence
\[
 0\to q^*\OO_H(h)\otimes\OO_{D_H}(1)\to E
          \to q^*\OO_H(h)\to0 .
\]
Its quotient is the value at the center, and its kernel is the value
along the cotangent line dual to the specified direction. The map
$H^0(S,h)\otimes\OO\to E$ is surjective because $h$ is very ample.
For a curve $C=(s=0)\in B$, the fiber $T_C$ is the zero scheme
$Z(s^{[2]})$ on $D_H$. In the following formula all Chow groups
of supports carry rational coefficients.
\begin{proposition}\label{prop:step21-singularpull}
If $C\ne H$ is integral, then
\[
 i_C^![T]=\sum_{x\in H\cap C} i_x(C,H)[v_x]\quad\text{in } A_0(T_C).
\]
At a smooth $x\in C$ put $v_x=T_xC$. At a singular $x\in C$
one may take any line $v_x\subset T_xS$. The same summands pushed to
$J_C$ are the cycle $i_C^! A_\tau$ in its compactified Jacobian.
\end{proposition}
\begin{proof}
The universal incidence in $D_H\times B$ is the regular zero scheme
of the section of $E\boxtimes\OO_B(1)$ induced by evaluation.
Indeed it is a projective-line bundle over $D_H$ and has the expected
codimension two. Refined base change for this zero locus
\cite[\S6.2 and \S14.1]{Fulton} identifies
$i_C^![T]$ with the localized top Chern class on $Z(s^{[2]})$.
The quotient section is $s|_H$ and cuts $H\cap C$ with length
$m_x=i_x(C,H)$ at $x$.

Locally split the bundle sequence and intersect successively with the quotient section $s|_H$ and the kernel section $ds_x(v)$. The first Cartier pullback gives $\sum_xm_x[q^{-1}(x)]$. On $q^{-1}(x)=\PP(T_xS)$ the residual line bundle is $\OO(1)$ and the second localized class is
\[
 \begin{cases}
 [T_xC], & ds_x\ne0,\\
 c_1(\OO(1))\cap[\PP(T_xS)]=[v_x], & ds_x=0 .
 \end{cases}
\]
Changing the splitting adds a multiple of the first section to the second. Proper pushforward to $J_C$ commutes with $i_C^!$.
\end{proof}
For a singular base curve the preceding proposition computes the
slice of the proper incidence class. An identification with the
fixed singular threefold sheaf obstruction theory would require a
virtual comparison for that target.

\begin{corollary}\label{cor:step21-slice-kernel}
The collection of all smooth fixed-pencil slice classes does not
determine the conic pairing with a middle boundary cycle. More
precisely, fix a smooth $C_0\in B$ and set
$f=(\iota_{C_0})_*[J_{C_0}]$. Then
\[
 i_C^! f=0\quad(C\in B^{\rm sm}),\qquad
 \deg_J(A_D f)=5016
\]
for the quartic conic class $A_D$ of
Corollary~\ref{cor:mark20-conictest}.
\end{corollary}
\begin{proof}
For $C\ne C_0$ the supports are disjoint. For $C=C_0$ the
self-intersection is the top Chern class of
$T_{C_0}B\otimes\OO_{J_{C_0}}$, a trivial rank-three bundle.
The conic pairing was computed in Corollary~\ref{cor:mark20-conictest}.
In particular $A_\tau$ and $A_\tau+f$ have identical fixed-pencil
zero-cycles and their conic pairings differ by $5016$.
\end{proof}

\subsection{A correspondence from the fixed relative ruled component}
We now fix the original pencil once and for all. Denote its smooth base
curve by $H=C_*\subset S$, its complement by $W_H$, and its rigid
ruled surface by $R_H\subset W_H$. The letter $H$ in this subsection
refers to this curve. The auxiliary curve in the definition
$\tau=j_*\varpi^*h$ can be taken to be $H$: this substitution uses
linear equivalence in $A^1(S)$, before any numerical pairing.

Write $P_H=\Pic^{-2}(H)$, identified by boundary evaluation with the
proper relative sheaf component on $W_H$ in Theorem~\ref{thm:step21-fixedslice}.
Let $k_H:P_H\hookrightarrow J$ be the inclusion.
An object of $P_H$ represents the surface sheaf
$j_{H*}b_H^*L$ with $L\in\Pic^{-2}(H)$.
There is a correspondence
\[
\begin{tikzcd}[column sep=large]
 P_H & H\arrow[l,"\delta_H"']&
 T_H\arrow[l,"q_H"']\arrow[r,"a_H"]&J .
\end{tikzcd}
\]
Here $\delta_H(x)=\OO_H(-2x)$ and a point of $T_H$ is
$(x,[v],C)$, where $x\in H$, $[v]\in\PP(T_xS)$, and the
plane quartic $C$ contains the length-two tangent subscheme $Z(x,v)$.
The map $q_H$ forgets $[v],C$; $a_H$ takes the value $[I_{Z(x,v)/C}]$.
The incoming sheaf is on the fixed $W_H$. The outgoing sheaf is on
the K3, with its supporting curve allowed to vary in $|h|$.

\begin{theorem}[Correspondence with the fixed relative virtual class]
\label{thm:step21-flag}
The proper correspondence defines a homomorphism
\[
 \mathcal K_H=(a_H)_*q_H^*\delta_H^!:A_i(P_H)_\QQ\longrightarrow A_i(J)_\QQ .
\]
At the fixed ruled component it satisfies the following Chow identity:
\begin{equation}\label{eq:step21-flag}
 \mathcal K_H([P_H]^{\vir})=A_\tau\quad\text{in }A_3(J)_\QQ .
\end{equation}
For every $D\in A^2(S^{[2]})_\QQ$, the two-contact comparison is
\begin{equation}\label{eq:step21-flag-eval}
 \deg_J(A_D A_g)=
 -50\deg_J\big(A_D\,(k_H)_*[P_H]^{\vir}\big)
 +10\deg_J\big(A_D\,\mathcal K_H([P_H]^{\vir})\big).
\end{equation}
\end{theorem}
\begin{proof}
After a translation on the Picard varieties $\delta_H$ is the
Abel map $x\mapsto\OO_H(x)$ followed by multiplication by $-2$.
The latter is \'etale, and the Abel map is an immersion because
$\omega_H=\OO_H(1)$ is basepoint free. If $2x\sim2y$ for $x\ne y$,
these two disjoint divisors span a basepoint-free pencil of degree
two, contradicting the canonical plane embedding.

Thus $\delta_H$ is
proper, unramified and universally injective over $\CC$, hence a
closed immersion. Its source and target are smooth of dimensions
one and three, so it is a regular immersion of codimension two.
The bundle maps
\[
 T_H\longrightarrow D_H=\PP(T_S|_H)\longrightarrow H
\]
are both projective-line bundles by Lemma~\ref{lem:mark20-abel2}.
In particular $q_H$ is smooth of relative dimension two.
On the ruled sheaf component $[P_H]^{\vir}=[P_H]$; hence
\[
 \delta_H^![P_H]^{\vir}=[H],\qquad q_H^*[H]=[T_H],\qquad
 (a_H)_*[T_H]=A_\tau .
\]
The last identity is Proposition~\ref{prop:mark20-moving} applied
with this representative of $h$. It includes every curve in the
proper flag family, including singular outgoing plane quartics.
Equation~\eqref{eq:step21-flag-eval} follows from
Proposition~\ref{prop:mark20-twoprobe} and
$(k_H)_*[P_H]^{\vir}=f$.
\end{proof}
Ordinary matching with the quartic uses $k_H$ and evaluates both sheaves on $H$. The insertion also permits an outgoing support $C\ne H$. Its contribution to the pencil degeneration requires a virtual comparison carrying this additional boundary datum.

\begin{proposition}[The zero-cycle of the inserted fixed component]
\label{prop:step21-cap}
In Theorem~\ref{thm:step21-flag}, set $H=C_*$ and write
$P_H=\Pic^{-2}(H)$, $\delta_H(x)=\OO_H(-2x)$ and
$\xi_H=(\delta_H)_*(c_1(\omega_H)\cap[H])$. Set
$D_H=\PP(T_S|_H)$, $k:D_H\hookrightarrow X_2=S^{[2]}$ and
$\pi:D_H\to H$. Write
\begin{gather*}
 T_H=U_2\times_{X_2}D_H,\qquad r:T_H\to D_H,\qquad q_H=\pi r,\\
 a_H:T_H\longrightarrow U_2\xrightarrow{a}J,\qquad
 m_H=\delta_H q_H:T_H\to P_H .
\end{gather*}
For a codimension-three class $A$ on $J$ define the proper
correspondence
\[
 \mathcal R_H(A)=(m_H)_*\big(a_H^*A\cap[T_H]\big)\in A_0(P_H)_\QQ .
\]
If $A=e_*[M]^{\vir}$ for a proper relative evaluation, the argument
of $(m_H)_*$ is the pushforward to $T_H$ of
$\Delta_J^!([M]^{\vir}\times[T_H])$ from $M\times_J T_H$.
The subsequent matching with the fixed ruled component
$\mathcal P_H\xrightarrow{\sim}P_H$ leaves this zero-cycle unchanged.
For every $D\in A^2(X_2)_\QQ$ one has
\begin{equation}\label{eq:step21-fixedcap}
 \mathcal R_H(A_D)
 =-\tfrac12\!\left(\int_{X_2}D\tau\right)\xi_H .
\end{equation}
In particular the quartic conic cycle gives
$\mathcal R_H(A_D)=-6640\,\xi_H$ in the Chow group of the fixed
Picard variety.
\end{proposition}
\begin{proof}
On the smooth $J$ we identify cycles and operational classes by
capping with $[J]$. Refined diagonal intersection and proper
pushforward give the stated base change description of $a_H^*A$.
The evaluation of the fixed smooth component is an isomorphism and
its virtual class is fundamental by Theorem~\ref{thm:step21-fixedslice}.
The refined pullback of an exterior product with this class by
$\Delta_{P_H}$ equals the original class, by the identity
correspondence.

For the closed immersion $a$ let $N=N_{U_2/J}$. By its
self-intersection formula and $p_*c_1(N)=-2$ from
Lemma~\ref{lem:mark20-abel2}, one has on $T_H$
\[
 a_H^*A_D=c_1(N|_{T_H})r^*k^*D,qquad
 r_*\big(c_1(N|_{T_H})\cap[T_H]\big)=-2[D_H].
\]
The second equality follows from the degree on each projective-line
fiber. Hence
\[
 \mathcal R_H(A_D)=-2(\delta_H)_*\pi_*(k^*D\cap[D_H]).
\]
Let $j:\mathcal E=\PP(T_S)\hookrightarrow X_2$ be the full exceptional
divisor, $\varpi:\mathcal E\to S$ its projection, and $i:H\hookrightarrow S$.
Flatness of $\varpi$ and refined base change in its Cartesian square
over $i$ give
\[
 \pi_*(k^*D\cap[D_H])=i^!\varpi_*(j^*D\cap[\mathcal E]).
\]
The last pushforward on $S$ is a divisor. Since $\Pic(S)=\mathbb Zh$,
write it as $d h\cap[S]$. Computing its intersection with $h$ gives
\[
 4d=\int_{\mathcal E}j^*D\cdot\varpi^*h
    =\int_{X_2}D\cdot j_*(\varpi^*h)=\int_{X_2}D\tau .
\]
Now $i^*h=c_1(\omega_H)$ by adjunction, proving
\eqref{eq:step21-fixedcap}. For the conic cycle,
Corollary~\ref{cor:mark20-conictest} gives $\int D\tau=13280$.
\end{proof}
The refined intersection $\mathcal R_H$ uses both maps from $T_H$, to the fixed Picard component and to the full Jacobian. To obtain it from a threefold invariant one must extend the insertion through the total pencil degeneration.

\begin{proposition}[Supports at the ends of an ordinary bubble chain]
\label{prop:step21-bubble}
Consider a chain of the ordinary bubbles for the internal K3 of the
second degeneration, with $\Pic(S)=\mathbb Zh$, and fix the
component identifications with $S$.
Assume the branch sheaves are pure, their derived restrictions match
at the intervening joins, and the two end values are pure curve
sheaves with first Chern class $h$. Then their fundamental support
curves coincide on $S$. In particular the pairs of boundary values
with supporting curves $H,C$ and $C\ne H$ in
Theorem~\ref{thm:step21-flag} do not arise by inserting such a chain.
\end{proposition}
\begin{proof}
Orient each bubble as in \textup{(\companionref{eq:generalcomponentclasses})}.
On $\Delta_j=\PP_S(\OO\oplus\OO(h))$ write its effective divisorial
cycle as $u_jp_j^*h+s_jD_{j,-}$, and the successive boundary
first Chern classes as $b_jh$. Intersecting with a general ruling
and restricting the determinant to the two sections give
\[
  s_j\ge0,\qquad b_j=u_j,\qquad b_{j-1}=u_j-s_j.
\]
These identities use only the bubble and its normal bundles, so apply
equally to the internal K3 of this degeneration. Thus
$0=b_{\rm last}-b_{\rm first}=\sum_j s_j$ forces $s_j=0$,
$u_j=1$ on each bubble. By $p_{j*}\OO_{\Delta_j}=\OO_S$, every
effective representative of $p_j^*h$ is $p_j^{-1}C_j$ for an
integral $C_j\in|h|$.

Neither section is a component of this divisor:
the local equation of either section acts injectively on the pure branch sheaf.
Its restriction is supported on $C_j$, and its first Chern class
$h$ makes $C_j$ its fundamental curve. Matching identifies these
curves across consecutive joins, hence at the two ends.
\end{proof}

\subsection{The ruled-sheaf obstruction theory at singular members}
The incidence calculation can be compared with the threefold sheaf
obstruction theory over a larger open set than the smooth members.
Write $B=|h|=(\PP^3)^\vee$ and
$\nu:S\to B$, $x\mapsto \mathbb T_xS$, for the Gauss morphism.
Choose the fixed smooth hyperplane curve $H$ generally if necessary.
In particular $H$ is not contained in the locus where the second
fundamental form $\mathrm{II}$ of $S\subset\PP^3$ is zero.
Write $\pi:T\to B$ for the incidence of
Proposition~\ref{prop:mark20-moving}.

\begin{proposition}\label{prop:step21-singular}
For $t=(x,[v],C)\in T$ put $I_t=I_{Z/C}$.
The noninvertible locus of the $I_t$ has underlying set
\[
 D=\{(x,[v],C):x\in H,\ C=S\cap\mathbb T_xS,                   \mathrm{II}_x(v,v)=0\},\qquad \dim D\le1.
\]
On $T\setminus D$, the moving ruled sheaf $G_t=j_{C*}b_C^*I_t$ is
perfect on $W_C$, also for singular $C$. Its scalar-rigidified
fixed-target deformation and obstruction spaces are
$H^1(C,\OO_C)$ and $0$. In the moving relative Picard chart the
collision incidence is regular of codimension three and its refined
virtual class equals $[T\setminus D]$.
\end{proposition}
\begin{proof}
Take formal parameters $(s,t')$ on $S$ with
$I_{Z/S}=(s,(t')^2)$ and $v=\partial/\partial t'$ at $x$.
The equation of $C$ can be written
$g=as+b(t')^2$. On $C$, the two generators of $I_t$ have relations
$(-(t')^2,s)$ and $(a,b)$. Indeed for any relation lift its
coefficients to the regular local ring, subtract a multiple of
$(a,b)$, and use the Koszul relation on $s,(t')^2$.
Nakayama now gives
\[
 \dim_{\CC} I_t/\mathfrak m I_t
    =2-\operatorname{rank}\begin{pmatrix}0&a(x)\\0&b(x)\end{pmatrix}.
\]
A principal ideal of finite colength on this Cohen--Macaulay curve is
generated by a nonzerodivisor. The condition $a(x)=0$ is exactly
$dg(x)=0$, i.e. the plane cutting $C$ is $\mathbb T_xS$.
In this case $b(x)=0$ means the quadratic term of $g$ vanishes
on $v$, equivalently $\mathrm{II}_x(v,v)=0$.

For the Gauss map $\gamma$ one has $\gamma^*\OO_B(1)=\OO_S(3h)$ by the four cubic partials. Ampleness gives finite fibers and separability gives generic rank two, so $\mathrm{II}$ is generically nonzero. Along the general curve $H$, the equation $\mathrm{II}_x(v,v)=0$ has at most two directions except at finitely many $x$, where the fiber has dimension at most one. Thus $\dim D\le1$.

The ruled surface $R_C=P|_C$ is Cartier in $W_C$: in the
chart $zv=g$ its equation is $v=0$, and $v$ is a nonzerodivisor.
Away from the exceptional locus this follows from the equation of
$C$ on $S$. For $t\notin D$, the line bundle $L_R=b_C^*I_t$
therefore has pushforward of projective dimension one.
The local-to-global spectral sequence has rows
\[
 H^p(R_C,\OO_{R_C}),\qquad H^p(R_C,N_{R_C/W_C}).
\]
The degree of $N_{R_C/W_C}$ on each ruling fiber is $-1$, by
the exceptional-section calculation of
Proposition~\ref{prop:rigidruledfull}, valid on the projective
bundle over the integral curve $C$. Hence
$\mathbf Rb_{C*}N_{R_C/W_C}=0$ and
$\mathbf Rb_{C*}\OO_{R_C}=\OO_C$. The first row gives the
asserted Ext spaces. The relative Picard scheme of line bundles is
smooth over $B$ since obstructions to lifting a line bundle on a
curve lie in $H^2(C,\OO_C)=0$.

Pullback to $R_C$ and Cartier
pushforward induce an isomorphism on relative tangent spaces by
the first row; the construction is \'etale on the scalar-rigidified
sheaf chart. Its evaluation is restriction to the ruling section,
the inclusion of the line-bundle locus into $J_{-4}$.

Finally $T\hookrightarrow U_2\hookrightarrow J_{-4}$ is a closed
immersion: the first inclusion is the pullback along the smooth $p$
of $D_H\hookrightarrow S^{[2]}$ and the second is
Lemma~\ref{lem:mark20-abel2}. Restrict to the Picard locus. Both
source and target are smooth, of dimensions three and six, so the
immersion is regular. The vanishing of $\Ext^2$ and the regular
Gysin formula give the stated refined virtual class.
\end{proof}

\begin{corollary}\label{cor:step21-testpencil}
Choose a general line $\Lambda\subset B$ of moving base curves.
Then $T_\Lambda=T\times_B\Lambda$ is a smooth projective curve,
and its intersection class with the moving relative ruled-sheaf
chart has multiplicity one through all members of $\Lambda$.
If $\rho:J_{-4}\to B$ is support, the two evaluation classes satisfy
\[
    (e|_{T_\Lambda})_*[T_\Lambda]
       = c_1(\rho^*\OO_B(1))^2\cap A_\tau
            \quad\text{in } A_1(J_{-4})_\QQ .
\]
\end{corollary}
\begin{proof}
A general line misses the curve $\nu(H)$ and the point $[H]$.
Every pair $(x,C)$ with $C\in\Lambda$ and $x\in C\cap H$ thus
satisfies $dg_C(x)\ne0$. The unique double point in $C$ at $x$ has
direction $T_xC$, also in families by the relative tangent line.
The base line in $\PP^3$ of the pencil of planes $\Lambda$ misses $H$.
For each $x\in H$ exactly one plane of $\Lambda$ contains $x$.
Consequently $T_\Lambda$ is the graph of the linear projection
$H\to\Lambda$, of degree four.

This proof permits ramification of that
map and singularities of $C$ away from $x$.
The sheaf chart above is smooth over $B$, so its restriction to
$\Lambda$ is smooth of dimension four. The closed immersion of the
smooth graph has codimension three. Regular Gysin pullback on this
chart gives the graph's fundamental class.

On $T$ the two hyperplanes
defining $\Lambda\subset B$ cut out this graph scheme-theoretically
in the expected codimension. The projection formula for $e$ proves
the cycle equality. Here $\Lambda$ varies the \emph{base curve}
of the second degeneration; a fixed pencil on the quartic threefold
chooses a single $C_*\in B$.
\end{proof}

The exceptional set in Proposition~\ref{prop:step21-singular}
can contain genuinely nonperfect threefold sheaves. At an ordinary
node choose coordinates with $g=xy$ and tangent jet $Z$ having
ideal $(y,x^2)$ in $\CC[[x,y]]/(xy)$. The moving blowup chart
at the vertex has completed ring
$A=\CC[[x,y,z,v]]/(xy-zv)$; the ruled divisor is $v=0$.
The completed sheaf is
\[
       G^{\wedge}\simeq (y,x^2)_{A/(v)}
       \simeq A/(x,v)\oplus A/(y,v),
\]
because $(y)\cap(x^2)=0$ in $A/(v)$. The module $A/(x,v)$ has depth
two over the three-dimensional Cohen--Macaulay ring $A$. Finite
projective dimension would therefore give projective dimension one.
Its kernel $(x,v)$ in the free cover $A\twoheadrightarrow A/(x,v)$
would be free of rank one; the independent linear terms $x,v$
give two minimal generators.

This contradiction proves
nonperfectness at the vertex, despite the perfect boundary ideal
viewed on the smooth K3. In particular the length-two comparison
above extends across singular members of a general test pencil
without presupposing that every object over the full plane system
has the usual threefold perfect obstruction theory.

\begin{corollary}[Evaluation of the main incidence component]
\label{cor:step21-maincollision}
Let $P$ be a proper integral compactification of the open collision
incidence $T\setminus D$, with evaluation $e_P:P\to J_{-4}$
extending $e$. Then $e_{P*}[P]=A_\tau$. A virtual three-cycle on $P$
which restricts to the incidence class of
Proposition~\ref{prop:step21-singular} has the same pushforward.
\end{corollary}
\begin{proof}
The open graph is dense in $P$, and $e:T\hookrightarrow J_{-4}$ is
a closed immersion. Therefore, as closed sets,
$e_P(P)=\overline{e(T\setminus D)}=e(T)$.
The morphism $P\to e(T)$ has degree one on function fields, since
it is the identity on the common graph. Hence
$e_{P*}[P]=[e(T)]=A_\tau$.
As $\dim(P\setminus(T\setminus D))\le2$, localization identifies
$A_3(P)_\QQ$ with $A_3(T\setminus D)_\QQ$.
This proves the assertion about the virtual class. A proper matching
space for the full charge can include further irreducible components;
the assertion applies to the closure of the stated graph.
\end{proof}

\subsection{An exceptional-divisor restriction on the total pencil degeneration}
\label{sec:step22-family}
The fixed component used in Theorem~\ref{thm:step21-flag} has two
maps: restriction of a ruled surface sheaf to its base curve and an
outgoing map from the flag variety to $J_{-4}$. We extend the first
map over the degeneration. The order of intersection with the quartic
relative class determines which part of this construction can enter
an ordinary degeneration formula.

Write $H=C_*$ for the smooth pencil base curve and $s$ for the
parameter of $\mathscr Y$ in Proposition~\ref{prop:secondfamilyfull}.
Let $E\subset\mathscr Y$ be the exceptional divisor of the second
blowup. Its center is $H\times\AA^1_s$: blowing up $H\times\{0\}$
on $H\times\AA^1$ leaves that surface unchanged, and the chart
$x=sz$ identifies its strict transform with $z=u=0$.
Thus, with $N$ its rank-two normal bundle in $\mathscr D$, we have
\[
 H\times\AA^1_s\ \xleftarrow{\ \pi\ } E=\PP_{H\times\AA^1}(N)\ \xrightarrow{\ i\ }\ \mathscr Y .
\]
The central divisor $E_0$ lies on $W_H$ and misses the internal K3
$S_\infty$. Contractions of expansions at internal joins are
isomorphisms on a neighborhood of $E$.

\begin{proposition}[Restriction to the pencil base curve in families]
\label{prop:step22-index}
For a test scheme $T\to\AA^1_s$ and a family $\mathcal F$ of surface
sheaves perfect on a neighborhood of $E_T$, the index
\[
 \mathcal V_E(\mathcal F)=\mathbf R\pi_{T*}\mathbf L i_T^*\mathcal F
 \in\perf(H\times T),\qquad
 d_E(\mathcal F)=[\det\mathcal V_E(\mathcal F)]\in\Pic_{H/\CC}(T)
\]
commutes with arbitrary change of $T$. It defines a morphism of the
corresponding moduli functor to the Picard scheme, also after internal
expansion and scalar rigidification. On a smooth threefold fiber
$Y_s=Y$, let $j_b:F_b\hookrightarrow Y$ be a pencil member and
$k_b:H\to F_b$ the exceptional section. For a perfect complex $V$ on
$F_b$ one has
\begin{equation}\label{eq:step22-indexspec}
 \mathcal V_E(j_{b*}V)\simeq\mathbf Lk_b^*V .
\end{equation}
On the fixed ruled component the map $d_E$ is $L\mapsto L$ on
$P_H=\Pic^{-2}(H)$. The same formula holds for sheaves on internal
expanded targets which restrict near $E_0$ to the corresponding ruled
sheaf.
\end{proposition}
\begin{proof}
The projection $\pi_T$ is a smooth projective-line bundle.
The derived pullback to $E_T$ is perfect. Proper perfect
pushforward and flatness of $\pi_T$ give a perfect index and
derived base change \cite[Tag~0A1G]{StacksFam22}.
Formation of the determinant commutes with pullback. Tensoring
$\mathcal F$ by a line bundle $L_T$ from the test base multiplies
its determinant by the base line bundle $L_T^{\otimes\rk\mathcal V_E}$,
which has trivial class in $\Pic_{H/\CC}(T)$ (componentwise in $T$).
Every internal expansion isomorphism preserves the projection and
restriction near $E_T$.

On $Y_s$, $F_b$ and $E_s$ are Cartier divisors with intersection a
section of $E_s\to H$. The equation of $F_b$ is regular on $E_s$,
so the Koszul resolution shows that their intersection square is
Tor-independent. Writing $l_b:H\to E_s$ for the other inclusion gives
\[
 \mathbf L i_s^*(j_{b*}V)=l_{b*}\mathbf L k_b^*V,\qquad
 \pi_s l_b=\id_H .
\]
The argument does not require that the support of $V$ avoid $H$.
On $W_H$ the proper transform of $P|_H$ meets $E_0$ in a section.
In the chart $zv=g$ of Proposition~\ref{prop:mark20-moving}, $g$
cuts out the smooth $H$, the exceptional divisor has equation $z=0$
and the proper transform has equation $v=0$. The two divisors meet
transversely. For $G=j_{H*}b_H^*L$ the same Koszul calculation yields
$\mathcal V_E(G)=L$. The internal node is disjoint from this
calculation.
\end{proof}

The morphism in Proposition~\ref{prop:step22-index} makes the
fixed-curve incidence $d_E(\mathcal F)=\OO_H(-2x)$ a condition on
the total family in the degree $-2$ sector. On the general smooth
pencil member, sheaves with finite point defects illustrate what this
condition measures. If $V=L\otimes I_Z$ on a smooth $F_b$, then
\[
 \det\mathbf L k_b^*V=(\det V)|_H=L|_H
\]
for every zero-dimensional $Z\subset F_b$. Indeed $I_Z$ is perfect
on the smooth surface, and its determinant is canonically trivial
off $Z$ and extends trivially by codimension two. Points of $Z$ lying
on $H$ can affect the derived restriction while leaving its determinant
unchanged.

\begin{proposition}[Order of the fixed-component intersection]
\label{prop:step22-order}
Let $G$ be a proper marked matching space for the fixed pencil with
virtual dimension zero and determinant map $d_E:G\to P_H$. On the
space with extra flag $G\times_{P_H}T_H$, the refined pullback of the
fixed matched class along $T_H\to P_H$ is zero. In contrast, applying
the same pullback to the fixed relative component before matching
with a quartic cycle gives the second, generally nonzero, term of
\eqref{eq:step21-flag-eval}.
\end{proposition}
\begin{proof}
Factor $T_H\to P_H$ as
$T_H\xrightarrow{q_H}H\xrightarrow{\delta_H}P_H$.
The codimension of $\delta_H$ is two and the relative dimension of
$q_H$ is two. With $G_H=G\times_{P_H}H$ their base changes act as
\[
 A_0(G)\xrightarrow{\delta_H^!} A_{-2}(G_H)
       \xrightarrow{q_H^*} A_0(G\times_{P_H}T_H).
\]
For the scalar-rigidified Deligne--Mumford spaces considered here,
$A_{-2}(G_H)_\QQ=0$. On $[P_H]^{\vir}$ the intermediate class has
dimension $3-2=1$, so that
$q_H^*\delta_H^![P_H]^{\vir}=[T_H]$.
The latter class can be matched with the quartic on its outgoing map
$a_H:T_H\to J$. Theorem~\ref{thm:step21-flag} gives
$a_{H*}[T_H]=A_\tau$, and
Proposition~\ref{prop:step21-cap} evaluates its conic pairing as
$-26560$.
\end{proof}
In particular positive-codimension ordinary insertions pulled back
from $P_H$ cap a fixed matched virtual zero-cycle to zero, while an
insertion of codimension zero on the connected smooth $P_H$ acts by
a scalar. The fixed-component collision and Picard-fiber functionals
on the Abel cycles are linearly independent by
Proposition~\ref{prop:mark20-twoprobe}. The extension of the incoming
determinant therefore specifies only one of the two maps needed in a
total-family comparison.

\begin{proposition}[Specialization of vertical sheaf supports on the quartic root]
\label{prop:step22-support}
Let $R$ be a discrete valuation ring mapping to the smoothing, and let
$\mathcal F$ be an $R$-flat coherent sheaf on an internal expansion of
$\mathscr Y_R$. Suppose the geometric generic fiber is stable and pure
of class $r F_{\mathrm{fib}}$, $r\ge1$. Every two-dimensional
irreducible component of the specialized support in $Q\setminus S$
is contained in the limit pencil member $F_b^Q\subset Q$.

If that
member is $S$ the specialized support is empty on $Q\setminus S$.
Otherwise every two-dimensional root component there has closure
$F_b^Q$. If the whole root branch is pure of dimension two, has no
component supported on $S$, and has a nonzero pure curve restriction
on $S$, its fundamental boundary curve is supported on $H$ in every
such multiplicity.
\end{proposition}
\begin{proof}
After finite base change the stable vertical support lies on a single generic pencil fiber by Theorem~\ref{thm:outerbitrarymultiplicity}. Properness of the pencil extends its parameter to a primitive pair $[u:v]\in\PP^1(R)$. Put $L=us_0+vs_1$ and pull it back from $Q$. On $L\ne0$, flatness over $R$ gives, for a uniformizer $\pi$,
\[
 \mathcal F|_{L\ne0}\hookrightarrow\mathcal F[\pi^{-1}]|_{L\ne0}=0 .
\]
The special root support away from $S$ therefore lies in $V(\bar u s_0+\bar v s_1)\setminus S$. All pencil members are integral surfaces. Under the last hypotheses the root support is this member; for a member distinct from $S$ its intersection with $S$ is $V(s_0,s_1)=H$. Purity and nonvanishing force its boundary curve to have support $H$.
\end{proof}
The restriction to general-fiber limits in this proposition holds
before choosing a compactification or proving a virtual degeneration
formula. In particular, a point of the flag with outgoing support
$C\ne H$ and a root surface distinct from $S$ cannot occur as that
boundary of an ordinary vertical limit on the fixed pencil.

\subsection{Lifting the incidence by varying the pencil center}
The preceding virtual-dimension computation suggests varying the
base curve itself. Put $B=|h|$ and let
$\mathcal C\subset S\times B$ be the universal plane section.
Here $H\in|h|$ is a fixed auxiliary smooth curve as in
Proposition~\ref{prop:step21-singular}. The parameter $b\in B$ chooses
the possibly different base curve $C_b$ of the pencil on $Q$.

Keep $s$ for its threefold degeneration parameter.
Let $J^{\rm lin}\subset J_{-4}$ be the locus of line bundles on
the integral support curves. It represents the degree $-2$ relative
Picard functor of $\mathcal C/B$ after the usual sheafification.

\begin{theorem}[A total-family incidence for moving pencils]
\label{thm:step22-universal}
There is a flat projective family of threefolds
$\mathscr Y_B\to B\times\AA^1_s$ restricting to the pencil
degeneration with center $C_b$ for every $b$. On the part of
the sheaf functor perfect near the exceptional divisor, in the
component where the index determinant has degree $-2$, there is a
morphism
\[
       d_B:\mathfrak M_B\longrightarrow J^{\rm lin}.
\]
On $s\ne0$ it sends $j_*V$ on a smooth pencil surface to
$[\det V|_{C_b}]$ for perfect $V$. On the marked ruled line-bundle
component it is boundary evaluation for every integral $C_b$.
Write $T^{\rm lin}=T\cap J^{\rm lin}$ for the collision incidence.
Then $T^{\rm lin}\hookrightarrow J^{\rm lin}$ is a regular
immersion of codimension three. Its base change along $d_B$ defines
a relative insertion over the full smoothing parameter.

Suppose a chart of $\mathfrak M_B$ carries a relative virtual
class of virtual dimension zero over $B$ on each threefold fiber,
compatible with specialization in $s$ and with marked gluing.
At a matching with the ruled line-bundle component, assume the
quartic factors are induced by the same relative quartic chart
$M_Q$ over varying $b$, with boundary evaluation $e_Q:M_Q\to J_{-4}$.
The marked class after this insertion is the restriction to
$T^{\rm lin}$ of
\begin{equation}\label{eq:step22-univmatch}
       \Delta_{J_{-4}}^!
          \big([M_Q]^{\vir}\times[T]\big)
          \qquad\text{on }M_Q\times_{J_{-4}}T.
\end{equation}
At any singular plane quartic where the collision ideal is invertible,
the same comparison applies. No change of boundary support is required
at a gluing in this moving family.
\end{theorem}
\begin{proof}
Blow up $\mathscr D\times B$ along the relative Cartier curve
$\mathcal C\times\AA^1_s$ in its distinguished surface family.
More precisely, the strict transform of $S\times\AA^1_s$ in
$\mathscr D$ is $S\times\AA^1_s$, on which the center has equation
$g_b=0$, the universal plane equation. In the affine normal chart
the blowup has equations $zv=g_b$ and $z=g_bu$ as in
Proposition~\ref{prop:secondfamilyfull}. Both equations remain
nonzerodivisors after restriction to each geometric base fiber.

Near the center these are relative effective Cartier hypersurfaces in flat charts, proving flatness there. Away from it the family is the original flat degeneration. For singular $C_b$ this constructs the flat blowup fiber and makes no assertion of semistability of its total target. The center is a regular embedding with
exceptional divisor a projective-line bundle over
$\mathcal C\times\AA^1_s$. It misses all internal joins.

The
index proof of Proposition~\ref{prop:step22-index}, applied to this
bundle, yields the line bundle
$\det\mathbf R\pi_*\mathbf L i^*\mathcal F$ on $\mathcal C_T$
functorially in $T$. On the generic pencil surfaces the Koszul square
with the exceptional divisor gives its stated value. On $W_{C_b}$ the
ruled surface $P|_{C_b}$ is Cartier.

At its exceptional section the
completed chart is $zv=g_b$, with $v=0$ on the ruled divisor and
$z=0$ on the exceptional divisor. On the ring of the latter $v$ is
a polynomial coordinate, so is regular, also if $C_b$ is singular.
For a line bundle $L'$ this gives
$\mathbf R\pi_*\mathbf L i^*(j_*b_{C_b}^*L')=L'$.
The same identities hold in families by these charts and base change.

The line-bundle locus of the ruled sheaves is identified in families
with $J^{\rm lin}$ on the charts under consideration. Indeed their
self-Ext spaces in degrees one and two are $H^1(C_b,\OO_{C_b})$ and
$0$ by Proposition~\ref{prop:step21-singular}. Their relative
obstruction theory over $B$ is the smooth one of relative
dimension three. In $J^{\rm lin}$ we have
$T^{\rm lin}=T\setminus D$, where $D$ is the noninvertible
locus of that proposition.

The closed immersion of this smooth
threefold into the smooth sixfold $J^{\rm lin}$ is regular.
With $i:T^{\rm lin}\hookrightarrow J^{\rm lin}$, the pullback
$i^!$ commutes with refined specialization in $s$ by commutativity
of the two regular Gysin operations \cite[Theorem~6.4]{Fulton}.
On a fixed $s$-fiber the virtual dimension before insertion is
$0+\dim B=3$; after insertion it is zero.

On the marked sector the moving relative pair functors match when
their boundary curves coincide. Over $J^{\rm lin}$ the ruled space
is $J^{\rm lin}$ itself, and the curve parameter of a matching is
forced to be $b=\rho(e_Q)$. Consequently the product on that sector
is $M_Q\times_J J^{\rm lin}$. By the assumed compatibility of
obstruction theories, its class is the diagonal pullback of
$[M_Q]^{\vir}\times[J^{\rm lin}]$. The restriction of $d_B$
is the ruled evaluation. Composing with $i^!$ gives
\[
 i^!\Delta_{J^{\rm lin}}^!
       \big([M_Q|_{J^{\rm lin}}]^{\vir}\times[J^{\rm lin}]\big)
  =\Delta_{J^{\rm lin}}^!
       \big([M_Q|_{J^{\rm lin}}]^{\vir}\times[T^{\rm lin}]\big).
\]
This is \eqref{eq:step22-univmatch} on the displayed open.
\end{proof}
The missing boundary in this calculation on $T$ is contained in
$D$, of dimension at most one for a general $H$. Knowing a virtual
zero-cycle on the complement does not determine its degree on a
proper compactification. Localization leaves the possible correction
in the image of
\[
 A_0(M_Q\times_J D)_\QQ\longrightarrow
 A_0(M_Q\times_J T)_\QQ .
\]
A statement about the full degree must also account for other
marked sectors and any other proper factors. Moreover $B$ moves
the pencils: a generating series of the inserted invariants over
$B$ needs its own comparison with K3 period or Hilbert-scheme series.

\begin{proposition}[Avoidance of the noninvertible jets for quartic conics]
\label{prop:step22-conicavoid}
For a general quartic threefold $Q$, general hyperplane
$S=Q\cap\Pi_0$, and general plane section $H=S\cap\Pi_1$ of this
quartic K3, the support of the proper quartic conic boundary
cycle of Corollary~\ref{cor:mark20-conictest} is disjoint from the
noninvertible locus $D\subset T$. Consequently its entire proper
refined intersection with $T$ is supported in $T^{\rm lin}$.
The degree there computed using its Abel cycle is $-26560$.
\end{proposition}
\begin{proof}
Fix a hyperplane $\Pi_0\simeq\PP^3$ in $\PP^4$ and a plane
$\Pi_1\subset\Pi_0$. Let $K$ be a plane conic on $Q$, with plane
$\Pi\not\subset\Pi_0$ and $\ell=\Pi\cap\Pi_0$. On the general
quartic section $S$ there are no lines, so $K$ has no component in
$\Pi_0$ and $\ell$ is not contained in $K$.
If its boundary value equals $[I_{Z/C}]\in D$, the closed immersion
of Lemma~\ref{lem:mark20-abel2} forces $K\cap\ell=Z$, a double
point centered at some $x\in\ell\cap\Pi_1$ with direction
$T_x\ell$.

By Proposition~\ref{prop:step21-singular}, $\mathrm{II}_{S,x}(T_x\ell,T_x\ell)=0$. Write $Q|_\Pi=q_Kr$, with $q_K,r$ quadratic. Since $\operatorname{ord}_x(q_K|_\ell)=2$, the required order of contact $\ge3$ is equivalent to $r(x)=0$, including when $r|_\ell=0$.

Before choosing $Q$, the dimensions of the flag parameters $(\ell,\Pi,x)$ are
\[
 \begin{array}{c|ccc}
      &\ell&\Pi\supset\ell&x\\ \hline
 \ell\not\subset\Pi_1&4&2&0\\
 \ell\subset\Pi_1&2&2&1 .
 \end{array}
\]
Here $x=\ell\cap\Pi_1$ in the first row. Restriction of plane quadrics to $\ell$ is onto; the required double roots span a line in the three-dimensional target. Hence the admissible conics form an open in $\PP^3$, including the allowed reducible and doubled conics. The total parameter dimension is at most $9$. Containment $Q|_\Pi=q_Kr$ imposes $15-6=9$ independent conditions in $|\OO_{\PP^4}(4)|$. The residual quadratic is arbitrary, so $r(x)=0$ imposes one more.
If $\mathcal D$ denotes this incidence in the quartic parameter space, then
\[
 \dim\overline{\operatorname{pr}(\mathcal D)}
 \le 9+(69-10)=68 < \dim|\OO_{\PP^4}(4)|.
\]
The same dimension estimate over the flag variety of $(\Pi_0,\Pi_1)$ proves the assertion for general $(Q,S,H)$.

The coincidence locus is proper because both input cycles have proper
models. Its support misses $D$, so its refined diagonal class, with
its excess multiplicity, is unchanged on passing to $T^{\rm lin}$.
Its degree is $-2\int_{S^{[2]}} D_{\rm conic}\tau=-26560$ by
Corollary~\ref{cor:mark20-conictest}; $D_{\rm conic}$ here denotes
the conic surface pushed to $S^{[2]}$.

It follows that the entire marked conic intersection lies where the boundary ideal is a line bundle. We use this support assertion in Theorem~\ref{thm:rev27-moving} to construct the obstruction theory along the full incidence on the noncentral family.
\end{proof}

\begin{corollary}[A horizontal conic family with the same inserted degree]
\label{cor:step22-horizontalconic}
For $Q,S,H$ as in Proposition~\ref{prop:step22-conicavoid}, let
$N$ parametrize a plane conic $K\subset Q$ and a hyperplane surface
$F\subset Q$ containing $K$. There is a family of stable surface
sheaves on the $s\ne0$ fibers of the moving-pencil family, parametrized
by $N\times\mathbb G_m$. The derived restriction to the exceptional
divisor, pushed down to the base curve, is the rank-one ideal
$I_{K\cap S/C_b}$, where $C_b=F\cap S$. It defines a proper
morphism $e_N:N\to J_{-4}$, independent of $s\ne0$. One has
\[
 (e_N)_*[N]=A_{D_{\rm conic}},\qquad
 \deg\Delta_J^!\big([N]\times[T]\big)=-26560 ,
\]
where the fiber product uses $e_N$ and the embedding $T\hookrightarrow J$.
The ordinary scalar-rigidified vertical sheaf obstruction theory over
$B$ has virtual multiplicity one on the dense open locus where the
quartic pencil member and the conic are smooth and the pencil meets
the conic Noether--Lefschetz divisor transversely.
\end{corollary}
\begin{proof}
Write $\Pi_K$ for the plane of $K$ and $\ell_K=\Pi_K\cap\Pi_0$.
Since $S$ contains no conic, $F\ne S$. Thus $F\cap S=C_b$
determines a plane member of $|h|$ and $F$ belongs to the pencil
generated by $S$ and any lift of $C_b$. The blowup of $F$ along its
Cartier divisor $C_b$ is $F$ itself. Via its strict transform in
$Y_b=\operatorname{Bl}_{C_b}Q$, with inclusion $j_F$, form
the pure stable sheaf $j_{F*}I_{K/F}$.

Indeed every hyperplane surface
of the smooth quartic threefold is integral since $\Pic(Q)=\mathbb Z
H_Q$ with $H_Q$ primitive ample, and the ideal has rank one and is
torsion free. The universal conic and the universal hyperplane
surface are flat over $N$, so the universal ideal is flat. The pencil
degeneration is the constant blowup family over $\mathbb G_m$ for each
$b$, also after $b$ varies.

Every conic cuts $S$ in a finite flat subscheme of length two and
no component is contained there. The defining equation of $S$ on $Q$
acts injectively on $\OO_F$ and $\OO_K$ since both are Cohen--Macaulay
and the intersections are proper. By the ideal sequence
\[
 \mathbf L k_b^* I_{K/F}=I_{K\cap S/C_b}.
\]
Here $k_b$ denotes the inclusion into $F$ of its Cartier
divisor $C_b$. The Cartier square with the exceptional divisor
is Tor-independent, so the Koszul proof of
\eqref{eq:step22-indexspec} identifies the derived index with this
ideal, including on singular $C_b$. This identity may equally be
computed for the given sheaf before establishing perfectness:
the Koszul differential is multiplication by the exceptional section
equation, which acts injectively on the ideal.

Write $I=I_{K/F}$ and let $x=0$ cut out $C_b$. Where $I/xI$ is invertible, lift its generator to obtain $\OO_F\twoheadrightarrow I$. Since $x$ is $I$-regular, the kernel $K'$ satisfies $K'/xK'=0$, hence $K'=0$ by Nakayama. The surface pushforward is therefore perfect near the exceptional divisor and $e_N=d_B$ there. Proposition~\ref{prop:step22-conicavoid} excludes the other jets from the incidence.

A plane in $\Pi_0$ containing $K\cap S$ contains the linear span
$\ell_K$. Conversely, such a plane spans with $\Pi_K$ a unique
hyperplane in $\PP^4$. In families these inverse maps identify
\[
     N=\mathcal F_2(Q)\times_{S^{[2]}}U_2,
     \qquad e_N=a\circ\mathrm{pr}_{U_2}.
\]
Flat base change for the projective-line bundle $U_2\to S^{[2]}$
gives the asserted pushforward. The closed immersion $T\to J$
is regular of codimension three, by Lemma~\ref{lem:mark20-abel2}
and the smoothness of $T$. Over $\mathbb G_m$ the resulting Cartesian
incidence is $(N\times_J T)\times\mathbb G_m$. Its relative
refined cycle is the exterior product of the cycle on $N\times_J T$
with $[\mathbb G_m]$. The degree is the conic pairing in
Corollary~\ref{cor:mark20-conictest}.

For the generic virtual multiplicity, the general quartic cuts the
generated rank-nine restriction bundle on the parameter space of conics
in a smooth surface, by Bertini for vector bundles. Over smooth conics
in smooth $F$, the projection $N\to |H_Q|$ has injective
differential: its kernel is
$H^0(K,N_{K/F})=H^0(\PP^1,\OO(-2))=0$ by adjunction on the K3.

Thus each such component projects generically finitely to a
hypersurface in the four-dimensional hyperplane system. Projection
from the general point $[S]$ to the space $B$ of pencils through it
is generically finite on this hypersurface (and separable in
characteristic zero). Therefore on a dense open of these components
the pencil tangent is transverse to the conic locus.

Put $L=\OO_F(-K)$ at such a point. Riemann--Roch and the divisor
sequence give
$h^0(L^\vee)=1$, $h^1(L^\vee)=h^2(L^\vee)=0$:
the linear system of the smooth negative curve $K$ consists of $K$
alone. In every first-order deformation of $L^\vee$ the section lifts
and its flat zero divisor deforms $K$. Line-bundle deformations
along a pencil direction thus exist exactly when the conic deforms.
For the Cartier fiber $j:F\hookrightarrow Y_b$ with trivial normal
bundle the local-to-global self-Ext sequence reads
\[
 0\to H^1(\OO_F)\to\Ext^1_{Y_b}(G,G)
 \to H^0(\OO_F)\xrightarrow{\partial}H^2(\OO_F)
 \to\Ext^2_{Y_b}(G,G)\to H^1(\OO_F),
\]
where $G=j_*L$. The arrow from $\Ext^1$ records infinitesimal support
motion in the pencil. Indeed simplicity and the scalar action of
the pencil parameter recover the infinitesimal support section as
in Theorem~\ref{thm:outerbitrarymultiplicity}; its sheaf on a smooth
deformation of the fiber is a line bundle by flatness and Nakayama.

Transversality makes this support motion zero. Thus
$\partial:\CC\to\CC$ is an isomorphism and the two Ext spaces vanish.
On this open the relative sheaf chart is \'etale over $B$, $N$ maps
to it \'etale, and the induced class is fundamental. Theorem~\ref{thm:rev27-moving} treats the points of the insertion outside this transverse locus.
\end{proof}

\begin{proposition}[The conic sheaf across the central fiber]
\label{prop:step22-across}
Over the whole parameter scheme $N\times\AA^1_s$ of
Corollary~\ref{cor:step22-horizontalconic}, there exists a flat family
of surface sheaves on the moving degeneration. Its value on $s\ne0$
is $j_*I_{K/F}$ on the strict pencil fiber. Its value over $s=0$ has
quartic and ruled restrictions
\[
 I_{K/F},\qquad b_{C_b}^*I_{Z/C_b},\qquad Z=K\cap S,
\]
on $F$ and $P|_{C_b}$, respectively, glued by the identity on
$I_{Z/C_b}$. Its exceptional-divisor index as a derived object is
$I_{Z/C_b}$ on every fiber, including the noninvertible locus.
Thus the evaluation on this parameter scheme is $e_N$ independently
of $s$, with values in the compactified Jacobian. The refined incidence
with $T$ is a proper family of cycle classes whose specialized
degree on these parameters is $-26560$.
\end{proposition}
\begin{proof}
Let $\mathcal K\subset\mathcal F\subset Q\times N$ and
$\mathcal C_N=\mathcal F\cap(S\times N)$ be the universal conic,
its hyperplane section, and the pencil center. They are flat over $N$.
Take the proper transform of $\mathcal F\times\AA^1_s$ in the
deformation to the normal cone of $S$. If $x$ cuts out $S$, its two
charts have relations $x=sz$ and $s=xw$. They are the blowup charts of
\[
 p:\widetilde{\mathcal F}=
 \operatorname{Bl}_{\mathcal C_N\times\{0\}}
 (\mathcal F\times\AA^1_s)\longrightarrow\mathcal F\times\AA^1_s,
\]
since $x,s$ is a regular sequence on $\mathcal F\times\AA^1_s$.
In the first chart let $f$ be the equation of $\mathcal F$.
The second blowup center has ideal $(z,f)$: its plane equation
on $S$ lifts to $f$, and any two lifts differ by a multiple
of $x=sz$. Its restriction to the strict surface is $(z)$,
an invertible ideal.

Consequently the final strict surface is still
$\widetilde{\mathcal F}$. In the chart $f=za$ it is $a=0$, and
its intersection with the second exceptional divisor $E$ is the
section $a=z=0$ over $\mathcal C_N\times\AA^1_s$.
The proper transform of this surface in the opposite normal chart
is disjoint from the second center.

Write $V=I_{\mathcal K/\mathcal F}$ and let $\widetilde V$ be the
ordinary pullback of $V$ along contraction to $\mathcal F$.
For each of $M=\OO_{\mathcal F},\OO_{\mathcal K},V$, the sheaf $M$
is $N$-flat, and $x$ acts injectively on each $M_n$ along $x=0$.
For $\mathcal K_n$ this follows since its Cohen--Macaulay curve has
no component in $S$, including when the conic is nonreduced.
Pullback to the two blowup charts is given, respectively, by
\[
 M[s,z]/(x-sz)M[s,z],\qquad M[s,w]/(s-xw)M[s,w].
\]
At each geometric $(n,t)$ the two displayed multipliers act
injectively on $M_n[z]$ and $M_n[w]$. At $t=0$ use regularity of
$x$; at $t\ne0$ use the leading coefficient of $x-tz$ and the
constant coefficient of $t-xw$, respectively. The fiberwise
criterion for flatness applied to the two-term presentations
gives $N\times\AA^1_s$-flat cokernels and commutation with
base change. Pushing $\widetilde V$ forward along the closed
immersion of the surface preserves this flatness.

At $s=0$ the blowup has the two components $\mathcal F$ and the
projective-line bundle over $\mathcal C_N$. Their restrictions of
$\widetilde V$ are $V$ and the pullback of $V/xV$ locally.
Regularity of $x$ on $\OO_{\mathcal K}$ gives
\[
   0\longrightarrow V/xV\longrightarrow \OO_{\mathcal C_N}
    \longrightarrow \OO_{\mathcal K\cap\mathcal C_N}
    \longrightarrow0.
\]
In the chart at their join the module is $M[w]/xwM[w]$,
with the two restrictions $M$ and $(M/xM)[w]$. Its map to their
direct sum identifies it with their fiber product over $M/xM$
by taking the constant coefficient in $w$. This gives the asserted
gluing on the central fiber.

Finally $z$ acts regularly on
$V[s,z]/(x-sz)V[s,z]$. Indeed reducing
$zv=(x-sz)v'$ modulo $z$ forces $v'\equiv0\pmod z$ by regularity
of $x$; cancellation of $z$ in the polynomial module finishes
the assertion. Derived intersection with $E$ is therefore
the ordinary restriction to its section, with value
$(V/xV)[s]$. Pushing along $E\to\mathcal C_N\times\AA^1_s$
returns this same ideal in degree zero. The classifying map to $J$ is
$e_N\circ\mathrm{pr}_N$. The incidence fiber product and its refined
class are the exterior products
\[
 (N\times_J T)\times\AA^1_s,\qquad
 \big(\Delta_J^!([N]\times[T])\big)\times[\AA^1_s].
\]
The first is proper over $\AA^1_s$, since $N,T$ are proper. Applying
the fiber pullback and Corollary~\ref{cor:step22-horizontalconic}
proves the degree statement. The calculation gives Chow specialization in the proper parameter space of the conic family. Theorem~\ref{thm:rev27-moving} below identifies its general-fiber input with the inserted sheaf virtual cycle.
\end{proof}

\begin{theorem}[Classification of the vertical conic charge]
\label{thm:step22-conicclassification}
Let $F$ be any integral quartic surface in $\PP^3$ over $\CC$.
The rank-one torsion-free sheaves on $F$ of Hilbert polynomial
$2m^2-2m+1$ are exactly the ideals $I_{K/F}$ of plane conics
$K\subset F$, including reducible and double conics.
The conic associated to the isomorphism class of the sheaf is unique.

This identifies families modulo tensoring by a line bundle from
the base. In particular, on the moving family of
pencil threefolds at $s\ne0$ (including singular base curves $C_b$), the scalar-rigidified moduli functor
of stable sheaves with vertical fundamental cycle one pencil fiber
and this polynomial for the hyperplane class on that fiber is $N$.
\end{theorem}
\begin{proof}
If $H^0(F,V)\ne0$, a nonzero section of the torsion-free rank-one
sheaf $V$ gives an injection $\OO_F\hookrightarrow V$. The Hilbert
polynomial of its cokernel would be
\[
        (2m^2-2m+1)-(2m^2+2)=-2m-1,
\]
which is negative for $m\gg0$. Thus $h^0(V)=0$, and
$h^2(V)=1+h^1(V)\ge1$. A quartic hypersurface is Gorenstein
with $\omega_F\simeq\OO_F$. By Serre duality
$\Hom_F(V,\OO_F)=H^2(F,V)^\vee\ne0$. A nonzero map to
$\OO_F$ is generically nonzero and has torsion kernel, hence is
injective. Its image is an ideal $I_{Z/F}$ with
$P_{\OO_Z}(m)=2m+1$.

To classify $Z$ as a subscheme of $\PP^3$, remove the maximal
zero-dimensional subsheaf of $\OO_Z$. This subsheaf is an ideal,
giving a pure degree-two curve $C$ and
\[
        1=\chi(\OO_Z)=\chi(\OO_C)+\operatorname{length}
                            \ker(\OO_Z\to\OO_C).
\]
An integral degree-two curve is a plane conic by the degree bound
for nondegenerate curves. A pure generically reduced curve is
reduced; for the union of two lines the normalization sequence
gives $\chi=1$ when the lines meet and $\chi=2$ when they are skew.
The remaining case has support a line $D$. Generically
$I_D^2\OO_C=0$ since the local algebra has length two.
Purity then gives $I_D^2\OO_C=0$ everywhere. Its nilradical
$I_D/I_C$ is torsion-free of rank one on the smooth line $D$,
so for some $a\in\mathbb Z$ we have
\[
 \OO_D(-1)^{\oplus2}=I_D/I_D^2\twoheadrightarrow
 I_D/I_C\simeq\OO_D(a),\qquad
 \chi(\OO_C)=1+(a+1).
\]
Consequently $a\ge-1$. Equality supplies a nonzero constant
relation between the two linear equations generating $I_D$:
a plane equation belongs to $I_C$. This case is the doubled line
in that plane. In every case $\chi(\OO_C)\ge1$ and equality requires
a plane conic. Therefore $Z=C=K$ is such a conic.

Conversely $I_{K/F}$ has the required polynomial and is rank-one torsion-free. Each proper subsheaf of equal dimension has a nonzero lower-dimensional quotient, giving stability for every ample polarization. The ideal sequence gives
\[
 H^0(\OO_F)\xrightarrow{\sim}H^0(\OO_K)=\CC,\qquad
 H^1(I_{K/F})=0,\qquad h^2(I_{K/F})=1 ,
\]
using $H^1(\OO_F)=H^1(\OO_K)=0$ and $h^2(\OO_F)=1$. Thus the inclusion in $\OO_F$ recovers $K$ uniquely up to scalar by duality.

For a flat family $p_T:F_T\to T$ of quartics and a $T$-flat family $V_T$, trivialize locally the base line factor in $\omega_{F_T/T}$. The proper pushforward is perfect, as is seen after pushing to the smooth ambient projective space over $T$. Its only fiber cohomology has rank one in degree two. Hence $Rp_{T*}V_T=L[-2]$ for the line bundle $L=R^2p_{T*}V_T$, and relative duality gives
\[
  Rp_{T*}\mathbf R\mathcal Hom(V_T,\omega_{F_T/T})[2]
        =\mathbf R\mathcal Hom(Rp_{T*}V_T,\OO_T).
\]
Thus $p_{T*}\mathcal Hom(V_T,\OO_{F_T})$ is invertible and commutes with base change. A local generator defines
\[
 0\to V_T\to\OO_{F_T}\to\OO_{K_T}\to0 .
\]
Fiberwise injectivity gives exactness and flatness of $K_T/T$. The image is unchanged on rescaling the generator; its ideal recovers the family up to a line factor from $T$.

On a smooth quartic threefold each hyperplane section is integral:
the Picard group is generated by the ample hyperplane class and
the positive classes of two effective divisors cannot sum to that
generator. A vertical stable sheaf of fundamental cycle one fiber
comes from a rank-one torsion-free sheaf on that fiber.

Here one can recover the supporting pencil member even over a nonreduced test scheme $T$ and a singular center $C_b$. Locally on the base of pencils, embed the incidence threefold by $j:Y_T\hookrightarrow Q\times\PP^1\times T$. A flat family $G$ of stable sheaves has perfect pushforward $j_*G$ on this smooth relative ambient space, by lifting fiberwise finite locally free resolutions. For the projection $p$ to $T$, put
\[
 V=\RR p_*\RR\mathcal Hom(j_*G,j_*G),\qquad A=\cofib(\OO_T\xrightarrow{1}V).
\]
Perfect proper pushforward gives base change for $V$. Since a closed immersion preserves ordinary Hom, simplicity on geometric fibers implies $H^{\le0}(A\otimes^{\LL}k(t))=0$. A minimal free representative over the local ring of $t$ starts in degree one. In particular the only endomorphisms of $G$ after any base change are scalar locally on the base.

Locally on $T$, properness places the image of $\Supp G$ in an affine pencil chart with coordinate $z$. The scalar endomorphism calculation gives
\[
 z\,\mathrm{id}_G=b'\,\mathrm{id}_G,\qquad b'\in\OO_T ,
\]
so the graph of $b'$ is the unique supporting member, also over nonreduced $T$. The preceding construction recovers its conic and is inverse to taking the relative ideal. Since the general $S$ contains no conic, the parameter scheme is $N$.
\end{proof}

\begin{corollary}[The fixed-pencil virtual degree in the conic charge]
\label{cor:step22-conicpencilnumber}
For general $Q,S$ and a general pencil center $b$, the stable vertical
moduli in Theorem~\ref{thm:step22-conicclassification} is a reduced
finite scheme. Its threefold deformation and obstruction spaces
in degrees one and two vanish. Its virtual degree is $5016$.
\end{corollary}
\begin{proof}
Consider the incidence of triples $(K,F_{\rm lin},Q')$, where $K$
is a plane conic, $F_{\rm lin}$ is a hyperplane containing it,
and $Q'$ is a quartic hypersurface containing $K$. The first
two entries form a projective-line bundle over the projective
bundle of all plane conics. Restriction of quartics to every conic
surjects onto the nine-dimensional $H^0(K,\OO_K(4))$, by the
plane conic resolution.

The incidence is thus a projective bundle
over this parameter space, irreducible of dimension $72$.
The union where $K$ is singular or $F_{\rm lin}\cap Q'$ is singular
is closed proper in this incidence. Indeed for a fixed smooth conic
on a hyperplane $\PP^3$ a quartic through it has equation
$uA_3+qB_2$ for $I_K=(u,q)$, and is smooth along $K$ whenever
$A_3,B_2$ have no common zero there. Away from $K$ apply Bertini.

For a general $Q$ the bad locus in its incidence fiber has dimension
at most $71-69=2$. Since $N$ is smooth of pure dimension three by
the generated rank-nine bundle description of
Corollary~\ref{cor:step22-horizontalconic}, this proves density of
the smooth conic and smooth quartic-surface locus on each component.
Projection from the general $[S]$ to $B$ is generically finite
and separable on the hyperplane image of each component as proved
in that corollary.

The bad loci and branch loci have images of
dimension at most two. For general $b$ the fiber $N_b$ is finite
\'etale and all its points have the vanishing of the two fixed-target
self-Ext groups proved there. A perfect obstruction theory with
zero tangent and obstruction spaces gives local virtual degree one.

Write $\rho:J_{-4}\to B$ for the support map and
$f=\rho^*[b]$. The identity $\rho e_N=\mathrm{pr}_B$ and the
projection formula give
\[
 \deg[N_b]^{\vir}=\deg[N_b]
    =\int_N e_N^* f
    =\int_{J_{-4}}(e_N)_*[N]\, f
    =\int_{J_{-4}} A_{D_{\rm conic}} f
    =5016,
\]
by Corollary~\ref{cor:mark20-conictest}. This is the fixed-pencil
vertical invariant in the stated charge; the quintic formula also
requires the other relative contributions.
\end{proof}

\begin{theorem}[The virtual intersection for the moving conic family]
\label{thm:rev27-moving}
Choose $Q,S,H$ as in Proposition~\ref{prop:step22-conicavoid}. Let $N$ be the proper conic moduli of Theorem~\ref{thm:step22-conicclassification} over $B=|h|$ at $s\ne0$, and put $J=J_{-4}$. Define
\[
 U=e_N^{-1}(J^{\rm lin}),\qquad
 Z=N\times_J T=U\times_{J^{\rm lin}} T^{\rm lin}.
\]
For the stable admissible outer-relative sheaves of Section~\ref{sec:outer-fiber-multiple} this charge uses the unexpanded target throughout $s\ne0$. On the open subscheme $U$ of the classical conic moduli, equipped with its universal ideal, the sheaf obstruction theory relative to $B$ identifies with $\mathbb L_{U/B}$. Write $[U]^{\vir}$ for the virtual pullback of $[B]$ to $U$, of dimension three. The complete incidence in this charge has a proper virtual zero-cycle. In $A_0(Z)_\QQ$ its class is
\[
 \Delta_{J^{\rm lin}}^!\big([U]^{\vir}\times[T^{\rm lin}]\big)
       =\Delta_J^!\big([N]\times[T]\big).
\]
Its degree is $-26\,560$. The family of Proposition~\ref{prop:step22-across} extends this class over the degeneration parameter as the exterior product of the displayed zero-cycle with $[\AA^1_s]$ on $Z\times\AA^1_s$.
\end{theorem}
\begin{proof}
\emph{Perfectness at every point of the insertion.}
Fix a pencil center $C_b$. Write $Y_b=\operatorname{Bl}_{C_b}Q$, let $D_b\subset Y_b$ be the member corresponding to $S$, and let $G=j_{F*}I_{K/F}$. The member $D_b$ is disjoint from $F$, since $F\ne S$ and the blowup separates the pencil members. The entire locus $Z$ lies in $U$ by Proposition~\ref{prop:step22-conicavoid}. This equality of fiber products holds scheme-theoretically, since restriction to an open containing every point of a scheme gives the whole scheme.

Away from the exceptional divisor the threefold is smooth. Along it, take the local equation $x$ of $C_b$ on $F$. The conic has no component in $S$, so $x$ is regular on $\OO_K$ and $I_{K/F}$. Hence
\[
 I_{K/F}/xI_{K/F}=I_{K\cap S/C_b}.
\]
On $U$ this is a line. The generator-lifting argument above shows $I_{K/F}$ is invertible along $C_b$; as $F$ is Cartier in $Y_b$, its pushforward is perfect there, also for singular $F,C_b$. Flatness of the universal family lifts the finite fiber resolutions.

On an outer expansion every root support is vertical by intersection with a fiber and an ample divisor. At the first join its equation acts both nilpotently and injectively on the local root sheaf, so the restriction is zero. By Theorem~\ref{thm:outerbitrarymultiplicity} an occupied bubble would have one interior slice $S\times\{z\}$. Stability makes its nilpotent transverse coordinate a scalar, hence zero. The resulting rank-one sheaf $V$ would satisfy
\[
 h\cdot c_1(V)=-2,\qquad c_1(V)\in\mathbb Zh,\qquad h^2=4,
\]
which is impossible. An empty bubble has a dilation stabilizer. All geometric points, and hence all families, belong to the open of unexpanded targets.

\emph{The obstruction morphism.}
Locally on $B$ the pencil is the hypersurface of bidegree $(H_Q,1)$ in $Q\times\PP^1$. Thus the threefolds are flat and Gorenstein over $B$, with $\omega_{Y_b}=\OO_{Y_b}(-D_b)$ by adjunction. The divisor sequence of the incidence hypersurface has first term $\OO_Q(-H_Q)\boxtimes\OO_{\PP^1}(-1)$. Its cohomology vanishes. As $H^{>0}(Q,\OO_Q)=0$, duality gives
\[
 R\Gamma(\OO_{Y_b})=\CC,\qquad
 H^i(Y_b,\OO_{Y_b}(-D_b))=
 \begin{cases}\CC&i=3,\\ 0&i\ne3.\end{cases}
\]
These dimensions are constant over $B$ and commute with base change. The relative anticanonical identification may involve a line bundle from the base. The construction below uses the divisor and its canonical trivialization on the sheaf determinant and hence glues without choosing any such identification.

For the universal family $G$ over $U$, let $p$ project to $U$ and write $K=\RR p_*\RR\mathcal Hom(G,G)$. In the derived moduli of perfect complexes rigidify the scalar unit. Take the derived pullback of the determinant morphism for line bundles trivialized along the outer divisor from its classical relative Picard space, as in Theorem~\companionref{thm:zeroboundaryPOT}.

This operation preserves the classical sheaf moduli on $U$. The determinant is canonically trivialized on the divisor because the restricted complex is zero. The classical relative Picard space is \'etale locally discrete over the pencil base: the tangent and obstruction groups $H^1(\OO(-D_b))$ and $H^2(\OO(-D_b))$ vanish, as does its infinitesimal automorphism group.

The sheaf, purity and stability conditions are open on the flat perfect locus. The fiber class used here is specified locally by Hilbert polynomials: the mixed intersection of the divisorial cycle with the hyperplane and a pencil fiber is zero, and its hyperplane-square degree is four.

A surface dominating the pencil has strictly positive mixed intersection. Since every fiber is integral of degree four, these conditions force a single fiber with multiplicity one. The remaining coefficients of the Hilbert polynomial select the conic charge. They are locally constant in flat families.

The differential of this derived pullback uses the unit and the trace, the latter being the differential of the determinant \cite[Proposition~3.2]{STV}. Its tangent complex over $B$ is
\[
 \fib\!\left(\cofib(\OO_U\xrightarrow{1}K)
        \longrightarrow \RR p_*\OO(-D)\right)[1].
\]
Here $D$ is the relative outer fiber pulled back to the universal family of targets. The arrow is the twisted trace, using $G(-D)\simeq G$. Scalar automorphisms act trivially on the framed determinant since $G$ has rank zero. At a geometric point the top trace in degree three is dual, up to a one-dimensional factor, to the unit $\CC\to\Hom(G,G)$ and is an isomorphism. Duality for a perfect sheaf on the Gorenstein threefold also identifies the middle Ext spaces as dual up to a line.

Thus the tangent complex has only $\Ext^1(G,G)$ and $\Ext^2(G,G)$ in degrees $0,1$, of equal dimensions. Perfect proper pushforward and base change give the same Tor-amplitude in families. The truncation morphism of this derived moduli problem gives the relative perfect obstruction theory $E^\bullet\to\mathbb L_{U/B}$ of rank zero. This is the construction in Section~\companionref{sec:zerorelative}, applied to the Gorenstein incidence threefold with the cohomology just calculated.

\emph{Comparison of complexes.}
The space of conics on a general $Q$ is a smooth surface, cut out by a regular section of the generated rank-nine bundle of quartic restrictions. The choice of hyperplane containing the conic plane is a projective line. Consequently $U$ is open in a smooth threefold $N$, and $B$ is smooth of dimension three. The cotangent complex
\[
 \mathbb L_{U/B}=[\Omega_B|_U\longrightarrow\Omega_U]
\]
has rank zero and amplitude $[-1,0]$. The cone of the obstruction morphism has amplitude $[-2,0]$. Its cohomology in degrees $-1,0$ vanishes by the defining surjectivity and isomorphism of an obstruction theory. Locally split the two terminal surjections in a three-term vector-bundle representative of the cone. The remaining bundle in degree $-2$ has rank zero. Hence $E^\bullet\simeq\mathbb L_{U/B}$.

Factor $U\to B$ as $\Gamma_{U/B}:U\hookrightarrow U\times B$ followed by projection. The graph is regular, and $E^\bullet\simeq\mathbb L_{U/B}$ gives
\[
 [U/B]^{\vir}([B])=\Gamma_{U/B}^{!}[U\times B]=[U].
\]
The fiber product $Z$ is proper and contained over $J^{\rm lin}$. Restriction to that open commutes with the refined diagonal pullback \cite[\S6.2]{Fulton}, so the inserted degree is the one in Corollary~\ref{cor:step22-horizontalconic}. Over $\AA^1_s$ the evaluation extends as $e_N\circ\mathrm{pr}_N$ by Proposition~\ref{prop:step22-across}. Flat base change takes the intersection class on $Z$ to its product with $[\AA^1_s]$; restriction to $s=0$ gives the asserted specialization.

\end{proof}

This comparison gives the inserted virtual degree for the moving noncentral pencils. Its specialization in the explicit conic family is a class on $Z$. To identify the contributions of all marked expansions in a relative sheaf compactification one needs a proper family of those expansions and compatibility of its obstruction theory with the sheaf gluing morphism in Theorem~\ref{thm:step22-universal}.

\subsection{The two curve conditions at arbitrary contact length}
Let $D_k\in|kh|$ be smooth. The two complementary geometries prescribe
$D_1=C_*$ and $D_5=\Gamma$. For $n\ge0$ set
\[
 H_n(k)=[D_k^{[n]}]\in H^{2n}(S^{[n]},\QQ),\qquad
 \mathfrak q_\lambda(h)=
 \prod_{m\ge1}\mathfrak q_m(h)^{a_m}|0\rangle ,
 \quad \lambda=(1^{a_1}2^{a_2}\cdots)\vdash n .
\]
Here $|0\rangle=1\in H^0(S^{[0]})$. In the convention of
\cite[\S2 and Theorem~3.17]{LehnMark20}, $\mathfrak q_m(h)$ is
the correspondence of nested subschemes $Z\subset Z'$ with
$\length(I_Z/I_{Z'})=m$ and with this quotient supported at one
point, integrated against $h$ at that point.
Put $z_\lambda=\prod_m m^{a_m}a_m!$ and $\ell(\lambda)=\sum_m a_m$.

\begin{proposition}\label{prop:mark20-allcontacts}
For all $n$ the following cohomological comparison holds:
\begin{equation}\label{eq:mark20-allcontacts}
 H_n(k)=\sum_{\lambda\vdash n}
 \frac{(-1)^{n-\ell(\lambda)}k^{\ell(\lambda)}}{z_\lambda}
 \mathfrak q_\lambda(h).
\end{equation}
In particular the linear operator on the span of the
$\mathfrak q_\lambda(h)$ defined by
$\mathfrak q_\lambda(h)\mapsto5^{\ell(\lambda)}\mathfrak q_\lambda(h)$
takes $H_n(1)$ to $H_n(5)$. For $n>0$ every pairing of $H_n(5)$ with
a test cycle is consequently a linear combination of its pairings
with the $n$ classes
\[
  \sum_{\substack{\lambda\vdash n\\\ell(\lambda)=j}}
   \frac{(-1)^{n-j}}{z_\lambda}\mathfrak q_\lambda(h),
  \qquad 1\le j\le n.
\]
\end{proposition}
\begin{proof}
The section cutting out $D_k$ induces a section of $\OO_S(kh)^{[n]}$
with zero scheme $\operatorname{Sym}^nD_k$. This zero scheme is
smooth of codimension $n$ in $S^{[n]}$. Thus
$H_n(k)=c_n(\OO_S(kh)^{[n]})$. By
\cite[Theorem~3.17, Theorem~4.6 and Remark~4.7]{LehnMark20},
\[
 \sum_{n\ge0}H_n(k)t^n =
 \exp\!\left(\sum_{m\ge1}
  \frac{(-1)^{m-1} kt^m}{m}\mathfrak q_m(h)\right)|0\rangle.
\]
The positive creation operators commute. The monomial indexed by
$(a_m)$ has coefficient
$\prod_m (-1)^{(m-1)a_m}k^{a_m}/(m^{a_m}a_m!)$.
The Nakajima basis theorem gives linear independence of the listed
monomials (extend $h$ to a basis of $H^*(S,\QQ)$).
Collecting terms gives \eqref{eq:mark20-allcontacts} and the stated
operator. Linearity of integration gives the pairing assertion.
\end{proof}
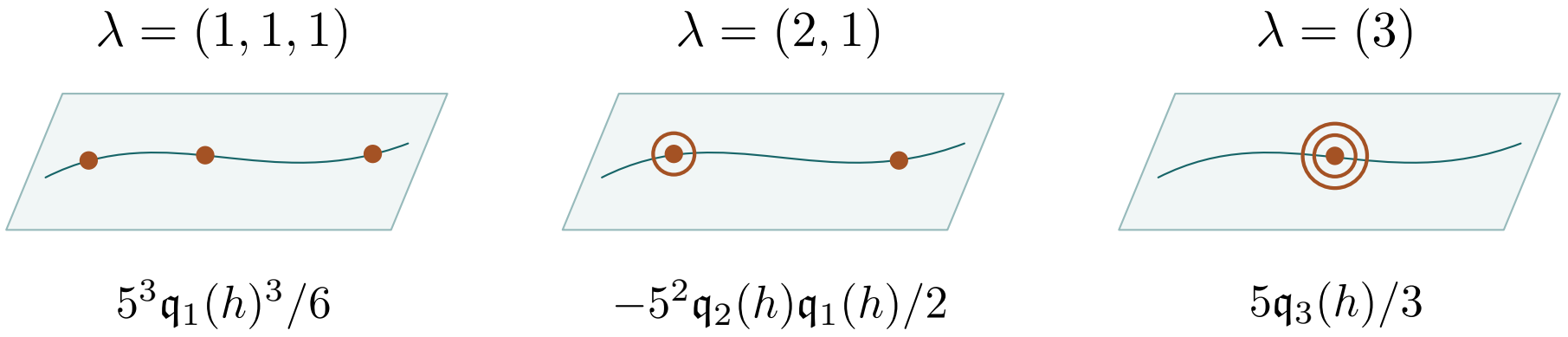
\begin{figure}[htbp]\centering
\begin{tikzpicture}[font=\footnotesize,line cap=round,line join=round,x=1cm,y=1cm]
\foreach \x/\lab in {0/{(1,1,1)},4.25/{(2,1)},8.50/{(3)}}{
 \begin{scope}[shift={(\x,0)}]
 \node[font=\small] at (1.74,2.05) {$\lambda=\lab$};
 \path[qIIteal!45,draw,fill=qIIteal!6] (.08,.55)--(.51,1.59)--(3.45,1.59)--(3.02,.55)--cycle;
 \draw[qIIteal,thin] (.38,.95)..controls (1.33,1.43) and (2.07,.80)..(3.15,1.21);
 \end{scope}}
\foreach \x/\y in {.71/1.08,1.60/1.12,2.88/1.13,5.18/1.13,6.90/1.08,10.23/1.115}{\fill[qIIclay] (\x,\y) circle(2pt);}
\draw[qIIclay,line width=.8pt] (5.18,1.13) circle(4.5pt);
\draw[qIIclay,line width=.8pt] (10.23,1.115) circle(4.5pt) (10.23,1.115) circle(7pt);
\node at (1.74,-.02) {$5^3\mathfrak q_1(h)^3/6$};
\node at (5.99,-.02) {$-5^2\mathfrak q_2(h)\mathfrak q_1(h)/2$};
\node at (10.24,-.02) {$5\mathfrak q_3(h)/3$};
\end{tikzpicture}
\caption{The three creation-operator monomials for a contact of length three in \eqref{eq:mark20-allcontacts}. Each ring denotes one additional unit of length at the same support point. The cycle of the curve in $|5h|$ supplies a factor of five for each part of the partition. The coefficients below include the signs and symmetry factors.}\label{fig:II-contact-partitions}
\end{figure}

In length three, with $Q_m=\mathfrak q_m(h)$, this reads
\[
 H_3(5)=\big(125Q_1^3/6-25Q_1Q_2/2+5Q_3/3\big)|0\rangle.
\]
Here $t$ records the length of a boundary subscheme on the fixed K3.
The DT variable $q$ records the total normalized threefold charge.
To substitute \eqref{eq:mark20-allcontacts} into the degeneration
formula one must identify the inserted relative virtual classes with
the indicated correspondences and their excess classes, including
the images of coincident contacts in the boundary moduli space.
\medskip

The cohomological comparison in \eqref{eq:mark20-allcontacts} gives a
numerical test for a proposed Chern class insertion. Put $L=\OO_S(h)$
and $\Xi(q)=\prod_{m\ge1}(1-q^m)$. In the following fixed-surface series the exponent of $q$ is the Hilbert-scheme length. For $k\ge0$ define
\[
 A_k(q)=\Xi(q)^{24}\sum_{n\ge2}q^n
 \int_{S^{[n]}}c_2((L^{\otimes k})^{[n]})\,
                         c_{2n-2}(T_{S^{[n]}}),\qquad
 B(q)=A_5(q)-25A_1(q).
\]
In length two, $c_2((L^{\otimes k})^{[2]})=H_2(k)$ for $k\ge1$.
In higher length the insertion has fixed codimension two.

\begin{proposition}\label{prop:mark20-k3bracket}
Write $\sigma_j(n)=\sum_{d\mid n}d^j$. For every $n\ge1$ and $k\ge0$,
\begin{align}
 [q^n]A_k
  &= -5(n-2)\bigl(\sigma_3(n)-n\sigma_1(n)\bigr)
          +2k^2(n-1)\sigma_1(n),\label{eq:mark20-an}\\
 {[q^n]B}
  &=120(n-2)\bigl(\sigma_3(n)-n\sigma_1(n)\bigr).
       \label{eq:mark20-bn}
\end{align}
In particular, $[q^2]B=0$ and $[q^n]B>0$ for every $n\ge3$.
Both series are quasimodular of mixed weight at most six.
\end{proposition}
\begin{proof}
For a smooth projective surface $Y$ and a line bundle $M$ on $Y$, set
\[
 V(Y,M;u)=\sum_{n\ge0}q^n\int_{Y^{[n]}}c_u(M^{[n]})c(T_{Y^{[n]}}),
 \qquad c_u(V)=\sum_{j\ge0}u^jc_j(V).
\]
By \cite[Theorem~4.2]{EGLMark20},
$\log(V(Y,M;u)/V(Y,M;0))$ is linear in the four numbers
$c_1(Y)^2,c_2(Y),c_1(Y)c_1(M),c_1(M)^2$. For $Y=S$, the symplectic
form identifies $T_{S^{[n]}}$ with its dual. Hence
$c_{2n-1}(T_{S^{[n]}})=0$ rationally, so $[u]V(S,M;u)=0$.
G\"ottsche's Euler formula gives $V(S,M;0)=\Xi^{-24}$. It follows that
\[
 A_k=[u^2]\log\big(V(S,L^k;u)/V(S,L^k;0)\big).
\]
Write $p_d$ for $[u^2]$ of the same logarithm for
$(\PP^2,\OO(d))$, and $s$ for its value on $(\PP^1\times\PP^1,\OO)$.
The four characteristic numbers of these pairs are $(9,3,3d,d^2)$
and $(8,4,0,0)$, respectively. Since $(S,L^k)$ has vector
$(0,24,0,4k^2)$, linearity gives
\begin{equation}\label{eq:mark20-egl}
 A_k=A_0+k^2 U,\qquad
 A_0=-16p_0+18s,\qquad U=2(p_2-2p_1+p_0).
\end{equation}

Here are the finite localization sums determining $A_0,U$.
At a torus fixed point of $Y$ write $a,b$ for its tangent weights
and $f$ for the fiber weight of $M$. For a partition $\lambda$ and
a box $(i,j)$, indexed from zero, put $r=\lambda_i-j-1$,
$t=\lambda'_j-i-1$. The tautological weight of the monomial
$x^jy^i$ is $v=f-ja-ib$. The two tangent weights are
$w_1=(r+1)a-tb$, $w_2=(t+1)b-ra$
by \cite[Theorem~2.11, rank-one case]{NYMark20}.
Form the power series
\[
 F_{a,b,f}=\sum_{\lambda}q^{|\lambda|}
 \prod_{(i,j)\in\lambda}
 (1+uv)(1+z/w_1)(1+z/w_2),\qquad
 W_{ij}=[u^iz^j]\prod_{y\in Y^T}F_{a_y,b_y,f_y}.
\]
Dividing by the tangent Euler class in the localization formula shows
that $W_{jj}=[u^j]V(Y,M;u)$. Thus for each of the four toric pairs
the required logarithmic coefficient is
$W_{22}/W_{00}-\tfrac12(W_{11}/W_{00})^2$.
For $\PP^2$ choose coordinate weights $(e_0,e_1,e_2)=(0,17,239)$.
At the $i$th vertex the triple is $(e_j-e_i,e_l-e_i,-de_i)$ with
$\{i,j,l\}=\{0,1,2\}$. For $\PP^1\times\PP^1$ take
$(\pm7,\pm191,0)$ independently. In the above partition sums retain
$|\lambda|\le6$, discard $u^3,z^3,q^7$, and apply
\eqref{eq:mark20-egl}. The resulting rows are
\[
 \begin{array}{c|rrrrrrr}
 n&0&1&2&3&4&5&6\\\hline
 {[q^n]A_0}&0&0&0&-80&-450&-1440&-3600\\
 {[q^n]U}  &0&0&6&16&42&48&120 .
 \end{array}
\]

For any vector bundle $E$ one has
$c_2(E)=\tfrac12\operatorname{ch}_1(E)^2-\operatorname{ch}_2(E)$.
Apply \cite[Theorem~1.2, equations~(1.1)--(1.4)]{QinMark20} on the
K3 surface to the two summands. Their respective weight bounds are
$(1+2)+(1+2)=6$ and $2+2=4$. Thus $A_0$ and $U=A_1-A_0$ lie in the
span of $\mathcal E=(1,E_2,E_2^2,E_4,E_2^3,E_2E_4,E_6)$, with the
Eisenstein series normalized to constant term one. The Fourier matrix
$([q^n]\mathcal E_j)_{0\le n,j\le6}$ has determinant
$2^{40}3^{16}5^3\cdot7$, using
\[
  (E_2,E_4,E_6)=(1,1,1)+
    \sum_{n\ge1}q^n(-24\sigma_1(n),240\sigma_3(n),-504\sigma_5(n)).
\]
Put $D=q\,d/dq$. The quasimodular expressions
\[
 -\tfrac1{48}D E_4-\tfrac{5}{24}D^2 E_2+
 \tfrac{1}{24}(E_4-1)+\tfrac5{12}DE_2,
 \qquad \tfrac{1}{12}(E_2-1-DE_2)
\]
belong to the same span by Ramanujan's differential identities. Their
first seven coefficients are the two displayed rows, hence they equal
$A_0$ and $U$, respectively. Taking coefficients gives
\eqref{eq:mark20-an}. Substitution of $k=5,1$ gives
\eqref{eq:mark20-bn}. Pairing the divisors $d$ and $n/d$ gives
\[
 2(\sigma_3(n)-n\sigma_1(n))
  =\sum_{d\mid n}(d+n/d)(d-n/d)^2>0\qquad(n>1).
\]
\end{proof}
Proposition~\ref{prop:mark20-k3bracket} evaluates an ordinary tangent
Chern number on the fixed middle K3 surface. Its use in a coefficient
of the quintic partition function requires identifying the relative
virtual insertion and its multiplicity on the corresponding contact
stratum.

The passage from a contact subscheme to its boundary sheaf changes at
length three. On a smooth plane quartic $C$ the fiber of
$b:C^{[3]}\to\Pic^{-3}(C)$, $Z\mapsto\OO_C(-Z)$, over
$\omega_C^{-1}(p)$ is $|\omega_C(-p)|$. Since the canonical
series embeds $C$ as a plane quartic, evaluation at $p$ is onto and
$h^0(\omega_C(-p))=3-1=2$. Hence
\[
      b^{-1}(\omega_C^{-1}(p))\simeq\PP^1.
\]
Already over smooth curves the degree-three Abel map contracts
families. Transferring the Hilbert-scheme formula to the Jacobian in
this length requires the corresponding pushforward along the Abel map.

\section{One point on a surface carrying a line defect}
At a singular point of the hyperplane surface lying on the line, the length-one quotients of its line hull form a $\PP^1$. We use the obstruction theory of the kernels.
\subsection{A point added to a line defect}\label{sec:line-point-verification}
A length-one quotient changes the third Chern character of its kernel by $-[\mathrm{pt}]$. Write $h=H|_T$ and let $e(G)$ denote the degree in the linear coefficient of the Hilbert polynomial on the primitive hyperplane section $T$.

\begin{proposition}\label{prop:line-point-hull-d}
Let $T$ be a hyperplane section of a smooth quintic $X$. Every rank-one torsion-free sheaf $G$ on $T$ with $e(G)=-1$ and $\chi(G)=3$ fits into an exact sequence
\[
 0\longrightarrow G\longrightarrow I_{\ell/T}\xrightarrow{q} k(p)\longrightarrow0.
\]
The line $\ell\subset T$, the point $p\in T$, and the quotient up to scalar are unique. Every such kernel gives a stable pure surface sheaf on $X$.
\end{proposition}
\begin{proof}
The finite Gauss map shows that $T$ has isolated singularities, and the integral hypersurface $T$ is normal. Put $R=G^{\vee\vee}$ and $s=\length(R/G)$. Then $e(R)=-1$ and $\chi(R)=3+s$. Restricting $R^\vee$ to a general smooth plane quintic section, its degree there is one. Since $H^0(R^\vee(-h))=0$, this gives $h^0(R^\vee)\le1$.

Applying the same restriction to $R^\vee(h)$ gives $h^0(R^\vee(h))\le4$: a line bundle of degree six on a smooth plane quintic has at most three sections by Clifford's inequality and its equality case. Since $h^0(R)=0$, duality gives $3+s=\chi(R)\le4$. If $s=1$, equality forces a section of $R^\vee$, whose divisor is a line $\ell$, so $R=I_{\ell/T}$. Its ideal is reflexive by the depth lemma applied to $\OO_T\twoheadrightarrow\OO_\ell$.

Suppose $s=0$. The rank-one Hilbert polynomial gives $\chi(R(h))=\chi(R)+e(R)=2$, and $h^2(R(h))=h^0(R^\vee)\le1$. A nonzero section of $R(h)$ dualizes to an injection $R^\vee(-h)\hookrightarrow\OO_T$. Its image cuts out a locally Cohen--Macaulay curve $A$ by the depth lemma. Its degree is four. The identities $\chi(R^\vee)=\chi(R)+e(R)=2$ and $e(R^\vee)=1$ give $\chi(R^\vee(-h))=6$. Hence
\[
 \chi(\OO_A)=5-6=-1,\qquad p_a(A)=2.
\]
A locally Cohen--Macaulay plane quartic has arithmetic genus three, while a nonplanar locally Cohen--Macaulay quartic in $\PP^3$ has arithmetic genus at most one by \cite[Theorem~2.2]{HLSspace}. This excludes $s=0$. Uniqueness follows because $G$ determines its hull and $\Hom_T(I_{\ell/T},\OO_T)$ has dimension one. Conversely the kernel of any indicated quotient is torsion-free of rank one on the integral surface $T$. Purity and Gieseker stability follow.
\end{proof}

\begin{proposition}\label{prop:linepoint-blowup-d}
Fix a line $\ell$ on a smooth quintic and write $V=H^0(\PP^4,I_\ell(1))$. The functor of triples $(T,p,q)$ for this line, with $q:I_{\ell/T}\twoheadrightarrow k(p)$ taken up to scalar (and an invertible quotient along a section for families), is a smooth projective scheme
\[
 \PP_{\operatorname{Bl}_\ell X}(K)
 \longrightarrow \operatorname{Bl}_\ell X,
 \qquad 0\longrightarrow K\longrightarrow V\otimes\OO
 \longrightarrow b^*\OO_X(H)(-E)\longrightarrow0 ,
\]
where $b$ denotes blowdown, $E$ is its exceptional divisor, and $\PP(K)$ in this display parametrizes one-dimensional subspaces of $K$. The map is a $\PP^1$-bundle.
\end{proposition}
\begin{proof}
Locally $I_{\ell/X}=(u,v)$ is a regular sequence. In degree $m$, a homogeneous relation on $(u,v)$ reduces modulo $u$ to $a_m v^m=0$, hence $a_m\in(u)$. Subtract a multiple of $vU-uV'$ and divide by $U$; induction gives
\[
 \operatorname{Sym}I_{\ell/X}=\operatorname{Rees}I_{\ell/X}
        =\OO_X[U,V']/(vU-uV').
\]
Consequently the length-one Quot scheme is $\operatorname{Bl}_\ell X$, with quotient $b^*I_{\ell/X}\twoheadrightarrow\OO(-E)$. The three linear equations of $\ell$ generate $I_{\ell/X}(H)$, so $K$ has rank two. For every base, descent of the quotient to $I_{\ell/T}=I_{\ell/X}/(f)$ is equivalent to $q(f)=0$. This is the incidence $\PP(K)$.

At a geometric point $p\in\ell$ the fiber of $I_{\ell/T}$ is
\[
 \frac{I_{\ell/X}}{\mathfrak m_p I_{\ell/X}+(f)}\otimes_{\OO_{X,p}}k(p).
\]
Since $\ell\subset X$ is a regular immersion, its dimension is one if $d(f|_X)_p\ne0$ and two if $d(f|_X)_p=0$. The second case is the unique tangent hyperplane $T_pX$ at $p$, which contains $\ell$. The extra $\PP^1$ of quotients over this locus is included by the construction.
\end{proof}

\begin{proposition}[The whole generic scheme in the mixed charge]\label{prop:linepoint-whole-generic}
On a general smooth quintic the scalar-rigidified moduli of surface sheaves in Proposition~\ref{prop:line-point-hull-d} is the disjoint union, over its lines, of the projective-line bundles in Proposition~\ref{prop:linepoint-blowup-d}. In particular the moduli is smooth of dimension four.
\end{proposition}
\begin{proof}
Let $F=i_*G$, $R=I_{\ell/T}$ and $R/G=k(p)$. The subsheaf $R$ of $\OO_T$ is maximal Cohen--Macaulay. One has $H^1(T,R^\vee)=0$: its Euler characteristic is $3$, $h^0(R^\vee)=1$, and duality gives $h^2(R^\vee)=h^0(R(h))=2$. Thus $\Ext_T^1(R,\OO_T)=0$. The map
\[
 \Ext^2_T(k(p),\OO_T)\longrightarrow\Ext^2_T(R,\OO_T)
\]
is dual to evaluation $H^0(R(h))\to k(p)(h)$. This is surjective even when $T$ is singular at $p$, since the two linear equations of the line in its ambient hyperplane generate $I_{\ell/T}(h)$. The exact sequence for $G\subset R$ and Cartier adjunction consequently give
\[
 (\ext_X^0,\ext_X^1,\ext_X^2)(F,\OO_X(-H))=(0,1,0).
\]
The nonsplit extension with middle term $I_{Z/X}$ is obtained by taking the inverse image of $G$ under $I_{\ell/X}\twoheadrightarrow i_*R$. In particular $I_{\ell/X}/I_{Z/X}=k(p)$. To recover the embedding of this middle term into $\OO_X$, note that
\[
 \Ext_X^1(I_{\ell/X},\OO_X)
 \hookrightarrow\Ext_X^2(\OO_\ell,\OO_X)
 =H^0(\ell,\det N_{\ell/X})=0 ,
\]
where $\det N_{\ell/X}=\OO_\ell(-2)$ by adjunction. Dualizing $I_{Z/X}\subset I_{\ell/X}$ identifies the groups in degrees $0,1$ for both ideals, since $\Ext_X^{<3}(k(p),\OO_X)=0$. It follows that $\Hom_X(I_{Z/X},\OO_X)=\CC$ and $\Ext_X^1(I_{Z/X},\OO_X)=0$.

Over any local Artin algebra the generator of the first extension and the embedding of its middle term therefore lift, uniquely up to units, by perfect base change and the exact sequence of a small extension. The fiberwise criterion for flatness recovers a flat embedded family $Z_A$, and the composite from $\OO(-H)$ recovers its containing hyperplane $T_A$.

The evaluation sequence of $I_{Z/X}(H)\subset I_{\ell/X}(H)$ gives
\[
 (h^0,h^1)(I_{\ell/X}(H))=(3,0),\qquad
 (h^0,h^1)(I_{Z/X}(H))=(2,0).
\]
The first evaluation is onto the length-one quotient by global generation. Base change gives a unique relative plane $\Pi_A$ containing $Z_A$. On the closed fiber
\[
 0\to I_{Z/\Pi}\to I_{\ell/\Pi}=\OO_\Pi(-1)\to k(p)\to0,\qquad
 I_{Z/\Pi}=I_{p/\Pi}(-1).
\]

Put $J=I_{Z_A/\Pi_A}$ and $D=\det J$. Near $p$ lift two generators and their one relation, using flatness and Nakayama at both steps. In local trivializations the determinant map is
\[
 0\longrightarrow\OO\xrightarrow{(a,b)^t}\OO^2
 \xrightarrow{(-b,a)}D\longrightarrow D|_{V(a,b)}\longrightarrow0.
\]
The pair $(a,b)$ reduces to parameters of the plane. It is a relative regular sequence cutting out a finite flat section $p_A$. Away from $p$, $J$ is invertible. Thus $J=I_{p_A/\Pi_A}\otimes D$. Since $\mathcal Hom(I_{p_A/\Pi_A},\OO_{\Pi_A})=\OO_{\Pi_A}$, the inclusion $J\to\OO_{\Pi_A}$ extends uniquely to $D\hookrightarrow\OO_{\Pi_A}$. Its fiber is $\OO_\Pi(-1)\hookrightarrow\OO_\Pi$, hence its zero divisor is a flat line $\ell_A\subset Z_A$.

This recovers $(\ell_A,T_A,I_{\ell_A/X_A}\twoheadrightarrow I_{\ell_A/X_A}/I_{Z_A/X_A})$ up to the stipulated units. The incidence morphism is proper, bijective on geometric points and formally \'etale by this inverse, hence an isomorphism. The line scheme is \'etale by Theorem~\ref{thm:direct-line}.

\end{proof}

\begin{theorem}[The mixed line--point virtual degree by specialization]\label{thm:mixedline-special}
With the degeneration and sign conventions of Theorem~\ref{thm:direct-line}, write $\overline M_{\ell,p}$ for the closure of the mixed-charge generic moduli in the projective relative Simpson moduli. The degree of the specialization of its intrinsic virtual zero-cycle can be computed on the two root components. Per horizontal sheet of the defect line its contributions are
\[
\begin{array}{c|rr|r}
 &Q\text{ root}&B\text{ root}&\text{sum}\\\hline
 D_Q&-158&-238&-396\\
 D_B&-160&-236&-396 .
\end{array}
\]
In particular
\[
 [q^{37/120}]h^{\mathrm{DT}}_{1,1}=(-396)|D_Q|+(-396)|D_B|
 =-1\,138\,500 .
\]
\end{theorem}
\begin{proof}
Theorem~\ref{thm:direct-line} extends the line on each \'etale sheet to a flat curve $\mathcal L\subset\mathcal X$ with smooth total space. On the quartic sheet $\mathcal L$ is the blowup of the line family at $(p,0)$, with
\[
  \mathcal L_0=L_Q\cup E_p,\qquad L_Q\cap S=E_p\cap S=\{p\}
\]
scheme-theoretically. Thus both contacts have order one.

At their crossing, take smooth coordinates $u,v$ along $S$ cutting out $\mathcal L$. Indeed, the two branch tangents of $\mathcal L$ at the node project nontrivially onto the $x,y$ axes of $xy=t$ and their tangent plane meets $T_pS$ trivially. Thus the local coordinates on the total space can be chosen as $(x,y,u,v)$, with $\mathcal L=(u=v=0)$ and $t=xy$.

At every other intersection with a root, the central divisor of $\mathcal L$ is smooth and reduced, so $\mathcal L$ meets that root transversely. On the plane sheet $\mathcal L$ meets only $B$, along a smooth line $L_B$. Hence
\[
 V=\operatorname{Bl}_{\mathcal L}\mathcal X\longrightarrow
 \operatorname{Spec} R
\]
is a regular semistable family: at the node its two charts are obtained by blowing up $(u,v)$ in $xy=t$, and at a smooth point of a root the center is $(u,v)$ in coordinates $(u,v,z,t)$. On the quartic sheet the central components and their double surface are
\[
 (V_Q,V_B,D)=(\operatorname{Bl}_{L_Q}Q,\operatorname{Bl}_{E_p}B,
                        \operatorname{Bl}_pS).
\]
On the plane sheet they are $(Q,\operatorname{Bl}_{L_B}B,S)$.

On $V_\eta$ let $K$ be the rank-two kernel of Proposition~\ref{prop:linepoint-blowup-d} and $a=c_1(\OO_{\PP(K)}(1))$, with $\PP(K)$ parametrizing lines. The generic virtual class is $c_4(\Omega^1_{\PP(K)})\cap[\PP(K)]$ by smoothness and the symmetric sheaf obstruction theory. If $b:\PP(K)\to V_\eta$ is projection, its cotangent sequence and relative Euler sequence give the following equality in the Chow group:
\[
 c_4(\Omega^1_{\PP(K)})
  =b^*c_3(\Omega^1_{V_\eta})(-2a-b^*c_1(K))
  =-2a\,b^*c_3(\Omega^1_{V_\eta}) .
\]
The last term from $V_\eta$ alone has codimension four on a threefold. Let $E=\Omega^1_{V/R}(\log)$ denote the relative logarithmic cotangent bundle, formed by quotienting the absolute logarithmic differentials by $d\log t$. It restricts to $\Omega^1_{V_\eta}$ and to $\Omega^1_{V_i}(\log D)$ on the two roots.

Write $\Lambda=\PP H^0(I_{L/\PP^4}(1))\simeq\PP^2$ and let $U\subset\Lambda$ be the good hyperplanes of Theorem~\ref{thm:direct-line}. Away from $L$, a point $z$ defines the pencil
\[
  \Lambda_z=\{P\in\Lambda:\langle L,z\rangle\subset P\}\simeq\PP^1 .
\]
On each root remove the join, the exceptional loci (including the inverse image of the original exceptional divisor of $B$), and the points for which $\Lambda_z\cap U=\varnothing$. The complement $V^\circ$ is dense on each component. On the quartic root a general plane through $L$ meets the quartic off $L$. On the plane root, varying $z\in\PP^3$ gives pencils sweeping a dense subset of $\Lambda$. When $L\subset\PP^3$ they all pass through the distinguished member $\PP^3$. On $V^\circ$ a general functional on the three equations through $L$ cuts a relative divisor of class $a$ in the pencil.

Here are details ensuring a specialization on $V^\circ$. Choose a relatively very ample $A$ on $V$ so positive that $E(A)$ is relatively generated, and choose sections over the trait lifting a sufficiently general tuple on the closed fiber; a larger twist gives surjectivity on sections. Write $z=c_1(A)$. The identity for a rank-three bundle
\[
 c_3(E)=c_3(E(A))-z c_2(E(A))+z^2c_1(E(A))-z^3
\]
represents $c_3(E)\cap[V]$ by a signed sum of finite horizontal cycles specializing in $V^\circ$. For $1\le i\le3$ use the degeneracy locus of $4-i$ general sections of $E(A)$. The first rank drop has codimension $i$ and the next, when possible, has codimension $2(i+1)>3$. After intersecting with $3-i$ general representatives of $z$, global generation therefore gives points off the excluded loci, all of which have dimension at most two on the special fiber. The same dimension count treats $z^3$.

The strata on the smooth open in each root are transverse by generic smoothness in characteristic zero. The local determinantal equations have full-rank Jacobian on that fiber, so the intersections are \'etale over the trait there. Properness and fiber dimension then make the resulting one-dimensional schemes finite and flat, including multiplicities if one works over a finite extension. Their signed central terms on $V_i$ have degree $\int_{V_i}c_3(\Omega^1_{V_i}(\log D))$.

Choose the functional defining $a$ so that its point on each of the finitely many special pencils lies in $U$. It cuts the relative pencils in sections. By Theorem~\ref{thm:direct-line} the transformed defect is Cartier on their regular total support $\mathcal T$. Write it as $\mathcal Z$. The marked point misses $\mathcal Z$ and the special joins, so the sheaf
\[
 i_*\ker\big(\OO_{\mathcal T}(-\mathcal Z)\longrightarrow
                      \OO_{p}\otimes\OO_{\mathcal T}(-\mathcal Z)\big)
\]
is flat. On smooth fibers it is the required mixed sheaf. On a special fiber its codimension-two modification changes only constant coefficients of the Hilbert polynomials, so the strict linear stability inequalities from Theorem~\ref{thm:direct-line} apply on both roots. This gives a family in $\overline M_{\ell,p}$. Taking closures of these cycles, with their signs and the factor $-2$, therefore computes Chow specialization of the generic virtual class into $A_0(\overline M_{\ell,p,0})_{\QQ}$.

Finally the logarithmic Gauss--Bonnet formula gives
\[
 -2\int_{V_i} c_3(\Omega^1_{V_i}(\log D))
       =2\big(e(V_i)-e(D)\big).
\]
One has $e(Q)=-56,e(B)=-96,e(S)=24$. Blowing up a smooth curve in a threefold adds $e(\ell)=2$ here, and blowing up a point in a surface adds one. Thus the quartic sheet gives $2(-54-25)+2(-94-25)=-396$ and the plane sheet gives $2(-56-24)+2(-94-24)=-396$. The stated counts of sheets from Theorem~\ref{thm:direct-line} finish the computation. At $e=-1$ and $\chi=3$ the charge formula gives $q^{37/120}$ in residue $1$.
\end{proof}

\subsection{The mixed charge in the opposite residue}
\label{sec:dual-line-point}
\begin{theorem}[Specialization after recording the point of the dual line defect]
\label{thm:dual-line-point-d}
For $e(G)=1$ and $\chi(G)=2$ the stable-sheaf virtual degree in
the first quintic degeneration is
\[
 (1600+1275)\bigl(2(-56-24)+2(-96-24)+4\bigr)
     =-1\,138\,500 .
\]
This is the coefficient at $q^{37/120}$ in residue $4$.
\end{theorem}
\begin{proof}
\emph{The full generic sheaf moduli.}
Write $R=G^{\vee\vee}$ and $k=\length(R/G)$ on its hyperplane
support. Then $e(R^\vee)=-1$ and $\chi(R^\vee)=\chi(R)+1=3+k$.
The line-hull bound gives $k\le1$. If $k=0$, $R^\vee$ would
be a reflexive sheaf of degree $-1$ and Euler characteristic three,
excluded in the proof of Proposition~\ref{prop:line-point-hull-d}.
It follows that $k=1$ and $R=(I_{L/T})^\vee$. Conversely every
kernel of a length-one quotient of such an $R$ is stable.
In what follows $R,G$ also denote the sheaves pushed forward to
the smooth threefold.

For one rigid line let $U_L$ denote the universal surface of
hyperplanes containing $L$, a hypersurface of
$X\times\mathbb P^2$. The projectivization $M_L^+$ of the universal
dual line ideal over $U_L$ classifies the quotient maps $R\to k(p)$
up to scalar. We show first that $M_L^+$ is smooth. Off the graph
$\Gamma_L=\{(p,T_pX):p\in L\}$ the surfaces are regular along $L$
or the point is off $L$, and the universal module is invertible.

The universal hypersurface is smooth off that graph: away from $L$
the equation varies surjectively in the hyperplane parameters, and
on $L$ at a nontangent hyperplane its differential on $X$ is nonzero.
At a graph point write the regular ideal of $L$ on $X$ as $(u,v)$
and its containing hyperplane equation as $au+bv=0$.

The differentials
of $u,v,a,b$ are independent on $X\times\mathbb P^2$ there.
Indeed the three linear equations defining $L\subset\mathbb P^4$
generate $I_{L/X}(H)$; at $p\in L$ the kernel of the evaluation on
its two-dimensional conormal fiber is the line spanned by $T_pX$.
Thus the two parameters of the hyperplane plane provide the two
independent derivatives of $(a,b)$.

With column relations, the module $(u,v)$ on $U_L$ and its dual
have square presentation matrices
\[
 P=\begin{pmatrix}-v&a\\u&b\end{pmatrix},\qquad P^{\mathsf t}.
\]
For the dual this follows by dualizing its projective-dimension-one
resolution on the smooth ambient space and applying Cartier
adjunction, with local line bundle factors trivialized.
The quotient is specified by a nonzero row vector $q$ modulo
scalars satisfying $qP^{\mathsf t}=0$.

These two equations have
independent derivatives at every $[q]$ above $\Gamma_L$, since
the four entries of $P$ have independent differentials. The
projectivization is consequently a smooth fourfold. Above
$\Gamma_L$ its fiber is $\mathbb P^1$ and elsewhere its map
to $U_L$ is an isomorphism.

We verify that this fourfold has the infinitesimal structure of the
moduli of the kernels. This matters at a point on the noninvertible
locus. Work at an exact sequence
\[
 0\longrightarrow G\xrightarrow{i}R\xrightarrow{q}S_p
     \xrightarrow{\delta}G[1],\qquad S_p=k(p).
\]
The arrow to $G[1]$ is the connecting map of the associated
triangle. The obstruction to lifting $q$ between given
first-order deformations of $R,S_p$ is the difference of the two
Yoneda products. The resulting map
\[
 q_*\oplus(-q^*):
 \Ext^1(R,R)\oplus\Ext^1(S_p,S_p)\longrightarrow\Ext^1(R,S_p)
\]
is surjective in this situation. In a minimal local free resolution
of $R$ on the threefold, the differential is a square matrix $P_0$
of size $s$ with $P_0(p)=0$. Here $s=2$ on the graph, by the
presentation above; off the graph $s=1$ and $P_0$ is a local
surface equation. Thus $\Ext^1(R,S_p)=\mathbb C^s$. Vary the
hyperplane and the point, represented by a tangent vector $v$
of $\mathbb P^2\times X$.

The product map on these variations is
$v\mapsto q\,dP_0(v)$, up to the sign convention for pullback:
this follows by lifting the row $q$ to the free resolution modulo
a square-zero parameter. At a graph point the four independent
differentials above give surjectivity onto $\mathbb C^2$. Elsewhere
smoothness of $U_L$ gives surjectivity onto $\mathbb C$.
The hyperplane variations are actual deformations of $R$ by
Proposition~\ref{prop:polar-incidence-functors}.

The obstruction for a first-order extension
$\xi\in\Ext^1(G,G)$ to carry a lift of $i$ into some deformation
of $R$ lies in $\Ext^2(S_p,R)$ and equals
$\theta=i[2]\xi[1]\delta$. Indeed the exact sequence obtained by
applying $\operatorname{Hom}(-,R)$ identifies its vanishing with
$i_*\xi=i^*\zeta$ for some $\zeta\in\Ext^1(R,R)$; this equality
is precisely the existence of a morphism of the two short exact
deformation sequences. By Calabi--Yau Serre duality $\theta$
pairs with $\Ext^1(R,S_p)$.

It annihilates both summands of the
displayed product map: for the first use $\delta q=0$, and for
the second use $qi=0$ after cyclic permutation in the trace.
Surjectivity proves $\theta=0$. A lifted map reducing to $i$
is injective with flat quotient by the local flatness criterion;
this quotient deforms $S_p$.

Thus every first-order deformation
of $G$ comes from a tangent vector of the flag fourfold.
For injectivity on tangents with $G$ fixed,
$\Ext^1(S_p,R)=0$ by $\depth_p R=2$, hence
$\Ext^1(R,R)\to\Ext^1(G,R)$ is injective. The target $R$ then
has trivial deformation, and $\operatorname{Hom}(G,R)=
\operatorname{Hom}(R,R)=\mathbb C$ by the same exact sequence.
A change in the lifted inclusion is scalar and preserves the flag.

The morphism from the disjoint union of the projective flag spaces
to the separated stable-sheaf moduli is proper and bijective on
geometric points by recovery of the hull. It is finite, and both
sides have local dimension four. The tangent calculation gives
embedding dimension four on the target; its local rings are regular.

A finite bijection of these schemes is birational in characteristic
zero and the target is normal. It is therefore an isomorphism.
In particular the symmetric sheaf virtual class on each component
is $c_4(T M_L^+)\cap[M_L^+]$.

\emph{Its central evaluation.} Project $\pi:M_L^+\to X$ to the
point $p$. Away from $L$ its fiber consists of the pencil of
containing hyperplanes. At $p\in L$ the base of the fiber is the
whole plane of such hyperplanes, all of whose points have one
preimage except $T_pX$, which has a projective line of preimages.

The pushforward of constructible functions is consequently
$\pi_*\mathbf 1_{M_L^+}=2\mathbf1_X+2\mathbf1_L$.
The proper naturality and smooth normalization of the
Chern--Schwartz--MacPherson transformation \cite[\S1]{KennedyCSM}, taking its zero-dimensional
component, give in $A_0(X)$
\[
 \pi_*[M_L^+]^\vir
 =2c_3(TX)\cap[X]+2j_*\bigl(c_1(TL)\cap[L]\bigr),\qquad j:L\hookrightarrow X.
\]
The algebraic transformation applies over the algebraic closure of
the characteristic-zero fraction field. The Chow identity descends
with rational coefficients: a rational equivalence of finite type
over the algebraic closure is defined over some finite extension,
and restriction followed by pushforward multiplies by its degree.

Take the closure of the graph of $\pi_\eta$ in the product over the
trait of the projective relative Simpson space and $\mathcal X$ for
each horizontal sheet of Theorem~\ref{thm:direct-line}. It is proper
and is the surface-sheaf moduli on that generic sheet. Chow
specialization of its virtual zero-cycle commutes with its proper
projection to $\mathcal X$.

The relative logarithmic tangent bundle
of $\mathcal X$ specializes to
$T_Q(-\log S)$ and $T_B(-\log S)$ on the two roots. The corresponding
bundle of the horizontal line family specializes to $T_{L_i}(-D_i)$
on each component $L_i$, where $D_i$ consists of its joining points.
Thus the degree of the pushed-forward special zero-cycle on this
sheet is
\[
 2\bigl(e(Q)-e(S)\bigr)+2\bigl(e(B)-e(S)\bigr)
       +2\sum_{i}\bigl(2-\#D_i\bigr).
\]
Indeed Chow specialization commutes with Chern classes of a bundle
on the total space and logarithmic Gauss--Bonnet evaluates each
degree. For a quartic-root line there are two projective lines with
one joint on each, so the final sum before its factor $2$ is
$1+1=2$. For a plane-root line it is $2-0=2$.

Consequently the
degree on either sheet is $-396$. Summing the $1600$ and $1275$
special sheets of Theorem~\ref{thm:direct-line} gives the formula.
This argument evaluates the specialization of the sheaf class after
proper pushforward to the position of the point; it does not
identify a full product of relative virtual classes in this charge.
\end{proof}

\subsection{The split-line summands with one further point}
\label{sec:mixed-split-marked}
Work over a strictly henselian trait $A$. Fix a finite \'etale line sheet in Theorem~\ref{thm:direct-line} with
\[
 C\subset\XX_A,\qquad C_0=L\cup f_p,\qquad p\in S .
\]
The point $z$ of the additional quotient may specialize to $p$.

\begin{proposition}\label{prop:mixed-split-marked}
For the negative line ideal $I_{C/T}$ with one further punctual quotient the complete smooth-fiber sheaf component over this sheet admits a proper expanded compactification $M_C$. Its two marked virtual zero-cycles have degrees $-158$ and $-238$. Their contribution, summed over the split-line sheets, is $-633\,600$ in residue $1$. This is a partial contribution to the coefficient at $q^{37/120}$.
\end{proposition}
\begin{proof}
Write $V$ for the rank-three free $A$-module of hyperplane equations through the moving line. The evaluation
\[
 V\otimes\OO_{\XX_A}(-H)\longrightarrow I_{C/\XX_A}
\]
is surjective. The regular immersion $C$ has codimension two on the total space and restricts to a regular immersion on either root. On $Q$ use the linear equations of $L$. On $B$ the three equations pull back the ideal $(a,b,c)$ of $p\in\PP^3$. With $\Gamma=(a,b)$ their pullbacks on its two blowup charts are
\[
  (a,a(b/a),c)=(a,c),\qquad (b(a/b),b,c)=(b,c),
\]
the ideals of $f_p$.

At the joint the two generators $u,v$ of $I_C$ restrict to a regular system of parameters on $S$ at $p$, by Lemma~\ref{lem:linedefect-target-rigidity}. The images of $V$ span this cotangent space, so Nakayama's lemma gives the assertion there. Flatness of $C$ and Nakayama give the relative assertion.

Let $b:Y=\operatorname{Bl}_{C}\XX_A\to\XX_A$ and let $E$ be the exceptional divisor. A regular sequence of length two is of linear type, by the elementary symmetric-algebra calculation in Proposition~\ref{prop:linepoint-blowup-d}. Thus $Y$ represents a point with a one-dimensional quotient of the ideal of $C$. Define
\[
 0\to K\to V\otimes\OO_Y\to b^*\OO(H)(-E)\to0,
 \qquad P=\PP_Y(K),
\]
where $\PP$ parametrizes lines in $K$. The section of $V$ defining $T$ lies in $K$ precisely when the point quotient descends to $I_{C/T}$. This equality is valid over arbitrary test rings. On the closed fiber $Y$ has roots $Y_Q=\operatorname{Bl}_L Q$, $Y_B=\operatorname{Bl}_{f_p}B$ and common surface $D=\operatorname{Bl}_p S$: near the joint blow up $(u,v)$ in $xy=t$. In each blowup chart $xy=t$ is still the equation of the family. Thus $P$ is regular and semistable, with roots $P_i=\PP_{Y_i}(K|_{Y_i})$ and crossing $D_P=\PP_D(K|_D)$.

To construct the sheaf over $P$ use the section $z$ furnished by $P\to Y\to\XX_A$. Refine the expanded target at $z$ by Lemma~\ref{lem:one-defect-universal}, with $ab=t$, $z'(P)$ in the interior of the bubble at a crossing. Pull back $T$ and $C$ to this refinement. Near the crossing their ideals are a Cartier equation and a regular sequence $(u,v)$ pulled back from coordinates transverse to the branches.

Consequently these pullbacks remain relatively Cartier and relatively regular of codimension two, respectively. On the new component the curve gains the ruling above $p$; if $z$ contracts to a different point the same ruling is added over $p$ away from $z'$. The pulled-back ideal together with its point quotient gives a kernel $G$ flat over the parameter, agreeing with the line hull near all refined joins.

The smooth-target point has finite projective dimension and the hull has a square locally free resolution, hence $G$ is perfect. Conversely a quotient of this pulled-back ideal at the interior section is exactly the line quotient on the contracted target: its evaluation on the section is the pullback of the two-term presentation of $(u,v)$. The equations $q(T)=0$ descend. The point-refinement inverse and its changes of trivialization in Lemma~\ref{lem:one-defect-universal} prove these assertions over nilpotent test rings as well.

We verify that forgetting the presentation gives the local deformation theory of $G$ as a sheaf. Fix an expanded target $W$ at a geometric point of this construction. The contracted sheaf hull pulls back to $R=I_{C'/T'}$ on $W$, where $C',T'$ denote the pulled-back curve and support. The recovery lemma for the polar line ideal applies also on the transparent refinement.

Indeed the normal bundle of $C'$ pulls back from $C_0$; on the new ruling it is the constant rank-two bundle, and normalization cancels its sections against one copy of the matching fiber. Hence $H^0(N_{C'/W})=H^1(N_{C'/W})=0$. Likewise $H^0(\omega_{T'}|_{C'})$ and the restriction from $H^0(\omega_{T'})$ are unchanged, since $H$ pulls back, has degree zero on the extra ruling, and the complete target has trivial dualizing bundle.

Pullback cancels the inserted $\PP^1$-bundle in the normalization sequences for $\OO_{T'}$, $\OO_W(H)$, giving the same section and vanishing groups. The square presentation and ideal-sequence proof of Lemma~\ref{lem:lineideal-recovery-on-expanded-support} therefore applies: every fixed-target hull deformation is the hyperplane variation through the uniquely lifted curve. In particular these hyperplane variations are actual deformations of $R$.

At $0\to G\xrightarrow{i} R\xrightarrow{q} k(z')\to0$ it suffices for the first-order deformation of every kernel on the fixed target to lift to this flag that
\[
 \Ext_W^1(R,R)\oplus\Ext_W^1(k(z'),k(z'))
 \xrightarrow{q_*\oplus(-q^*)}\Ext_W^1(R,k(z'))
 \quad\text{be surjective}.
\]
Here is the check, including the noninvertible hull. The point $z'$ is smooth on the threefold $W$. Off $C'$ the hull is locally the hypersurface module; the three generators from $V$ have no common zero, so the incidence equation $T(z')=0$ has a nonzero derivative in a hyperplane parameter. On $C'$ use the local matrix with columns $(-v,u)$ and $(a,b)$ for $T'=au+bv$.

At a fixed $z'$ where $(a,b)\ne0$ the supporting surface is smooth and motion of $z'$ yields the required one-dimensional image. Where $(a,b)=0$, variations of $u,v$ give two independent cotangent directions, and variation of the hyperplane at this point maps surjectively to the conormal fiber spanned by $(a,b)$.

This last surjection follows from the generation by $V$ proved above, also on the added ruling by pullback. Thus the differential of the matrix is onto $\mathrm{Mat}_2$, and its composition with the nonzero row $q$ is onto the two-dimensional $\Ext_W^1(R,k(z'))$. These computations take place at the smooth interior point even when its contraction equals $p$.

Let $k=k(z')$ and let $\delta:k\to G[1]$ be the boundary map. The obstruction to lifting $i$ along $\xi\in\Ext^1_W(G,G)$ to a first-order deformation of $R$ is
\[
  i[2]\circ\xi[1]\circ\delta\in\Ext^2_W(k,R)\simeq\Ext^1_W(R,k)^\vee .
\]
The duality holds for these perfect objects on the proper Gorenstein Calabi--Yau threefold $W$. Since $\delta q=0$ and $qi=0$, cyclicity of trace kills the pairings with $q_*\Ext^1_W(R,R)$ and $q^*\Ext^1_W(k,k)$. Surjectivity above kills the obstruction.

Conversely, $\operatorname{depth}_{z'} R=2$ gives $\Hom_W(k,R)=\Ext^1_W(k,R)=0$. Applying $\Hom_W(-,R)$ to the flag gives
\[
 \Hom_W(G,R)=\Hom_W(R,R)=\CC,\qquad
 \Ext^1_W(R,R)\hookrightarrow\Ext^1_W(G,R).
\]
Thus flags and kernels have the same fixed-target tangents modulo scalars. Recovery on geometric points uses the annihilator for the support, the double dual on that support for $R$, and $R/G$ for the marked point.

Work on a framed expansion chart $B$. Write
$g:P^{\rm fr}\to M^{\rm fr}$ for the morphism from flags to
scalar-rigidified kernels on this chart. The flag space is smooth over
$B$. At the refined join this follows on the blowup chart
$xy=t$: up to the freely framed ruling rescaling the two target
parameters pull back to $a=x(z)$, $b=y(z)$, independent coordinates
on the blowup; $P$ is smooth over that blowup. Away from a join use
the ordinary chart. For $p\mapsto m\mapsto b$, the fixed-target
tangent calculation and smoothness give
\[
\begin{array}{ccccccccc}
0&\to&T_p(P^{\rm fr}_b)&\to&T_pP^{\rm fr}&\to&T_b B&\to&0\\
 &&\downarrow\wr&&\downarrow dg&&\Vert\\
0&\to&T_m(M^{\rm fr}_b)&\to&T_mM^{\rm fr}&\to&T_b B&\to&0 .
\end{array}
\]
Surjectivity at the bottom right follows by following the top row.
Thus $dg$ is an isomorphism. At geometric points the support, double
dual and quotient recover the flag; the framed morphism is representable
and injective on geometric fibers after scalar rigidification. It is
locally quasi-finite. On an \'etale chart, at a closed point the image
closure of the smooth flag component through $p$ therefore has the
dimension of that component. Since $B$ is smooth this gives
\[
 \dim_pP^{\rm fr}\le\dim_m M^{\rm fr}\le
 \dim T_mM^{\rm fr}=\dim T_pP^{\rm fr}=\dim_pP^{\rm fr}.
\]
The local ring at $m$ is regular. Thus $g$ is \'etale there and,
being universally injective, is an open immersion. These conclusions
hold on the completed smooth base trait as well, using the same
regular expansion charts before completion and base change.
The new bubble dilation fixes the interior support of the quotient
only at $\lambda=1$; the compatible branch scalars coincide.
For the pullback nef polarization on the inserted ruling
$\rho:T_b\to C_b$ one has
\[
 P_{\rho^*E}=P_E=B_e,\qquad
 P'_{\le e+1}=P_{\le e}-B_e+P_{\rho^*E}=P_{\le e}.
\]
The two refined cuts therefore inherit the two original strict linear
bounds. They persist under a sufficiently small ample perturbation
as in Theorem~\ref{thm:one-defect-global}. The branch kernels
are rank-one torsion-free; the point quotient subtracts only $1$
from the appropriate constant term. This identifies the proper $P$
with an open expanded stable sheaf locus.

The two normalized marked components are $P_Q,P_B$. On either one the unmarked parameter cuts out $D_P$ simply. The smooth map to the marked expansion stack consequently has relative tangent $T_{P_i}(-\log D_P)$. By relative Serre duality and the cotangent triangle as in Proposition~\ref{prop:one-defect-positive-marked},
\[
 [P_i]^{\vir}=c_4(\Omega^1_{P_i}(\log D_P))\cap[P_i],\qquad
 \deg[P_i]^{\vir}=2\bigl(e(Y_i)-e(D)\bigr).
\]
The factor two is the Euler characteristic of the projective-line fiber. Here $e(Y_Q)=-54$, $e(Y_B)=-94$ and $e(D)=25$. The claimed marked degrees and their sum over $1600$ sheets follow. Both are intrinsic classes on proper marked spaces. The proof does not require independently proper unrestricted relative factors.
\end{proof}

\begin{theorem}[The complete mixed line and point degree by marked specialization]\label{thm:mixed-all-marked}
For a general first degeneration both mixed charges, $(e,\chi)=(-1,3)$ and $(1,2)$, admit proper expanded surface-sheaf spaces with complete respective generic fibers. Their marked degrees per split-line sheet are $(-158,-238)$; per four-secant sheet the degrees are $(-160,-114,-122)$ in the order positive with point on the quartic, positive with point on the blowup, and zero boundary. Consequently in both residues $1,4$ the marked degeneration gives
\[
 [q^{37/120}]h^{\mathrm{DT}}_{1,\mu}
       =1600(-158-238)+1275(-160-114-122)=-1\,138\,500 .
\]
\end{theorem}
\begin{proof}
Fix a horizontal line sheet. On the minimal universal surface $\mathcal T\subset\mathcal W$ over the proper hull space $J$ of Theorem~\ref{thm:expanded-line-compact}, put
\[
 R=I_{C/\mathcal T}\ \text{or}\ I_{C/\mathcal T}^{\vee},\qquad
 P=\operatorname{Proj}_{\mathcal T}\operatorname{Sym} R .
\]
Both hulls are flat in families and perfect on $\mathcal W$ (Corollary~\ref{cor:line-opposite-expanded} for the dual). The projective $J$-space $P$ represents one-dimensional quotients. Pull back the target to $P$ and refine the point to an interior section $z'$ by Lemma~\ref{lem:one-defect-universal}. The kernel of the tautological quotient along $z'$ is the required sheaf. At a split join the defect ideal is $(u,v)$ in $xy=t$. The four-secant defect misses the joins.

Thus the square presentations of the hull pull back to relative projective-dimension-one presentations, and the dual presentation pulls back as its transpose; the new support is the Cartier pullback. This gives precise base-change compatibility for the point quotient under refinement. Conversely the intermediate levels carrying a nontrivial graph or the chosen point are recovered from the sheaf: the geometric double dual on its Cartier support gives the polar hull, and the cokernel recovers the point.

The graph level is required already by the hull. A transparent level carrying only the point contracts on forgetting it. Stability and the intrinsic stabilizer follow by the strict cut computation at the end of Proposition~\ref{prop:mixed-split-marked}, which uses only stability of the chosen polar hull and the pullback rule on a transparent component.

On the inserted ruled surface $\rho:T_b\to C_b$ the
boundary module and the pulled-back hull are
\[
 E=\begin{cases}
 I_{p/C_b}\ \text{or }I_{p/C_b}^{\vee},&\text{split line},\\
 \OO_{C_b},&\text{nonempty secant join},
 \end{cases}
 \qquad R|_{T_b}=\rho^*E.
\]
The plane quartic is integral since $\Pic(S)=\mathbb Zh$.
A torsion-free rank-one $E$ embeds locally into the function
field of $C_b$; on the affine-line bundle charts flat pullback
preserves this injection into the function field of $T_b$.
This includes a singular point of $C_b$. The split ideal is given
locally by the pulled-back generators $(u,v)$, and its dual by their
transposed presentation (flat pullback commutes with dualization of a
finitely presented module).

Let $H_0$ denote the old polarization pulled back to the refinement.
By $R\rho_*\OO_{T_b}=\OO_{C_b}$ and the projection formula,
$P^{H_0}_{\rho^*E}=P_E=B_e$, and $[m^2]B_e=0$.
For the hull, using the partial-polynomial notation of
\textup{(\companionref{eq:cutdata})}, insertion gives
\[
\begin{aligned}
 P'&=P+B_e-B_e=P, & A'&=A,\\
 P'_{\le e}&=P_{\le e},&
 P'_{\le e+1}&=P_{\le e}-B_e+B_e=P_{\le e},\\
 A'_e&=A'_{e+1}=A_e,&
 [m]D'_e&=[m]D'_{e+1}=[m]D_e .
\end{aligned}
\]
The old polar cut bound is $0<[m]D_e<[m]B_e$ for the split and
ordinary secant joins, and uses $d_1=0$ at the graph join.
For the new interval alone the residual linear coefficient is
$[m](B_e-2B_e)=-[m]B_e<0$. All other cuts have their previous
linear coefficients. If $H_\epsilon=H_0+\epsilon H_{\rm amp}$ for a
chartwise relatively ample $H_{\rm amp}$, the coefficients in
\textup{(\companionref{eq:cutdata})} are continuous at $0$ since $A(0)>0$.

For $0<\epsilon\ll1$ rational, all the strict cut inequalities
persist on the specified chart. Taking the kernel of a length-one
quotient leaves $[m^2]P$ and $[m]P$ unchanged. The rank-one
torsion-free branches therefore satisfy Theorem~\companionref{thm:cutcriterion},
so the kernel is stable.

If a target automorphism fixes the kernel it
fixes the recovered hull and the support of its quotient. On the new
level the interior ruling coordinate $z'\ne0,\infty$ satisfies
$\lambda z'=z'$ only for $\lambda=1$; the remaining dilations
already fix the contracted hull only trivially. The pre-rigidification
stabilizer is the common scalar $\mathbb G_m$.

We verify target-relative smoothness. Work in an ordered expansion chart with the tori framed, and first fix its target parameters. The hull deformations are exactly the support variations of the polar family on this target. This holds on its minimal chart by Theorem~\ref{thm:expanded-line-compact}, and after a further transparent refinement by the normal-bundle and cohomology cancellation in Proposition~\ref{prop:mixed-split-marked}.

In the four-secant case the defect is disjoint from the extra level, so its normal-sheaf groups are literally unchanged. For the dual use the transposed presentation and its involution on hull deformations. We describe the fixed-target Jacobian of $P$.

At an interior point away from the defect the equation is the incidence equation $T(z')=0$. Along the defect in the smooth part of a root take $I_{C/W}=(u,v)$ and $T'=au+bv$. At a point with $a=b=0$, the two sign charts impose the two entries of $qP_0=0$ or $qP_0^{\mathsf t}=0$ for a nonzero row $q$, where
\[
  P_0=\begin{pmatrix}-v&a\\u&b\end{pmatrix}.
\]
It suffices that $u,v,a,b$ have independent differentials in the combined point and fixed-target hull parameters. Precisely then these quotient equations have independent derivatives for every $[q]\in\PP^1$. Off that locus one has the regular hypersurface incidence test. We check the tests on every root, including the base locus of the secant pencil.

For a split sheet the three linear equations generate the ideal of the curve on each root and, at the contact, on $S$. Hence the computations of Proposition~\ref{prop:mixed-split-marked} give the asserted four independent differentials for either presentation.

If the point is on a transparent bubble created at the contact, then $u,v$ are coordinates transverse to its ruling and both hyperplane parameters are still available for $a,b$. More generally the universal support incidence restricted to $S$ is smooth in the same quotient charts; away from $p$ some plane variation has nonzero value, and at $p$ the equations generate its cotangent space.

For a four-secant sheet write $L'\subset B$ for the strict transform. On the ordinary positive support the quartic root has no defect. The system there is the plane of hyperplanes through the projected line. Its base on $Q$ consists of the four points $p_j=L\cap S$. At each such point every containing hyperplane cuts $Q$ smoothly since $L$ is transverse to $S$ in $\PP^3$ and is therefore transverse to $Q$ in $\PP^4$.

Every section of $S$ by a plane containing $L$ is likewise smooth at those base points. Outside the four base points the plane system moves the incidence off the point. For the plane root, let $l_1,l_2$ be plane equations generating the ideal of $L$ in $\PP^3$. Both vanish along $L'$ and the four exceptional fibers $f_j$. Away from the four crossings $w_j=L'\cap f_j$ their common zero is locally the smooth curve $L'$ or $f_j$, with ideal generated by $l_1,l_2$.

For $f_j\setminus L'$ this follows from the blowup along $\Gamma$: the derivative along $\Gamma$ projects to a nonzero direction normal to $L$, and the ruling direction specifies a second independent direction except at the unique normal ray of $L$ itself. Thus the ordinary pencil incidence away from the defect is smooth there.

At $w_j$ an explicit chart also treats the zero-support and graph profiles. Use a local parameter $z$ along $L$ chosen as $f_4|_L$ using simple contact, extend $z=f_4$ via a local trivialization and set $u=f_5/f_4-u_0$, where $u_0$ is a linear extension of $f_5|_L/f_4|_L$. This is regular on the chart $B\setminus S$ and equals zero on $L'$. Write $U=l_1$ (after a change of basis), $l_2=B_1$, and
\[
 f_5-f_4u_0=dU+c B_1 .
\]
Here $c$ is a unit at $p_j$ after changing basis: $df_5,df_4$ are independent along $\Gamma$, and $f_5-f_4u_0$ vanishes on $L$. Solving the blowup equation gives
\[
 B_1=A(z,U,u)U+B(z,U,u)zu,\qquad B(w_j)\ne0 .
\]
Consequently $(z,U,u)$ are regular coordinates on $B$ at $w_j$, $L'=(U,u)$ and $f_j=(z,U)$. On the ordinary/graph-boundary plane pencil, the combination with vanishing linear term is unique. Its coefficient of $U$ varies with derivative a unit in the pencil parameter, while its coefficient of $u$ has derivative a unit in the $z$ direction.

The four differentials of $U,u,a,b$ are therefore independent there. All other plane members have nonzero differential along $L'$ at $w_j$. Away from the $w_j$ the generating pencil spans the conormal along $L'$, so no ordinary plane in the pencil is singular there along $L'$.

The $t$-chart of \textup{(\companionref{eq:vacuumtchart})} on $B\setminus S$ has equation $u+\lambda l_1+\mu l_2=0$ for its affine plane of supports through $L'$. At $w_j$ it is smooth, since the coefficient of $du$ is one. At a point of $f_j\setminus L'$ there is no such surface, since $l_1=l_2=0$ and $u\ne0$. Off the base, the affine plane of support parameters moves the incidence off the point. Along $L'\setminus\{w_j\}$ both $l_1,l_2$ generate the normal fiber, so there is exactly one tangent support and the two parameters give the two required derivatives.

At the intersection where the graph level appears, the support on $B$ is the plane pencil and the calculation at $w_j$ above still applies: even with both target parameters fixed at zero the plane-pencil parameter varies. On this fixed graph target the support section is $\zeta+\ell=0$ with the linear form $\ell$ constrained by $\ell|_{L'}=0$.

Its variations include the full pencil direction even at fixed smoothing parameters: in a frame fixing the two target parameters the remaining overall coefficient of $\ell$ is the dilation coordinate. Thus fixing the target imposes no extra plane-pencil condition on the preceding Jacobian test. The graph on the ruled level is a section, hence smooth in its interior.

At its active join the incidence is the plane-pencil incidence on $S$, already verified to be smooth; the other join misses the support. For a further point refinement, the new interior condition is precisely this same incidence on $S$, with an independent invertible ruling parameter. The four-secant defect misses the joining divisor throughout. This proves all the asserted fixed-target Jacobian tests.

These tests imply smoothness over the expansion charts themselves. After framing, the hull chart is smooth over its target chart: for the split hull it is the plane bundle before the point refinement, and for the secant hull the blowup has coordinates $s,v,\lambda$ with $sv=t$ at the crossing. Passage to the target chart records $s,v$ up to the freely framed torus, giving a smooth morphism there.

At a new point level replace the one target parameter $d$ by $a_0,b_0$ with $d=a_0b_0$. Adding this coordinate pair to the hull chart is base change of that smooth morphism. The lifted interior point moves in the smooth locus of the universal expanded threefold; if this chart has a new level its interior ruling coordinate gives a free framing, with gauge $1$ at the point.

The equations for lying on the surface or for its one-dimensional quotient are exactly the fixed-target incidences just checked. The corresponding invertible Jacobian minors give smoothness over the target chart by the relative Jacobian criterion, also when several target parameters vanish.

The unmarked smoothing parameters pull back to independent smooth divisors. In particular the normalized marked components are smooth with simple normal crossings boundary given by the remaining parameters. Their interiors are the unexpanded ordinary and zero-boundary incidences.

The objects so constructed are the kernels as sheaves in families. Indeed the same Jacobian tests show that the product map on $\Ext^1$ in the proof of Proposition~\ref{prop:mixed-split-marked} is onto at every point, for both signs: outside the defect the derivative maps to the rank-one target, and at a two-generator hull it is the row times the derivative of the displayed matrix or its transpose.

The Calabi--Yau trace and depth argument there now identifies all fixed-target kernel tangents with flag tangents. The smoothness on framed charts, quasi-finiteness and local dimension argument of that proposition gives a representable open immersion of the kernels in the scalar-rigidified stack of perfect sheaves on expansions.

This proves that the obstruction bundle on each proper normalized marked fourfold is $\Omega^1(\log D_{\rm mark})$. It gives virtual degree equal to the Euler characteristic of the complement of the remaining boundary. The generic identifications are those of Proposition~\ref{prop:linepoint-whole-generic} and Theorem~\ref{thm:dual-line-point-d}.

Here are the signed marked counts. Use $Q^\circ=Q\setminus S$, $B^\circ=B\setminus S$. For a split sheet, away from the respective curve the fiber over the position is $\PP^1$. On $L\setminus\{p\}$ and $f_p\setminus\{p\}$ each position admits the plane of support hyperplanes, with one tangent member having an extra $\PP^1$ of quotients in place of a point. Its total fiber Euler characteristic is $4$. Thus
\[
 2e(Q^\circ)+2e(\AA^1)=-158,\qquad
 2e(B^\circ)+2e(\AA^1)=-238
\]
also for the dual.

For a secant sheet in the positive ordinary stratum, a point on $Q^\circ$ determines a pencil of hyperplanes through the projected line and the point, none equal to $\PP^3$. Its Euler characteristic yields $2e(Q^\circ)=-160$. A point on $B^\circ$ off the inverse image of the projected line determines a unique plane, and the hyperplanes over that plane are an affine line.

Over the inverse image the plane parameter ranges in $\PP^1$, with the same affine-line fibers. This base locus is $L'\cup\bigcup_j(f_j\setminus S)$, of Euler characteristic $2+4-4=2$. At each $w_j$ the one tangent plane acquires a quotient $\PP^1$, for each parameter of the affine line of lifts. Thus the second positive-marked degree is
\[
 e(B^\circ)+2+4=-114.
\]
In the zero-boundary unexpanded stratum the two affine support parameters solve $u+\lambda l_1+\mu l_2=0$. The point must lie in $B^\circ$. At each $p_j$ the exceptional fiber meets $S$ in precisely one point, distinct from $w_j$ since $L'\cap S=\varnothing$. Thus $f_j\cap B^\circ$ has Euler characteristic one. Off the same base locus the fiber is an affine line.

On a ruling $f_j$ outside $L'$ there is no solution. Over $L'$ the fiber is an affine plane of choices with a single tangent support except at the four points $w_j$, where all the supports are smooth. Each tangent support adds one to the Euler characteristic of the quotient fiber. Hence this marked degree is
\[
 (e(B^\circ)-2)+e(L')+e(L'\setminus\{w_1,w_2,w_3,w_4\})
 =-122+2-2=-122.
\]
This gives $-160-114-122=-396$. The split-sheet degrees sum to the same number. In all counts the one-dimensional quotients are parameters for actual kernel sheaves, and the virtual multiplicity comes from their intrinsic marked obstruction theory. No numerical curve or pair invariant enters the calculation.
\end{proof}

For either mixed sign cut the sheaf on each positive marked piece at its joining K3. The one-ended refinement keeps the point in the smooth interior of its half, so both derived boundary restrictions equal the polar value $E\in J_\chi$.

Each active branch has generic rank one on an integral support. Its sheaf endomorphism algebra is finite birational over the projective support, hence has only constant global functions. Agreement across a nonempty join identifies the branch scalars. The interior point or graph support fixes every possible dilation. Thus the hypotheses on simplicity and stabilizers in Theorem~\ref{thm:point-open-diagonal} hold.

Take finite-type neighborhoods $U_1,U_2$ of the two cuts in the one-ended simple perfect sheaf stacks with this stable boundary as in that theorem. Properness of the marked locus gives a finite cover by such neighborhoods. Normalization glues the two kernels together with their boundary isomorphism; cutting in families inverts this construction. With the same scalar and determinant reductions one has on each proper marked locus
\[
 j_\alpha:P^\dagger_\alpha\hookrightarrow U_1\times_{J_\chi}U_2,\qquad
 [P^\dagger_\alpha]^{\vir}
   =j_\alpha^*\Delta_{J_\chi}^!
       ([U_1]^{\vir}\times[U_2]^{\vir}).
\]

The deformation complex before reduction is the fiber of the difference of boundary restrictions. The reduced virtual pullback is Theorem~\ref{thm:point-open-diagonal}. At an empty join the occupied sheaf misses its outer divisor, giving the one-sided obstruction theory of Theorem~\companionref{thm:zeroboundaryPOT}. Proper pushforwards apply to the marked loci in this formula.

\subsection{A mixed point arriving at a joining divisor}
At a joining divisor the derived restriction of a length-one kernel can have cohomology in arbitrarily negative degrees.

\begin{proposition}\label{prop:mixed-kernel-at-join}
Let $A=\mathbb C[[x,y,u,v]]/(xy)$, let $B=A/(x,y)$ be the local ring of the joining K3, and let $k=A/(x,y,u,v)$. Suppose $R$ has projective dimension at most one over $A$ and $0\to G\to R\to k\to0$ is exact. Then
\[
 H^{-i}(G\otimes_A^{\mathbf L} B)
   \cong k^{\oplus2}\quad (i\ge2).
\]
In particular, taking this kernel on the contracted fiber cannot extend the mixed moduli by perfect sheaves with perfect boundary values.
\end{proposition}
\begin{proof}
A free resolution of $B$ starts with $A^2\xrightarrow{(x\ y)} A$ and continues with the alternating diagonal matrices $\operatorname{diag}(y,x)$ and $\operatorname{diag}(x,y)$. Indeed $(x)\cap(y)=0$ in $A$, $\operatorname{ann}(x)=(y)$ and $\operatorname{ann}(y)=(x)$. These equalities give the kernel of the first arrow and prove exactness at each subsequent term. After tensoring with $k$ every differential is zero. Thus $\operatorname{Tor}_j^A(k,B)=k^2$ for $j\ge1$. For $i\ge2$, the exact segment
\[
 0=\operatorname{Tor}_{i+1}^A(R,B)
 \longrightarrow\operatorname{Tor}_{i+1}^A(k,B)
 \longrightarrow\operatorname{Tor}_i^A(G,B)
 \longrightarrow\operatorname{Tor}_i^A(R,B)=0
\]
gives the assertion.
\end{proof}

For example, at a split line choose $I_{L\cup f_p}=(u,v)$ as in Lemma~\ref{lem:linedefect-target-rigidity}. If the supporting surface is cut out by a nonzerodivisor $au+bv$, then $R=(u,v)/(au+bv)$ has the two-term presentation in Lemma~\ref{lem:lineideal-recovery-on-expanded-support}; the proposition applies to every length-one quotient at this point. On a refinement placing the point in the interior of a bubble, the quotient restricts to zero near the two new joins.

Thus the refined kernel and hull have the same derived restrictions at the new joins. Theorem~\ref{thm:mixed-all-marked} identifies their deformation-theoretic comparison for the two primitive line hulls. The proposition shows why the Quot space on the contracted hull alone cannot supply this identification.

\begin{figure}[htbp]\centering
\begin{tikzpicture}[x=1cm,y=1cm,font=\footnotesize,line join=round,line cap=round]
\node[font=\small\bfseries] at (2.35,2.87) {A point at the join};
\node[font=\small\bfseries] at (9.16,2.87) {A point on the refinement};
\path[draw=qIIblue,fill=qIIblue!7] (.25,.95)--(.73,1.97)--(2.68,1.97)--(2.20,.95)--cycle;
\path[draw=qIIblue,fill=qIIblue!7] (2.20,.95)--(2.68,1.97)--(4.65,1.97)--(4.17,.95)--cycle;
\draw[qIIteal,line width=1.5pt] (2.20,.95)--(2.68,1.97);
\fill[qIIclay] (2.43,1.44) circle(2.9pt);
\draw[-{Stealth[length=5pt]},qIIblue] (6.61,1.47)--(5.05,1.47);
\node[text=qIIblue] at (5.82,1.84) {$c$};
\path[draw=qIIblue,fill=qIIblue!7] (7.00,.95)--(7.48,1.97)--(8.52,1.97)--(8.04,.95)--cycle;
\path[draw=qIIteal,fill=qIIteal!11] (8.04,.95)--(8.52,1.97)--(10.19,1.97)--(9.71,.95)--cycle;
\path[draw=qIIblue,fill=qIIblue!7] (9.71,.95)--(10.19,1.97)--(11.38,1.97)--(10.90,.95)--cycle;
\foreach \y in {1.16,1.46,1.77}{\draw[qIIteal!42] ({8.04+.48*(\y-.95)/1.02},\y)--({9.71+.48*(\y-.95)/1.02},\y);}
\foreach \x in {8.04,9.71}{\draw[qIIteal,line width=1.35pt] (\x,.95)--({\x+.48},1.97);}
\fill[qIIclay] (9.04,1.46) circle(2.9pt);
\node[text=qIIclay,above] at (9.04,1.53) {$z$};
\node at (2.33,.38) {$H^{-i}(\mathbf L i^*G)=k^2,\quad i\ge2$};
\node at (9.12,.38) {$\mathbf L i_\pm^*k(z)=0$};
\end{tikzpicture}
\caption{A quotient at a joining divisor and its refinement in Theorem~\ref{thm:mixed-all-marked}. The point in the bubble interior misses both new joins, so the lifted hull and kernel have identical derived restrictions there. On the contracted target the kernel has the higher Tor groups of Proposition~\ref{prop:mixed-kernel-at-join}.}\label{fig:II-mixed-kernel}
\end{figure}
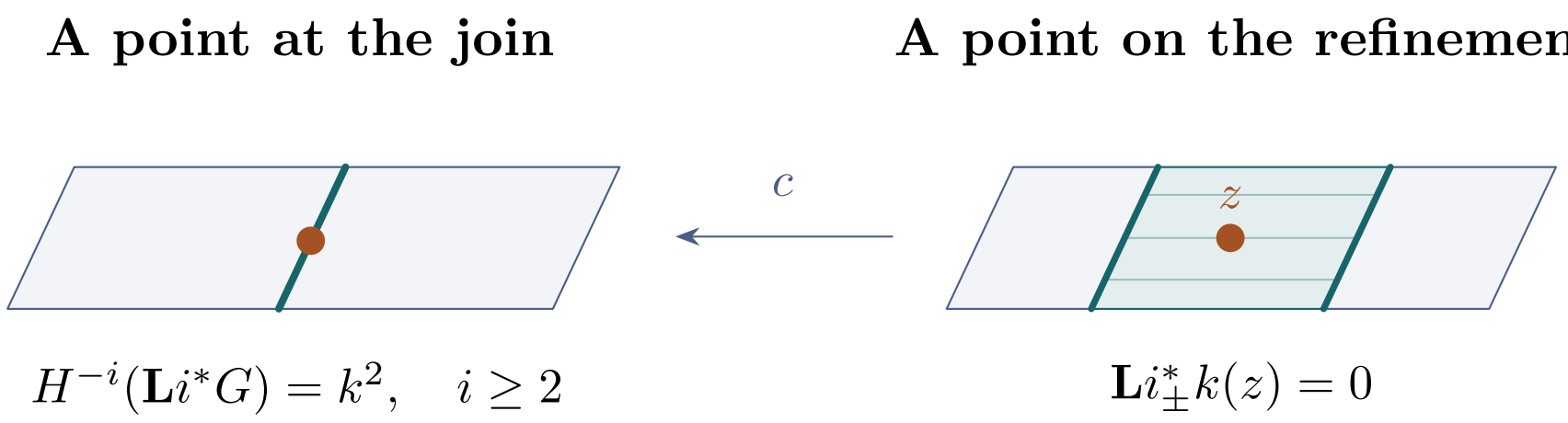

\subsection{The next neutral charge and higher divisor classes}
\begin{proposition}[Additional hulls at the third neutral defect]\label{prop:third-neutral-additional-hulls}
On a general smooth quintic threefold the class
$\ch(\OO_T)-3[\mathrm{pt}]$ contains stable sheaves of the form
$i_*\OO_T(L_1-L_2)$, where $L_1,L_2$ are disjoint lines on a smooth hyperplane section $T$. In particular the locus of length-three point ideals does not exhaust this charge.
\end{proposition}
\begin{proof}
Ordered intersecting lines in $\PP^4$ depend on $4+3+3=10$ parameters. For their union $C$ the plane conic sequence gives
\[
 H^0(\PP^4,\OO(5))\twoheadrightarrow H^0(C,\OO_C(5)),\qquad
 h^0(C,\OO_C(5))=6+6-1=11 .
\]
Thus a general quintic contains no such pair.

The general quintic has disjoint lines. For an ordered skew pair $L_1,L_2$ write $P=\langle L_1,L_2\rangle\simeq\PP^3$. The triples $(L_1,L_2,X)$ with $L_i\subset X$ and $L_1,L_2$ skew form an open set in a projective bundle over the irreducible locus of skew pairs: the restrictions to their union are surjective in degree five. This last assertion follows in coordinates in $P$ from $L_1=(x_0=x_1=0)$ and $L_2=(x_2=x_3=0)$, for which a sum of two binary quintics realizes arbitrary restrictions.

The general quintic through both lines has smooth $T=X\cap P$: Bertini applies off the two lines, and on either line the two normal derivatives are base-point-free quartics with no common zero for a general section. The bad locus is thus proper in the irreducible incidence space.

The previously computed finite \'etale defect scheme has more than one point on a general quintic, so the incidence dominates the space of quintics and is generically finite. It follows that for a general $X$ the spans of the indicated skew pairs cut smooth hyperplane sections.

Adjunction on $T$ gives $K_T=h$, $hL_i=1$ and $L_i^2=-3$. For $D=L_1-L_2$,
\[
 hD=0,\qquad D^2=-6,\qquad
 \chi(\OO_T(D))=5+\tfrac12D(D-h)=2 .
\]
Both $\OO_T(D)$ and the ideal $I_Z$ of three points are stable of the stated numerical class. Their hulls differ since $D^2\ne0$. Even within the ideal locus the fibers over $Z$ vary:
\[
 \{P:Z\subset P\}\simeq\PP^{3-\dim\langle Z\rangle}
  =\begin{cases}\PP^2& Z\text{ collinear},\\ \PP^1&\text{otherwise}.\end{cases}
\]
\end{proof}

\begin{proposition}[Additional hulls in the next conic charge]\label{prop:additional-conic-hulls}
The primitive charge $(e,\chi)=(-2,3)$ on a general quintic contains both length-one modifications of conic ideals and locally free sheaves on their supporting surfaces with trivial punctual defect. In particular the compactification of conic ideals with one added point computes at most one locus of that charge.
\end{proposition}
\begin{proof}
For the first assertion take a smooth conic $C\subset T$ and modify $I_{C/T}$ by a quotient of length one. Its degree is $-2$ and the ideal sequence gives $\chi(I_{C/T})=5-1=4$, hence its kernel has Euler characteristic three. For the second assertion choose two disjoint lines $L_1,L_2$ on a smooth hyperplane section $T$ as in Proposition~\ref{prop:third-neutral-additional-hulls}. Write $h=H|_T$. For $D=-L_1-L_2$, adjunction gives $K_T=h$, $hD=-2$ and $D^2=-6$. Therefore
\[
 \chi(\OO_T(D))=5+\tfrac12(D^2-hD)=3.
\]
The sheaf $i_*\OO_T(D)$ is stable by rank one on the integral support. Its hull quotient has length zero, whereas the kernel in the first construction has hull quotient of length one. These are distinct sheaves with the same charge.
\end{proof}

\begin{proposition}[Additional line bundles in the conic--point character]\label{prop:conicpoint-competitor18}
On a general smooth quintic both primitive charges $(e,\chi)=(-2,3)$ and $(2,1)$ contain invertible sheaves on smooth hyperplane supports. They also contain punctual modifications of the two conic hulls in Lemma~\ref{lem:wholeconic-18}. Consequently a marked construction restricted to one point on those hulls cannot compute the entire coefficient in either charge.
\end{proposition}
\begin{proof}
Let $C\subset P=\PP^3\subset\PP^4$ be a twisted cubic. The family of such embedded curves is irreducible: choose $P$ and take the orbit under $\mathrm{PGL}(P)$ of the cubic Veronese embedding of $\PP^1$. The degree-five restriction to $C$ is surjective, since products of five binary cubics span all binary forms of degree fifteen. Therefore the incidence of $(C,X)$ with $C\subset X$ is a projective bundle on this family.

At the example of \cite[Appendix, (6.1)--(6.4)]{KatzCubic} the normal sequence gives
\[
 N_{C/X}\simeq\OO_{\PP^1}(-1)^{\oplus2},\qquad
 H^0(N_{C/\PP^4})\longrightarrow H^0(\OO_C(5H))
            \longrightarrow H^1(N_{C/X})=0 .
\]
The incidence is smooth over the space of quintics at this pair by the deformation criterion for a regular immersion. It therefore dominates. In the span $P$ the three quadratic generators of $I_C$ show that $I_C(5)$ generates $N^\vee_{C/P}(5)$ along $C$. A general section of this rank-two bundle on the curve has no zero. Bertini off $C$ then gives a smooth $T=X\cap P$. This dense open of the irreducible incidence still dominates.

On such a smooth $T$, $K_T=h$, $hC=3$, $h^2=5$, and $C^2=-2-hC=-5$. Set $D=C-h$. Then
\[
 hD=-2,\qquad D^2=-6,\qquad
 \chi(\OO_T(D))=5+\tfrac12(D^2-hD)=3,\qquad
 \chi(\OO_T(-D))=1 .
\]
Both pushforwards are stable by integrality and rank one. Each equals its surface double dual. In the same characters a length-one modification of a conic ideal or its dual is strictly contained in its surface double dual. This distinguishes the sheaves in the moduli without making any comparison of numerical virtual degrees.
\end{proof}

For every $a\in\mathbb Z$, tensoring with $\OO_X(aH)$ gives an isomorphism on the smooth-fiber sheaf moduli with character
\[
 (0,D,\beta,n)\longmapsto
 (0,D,\beta+aHD,n+aH\beta+\tfrac{a^2}{2}H^2D).
\]
Indeed $P_{F(a)}(m)=P_F(m+a)$ and $R\!\operatorname{Hom}(F(a),F(a))\simeq R\!\operatorname{Hom}(F,F)$ compatibly with trace. The virtual degree and Fourier index are constant on each twist orbit (Proposition~\ref{prop:generalspectralflow}). For $r\ge2$, Theorem~\companionref{thm:matrixexpandedallr} gives higher-multiplicity limits and their point modifications. At $r(1-e^{-H})$ the sheaf $\OO_T^{\oplus r}$ has scalar-rigidified stabilizer $\mathrm{PGL}_r$ (Proposition~\companionref{prop:determultiple}).
\begin{proposition}[The virtual contribution of two skew lines]\label{prop:twoskew-component}
Let $X$ be a general quintic. For distinct lines $L_1,L_2\subset X$, let $T=X\cap\langle L_1,L_2\rangle$ with inclusion $i:T\hookrightarrow X$. Fix nonzero integers $a,b$. Each sheaf
\[
                 F_{a,b}=i_*\OO_T(aL_1+bL_2)
\]
is an isolated reduced point of its stable-sheaf moduli and contributes $+1$ to its virtual degree. For $(a,b)=(1,-1)$ these points contribute $8\,262\,750$ at $(\mu,x)=(0,17/24)$. For $(a,b)=(-1,-1)$ they contribute $4\,131\,375$ at $(\mu,x)=(2,133/120)$. Each number is the degree of the Chow specialization of the indicated union of points in the projective relative Simpson space.
\end{proposition}

\begin{proof}
\emph{The support and its two lines.}
Proposition~\ref{prop:third-neutral-additional-hulls} gives disjoint lines and smooth sections $T$ for a general $X$. Such a $T$ contains just the two lines spanning its hyperplane. Indeed the space of a hyperplane in $\PP^4$ and three ordered skew lines in it has dimension $4+3\cdot4=16$. Three skew lines in $\PP^3$ lie in one ruling of a smooth quadric $V$. Their union imposes $18$ independent conditions on quintics. To see surjectivity of restriction, first restrict to $V$, then use
\[
 H^1(\PP^3,\OO(3))=0,\qquad
 H^1(V,\OO_V(5,5)\otimes\OO_V(-3,0))=H^1(V,\OO_V(2,5))=0 .
\]
The incidence of triples therefore has image of codimension at least two in the space of quintic equations.

\emph{First-order deformations.}
Set $h=H|_T$ and $V_0=H^0(T,\OO_T(h))$. Adjunction gives $K_T=N_{T/X}=\OO_T(h)$, so $H^2(T,\OO_T)=V_0^\vee$. A deformation of $i_*M$ over the dual numbers, for a line bundle $M$ on $T$, is a deformation of the Cartier divisor together with a line bundle on it. Locally the sheaf has a free resolution by a single equation.

Flat lifting of its generator and relation gives a presentation by a lifted equation. Its zeroth Fitting ideal recovers the deformed divisor, and the lifted modules glue as line bundles on it. For a fixed divisor the tangent space of its Picard scheme is $H^1(T,\OO_T)=0$.

Write $o_j:V_0\to V_0^\vee$ for the obstruction to lifting $\OO_T(L_j)$ with the divisor, and $r_j:V_0\to H^0(L_j,\OO_{L_j}(1))$ for restriction. We claim
\[
   \ker o_j=\ker r_j,qquad \operatorname{im}o_j=(\ker r_j)^\perp .
\]
The sequence $0\to\OO_T\to\OO_T(L_j)\to\OO_{L_j}(-3)\to0$ gives $H^1(T,\OO_T(L_j))=0$. Indeed the connecting map $H^1(\OO_{L_j}(-3))\to H^2(\OO_T)$ is dual to the surjection $r_j$. Thus whenever the line bundle lifts its canonical section lifts with it. The zero divisor of this lifted section is a flat embedded deformation of $L_j$ in $X$.

The line scheme of a general quintic is reduced and zero-dimensional by the calculation in Theorem~\ref{thm:direct-line}. Its tangent space is zero. The lifted hyperplane therefore still contains $L_j$. Conversely any first-order hyperplane containing $L_j$ carries the constant line as a relative Cartier divisor. This proves the kernel formula.

Write $s_j$ for the canonical section. If $1+\epsilon c_{\alpha\beta\gamma}$ is the obstruction cocycle for the transition functions, $s_jc_{\alpha\beta\gamma}$ is the coboundary of the local section discrepancies. Hence
\[
 \operatorname{im}o_j\subset
 \ker\!\left(H^2(\OO_T)\xrightarrow{\,s_j\,}H^2(\OO_T(L_j))\right)
       =(\ker r_j)^\perp .
\]
The last equality is Serre duality. Both sides have dimension two by the kernel formula.

The lines span the hyperplane, so $(r_1,r_2)$ identifies the four-dimensional $V_0$ with $H^0(\OO_{L_1}(1))\oplus H^0(\OO_{L_2}(1))$. Hence $\ker r_1\cap\ker r_2=0$ and $\operatorname{im}o_1\cap\operatorname{im}o_2=0$. Tensor products add the obstruction cocycles. In particular
\[
   \ker(a o_1+b o_2)=0\qquad(ab\ne0).
\]
This eliminates every first-order support deformation of $F_{a,b}$. The tangent space along the fixed support was zero, so $\Ext^1_X(F_{a,b},F_{a,b})=0$. By Serre duality $\Ext^2_X(F_{a,b},F_{a,b})=0$ as well. Stability is rank-one stability on the integral surface. Each moduli point is consequently open and closed and has virtual degree one.

\emph{Specialization and charges.}
The horizontal line scheme in Theorem~\ref{thm:direct-line} has $N=1600+1275=2875$ geometric sheets. After a splitting base change, label its ordered pairs of distinct sheets. On the generic fiber every pair gives the sheaf above. Properness of the projective relative Simpson space extends the resulting moduli points over a trait after base change.

For a section the specialization of its fundamental zero-cycle has degree one. Proper pushforward and additivity therefore preserve the sum of the point contributions, including when several sections specialize to the same closed point. This calculation uses Chow specialization in the Simpson space.

Distinct unordered pairs determine distinct supports. Since $(L_1-L_2)^2=-6$, reversing the order changes $\OO_T(L_1-L_2)$, whereas $\OO_T(-L_1-L_2)$ is unchanged. For $D=aL_1+bL_2$, with $L_i^2=-3$, $hL_i=1$, $L_1L_2=0$, Riemann--Roch gives
\[
 hD=a+b,\qquad \chi(\OO_T(D))=5-\tfrac{3(a^2+b^2)+a+b}{2}.
\]
Thus
\[
\begin{array}{c|rr|rr|r}
(a,b)&hD&\chi(\OO_T(D))&\mu&x&\text{degree}\\ \hline
(1,-1)&0&2&0&17/24&N(N-1)\\
(-1,-1)&-2&3&2&133/120&N(N-1)/2 .
\end{array}
\]

These are virtual degrees of specified open-and-closed loci in their respective characters.
\end{proof}

\begin{corollary}[Point modifications of the two-line sheaves]\label{cor:twoskew-all-lengths}
Fix one of the line bundles $L=\OO_T(aL_1+bL_2)$ of Proposition~\ref{prop:twoskew-component}, and write $x_L$ for the exponent of $i_*L$. For each $k\ge0$ the sheaves $i_*(L\otimes I_{Z/T})$ with $\length Z=k$ form a smooth open-and-closed component of their moduli, isomorphic to $T^{[k]}$. Their contribution to the series in that residue is
\[
                  q^{x_L}\prod_{m\ge1}(1-q^m)^{-55}.
\]
\end{corollary}
\begin{proof}
Put $G=L\otimes I_{Z/T}$. On the smooth surface $T$ one has $\det G=L$ and $\mathcal End(G)=\OO_T$. For the latter identity an endomorphism is multiplication by a rational function, regular off $Z$ and hence everywhere. The spectral sequence for the Cartier immersion $i:T\hookrightarrow X$ gives
\[
 0\to\Ext^1_T(G,G)\to\Ext^1_X(i_*G,i_*G)
   \to H^0(\OO_T(h))\xrightarrow{\partial_G}\Ext^2_T(G,G).
\]
The middle arrow records the first-order action of a local equation of $T$ on a deformation. On an open set where the equation is $t$, an element $v\in H^0(\OO_T(h))$ acts by multiplication by $v$ on $G$. A deformation inducing this action is a module on the first-order divisor with equation $t-\epsilon v$.

Indeed $t$ acts on the deformed sheaf by $\epsilon v$, since the latter takes values in the square-zero submodule identified with $G$ by flatness. These local equations glue to the embedded first-order deformation of $T$ corresponding to $v$.

A flat family of $G$ on this first-order surface is perfect there, since the surface is smooth over the dual numbers and $G$ has a finite locally free resolution on the special fiber. Taking its determinant extends $L$ to that surface. The kernel calculation of Proposition~\ref{prop:twoskew-component} forces $v=0$. Thus the injection on the left in the displayed sequence is an isomorphism.

Riemann--Roch on $T$ gives $\chi(G,G)=\chi(\OO_T)-2k=5-2k$. By Serre duality $\ext^2_T(G,G)=h^0(K_T)=4$, while $\ext^0_T(G,G)=1$. The tangent dimension on $X$ is therefore $2k$. The universal ideal defines a proper morphism $T^{[k]}\to M_X(v_k)$ for $v_k=\ch(i_*L)-k[\mathrm{pt}]$. Every geometric point recovers $L=G^{\vee\vee}$ and $Z$ from $G\subset L$, so the morphism is finite. Fogarty's theorem \cite[Theorem~2.4]{Fogarty} and the tangent calculation give at every image point
\[
 2k=\dim T^{[k]}\le\dim\OO_{M_X(v_k),[i_*G]}
       \le\ext^1_X(i_*G,i_*G)=2k.
\]
The target is regular there. The finite map is birational to its open-and-closed image by geometric injectivity in characteristic zero. Normality gives the required isomorphism.

The Behrend sign is $(-1)^{2k}=1$. Adjunction gives $e(T)=11\int_T h^2=55$ from $c(T_T)=(1+h)^4/(1+5h)$. Hence \cite[Theorem~0.1]{GottscheBetti} gives
\[
 \sum_{k\ge0}q^{x_L+k}\deg[T^{[k]}]^\vir
       =q^{x_L}\sum_{k\ge0}e(T^{[k]})q^k
       =q^{x_L}\prod_{m\ge1}(1-q^m)^{-55}.
\]
Chow specialization in the proper relative Simpson space preserves each degree as in Proposition~\ref{prop:twoskew-component}.

\end{proof}

\section{Further coefficients from the complete moduli on the smooth quintic}\label{sec:smooth-evaluations}
Fix $v=1-e^{-H}-3[\mathrm{pt}]$ and $T=X\cap P$ for a hyperplane $P$. Besides the ordered pairs of lines in Proposition~\ref{prop:twoskew-component}, the same character has a component of length-three ideals.
\begin{theorem}\label{th:seven-total}
Let $\mathcal I=\{(Z,P): Z\subset X\cap P,\ \length Z=3\}$. There is a decomposition of schemes
\[
 M\simeq \mathcal I\ \sqcup \coprod_{\substack{L_1,L_2\subset X\text{ lines}\\L_1\ne L_2}}\Spec\CC.
\]
On $\mathcal I$ the associated sheaf is $i_*I_{Z/T}$. The point indexed by $(L_1,L_2)$ corresponds to $i_*\OO_T(L_1-L_2)$ on the section of the hyperplane spanned by the two lines. In particular
\[
             \deg[M]^{\mathrm{vir}}=-2\,445\,625+8\,262\,750
                                     =5\,817\,125.
\]
\end{theorem}

With the charge convention of the introduction $q_1=-H\ch_2(v)=5/2$, $q_0=\int_X\ch_3(v)+25/12=-1/12$, and $\mu=q_1-5/2=0$. Thus $q_1^2/10-q_0=17/24$. Equivalently this is the $q^3$ term after removing $q^{-55/24}$ from $Z^{\mathrm{GSY}}_0$ in \cite[(3.10)]{GSY}. The geometric proof uses the smooth threefold $X$ throughout.

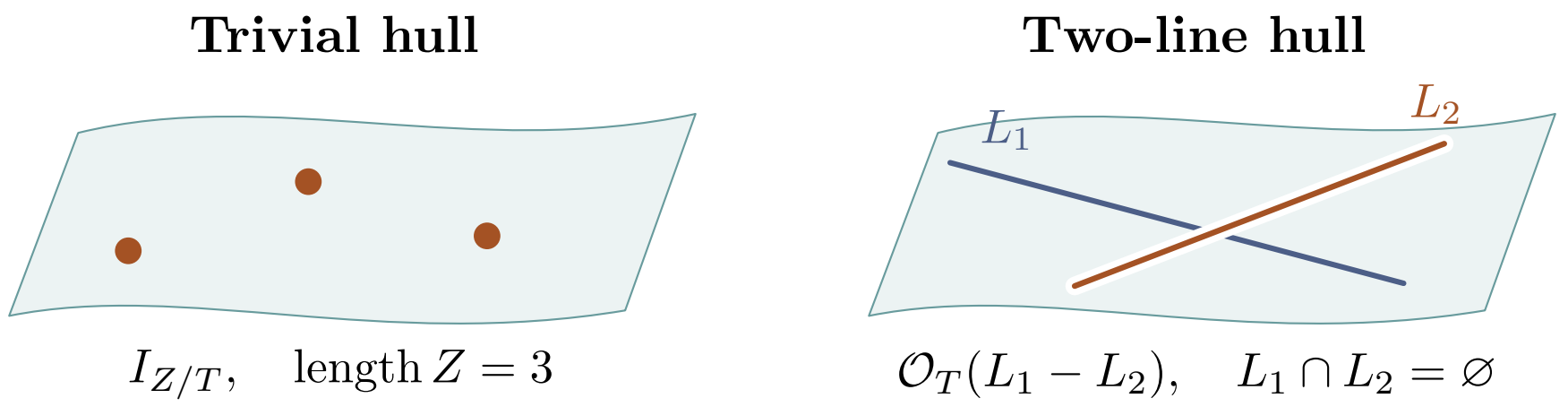
\begin{figure}[htbp]\centering
\begin{tikzpicture}[x=1cm,y=1cm,font=\footnotesize,line cap=round,line join=round]
\foreach \x/\lab in {0/{Trivial hull},6.35/{Two-line hull}}{
 \begin{scope}[shift={(\x,0)}]
 \node[font=\small\bfseries] at (2.65,2.91) {\lab};
 \path[draw=qIIteal!65,fill=qIIteal!8]
 (.25,.83)..controls (1.55,1.09) and (2.99,.57)..(4.80,.87)
 --(5.32,2.32)..controls (3.68,1.97) and (2.19,2.53)..(.76,2.18)--cycle;
 \end{scope}}
\foreach \x/\y in {1.13/1.31,2.46/1.82,3.78/1.42}{\fill[qIIclay] (\x,\y) circle(2.8pt);}
\node at (2.70,.40) {$I_{Z/T},\quad\length Z=3$};
\draw[qIIblue,line width=1.15pt] (7.20,1.96)--(10.55,1.07);
\draw[white,line width=3.9pt] (8.12,1.05)--(10.85,2.10);
\draw[qIIclay,line width=1.2pt] (8.12,1.05)--(10.85,2.10);
\node[above,text=qIIblue] at (7.62,1.91) {$L_1$};
\node[above,text=qIIclay] at (10.79,2.10) {$L_2$};
\node at (9.02,.40) {$\OO_T(L_1-L_2),\quad L_1\cap L_2=\varnothing$};
\end{tikzpicture}
\caption{The two types of sheaf in the neutral three-point charge. For the trivial hull the picture shows three distinct points; the incidence also contains nonreduced schemes. For the ordered skew lines the support is the section cut out by their span. The crossing is only in projection.}\label{fig:II-neutral-triple}
\end{figure}

\subsection{Classification of the closed points}
\begin{proposition}[Classification in character $v$]\label{prop:row7-hulls}
Let $X\subset\mathbb P^4_{\mathbb C}$ be a general smooth quintic and let $H$ denote its hyperplane class. The stable pure two-dimensional sheaves on $X$ with support class $H$ and Hilbert polynomial
\[
                       P_F(m)=\frac52m(m-1)+2
\]
are exactly the following sheaves. Write $i:T\hookrightarrow X$ for a hyperplane section.
\begin{enumerate}
\item $i_*I_{Z/T}$, where $Z\subset T$ has length three.
\item $i_*\mathcal O_T(L_1-L_2)$, where $L_1,L_2$ are distinct lines on $X$ and $T=X\cap\langle L_1,L_2\rangle$.
\end{enumerate}
The two cases are disjoint. The description includes sheaves supported on singular hyperplane sections.
\end{proposition}
\begin{proof}
Since $\operatorname{Pic}(X)=\mathbb ZH$ and $H$ is primitive ample, every $T\in|H|$ is integral. For the Gauss map $g:X\to(\PP^4)^\vee$ we have
\[
 g^*\OO(1)=\OO_X(4),\qquad \operatorname{Sing}(T)=g^{-1}([T])
          \quad\text{on closed points}.
\]
Thus $g$ has finite fibers. The hypersurface $T$ satisfies $R_1$ and $S_2$, hence is normal.

The support cycle of a sheaf in the statement is an integral $T$ with generic multiplicity one. Multiplication by a local equation of $T$ has image of dimension at most one and hence vanishes by purity. Consequently $F=i_*G$ for a rank-one torsion-free sheaf $G$ on $T$.

A nonzero section of $G$ would inject $\mathcal O_T$ into $G$. The polynomial of the cokernel would be $P_G-P_{\mathcal O_T}=-3$, which is impossible. Hence $h^2(G)=2+h^1(G)\ge2$. The dualizing sheaf of $T$ is $\mathcal O_T(1)$. By Serre duality there is an injection $G\hookrightarrow\mathcal O_T(1)$, whose image is $I_{W/T}(1)$ for a subscheme $W$ with
\[
 P_{\mathcal O_W}(m)=P_{\mathcal O_T}(m)-P_G(m-1)=5m-2.
\]
Divide $\mathcal O_W$ by its maximal zero-dimensional subsheaf, of length $b$, to obtain $\mathcal O_C$. Thus $C\subset T\subset\mathbb P^3$ is locally Cohen--Macaulay, $\deg C=5$ and $p_a(C)=3+b$.

If $C$ lies in a plane $\Pi$, the plane quintic $T\cap\Pi$ contains $C$ and they agree. Indeed they have the same degree and the Cohen--Macaulay curve $T\cap\Pi$ has no subsheaf supported at finitely many points. Hence $I_{C/T}(1)=\mathcal O_T$. The inclusion $I_{W/T}(1)\subset I_{C/T}(1)$ identifies $G$ with an ideal of colength three.

Suppose $C$ is nonplanar. The genus bound and the bound on $h^1(I_C(m))$ in \cite[Theorem~9.1]{HS} give $3+b\le3$ and $h^1(I_C(m))=0$ for every $m$. Therefore $b=0$, $C=W$, and $C$ attains the extremal Rao function in \cite[Definition~9.2]{HS}. By \cite[Proposition~9.5]{HS} it contains a plane quartic $D$ in a plane $\Pi$. In degree five Proposition~9.5 of that reference excludes its low-degree alternatives. Both $[T\cap\Pi]-[D]$ and $[C]-[D]$ are effective one-cycles of degree one, hence lines $L_1,L_2$. In the Weil divisor group of the normal surface $T$ this gives
\[
                     [C]=[T\cap\Pi]-L_1+L_2.
\]

The quotient $\mathcal O_T/I_{C/T}$ is Cohen--Macaulay, so the depth lemma makes $I_{C/T}$ reflexive. The divisor description now gives $G=I_{C/T}(1)\simeq\mathcal O_T(L_1-L_2)$. If the two lines coincide then $\chi(G)=5$, contrary to $\chi(G)=2$. By Proposition~\ref{prop:twoskew-component}, the two distinct lines are skew and span the unique hyperplane containing them, and its section $T$ is smooth.

The ideal sequence gives the polynomial for the first type. For the second put $A=L_1-L_2$. Adjunction and disjointness give
\[
 L_j^2=-3,\quad hA=0,\quad A^2=-6,\qquad
 P_{\OO_T(A)}(m)-P_{\OO_T}(m)=\tfrac12 A(A-h)=-3 .
\]
Both types have integral rank-one support, hence are stable. Their hulls have Euler characteristics $5$ and $2$, respectively.
\end{proof}
\subsection{The scheme structure on the point-ideal locus}
\begin{proposition}\label{prop:row7-point-component}
Let $X\subset\PP^4$ be a general smooth quintic threefold, write
$H=c_1(\OO_X(1))$, and let
\[
 \mathcal I=\{(Z,T):Z\in\operatorname{Hilb}^3 X,\quad
                         T\in|H|,\quad Z\subset T\}.
\]
The morphism $(Z,T)\mapsto i_* I_{Z/T}$, $i:T\hookrightarrow X$, identifies
$\mathcal I$ with an open and closed subscheme of the stable sheaf moduli space
$M_X(\operatorname{ch}(\OO_T)-3[\mathrm{pt}])$.
It is smooth and projective of dimension ten. The contribution of this
subscheme to the Donaldson--Thomas invariant is its ordinary topological
Euler characteristic.
\end{proposition}

\begin{proof}
\emph{The incidence scheme.}
For a smooth threefold $Y$, the scheme $\operatorname{Hilb}^3Y$ is smooth
of dimension nine. This can be checked in completed coordinates at the
support. A curvilinear local subscheme of length $m\leq3$ is a complete
intersection and its tangent space has dimension $3m$. The remaining
local algebra is $D=\CC[x,y]/(x,y)^2$, embedded in a smooth threefold
with ideal $J=(z,x^2,xy,y^2)$.

The two plane syzygies show that the
images of $x^2,xy,y^2$ under a homomorphism $J\to D$ lie in $(x,y)D$.
Conversely any such images and any image of $z$ respect the relations:
the syzygies involving $z$ are generated by the Koszul relations
between $z$ and the three plane generators. Thus
\[
             \dim\Hom(J,D)=3+3\cdot2=9.
\]
This length-three subscheme is the limit in $z=0$ of
$(0,0),(t,0),(0,t)$, with ideals $(x(x-t),xy,y(y-t))$. Curvilinear
subschemes split along a smooth curve. Each local type is consequently
in the closure of the distinct-point locus. Disjoint supports contribute
the sum of their local tangent dimensions. The asserted regularity follows.

The family of hyperplanes $\PP^3\subset\PP^4$ is a projective-space
bundle over $(\PP^4)^\vee$. Its relative Hilbert scheme of subschemes of length
three is smooth of dimension $4+9$. The universal quintic cuts out
$\mathcal I$ as the zero scheme of a section of the rank-three tautological
bundle $\OO_{\PP^4}(5)^{[3]}$ on this smooth scheme.

The evaluation on every length-three subscheme is surjective: after
dehomogenizing by a hyperplane avoiding the support, affine linear
coordinates generate its three-dimensional algebra, and monomials of
degree at most two span this algebra. Indeed the degree filtration
strictly increases until it is the whole algebra. Homogenize to degree
five.

The universal zero scheme over the space of quintic equations is
therefore a vector bundle over the relative Hilbert scheme. Generic
smoothness gives smooth $\mathcal I$ of dimension $13-3=10$ for general $X$.
Every smooth $X$ has triples on its smooth hyperplane sections, so the
universal projection dominates. Projectivity follows from the incidence
description.

\emph{A tangent calculation on the threefold.}
For $(Z,T)\in \mathcal I$, set $A=\OO_X(-H)$, $B=I_{Z/X}$ and
$F=i_*I_{Z/T}$. There is an exact triangle
\[
          A\xrightarrow{s} B\xrightarrow{p} F\xrightarrow{a} A[1].
\]
The integral hyperplane section $T$ carries $I_{Z/T}$ as a torsion-free
sheaf of rank one, hence $F$ is stable. Integrality here follows from
$\operatorname{Pic}(X)=\mathbb Z H$ and the primitivity of $H$.
The tangent space to $\operatorname{Hilb}^3X$ at $Z$ is
$\Ext^1_X(B,B)$. Indeed $\Hom_X(B,\OO_X)=\CC$ and
$\Ext^1_X(B,\OO_X)=0$ by $0\to B\to\OO_X\to\OO_Z\to0$,
local duality in dimension three and $H^1(X,\OO_X)=0$.
Applying $\Hom_X(B,-)$ identifies $\Hom_X(B,\OO_Z)$ with
$\Ext^1_X(B,B)$.

The three incidence equations on the smooth $\operatorname{Hilb}^3X\times|H|$ have surjective differential, since their zero locus is smooth of codimension three. Modulo the hyperplane directions, this differential is
\[
 \Ext^1(B,B)=\Hom(B,\OO_Z)\longrightarrow
 H^1(B(H))=\frac{H^0(\OO_Z(H))}{\operatorname{im}H^0(\OO_X(H))},
 \qquad b\longmapsto[b(s)] .
\]
Indeed $[b(s)]$ is the obstruction to lifting $s:A\to B$ to the extension specified by $b$. It is the pullback $s^*b$.

We record the same chase for any incidence triangle $A\xrightarrow{s}B\xrightarrow{p}F\xrightarrow{a}A[1]$ on $X$, with $A=\OO_X(-H)$. Serre duality and the two long exact Ext sequences give
\[
\begin{tikzcd}[column sep=small]
 \Ext^1(B,B)\ar[r,"s^*"]\ar[d,"p_*"']&
 \Ext^1(A,B)\ar[d,"\sim","p_*"']&\\
 \Ext^1(B,F)\ar[r,"s^*"]&
 \Ext^1(A,F)\ar[r,"\partial"]&\Ext^2(F,F).
\end{tikzcd}
\]
The right vertical map is invertible since $\Ext^{1,2}(A,A)=0$. If the top arrow is onto, $\partial=0$, and therefore
\[
 \partial^\vee=\pm a_*:\Ext^1(F,F)\longrightarrow \Ext^2(F,A),
 \qquad \xi\longmapsto\pm a[1]\circ\xi
\]
vanishes. This holds in particular at each point of $\mathcal I$.

\emph{Recovering the ideal from a deformation of the sheaf.}
Let $R=\CC[\epsilon]/(\epsilon^2)$ and let $F_R$ represent $\xi$. Lifting a chain representative of $a$ gives obstruction $a[1]\circ\xi=0$ and ambiguity $\epsilon\Ext^1(F,A)$. By duality and the ideal sequence,
\[
 \Ext^1(F,A)\simeq H^2(I_{Z/T}(H))^\vee
     \simeq H^2(\OO_T(H))^\vee\simeq\CC,\qquad \omega_T=\OO_T(H).
\]
The extension is nonsplit since $B$ is torsion-free and $F$ is torsion. Units $1+\epsilon c$ therefore identify its lifts. The inclusion $B_R\to\OO_{X_R}$ also lifts uniquely up to units by $(\ext^0,\ext^1)(B,\OO_X)=(1,0)$. It is injective with flat cokernel by the fiber criterion. The recovered flag is
\[
 A_R\lhook\joinrel\longrightarrow B_R=I_{Z_R/X_R}
     \lhook\joinrel\longrightarrow\OO_{X_R},\qquad
 T_R=V(A_R\longrightarrow\OO_{X_R}),\qquad
 F_R=B_R/A_R .
\]
Thus $T_{(Z,T)}\mathcal I\cong\Ext^1(F,F)$. The same recovery on closed points and properness give a finite geometrically injective morphism $\mathcal I\to M_X$. At a closed image point $y$,
\[
 10=\dim\mathcal I\le\dim\OO_{M_X,y}\le\dim T_yM_X=10.
\]
The local ring is regular. The closed image contains its component locally and is open and closed. The finite map is birational in characteristic zero and its image is normal, hence it is an isomorphism. Finally $\nu_{\mathcal I}=(-1)^{10}=1$, so the virtual degree is $e(\mathcal I)$ by \cite[Theorem~4.18]{Behrend27}.
\end{proof}

\subsection{The numerical degree}
Write $e$ for compactly supported topological Euler characteristic. Let $C_3\subset X^{[3]}$ consist of length-three subschemes whose linear span is a line. All other triples span a plane. The hyperplanes through a line form $\PP^2$ and those through a plane form $\PP^1$. Stratifying the incidence by the evaluation rank gives
\[
                     e(\mathcal I)=2e(X^{[3]})+e(C_3).
\]
Every double spans a unique line. Over lines $\ell\not\subset X$ the relative Cartier divisor $X\cap\ell$ of degree five gives the residual isomorphism $Z\mapsto(X\cap\ell)-Z$ between length-two and length-three subdivisors. Over each of the $2875$ contained lines the two fibers are $\operatorname{Sym}^2\ell=\PP^2$ and $\operatorname{Sym}^3\ell=\PP^3$. Their Euler numbers differ by one, so
\[
                     e(C_3)=e(X^{[2]})+2875.
\]

Let $b=e(X)=-200$, by $c(T_X)=(1+H)^5/(1+5H)$ and $\int_XH^3=5$. At a smooth point of a threefold the punctual Hilbert schemes of lengths one, two and three have Euler characteristics $1,3,6$. Length two gives the space $\PP^2$ of tangent directions.

At length three the curvilinear locus has an affine two-space of second jets over each tangent direction; the remaining ideals have square-zero maximal ideal and are parametrized by the tangent two-planes. Each of these two strata has Euler characteristic three. The strata of the Hilbert--Chow map now give
\begin{align*}
 e(X^{[2]})&=\binom b2+3b=19\,500,\\
 e(X^{[3]})&=\binom b3+3b(b-1)+6b=-1\,234\,000,\\
 e(\mathcal I)&=2(-1\,234\,000)+19\,500+2875=-2\,445\,625.
\end{align*}

The point-ideal locus is smooth of dimension ten, and each two-line point is smooth of dimension zero. The Behrend function on both is $+1$. By the classification these exhaust the stable moduli. Adding their Euler characteristics proves Theorem~\ref{th:seven-total}.

This coefficient is displayed in \cite[(3.10)]{GSY}. In the geometric count of Gaiotto and Yin \cite[\S2.1]{GY} the configurations on a contained line led to a discrepancy of $2875$. Here each contained line contributes $e(\PP^3)-e(\PP^2)=1$ to the residual-divisor comparison. A formula in terms of curve-counting invariants appears in \cite[(5.9)]{AFKPS}.
\subsection{The stable sheaves in the eighth Fourier position}
Let $X\subset\PP^4$ be a general quintic, $H=c_1(\OO_X(1))$, and let $T$ be a hyperplane section. Write $h=H|_T$. In the representative $(h\cdot c_1(G),\chi(G))=(-2,3)$ of the charge, $F=i_*G$ has Hilbert polynomial
\[
              P_F(m)=\frac52m(m-1)-2m+3 .
\]
Thus $d=H\ch_2(F)=-9/2$, $n=\int_X\ch_3(F)=-7/6$, and the Fourier position is $(\mu,x)=(2,133/120)$.

\begin{proposition}[The three types of sheaves]\label{prop:row8-exhaustion}
The stable sheaves of the indicated character on $X$ consist of precisely the following three disjoint loci.
\begin{enumerate}
\item $i_*I_{W/T}$, where $W\subset T$ contains a smooth conic $K$ and $I_{K/T}/I_{W/T}$ has length one.
\item $i_*\OO_T(-L_1-L_2)$ for an unordered pair of distinct lines, on the section $T=X\cap\langle L_1,L_2\rangle$.
\item $i_*\OO_T(C-h)$ for a smooth twisted cubic $C\subset X$, on the section $T=X\cap\langle C\rangle$.
\end{enumerate}
In the last two cases $T$ is smooth. The first case includes singular members of the pencil of hyperplane sections containing $K$.
\end{proposition}

\begin{proof}
As in Proposition~\ref{prop:row7-hulls}, purity and the primitivity of $H$ give $F=i_*G$ on an integral normal $T$, with $G$ torsion free of rank one. Adjunction gives $\omega_T=\OO_T(h)$.

\emph{A morphism to $\OO_T$.}
Every nonzero map $G\to\OO_T$ embeds $G$ as $I_{W/T}$. Its quotient has polynomial $2m+2$. Dividing $\OO_W$ by its maximal zero-dimensional subsheaf, of length $b$, gives a locally Cohen--Macaulay curve $B$ of degree two and genus $b-1$. If $B$ is planar then $p_a(B)=0$ and $b=1$. For nonplanar $B$, \cite[Theorem~9.1(1)--(2)]{HS} gives $p_a(B)\le0$, with equality only for a plane conic. Thus $b=0$ and $p_a(B)=-1$.

On a general quintic the lines are disjoint and infinitesimally rigid. Indeed quintics generate $\operatorname{Sym}^5 U^\vee$ on $\operatorname{Gr}(2,5)$, whose rank is its base dimension, and two meeting lines have ten parameters and impose eleven independent conditions. Every line $L\subset X$ consequently has $N_{L/X}=\OO_{\PP^1}(-1)^{\oplus 2}$ by adjunction and $H^0(N_{L/X})=0$.

A locally Cohen--Macaulay double line has $I_L^2\OO_B=0$: the left side vanishes generically and $\OO_B$ has no finite submodule. The kernel of $\OO_B\to\OO_L$ is a rank-one torsion-free sheaf on $L$, hence $\OO_L(a)$, and is a quotient of $I_{L/X}/I_{L/X}^2=\OO_L(1)^{\oplus 2}$. Thus $a\ge1$ and $p_a(B)=-a-1\le-2$.

This excludes double lines in both cases. A reduced planar conic in $X$ must therefore be smooth, while the nonplanar $B$ is the disjoint union of two lines. In the latter case its spanned hyperplane cuts a smooth section for general $X$, as in the two-line calculation above, so $I_{B/T}=\OO_T(-L_1-L_2)$.

\emph{A section of $G(h)$.}
Suppose $\Hom_T(G,\OO_T)=0$. Serre duality gives $H^2(G(h))=0$, while $P_G(1)=1$ gives $h^0(G(h))\ge1$. Let $A=G(h)^{**}$, and choose a nonzero section of $G(h)$. Its composite with $G(h)\hookrightarrow A$ identifies $A$ with $\OO_T(C)$ for an effective Weil divisor $C$. The corresponding divisorial ideal is $I_{C/T}=A^\vee$. The ideal is reflexive and the quotient is Cohen--Macaulay of pure dimension one by the depth lemma. Its degree is $-2+5=3$.

Let $b=\operatorname{length}(A/G(h))$. On the Gorenstein surface $T$ the maximal Cohen--Macaulay module $A$ has $\mathcal Ext_T^{>0}(A,\OO_T)=0$. Serre duality and the ideal sequence give
\[
 3+b=\chi(A(-h))=\chi(I_{C/T}(2h))
     =10-(6+1-p_a(C))=3+p_a(C).
\]
If $C$ is planar, residual intersection in the plane section gives a plane conic $K$ with $C+K\sim h$ as Weil divisors. The conic is smooth by the preceding argument and $A(-h)=\OO_T(-K)=I_{K/T}$, which embeds $G$ in $\OO_T$. In the present case $C$ is therefore nonplanar.

The genus bound of \cite[Theorem~9.1]{HS} gives $p_a(C)\le0$. Hence $b=p_a(C)=0$. The same theorem bounds the Rao function by zero in every degree. Thus $C$ is arithmetically Cohen--Macaulay. In its span $\PP^3$, its Artinian reduction has $h$-vector $(1,2)$, as it has length three and $C$ is nondegenerate. The Hilbert--Burch resolution is
\[
0\longrightarrow \OO_{\PP^3}(-3)^2
 \longrightarrow \OO_{\PP^3}(-2)^3
 \longrightarrow I_{C/\PP^3}\longrightarrow0.
\]
Indeed the $h$-vector gives regularity two for $I_C$ and three generators in degree two. Its Hilbert series determines the two syzygies in degree three.

The height-two maximal minors of $3\times2$ matrices of linear forms give a flat family over an irreducible open in matrix space, by this resolution. Its general member is a smooth twisted cubic. Allowing its span to move gives dimension
\[
 4+\dim\operatorname{PGL}_4-\dim\operatorname{PGL}_2=16 .
\]
The nonsmooth members in its Hilbert-scheme closure have dimension at most $15$. For these ACM curves the resolution gives
\[
 H^1(\PP^3,I_C(5))=0,\qquad h^0(C,\OO_C(5))=3\cdot5+1=16 .
\]
Thus the nonsmooth incidence has codimension at least $16-15=1$ in the space of quintics.

For a fixed smooth twisted cubic, quintic equations in its span generate $I_C(5)$ and $N_{C/\PP^3}^\vee(5)$ by the same resolution. Thus the general containing quintic surface in that span is smooth, by Bertini off $C$ and by the rank-two conormal evaluation on $C$. The incidence of smooth twisted cubics and containing quintic threefolds is an irreducible projective bundle over the smooth-cubic parameter space.

It is generically finite over quintics. A proper closed condition on its fibers over cubics has image of positive codimension in the quintic parameter space. This proves the required smoothness of $T$ simultaneously for the finitely many cubics on a general $X$.

Conversely for $A=C-h$ and $A=-L_1-L_2$, adjunction gives
\[
 hA=-2,\qquad A^2=-6,\qquad \chi(\OO_T(A))=5+\tfrac12 A(A-h)=3 .
\]
Here $C^2=-5$ and $hC=3$. All three types are stable by integral rank-one support. Since $C$ spans its hyperplane, $\Hom_T(\OO_T(C-h),\OO_T)=H^0(I_{C/T}(h))=0$. The other two types have hull quotients of lengths one and zero, respectively.
\end{proof}

\begin{proposition}[The twisted-cubic contribution]\label{prop:row8-cubic}
The sheaves $i_*\OO_T(C-h)$ in Proposition~\ref{prop:row8-exhaustion} contribute $317\,206\,375$ to the smooth-quintic virtual degree.
\end{proposition}
\begin{proof}
The smooth cubic parameter space is smooth of dimension sixteen, and quintics generate the rank-sixteen bundle with fiber $H^0(\OO_C(5))$. On a general quintic their zero scheme is reduced, and its length is $317\,206\,375$ by \cite[Theorem~1.1]{ESrow}. In particular $H^0(C,N_{C/X})=0$. The span $\langle C\rangle$ is the unique containing hyperplane. In the sequence
\[
 0\longrightarrow\OO_T\longrightarrow\OO_T(C)
   \longrightarrow \OO_{\PP^1}(-5)\longrightarrow0
\]
the connecting map $H^1(\OO_{\PP^1}(-5))\to H^2(\OO_T)$ is dual to $H^0(T,\OO_T(h))\to H^0(C,\OO_C(h))$. This restriction is an isomorphism since the four cubics parametrizing $C$ form a basis for $H^0(\PP^1,\OO(3))$. Hence $h^0(\OO_T(C))=1$ and $H^1(\OO_T(C))=0$.

A first-order deformation of the pushforward of a line bundle on $T$ deforms its Cartier support and the line bundle there, by its local one-generator, one-relation presentation. Twisting by $H$ and lifting the section of $\OO_T(C)$ deforms the curve in $X$. Since $H^0(N_{C/X})=0$, this curve is constant.

Its deformed hyperplane is also constant since $C$ spans $\PP^3$. With the support fixed the line-bundle tangent space is $H^1(\OO_T)=0$. Therefore each sheaf gives a reduced isolated moduli point and contributes $+1$. The uniquely determined cubic section distinguishes the sheaves belonging to different cubics on the same surface. This proves the assertion.
\end{proof}

\begin{proposition}[A conic with a quotient of length one]\label{prop:split28-conic1}
Let $X$ be a general quintic and let $C\subset X$ be one of its smooth conics. Write $\Pi=\langle C\rangle$ and $X\cap\Pi=C\cup D$. Let $B=\operatorname{Bl}_C X$ and let $D'\subset B$ be the strict transform of $D$.
The surface sheaves
\[
 i_*I_{W/T},\qquad T=X\cap P,\quad C\subset W\subset T,\quad
                    \length(I_{C/X}/I_{W/X})=1
\]
form a smooth open and closed subscheme of the stable sheaf moduli of character $(e,\chi)=(-2,3)$. It is isomorphic to $\operatorname{Bl}_{D'} B$ and has virtual degree $198$. Summing over all conics gives $120\,631\,500$.
\end{proposition}

\begin{proof}
\emph{The Hilbert scheme of $W$ in the threefold.}
The conics on $X$ are isolated and reduced by Lemma~\ref{lem:wholeconic-18}. In particular $H^0(C,N_{C/X})=0$. For a fixed smooth conic the restriction of quintics containing $C$ to its plane fills the polynomials $qf_3$, where $q$ defines $C$ and $f_3$ is an arbitrary cubic.

Thus for a general equation containing $C$, the residual $D$ is a smooth plane cubic meeting $C$ transversely in six points. The incidence of smooth conics with quintics containing them is irreducible and generically finite over quintics. The closed bad locus has smaller dimension, so the preceding assertions hold for every conic on a general $X$.

A length-one quotient of $I_C$ is a line bundle on the graph of its support. Since $C\subset X$ is a regular immersion of codimension two,
\[
 \operatorname{Quot}^1_X(I_C)=\PP_X(I_C)
       =\operatorname{Proj}_X\operatorname{Rees}(I_C)=B
\]
by the regular-sequence calculation in Proposition~\ref{prop:linepoint-blowup-d}. Its kernel is $I_W$. Off $C$ the extra point moves freely. For a quotient supported at $p\in C$, choose parameters $(x,y,z)$ with $I_C=(x,y)$ and $q(x)=1$, $q(y)=0$. Then
\[
 I_W=(y,x^2,xz),\qquad
 \OO_W=\CC[[z]]\oplus\CC\bar x,\qquad z\bar x=\bar x^2=0.
\]
Given $b:I_W\to\OO_W$, put $b(x^2)=f(z)+f_0\bar x$ and $b(xz)=g(z)+g_0\bar x$. The syzygy gives $zf(z)=g(0)\bar x$, so $f=0$ and $g(0)=0$. The restriction of $b$ to the punctured curve therefore extends in the normal bundle by $x\mapsto g(z)/z$, $y\mapsto b(y)\bmod\bar x$. The extension is unique and regular since $C$ is smooth. Thus it is zero by $H^0(N_{C/X})=0$. The ideal sequence now gives
\[
 \Hom_X(I_W,\OO_W)=\Hom_X(I_W,k(p)),\qquad
 \dim\Hom_X(I_W,k(p))=3 .
\]
The proper geometrically injective map $B\to\operatorname{Hilb}X$ is finite. At an image point the dimension is at least three and the tangent dimension equals three. Thus the image is open and closed, regular, and the finite birational map to it is an isomorphism.

\emph{Imposing a hyperplane.}
Let $\pi:B\to X$ have exceptional divisor $E$, and let $V=H^0(I_{C/X}(H))$, of dimension two. Every hyperplane containing $W$ contains $\Pi$. Pulling back its equation and dividing by the exceptional divisor gives the pencil
\[
 V\otimes\OO_B\longrightarrow \pi^*\OO_X(H)(-E).
\]
The scheme of pairs $(W,P)$ is the incidence of zero evaluations in $B\times\mathbb P(V)$, by the tautological quotient parametrized by $B$. The sections generate the ideal of $D'$ after a local trivialization. Away from $C\cap D$ this follows from the base scheme $C\cup D$. At a transverse intersection the completed local equations in $X$ are
\[
 I_C=(x,z),\qquad I_D=(y,z),\qquad I_{C\cup D}=(xy,z).
\]
Indeed the nodal plane section is contained in a smooth local surface obtained by taking a nonsingular linear combination of its two defining equations; the two branches are transverse in that surface. In the chart $z=xu$ of $B$ the residual ideal is $(y,u)$, and in the chart $x=zv$ it is the unit ideal.

Thus the base is the smooth curve $D'\simeq D$. The incidence of the pencil is the blowup $Y_C=\operatorname{Bl}_{D'}B$, since locally the incidence equation for a pair of regular generators is the Rees relation after exchanging the two pencil coordinates. In particular $Y_C$ is smooth of dimension three and
\[
 e(Y_C)=e(X)+e(C)+e(D')=-200+2+0=-198 .
\]

\emph{Deformations as sheaves on $X$.}
Set $A=\OO_X(-H)$, $J=I_{W/X}$, $F=i_*I_{W/T}$. Since $T$ is integral and $I_{W/T}$ is torsion-free of rank one, $F$ is stable. The defining sequence is
\[
             A\xrightarrow{s}J\xrightarrow{p}F\xrightarrow{a}A[1].
\]
The sequence $0\to k(p)\to\OO_W\to\OO_C\to0$ gives $H^1(\OO_W)=H^1(\OO_W(H))=0$. By Serre duality and $H^1(\OO_X)=H^2(\OO_X)=0$ one has $\Ext^1_X(J,\OO_X)=0$ and $\Hom_X(J,\OO_X)=\CC$. In particular $T_W B=\Hom_X(J,\OO_W)=\Ext^1_X(J,J)$. Furthermore
\[
 \Ext^1_X(A,J)=H^1(J(H))=
 \operatorname{coker}\big(V\longrightarrow H^0((I_C/J)(H))\big),
\]
because $H^1(I_C(H))=0$ by surjectivity of restriction of linear forms to the plane conic.

The map $s^*:\Ext^1_X(J,J)\to\Ext^1_X(A,J)$ is surjective. Off the base $D'$ its target vanishes. On $D'$ its target is one-dimensional, and the differential at $W$ of the evaluation of $s\in V$ on the tautological quotient is nonzero. Indeed the two generating sections cut out $D'$ transversely, so the differential of any nonzero linear combination is nonzero.

To identify the differential with $s^*$, represent a Hilbert tangent by $b:J\to\OO_W$. The obstruction to lifting $s$ is the connecting image of $b(s)$ from the sequence $0\to J(H)\to\OO_X(H)\to\OO_W(H)\to0$. Since the family keeps $C$ fixed, this lies in the stated cokernel of evaluation at the quotient. This description is valid on singular hyperplane sections.

Apply the Ext square of Proposition~\ref{prop:row7-point-component} to $B=J$. Its top map is onto by the preceding evaluation calculation, so $a[1]\circ\xi=0$ for every $\xi\in\Ext^1(F,F)$. The remaining groups in the inverse construction are
\[
 \Ext^1(F,A)\simeq H^2(I_{W/T}(H))^\vee=\CC,\qquad
 (\ext^0,\ext^1)(J,\OO_X)=(1,0).
\]
The first equality uses $H^1(\OO_W(H))=0$ and $\omega_T=\OO_T(H)$. The extension and embedded ideal therefore lift to first order up to units, recovering the hyperplane by composition. Thus $T_{(W,P)}Y_C=\Ext^1(F,F)$.

The same one-dimensional extension space recovers $(W,P)$ from the closed point $F$, and the unmixed part of $W$ recovers $C$. Hence the proper map from $Y_C$ is finite and geometrically injective. The local rings at its image have dimension at least three and embedding dimension three.

They are regular, so the image is open and closed and the finite birational map is an isomorphism. Its Behrend function is $(-1)^3$. The degree of this component is $198$ by \cite[Theorem~4.18]{Behrend27}. Lemma~\ref{lem:threeconics-18} gives $258\,200+187\,850+163\,200=609\,250$ conics. The product gives the asserted contribution.
\end{proof}

\begin{theorem}[The full conic coefficient]\label{th:eight-total}
For a general smooth quintic,
\[
 [q^{133/120}]h^{\mathrm{DT}}_{1,2}=[q^{133/120}]h^{\mathrm{DT}}_{1,3}
       =441\,969\,250 .
\]
\end{theorem}
\begin{proof}
By Proposition~\ref{prop:row8-exhaustion} the three loci exhaust residue $2$ in this character. Their contributions are the conic degree of Proposition~\ref{prop:split28-conic1}, the two-line degree of Proposition~\ref{prop:twoskew-component} and the cubic degree of Proposition~\ref{prop:row8-cubic}. Additivity gives
\[
 120\,631\,500+4\,131\,375+317\,206\,375=441\,969\,250 .
\]
Theorem~\ref{thm:rank-zero-charge-conjugation} identifies the opposite residue. This equals the coefficient of the quintic elliptic genus in \cite[(3.10)]{GSY}.
\end{proof}
\subsection{A line with two point modifications and the elliptic quartics}
In charge $(\mu,x)=(1,157/120)$ the two types of sheaf arise from a line with a quotient of length two (possibly noncyclic) and from a smooth elliptic quartic.
\begin{proposition}[Line defects and elliptic quartics]\label{prop:s28-row9hulls}
In the primitive charge $(h\cdot c_1(G),\chi(G))=(-1,2)$ on a general quintic, every stable sheaf $F=i_*G$ belongs to one of the following disjoint families.
\begin{enumerate}
\item $i_*I_{W/T}$, where $I_{L/T}/I_{W/T}$ has length two for a line $L\subset X$;
\item $i_*\OO_T(C-h)$, where $C\subset X$ is a smooth elliptic quartic.
\end{enumerate}
Each sheaf of the second kind is an isolated reduced point, and their total virtual degree is $3\,718\,024\,750$.
\end{proposition}
\begin{proof}
The support $T$ is integral and normal as in Proposition~\ref{prop:row8-exhaustion}. A nonzero morphism $G\to\OO_T$ embeds $G$ as an ideal $I_{W/T}$. The quotient has polynomial $m+3$. Its unmixed part has degree one, hence is a line, and the remaining length is two.

If $\Hom_T(G,\OO_T)=0$, $\chi(G(h))=1$ and Serre duality gives a nonzero section of $G(h)$. As in the proof of Proposition~\ref{prop:row8-exhaustion}, the hull of $G(h)$ is $\OO_T(C)$ for a locally Cohen--Macaulay effective curve of degree four. If $b$ is the length of the hull quotient, then
\[
  2+b=\chi(I_{C/T}(2h))=10-(8+1-p_a(C))=1+p_a(C).
\]
If $C$ is a plane curve, its residual line in the plane section identifies the untwisted hull with $I_{L/T}$ and supplies a morphism to $\OO_T$. For nonplanar $C$, \cite[Theorem~9.1]{HS} gives $p_a(C)\le1$, so $b=0$, $p_a(C)=1$, and its Rao module vanishes. The Artinian reduction of the homogeneous coordinate ring has $h$-vector $(1,2,1)$; consequently the ideal contains two independent quadrics.

Two quadrics without a common factor cut out the entire curve, by degree and Cohen--Macaulayness. With a common plane factor $p$, write the quadrics as $p\ell_1,p\ell_2$. The residual scheme with ideal $I_C:p$ is nonempty since $C$ is nonplanar. Its structure sheaf, up to a twist, is a submodule of $\OO_C$, so is pure, and the ideal contains $\ell_1,\ell_2$. Thus the residual is a line.

The hyperplane residual exact sequence gives a plane subscheme $C\cap V(p)$ of degree three. Its unmixed part is a plane cubic $D\subset C$. On a general quintic every plane cubic is smooth: plane cubics have $6+9=15$ parameters and impose fifteen independent quintic conditions, whereas the singular locus in this parameter space has smaller dimension.

The effective cycle $[C]-[D]$ has degree one, hence is a line $L$. Both components have multiplicity one, so purity kills the zero-dimensional nilradical. Thus $C=D\cup L$ is reduced, and nonplanarity gives $L\not\subset\langle D\rangle$. If $L$ meets $D$, the parameter count and the number of quintic conditions are
\[
 6+9+1+3=19,\qquad h^0(\OO_D(5))+h^0(\OO_L(5))-1=15+6-1=20 .
\]
The conditions are independent: after restriction to $D$, multiply a quartic on $L$ by a hyperplane through $\langle D\rangle$ to prescribe the remaining values on $L$. Hence on a general quintic $D\cap L=\varnothing$. This gives $p_a(D\sqcup L)=0$, contrary to $p_a(C)=1$.

Complete intersections of two quadrics in a hyperplane belong to a smooth parameter space of dimension $4+\dim\operatorname{Gr}(2,10)=20$. The singular members have dimension at most nineteen, and every member imposes twenty independent conditions on quintics by the Koszul resolution.

General $X$ thus contains finitely many reduced Hilbert points of smooth elliptic quartics and no singular complete intersections. The same normal generation and Bertini argument as for the twisted cubics proves that each $T=X\cap\langle C\rangle$ is smooth. In particular the smooth curve count is $3\,718\,024\,750$ by \cite[Theorem~1.3]{ESrow}.

For each $C$ adjunction gives $C^2=-4$, $Ch=4$. Its line bundle $\OO_T(C-h)$ has degree $-1$ and Euler characteristic two, and has no nonzero map to $\OO_T$ since the curve spans the hyperplane. The restriction $H^0(\OO_T(h))\to H^0(\OO_C(h))$ is an isomorphism of four-dimensional spaces. In the divisor sequence this gives $h^0(\OO_T(C))=1$ and $H^1(\OO_T(C))=0$.

On twisting a first-order lift by $h$, the section lifts and cuts out an embedded first-order deformation of $C$. Embedded rigidity of $C$ fixes its hyperplane, and $H^1(\OO_T)=0$ fixes the bundle. Distinct quartics are recovered from the unique section; each gives an isolated reduced point and intrinsic virtual degree $+1$.
\end{proof}
\begin{lemma}[Two points on a residual plane quartic]\label{lem:split28-quartic-net}
Let $L$ be a line on the general quintic $X$. Projection from $L$ gives the family $\rho:B=\operatorname{Bl}_LX\to\PP^2$ of residual plane quartics. Its relative Hilbert scheme of two points is smooth and
\[
                  e\bigl(\operatorname{Hilb}^2(B/\PP^2)\bigr)=1575 .
\]
\end{lemma}
\begin{proof}
Write $L=\PP(A)$ for $\dim A=2$. The universal plane containing $L$ is $\PP(U)\to K=\PP^2$ with $U=(A\otimes\OO_K)\oplus\OO_K(-1)$ after a choice of complement to $A$. Projective bundles in this calculation parametrize lines. On $\PP(U)$ the pullback of the quintic equation is divisible by the section cutting out $\PP(A)\times K$, of class $\OO_{\PP(U)}(1)\otimes\OO_K(-1)$.

The residual equation thus has class $\OO_{\PP(U)}(4)\otimes\OO_K(1)$. It defines the graph of projection after blowing up $L$: indeed $\PP(U)=\operatorname{Bl}_L\PP^4$, and this is the strict transform of the smooth hypersurface containing $L$ to order one. Its fibers are the residual quartics, with their scheme structures.

All residual equations in this linear system arise. Expanding the quintic in the three coordinates complementary to $A$, their terms of degree $j+1$, $0\le j\le4$, supply the entire spaces
\[
          \operatorname{Sym}^{4-j}A^\vee\otimes H^0(K,\OO_K(j+1)).
\]
For general $X$ the result holds for every line on $X$, by irreducibility and generic finiteness of the incidence of lines and quintics.

Each length-two subscheme in a plane spans a unique line, also in families by the evaluation map on linear forms. Write $G=\operatorname{Gr}_K(2,U)$, with tautological subbundle $S$, and $P=\PP_G(\operatorname{Sym}^2 S^\vee)$. A point of $P$ specifies a relative line and a quadratic divisor on it. Division by its quadratic equation gives the exact sequence of vector bundles
\[
 0\longrightarrow \OO_P(-1)\otimes\operatorname{Sym}^2 S^\vee
 \longrightarrow\operatorname{Sym}^4 S^\vee\longrightarrow V\longrightarrow 0 .
\]
The relative Hilbert scheme is the zero scheme of the residual quartic in $E=V\otimes\OO_K(1)$. The equations generate every fiber of $E$: for each plane they include all quartics, which restrict onto any length-two divisor. Bertini gives a smooth fourfold, since $P$ is smooth of dimension six and $E$ has rank two. Its Euler characteristic is
\[
                   \int_P c_2(E)\,[c(TP)c(E)^{-1}]_4 .
\]

For reproducibility, choose diagonal coordinate weights $w_0,\ldots,w_4=0,1,2,3,4$ on $\CC^5$, with $A$ spanned by the first two coordinates. At a fixed point of $K$ indexed by $i=2,3,4$, choose an unordered pair $a,b$ from $\{w_0,w_1,w_i\}$ for the weights of $S$, and call the remaining weight $r$. Put
\[
  D=\{-2a,-a-b,-2b\},\qquad Q=\{-(4-j)a-jb:0\le j\le4\}.
\]
A fixed point over this $S$ is a choice of $t\in D$. Its tangent and bundle weight multisets are
\begin{align*}
 T&=\{w_j-w_i:2\le j\le4,\ j\ne i\}\cup\{r-a,r-b\}
                           \cup\{v-t:v\in D,\ v\ne t\},\\
 N&=\{u-w_i:u\in Q\setminus(t+D)\}.
\end{align*}
The contribution at this fixed point is
\[
   \frac{\prod_{n\in N}n}{\prod_{v\in T}v}\,
               [z^4]\,\frac{\prod_{v\in T}(1+vz)}
                                    {\prod_{n\in N}(1+nz)} .
\]
Summing first over the three pairs $a,b$ and the three choices of $t$ gives, for $i=2,3,4$, respectively $6846,-24537,19266$. Their sum is $1575$.
\end{proof}
\subsubsection{Two point modifications on a surface through a line}
Fix a line $L$ on a general quintic $X$, put $V=H^0(X,I_{L/X}(H))$, and set $P=\mathbb P(V)\cong\mathbb P^2$. The hyperplane section represented by $s\in P$ is denoted by $T_s$. For this line let
\[
 Q_L=\{(q:I_{L/X}\twoheadrightarrow D,[s]):\ \ell(D)=2,\ q(s)=0\}.
\]
Here $q(s)$ is evaluated in $H^0(D(H))$. Write $B=\ker(q)=I_{W/X}$ and $F=I_{W/T_s}$, regarded as a sheaf on $X$. Over $P$, taking the quotient by $s$ gives a proper morphism $j:Q_L\to M_X(v)$ for the charge of the line ideal with two point defects. We compute the virtual contribution of its image. All operations with point quotients concern surface sheaves through this morphism.

\begin{proposition}[The virtual contribution for a fixed line]\label{prop:row9-line2}
The morphism $j$ identifies $Q_L$ with an open and closed subscheme of the stable surface-sheaf moduli. Its Behrend function equals $+1$ and its virtual degree is $20\,678$.
\end{proposition}

\begin{proof}
\emph{The incidence map from the blowup.}
Let $Y=\operatorname{Bl}_L X$ with exceptional divisor $E$, and let $p:Y\to X$. A hyperplane through $L$ cuts out a divisor $Y_s$ of class $p^*H-E$. Let
\[
 \widetilde Q_L=\{(\zeta,s):s\in P,\ \zeta\subset Y_s,\ \ell(\zeta)=2\}.
\]
The tautological quotient $p^*I_{L/X}\to\mathcal O_Y(-E)$ on $Y$ defines
\[
 I_{L/X}\longrightarrow p_*\mathcal O_{\zeta}(-E).
\]
It is surjective. At a point with scheme-theoretic image of length two, $\zeta$ maps isomorphically to that image and the assertion follows from the tautological quotient. At a point where the image has length one, $\zeta$ is a length-two subscheme of the exceptional $\mathbb P^1$ and evaluation of $\mathcal O_{\mathbb P^1}(1)$ on $\zeta$ is an isomorphism on global sections. Surjectivity in families follows by Nakayama's lemma. The lifted hyperplane equation vanishes on $\zeta$, giving a proper morphism $\widetilde Q_L\to Q_L$.

Call $D$ cyclic if it can be generated locally on $X$ by one element. Its annihilator then cuts out a length-two subscheme $Z\subset X$, and $D$ is an invertible $\mathcal O_Z$-module. The quotient defines a section of $\mathbb P_X(I_{L/X})\to X$ over $Z$. Since $L$ is regularly embedded of codimension two, its ideal is of linear type and this projectivization equals $Y$.

The construction commutes with base change on this locus: locally the chosen generator identifies $D$ with the quotient by its annihilator. Thus the morphism is an isomorphism over the cyclic locus. The same argument without $s$ identifies the corresponding cyclic open of $\operatorname{Quot}^2_X(I_{L/X})$ with an open in $Y^{[2]}$.

A noncyclic length-two quotient is $D=k(x)^2$. It occurs precisely at $x\in L$, where it is the canonical quotient of $I_{L/X}$ by $\mathfrak m_xI_{L/X}$. Its kernel determines a unique $[s]$: the normal first jets give a surjection $V\to I_{L/X}\otimes k(x)(H)$ of ranks three and two. This follows also by restricting linear equations of $L$ in $\mathbb P^4$ to $X$.

Its one-dimensional kernel defines the hyperplane tangent to $X$ at $x$. Hence the noncyclic stratum is $L$. Its preimage is the relative Hilbert square of $E\to L$, with fiber $\operatorname{Hilb}^2\mathbb P^1=\mathbb P^2$. Consequently
\[
                        e(Q_L)=e(\widetilde Q_L)-2e(L).
\]

\emph{Smoothness on the cyclic locus.}
The blowups of the hyperplanes $\PP^3$ along $L$ form a smooth family over $P$. The Hilbert square of each blowup $Y$ is $\operatorname{Bl}_\Delta(Y\times Y)/\mathfrak S_2$, smooth of dimension six. Hence the relative Hilbert square is smooth of dimension $2+6=8$. On it $\widetilde Q_L$ is cut out by the two values of the strict transform of the quintic on the universal double.

On $\operatorname{Bl}_L\mathbb P^4$, the linear system $5H-E$ separates all length-two subschemes. Indeed the graph of projection from $L$ embeds this blowup in $\mathbb P^4\times\mathbb P^2$, and $5H-E$ is the restriction of $\mathcal O(4,1)$. Products of quartic equations with linear equations vanishing on $L$ give all sections needed for separation.

The corresponding evaluation bundle of rank two is therefore generated by quintics through $L$. Bertini gives smoothness and dimension six for $\widetilde Q_L$ for a general such quintic. The exceptional locus in the family of pairs $(X,L)$ has positive codimension; the finite general fibers over the quintics show that the conclusion holds simultaneously for their lines.

\emph{The local model at a noncyclic quotient.}
At such a point $x$ choose two independent normal linear equations $u,v\in V$ and a coordinate $z$ along $L$. In analytic coordinates the ideal is $(u,v)$ with relation $vu-uv=0$. Trivialize the length-two quotient by the ordered basis of images of $u,v$, denoted $e_1,e_2$. In the analytic germ of the Quot scheme the coordinates are the three $2\times2$ multiplication matrices $U,V',Z$ for $u,v,z$, subject to
\[
 [U,V']=[V',Z]=[Z,U]=0,\qquad V'e_1=Ue_2.
\]
Commutation gives the local-ring action; $V'e_1=Ue_2$ is its sole module relation. The prescribed ordered basis removes change of frame. Write
\[
 U=\alpha I+U_0,\quad V'=\beta I+V_0,\quad Z=\gamma I+Z_0,
 \qquad U_0,V_0,Z_0\in\mathfrak{sl}_2 .
\]
The relation is $\beta e_1-\alpha e_2=U_0e_2-V_0e_1$. It eliminates $\alpha,\beta$ and leaves $\gamma$ free. In a basis of $\mathfrak{sl}_2$ put $N=(U_0\,V_0\,Z_0)$. The bracket $\bigwedge^2\mathfrak{sl}_2\xrightarrow{\sim}\mathfrak{sl}_2$ identifies the three commutators with the three pairs of column wedges. Their entries generate the ideal of all $2\times2$ minors of $N$.

The evaluations of $u,v$ in $D(H)$ remain a basis near the chosen point. Thus the section $[s]$ is the unique kernel of the surjective evaluation $V\to H^0(D(H))$, and imposes no further local equations. The analytic germ of $Q_L$ is the product of a smooth curve with $\{\operatorname{rk}N\le1\}\subset\operatorname{Mat}_{3\times3}$. This determinantal scheme is the critical locus of $\det N$.

\emph{Comparison with the surface-sheaf moduli.}
Set $A=\mathcal O_X(-H)$ and write the defining triangle as $A\xrightarrow{s} B\to F\xrightarrow{a}A[1]$. The ideals $B$ deform locally in $\operatorname{Quot}^2_X(I_{L/X})$: by \cite[Theorems~2.1 and 2.2]{Ric22} this Quot scheme is the scheme-theoretic fiber of Hilbert--Chow and the Hilbert scheme is locally analytically its product with the curve-cycle space near $[L]$. That cycle space is a reduced point locally, since the scheme of lines on a general quintic is reduced and zero-dimensional. Moreover
\[
 H^1(W,\mathcal O_W)=H^1(W,\mathcal O_W(H))=0,\quad
 \operatorname{Hom}_X(B,\mathcal O_X)=\mathbb C,\quad
 \operatorname{Ext}^1_X(B,\mathcal O_X)=0 .
\]
Indeed $\mathcal O_W$ is an extension of $\mathcal O_L$ by $D$. The last two assertions follow from its ideal sequence, Serre duality on $X$, and $H^1(\mathcal O_X)=0$. These identities identify the tangent of the Hilbert scheme at $W$ with $\operatorname{Ext}^1_X(B,B)$.

On the cyclic locus the differential of the two equations $q(s)=0$ is onto: their zero locus in the smooth eightfold $\operatorname{Quot}^{2,\rm cyc}_X(I_L)\times P$ is smooth of dimension six. Since $H^1(I_L(H))=0$, dividing out the hyperplane directions gives
\[
 \Ext^1(B,B)\xrightarrow{s^*}
 \Ext^1(A,B)=H^1(B(H))=
 \operatorname{coker}\big(V\longrightarrow H^0(D(H))\big)
 \longrightarrow0 .
\]
For $b:B\to\OO_W$ this sends $b$ to the class of $b(s)$, by the ideal sequence. The square in Proposition~\ref{prop:row7-point-component} gives $a_*:\Ext^1(F,F)\to\Ext^2(F,A)$ equal to zero. Moreover
\[
 \Ext^1(F,A)=H^2(T_s,F(H))^\vee
    =H^2(\OO_{T_s}(H))^\vee=\CC,
\]
using $H^1(\OO_W(H))=0$. These identities and $(\ext^0,\ext^1)(B,\OO_X)=(1,0)$ give the tangent inverse for $(W,T_s)\mapsto F$ exactly as in that proposition.

At a noncyclic quotient the evaluation of $V$ spans all of $H^0(D(H))$. Hence $H^1(B(H))=0$ and the triangle gives
\[
 \operatorname{Ext}^2_X(F,A)=\operatorname{Ext}^1_X(A,F)^\vee=0 .
\]
For a small extension $R'\twoheadrightarrow R$ with kernel $J$, fix a lift $F_{R'}$ and use the two lifting sequences
\[
\begin{aligned}
 \Ext^1_X(F,A)\otimes J &\to \Ext^1_{X_{R'}}(F_{R'},A_{R'})
     \to \Ext^1_{X_R}(F_R,A_R)\to0,\\
 \Hom_X(B,\OO_X)\otimes J &\to \Hom_{X_{R'}}(B_{R'},\OO_{X_{R'}})
     \to \Hom_{X_R}(B_R,\OO_{X_R})\to0.
\end{aligned}
\]
First lift the extension class to obtain $B_{R'}$, then lift its inclusion. The leftmost spaces are copies of $J$. Units $1+J$ account for both ambiguities because the closed extension is nonsplit and the closed inclusion is nonzero. Fiberwise injectivity gives a flat ideal quotient.

The composite from $A_{R'}$ cuts out the hyperplane containing the recovered $W_{R'}$. This scheme contains the fixed $L_{R'}$ by the Hilbert--Chow identification above. Iteration gives the inverse on completed local deformation functors. This applies also to singular $T_s$.

The map $j$ is proper and injective on geometric points: the nonzero extension of $F$ by $A$ recovers $B$, its unique inclusion in $\mathcal O_X$ recovers $W$, and the composite recovers $T_s$. Thus $j$ is finite. At cyclic quotients its target has local dimension at least six and embedding dimension six by the tangent comparison, hence is regular.

The finite degree-one map is an isomorphism on an open neighborhood there by normality. At noncyclic quotients the completed local rings are isomorphic, so $j$ is also a local isomorphism there. Properness makes the image open and closed.

\emph{The sign and the degree.}
On the smooth sixfold locus $\nu=1$. At a noncyclic point the critical chart is $\AA^1\times\operatorname{Mat}_{3\times3}$ with function $\det N$. The Milnor fiber of the homogeneous determinant has the homotopy type of $\mathrm{SL}(3,\CC)$, which retracts onto $\mathrm{SU}(3)$ by polar decomposition (contract the positive Hermitian factor in a Milnor ball). A nonzero invariant vector field on this compact Lie group gives $e(\mathrm{SU}(3))=0$. Thus \cite[\S1.2, equation~(4)]{Behrend27} gives
\[
 \nu=(-1)^{10}\bigl(1-e(\mathrm{MF}_0(\det))\bigr)=1 .
\]

Project from $L$ to obtain $f:Y\to\mathbb P^2$. The surfaces $Y_s$ are inverse images of lines in this plane. For a length-two subscheme of $Y$ whose image under $f$ has length two, there is a unique line containing that image. When the image has length one the lines form a $\mathbb P^1$. Therefore
\[
        e(\widetilde Q_L)=e(Y^{[2]})+e(\operatorname{Hilb}^2(Y/\mathbb P^2)).
\]
Here $e(Y)=e(X)+e(L)=-198$, so
$e(Y^{[2]})=\binom{-198}{2}+3(-198)=19\,107$.
Lemma~\ref{lem:split28-quartic-net} gives the relative Hilbert-square Euler characteristic $1\,575$. Subtracting $2e(L)=4$ and integrating the Behrend function gives $19\,107+1575-4=20\,678$.
\end{proof}

\begin{theorem}[The full coefficient in the line charge with two point modifications]\label{th:nine-total}
For a general smooth quintic,
\[
   [q^{157/120}]h^{\mathrm{DT}}_{1,1}=[q^{157/120}]h^{\mathrm{DT}}_{1,4}
       =3\,777\,474\,000.
\]
\end{theorem}
\begin{proof}
Proposition~\ref{prop:s28-row9hulls} accounts for all stable sheaves of character $(e,\chi)=(-1,2)$. Each of the $2875$ line components has virtual degree $20\,678$ by Proposition~\ref{prop:row9-line2}. The other components are the reduced points associated to smooth elliptic quartics. Thus
\[
   2875\cdot20\,678+3\,718\,024\,750
       =59\,449\,250+3\,718\,024\,750=3\,777\,474\,000.
\]
Charge conjugation (Theorem~\ref{thm:rank-zero-charge-conjugation}) gives the opposite residue. Comparison with \cite[(3.10)]{GSY} gives the asserted match.
\end{proof}
\subsection{The four-point component at the next Fourier exponent}
At $q^{41/24}$ in residue zero there is a closed family of sheaves obtained from four points on a hyperplane section. Its degree provides one summand of the coefficient. The distinct family of two-line sheaves with a point is treated in Corollary~\ref{cor:twoskew-all-lengths}.
\begin{proposition}[Four points in a hyperplane section]\label{prop:split28-fourpoints}
For a general quintic $X$, the sheaves $i_*I_{Z/T}$ with $\length Z=4$ form an open and closed locus of the stable-sheaf moduli. Its virtual degree is $105\,770\,850$.
\end{proposition}
\begin{proof}
Write $J_4=\{(Z,P): Z\subset X\cap P,\ \length Z=4\}$. Suppose each local component of $Z$ has embedding dimension at most two. In a smooth formal threefold its ideal is $(z,I')$, where $I'$ is an ideal in a smooth surface. The surface Hilbert theorem and the free image of $z$ give tangent dimension $2n+n=3n$ for a component of length $n$. It is smoothable on that surface, so the relative Hilbert scheme of four points in the moving hyperplanes is smooth here. Monomials of degree at most three span every quotient algebra of length four. Thus quintics generate the rank-four evaluation bundle, and Bertini gives
\[
  \dim J_4=\dim|H|+3\cdot4-4=12
\]
on this open. As for three points, the surjective differential of evaluation and Calabi--Yau duality for
\[
          A=\OO_X(-H)\longrightarrow I_{Z/X}\longrightarrow F=i_*I_{Z/T},
\]
lifts the unique extension class across every first-order deformation of $F$. Thus the incidence map is an isomorphism onto its image in the sheaf moduli locally at these points.

The remaining local algebra is $\OO_{X,p}/\mathfrak m_p^2$. Its embedding dimension is three, and it spans its projective tangent hyperplane. On the open neighborhood where the restriction of linear forms has rank four, $J_4$ is isomorphic to the corresponding open subset of $X^{[4]}$. There $H^1(I_{Z/X}(H))=0$, hence $\Ext^2_X(F,A)=\Ext^1_X(A,F)^\vee=0$.

On every Artin thickening the extension class of $F$ by $A$ lifts successively, uniquely up to a unit since $\Ext^1_X(F,A)=\CC$. The recovered ideal embeds in $\OO_X$ since $\Ext^1_X(I_{Z/X},\OO_X)=0$. These constructions identify the formal incidence and sheaf moduli, also at the nonsmooth points. Properness and injectivity on geometric points give the asserted open-and-closed immersion, since it is finite and an isomorphism on completed local rings.

At $\mathfrak m_p^2$ the germ of $X^{[4]}$ depends only on the germ of a smooth threefold at a point. Its Behrend function has a constant value $c$ along this copy of $X$. Everywhere else it is $+1$. The degree-zero formula of \cite[Theorem~4.12]{BFpoints22} gives $e(X^{[4]},\nu)=(-1)^4e(X^{[4]})$, whence $(c-1)e(X)=0$. Thus $c=1$ and the full virtual degree on $J_4$ is $e(J_4)$.

Denote by $P_4$ the scheme of pairs $(\Pi,Z)$ where $\Pi$ is a plane and $Z\subset X\cap\Pi$ has length four. Let $A_4$ be the locus of collinear quadruples. For spans of projective dimension $a=1,2,3$, the Euler characteristics of the hyperplane fibers are $(3,2,1)$ and those of the plane fibers are $(3,1,0)$, respectively. Indeed they are $e(\PP^{3-a})$ and $e(\operatorname{Gr}(2-a,4-a))$, with an empty plane fiber when $a=3$. Hence
\[
 (3,2,1)=(1,1,1)+(3,1,0)-(1,0,0),\qquad
 e(J_4)=e(X^{[4]})+e(P_4)-e(A_4).
\] On noncontained lines length-four divisors are residual to a point. The lines through a point of $X$ form a $\PP^3$. Each contained line replaces the fiber $\PP^1$ by $\operatorname{Sym}^4\PP^1=\PP^4$, so $e(A_4)=4e(X)+3\cdot2875=7825$. The punctual Euler formula gives
\[
      e(X^{[4]})=[z^4]\prod_{m\ge1}(1-z^m)^{200m}=56\,923\,950.
\]

On $G=\operatorname{Gr}(3,5)$ let $S$ be tautological and put $P=\operatorname{Hilb}^4(\PP(S)/G)$. This smooth space has dimension $6+8=14$. The rank-four bundle $E=\OO(5)^{[4]}$ is generated by quintics and a general section cuts out $P_4$. Adjunction gives
\[
 e(P_4)=\int_P c_4(E)\,[c(TP)/c(E)]_{10}.
\]
Use ambient coordinate weights $0,1,11,121,1331$.

Choose three of these to form $I$ and choose a partition $\lambda_a$ at each $a\in I$, of total size four. For $I\setminus\{a\}=\{b,c\}$ in increasing order write $u=b-a$, $v=c-a$. Start the tangent multiset $T$ with all weights $d-a$ for $a\in I$, $d\notin I$.

For a box $(j,k)\in\lambda_a$, indexed from zero, write $A,D$ for its arm and leg. Its tangent and tautological weights are respectively
\[
  (A+1)u-Dv,\quad (D+1)v-Au\ ;\qquad -5a-ju-kv .
\]
Adjoin the first two to $T$ and the last to $N$, initially empty. These are the two boundary maps in $\Hom(I_{\lambda_a},\OO_{\lambda_a})$ and the quotient monomial. Localization gives
\[
 \sum_{I,(\lambda_a)}
 \frac{\prod_{n\in N}n}{\prod_{w\in T}w}
    [z^{10}]\frac{\prod_{w\in T}(1+wz)}{\prod_{n\in N}(1+nz)}
    =48\,854\,725 .
\]
Consequently $e(J_4)=56\,923\,950+48\,854\,725-7825=105\,770\,850$.
\end{proof}
\section{The arithmetic type of the two-degeneration correspondence}\label{sec:arithfull}
We compare the quartic and quintic discriminant forms, of orders four and five, and the periods of the quintic hyperplane system.

\subsection{The primitive forms and an explicit complement}
For a smooth quintic surface $T\subset\PP^3$, put $h=c_1(\OO_T(1))$. Adjunction and the normal sequence give
\[
 K_T=h,\quad h^2=5,\quad
 c(T_T)=(1+h)^4/(1+5h)=1-h+11h^2 .
\]
Thus $e(T)=55$, $q(T)=0$ and $p_g(T)=4$. It follows that $b_2=53$ and $(b_2^+,b_2^-)=(9,44)$.

For $\lambda\in H^2(T,\mathbb Z)$ put $b=\lambda\cdot h$ and $\lambda_\perp=\lambda-(b/5)h$. Wu's formula gives $\lambda^2\equiv b\bmod 2$. Hence $h^\perp(-1)$ is even of signature $(44,8)$, and its discriminant form is
\[
  -\lambda_\perp^2/2\equiv b^2/10-b/2
  \equiv 3\mu^2/5\pmod{\mathbb Z},\qquad \mu=-b\bmod 5 .
\]
The negative primitive quartic K3 lattice has signature $(19,2)$ and discriminant form $D_4$. The following lattice models realize the complement at the level of quadratic forms.

\begin{proposition}\label{prop:latticefull}
There is a primitive embedding between the lattice models
\[
 L_4=U^2\oplus E_8^2\oplus\langle4\rangle,
 \qquad L_5=U^8\oplus E_8^4\oplus A_4
\]
with complement
\[
 N=U^6\oplus E_8^2\oplus B_3,
 \qquad B_3=\left(\mathbb Z^3,
 \begin{pmatrix}2&0&1\\0&2&1\\1&1&6\end{pmatrix}\right).
\]
The complement has rank $31$, signature $(25,6)$, and finite quadratic module
\[
 D_N=-D_4\oplus D_5.
\]
The index of $L_4\oplus N$ in $L_5$ is $4$.
\end{proposition}
\begin{proof}
Embed the two $U$ and two $E_8$ summands directly.  Realize
\[
 A_4=\{(x_1,\ldots,x_5)\in\mathbb Z^5:\textstyle\sum x_i=0\}
\]
with its positive dot product.  The primitive vector $v=(1,1,-1,-1,0)$ has square $4$.  Its orthogonal complement is generated by
\[
 a=(1,-1,0,0,0),\quad b=(0,0,1,-1,0),\quad c=(1,0,1,0,-2).
\]
Indeed an orthogonal vector has $x_1+x_2=x_3+x_4=s$ and $x_5=-2s$, and equals $-x_2a-x_4b+sc$.  These pairings give precisely the displayed Gram matrix $G$.  Its determinant and inverse are
\[
 \det G=20,\qquad G^{-1}=\frac1{20}
 \begin{pmatrix}11&1&-2\\1&11&-2\\-2&-2&4\end{pmatrix}.
\]
The first dual basis vector $a^*=(11a+b-2c)/20$ generates the order-twenty discriminant group.  Put $t=5a^*$ and $u=8a^*$.  Their orders are $4$ and $5$, and
\[
 q(t)=\frac{55}{8}\equiv-\frac18,\qquad
 q(u)=\frac{88}{5}\equiv\frac35,\qquad (t,u)=22\equiv0\pmod{\mathbb Z}.
\]
This is the asserted orthogonal discriminant decomposition.  The remaining summands are unimodular.  The rank and signature follow from the direct sum, and
\[
 [L_5:L_4\oplus N]^2=\frac{4\cdot20}{5}=16.
\]
This gives the embedding of abstract lattices.
\end{proof}

\subsection{The vacuum numerator and the period pairing}
The quintic side has its own zero-flux normalization, which can be computed over the entire linear system, including its singular members. Let $P=|H|$, and write $\pi:\mathcal T\to P$ for the universal hyperplane section of $X$.
\begin{proposition}\label{prop:quintic-hodge-normalization}
The rank-four bundle of fiberwise canonical forms $\mathcal H=\pi_*\omega_{\mathcal T/P}$ is isomorphic to $T_P$. Its dual Euler class gives the neutral constant of the primitive numerator:
\[
 \int_P c_4(\mathcal H^\vee)=5
 =[q^0]\big(\eta^{55}h^{\mathrm{DT}}_{1,0}\big).
\]
\end{proposition}
\begin{proof}
On $X\times P$, $\mathcal T$ is the zero locus of the universal section of $\OO_X(H)\boxtimes\OO_P(1)$. Relative adjunction and $\omega_X=\OO_X$ therefore give
\[
0\longrightarrow\OO_{X\times P}
 \longrightarrow\OO_X(H)\boxtimes\OO_P(1)
 \longrightarrow i_*\omega_{\mathcal T/P}\longrightarrow0 .
\]
Pushing forward to $P$ and using $H^1(X,\OO_X)=0$ gives the Euler sequence
\[
 0\longrightarrow\OO_P\longrightarrow H^0(X,\OO_X(H))\otimes\OO_P(1)
 \longrightarrow \mathcal H\simeq T_P\longrightarrow0 .
\]
With $a=c_1(\OO_P(1))$ we obtain $\int_P c_4(\mathcal H^\vee)=\int_P 5a^4=5$.
Theorem~\ref{thm:vacuumvirtual} gives $[q^{-55/24}]h^{\mathrm{DT}}_{1,0}=5$.
The degree-zero hull bound excludes smaller exponents, so multiplication by
$\eta^{55}=q^{55/24}(1+O(q))$ gives the same constant term.

\end{proof}
\subsection{The compact period pairing for quintic sections}\label{sec:quintic-period-comparison-new}
Let $P=|H|$ and let $U\subset P$ parametrize the smooth sections $T$ of a general quintic $X$. For $c\in H^2(T,\mathbb Z)$ write
\[
 a=h\cdot c,\qquad v=c-\tfrac a5 h,\qquad m=-\tfrac12v^2,\qquad
 \mu=-a\bmod 5 .
\]
The Hodge condition on $c$ is the vanishing of its pairing with $H^{2,0}(T)$. We use the period construction of \cite[\S6.2]{FTperiod} for the primitive hyperplane class.

\begin{proposition}\label{prop:quintic-compact-hodge}
The complete hyperplane family of a general quintic threefold defines a numerical period series $\mathcal N_X\in M_{26}(\rho_5)$. For $m>0$ the coefficient of $q^m e_\mu$ is the normalized intersection of a compact four-dimensional period cycle with the Hodge cycle of vectors of norm $2m$ and residue $\mu$ for the sign-reversed primitive lattice. The constant term belongs to the zero coset and is the normalized degree of the top Chern class of the dual Hodge bundle.
\end{proposition}
\begin{proof}
For a fixed flag $p\in P_0\subset\PP^4$, quintics restrict onto the affine five-jets on $P_0$ at $p$. In three local variables the singularity condition $f(p)=df(p)=0$ has codimension four. For the Hessian $B$,
\[
 \operatorname{codim}\{\operatorname{corank}B\ge c\}=c(c+1)/2 .
\]
Corank zero gives $A_1$. In corank one, split off two squares and write the remaining jet as $a_3z^3+a_4z^4+a_5z^5$. Successive pure powers of $z$ vary each coefficient with the preceding ones fixed. Unless all three vanish, the type is $A_2,A_3$ or $A_4$. In corank two a binary cubic with distinct roots gives $D_4$; a repeated root imposes one further condition. The excluded jets thus have codimension at least
\[
 \min\{4+1+3,\ 4+3+1,\ 4+6\}=8>\dim\{(P_0,p):p\in P_0\}=7 .
\]
They occur on no hyperplane section of a general quintic.

At a singular section the vanishing lattices $R_p$ lie in disjoint Milnor balls. Picard--Lefschetz gives a local monodromy group contained in the finite product $\prod_p W(R_p)$, acting trivially on $(\bigoplus_p R_p)^\perp$; see \cite[\S4]{GreerQuinticPeriod}. Choose a neat finite-index subgroup $\Gamma$ of the discriminant kernel in the identity component of the real orthogonal group, and take the connected finite cover of $U$ defined by the inverse image of $\Gamma$ under monodromy. Normalize $P$ in this cover and resolve its boundary to normal crossings. For the resulting map $\widetilde P\to P$ of generic degree $d$, the boundary monodromies are trivial since the corresponding subgroups of $\Gamma$ are finite. The period map therefore extends to $\bar\phi:\widetilde P\to\Gamma\backslash\mathcal D$ by \cite[Theorem~9.5]{GriffithsQuinticPeriod}, where $\mathcal D$ is the period domain.

Let $\Theta_L^{(8)}$ denote the real degree-eight term of Garc\'ia's genus-one theta class in Hodge-cycle normalization, using all discriminant cosets. At $m=0$ one uses the top Chern class of the dual Hodge bundle in this normalization. Set
\begin{equation}\label{eq:quintic-compact-period}
 \mathcal N_X(\tau)=\frac1d\int_{\widetilde P}\bar\phi^*\Theta_L^{(8)}(\tau).
\end{equation}
The class is holomorphic modular of weight $\rk(L)/2$ by \cite[Theorems~1.2, 5.2 and 5.4]{GarciaQuinticPeriod}. Integration is a fixed linear functional, and the projection formula on a common refinement makes the normalized integral independent of the cover and resolution. Here $L=(h^\perp,-\,\cdot)$ has signature $(44,8)$ and rank $52$. The characteristic property of $h=K_T$ gives $c^2\equiv a\bmod2$ and $m\equiv3a^2/5\bmod\mathbb Z$. The finite Fourier phase is $e(-(44-8)/8)=-1$, giving \eqref{eq:Weil5new} in the coset basis.
\end{proof}

\begin{proposition}[The twisted-cubic contributions in all point lengths]\label{prop:quintic-cubic-period}
For every smooth twisted cubic $C\subset X$ let $T=X\cap\langle C\rangle$ and $D=C-h$. The sheaves with hull $\OO_T(kD)$, $k\ge1$, and arbitrary punctual quotient contribute to $\eta^{55}h^{\mathrm{DT}}_1$ the series
\begin{equation}\label{eq:quintic-cubic-period-sector}
 317\,206\,375\sum_{k\ge1} q^{17k^2/5}e_{\,2k\bmod5}.
\end{equation}
The marked Hodge intersections at these smooth sections give the same summand of $\mathcal N_X$.
\end{proposition}
\begin{proof}
By Proposition~\ref{prop:row8-cubic}, $T$ is smooth, $h^0(\OO_T(C))=1$ and $H^1(\OO_T(C))=0$. Hence a deformation extending $\OO_T(D)$ extends the section cutting out $C$ after twisting by $h$. The rigid cubic fixes its span, so the obstruction to extending the hull is injective. Explicitly,
\[
 o_D:H^0(N_{T/X})\longrightarrow H^2(\OO_T),\qquad
 v\longmapsto\operatorname{KS}(v)\mathbin{\lrcorner}c_1(D).
\]
Both spaces have dimension four. This is the derivative of the Hodge section, and $o_{kD}=k\,o_D$. Thus for $k\ne0$ the marked intersection is transverse of multiplicity one.

For $G=\OO_T(kD)\otimes I_Z$ a first-order support displacement must extend $\det G=\OO_T(kD)$ (Corollary~\ref{cor:twoskew-all-lengths}), so it vanishes. Put $n=\length Z$. Riemann--Roch and Serre duality on $T$ give
\[
 \begin{gathered}
 \chi(G,G)=5-2n,\qquad \dim\Hom_T(G,G)=1,\qquad \ext^2_T(G,G)=4,\\
 \ext^1_X(i_*G,i_*G)=\ext^1_T(G,G)=2n .
 \end{gathered}
\]
The proper injective map $T^{[n]}\to M_X$ has image of dimension $2n$, equal to the target tangent dimension along it. As in that corollary, its image is regular and the map identifies an open and closed component of degree $(-1)^{2n}e(T^{[n]})$.

Here $hD=-2$ and $D^2=-6$, giving $-(D+2h/5)^2/2=17/5$. By \cite{GottscheBetti},
\[
 \eta^{55}q^{-55/24}\sum_{n\ge0}e(T^{[n]})q^n
 =\prod_{j\ge1}(1-q^j)^{55}\prod_{j\ge1}(1-q^j)^{-55}=1 .
\]
The norm distinguishes the positive $k$. For fixed $k$, equal degrees exclude a nonzero hyperplane twist; on the same $T$, torsion-freeness of $\Pic(T)$ recovers $D$ from $kD$ and the unique section of $\OO_T(C)$ recovers $C$. Each of the $317\,206\,375$ cubics of \cite[Theorem~1.1]{ESrow} is therefore counted once in \eqref{eq:quintic-cubic-period-sector}.
\end{proof}

\begin{proposition}[The first neutral terms on the discriminant]\label{prop:quintic-two-boundary}
Let $s(T)$ be the number of singular points and let $M(T)$ be their total Milnor number, both zero on a smooth section. With compactly supported Euler integration on $P$,
\begin{equation}\label{eq:quintic-two-boundary}
 \begin{aligned}
 {}[q]\bigl(\eta^{55}h^{\mathrm{DT}}_{1,0}\bigr)
   &=-\int_P M(T)\,de=-1075,\\
 {}[q^2]\bigl(\eta^{55}h^{\mathrm{DT}}_{1,0}\bigr)
   &=\int_P\left(\frac{M(T)^2-3M(T)}2+s(T)\right)de=109\,650 .
 \end{aligned}
\end{equation}
\end{proposition}
\begin{proof}
The Milnor fiber of an isolated surface hypersurface singularity has Euler number $1+\mu_p$. Smoothing gives $u=e(T)=55-M(T)$. Distinct unordered pairs on $T$ have Euler number $\binom u2$. At a smooth or singular point the punctual doubles form respectively a $\PP^1$ or a $\PP^2$ of tangent directions. Thus
\[
 e(T^{[2]})=\binom u2+2(u-s)+3s=\tfrac12(u^2+3u)+s.
\]
By Proposition~\ref{prop:polar-incidence-functors} the three incidences are smooth of even dimension. Their degrees are $\int_P(1,u,e(T^{[2]}))\,de=(5,-800,58\,500)$. Removing the point factor means multiplying by $1-55q+1430q^2+O(q^3)$. With $u=55-M$ this gives
\[
 u-55=-M,\qquad e(T^{[2]})-55u+1430=\tfrac12(M^2-3M)+s .
\]
Integration gives \eqref{eq:quintic-two-boundary}.
\end{proof}

The cluster types in Figure~\ref{fig:quintic-milnor-moments} illustrate the dependence on $s(T)$ at fixed $M(T)$. For example, a section with just two nodes and a section with just one $A_2$ point both have $e(T)=53$.
\begin{figure}[htbp]
\centering
\begin{tikzpicture}[x=1cm,y=1cm,font=\small,
 root/.style={circle,draw=blue!55!black,fill=blue!12,inner sep=0pt,minimum size=5pt},
 bond/.style={line width=.8pt,draw=blue!65!black}]
 \foreach \x/\lab/\mn/\ss/\first/\second in {0/{A_1}/1/1/-1/0,2.55/{2A_1}/2/2/-2/1,5.10/{A_2}/2/1/-2/0,7.65/{D_4}/4/1/-4/3}{
  \node at (\x,1.15) {$\lab$};
  \node at (\x,-.05) {$(\mn,\ss)$};
  \node at (\x,-.60) {$\first$};
  \node at (\x,-1.12) {$\second$};}
 \node[root] at (0,.60) {};
 \node[root] at (2.25,.60) {}; \node[root] at (2.85,.60) {};
 \draw[bond] (4.80,.6)--(5.40,.6);
 \node[root] at (4.80,.6) {}; \node[root] at (5.40,.6) {};
 \draw[bond] (7.24,.49)--(8.06,.49) (7.65,.49)--(7.65,.87);
 \foreach \x/\y in {7.24/.49,7.65/.49,8.06/.49,7.65/.87} \node[root] at (\x,\y) {};
 \node[anchor=east] at (-.62,-.05) {$(M,s)$};
 \node[anchor=east] at (-.62,-.60) {$[q]$};
 \node[anchor=east] at (-.62,-1.12) {$[q^2]$};
 \draw[black!30] (-1.50,-.31)--(8.35,-.31);
\end{tikzpicture}
\caption{Neutral integrands in \eqref{eq:quintic-two-boundary} for four singularity configurations. Disconnected diagrams correspond to distinct singular points. The last two rows are the local values over a hyperplane of the indicated type, before Euler integration over the base.}\label{fig:quintic-milnor-moments}
\end{figure}

Suppose an integral $(1,1)$-class on a smooth section satisfies $h\xi=0$, $\xi^2=-2$. Ampleness gives $h^0(\pm\xi)=0$, so Riemann--Roch and duality give $h^0(h+\xi)\ge\chi(h+\xi)=4$. For a smooth plane section $C\in|h|$ restriction gives
\[
 H^0(T,h+\xi)\hookrightarrow H^0(C,(h+\xi)|_C),\qquad
 \deg((h+\xi)|_C)=5,\quad g(C)=6 .
\]
Clifford's bound is $h^0\le3$ for a special bundle of degree five; a nonspecial bundle of this degree has no sections. This excludes $\xi$. On the crepant resolution of an ADE section the same argument shows that one of $\pm\xi$ must be effective. Its degree against the pullback of $h$ is zero, so it is exceptional. Hence only singular sections contribute to the neutral period coefficient at $q$.

The first number in \eqref{eq:quintic-two-boundary} also has a Chern class expression on the marked discriminant. Let $g:X\to P$ be the Gauss map. Evaluation of relative dualizing forms at the marked point gives $g^*\mathcal H\twoheadrightarrow\OO_X(5H)$. Indeed $\omega_{\mathcal T/P}=\OO_{\mathcal T}(H)\otimes\pi^*\OO_P(1)$ and $g^*\OO_P(1)=\OO_X(4H)$. Put $E=g^*\mathcal H^\vee/\OO_X(-5H)$. By Proposition~\ref{prop:quintic-hodge-normalization},
\[
  c(E)=\frac{(1-4H)^5}{1-5H},\qquad
  \int_X c_3(E)=-215\int_XH^3=-1075.
\]
Identifying this class with the excess intersection of the marked root cycle would require its normal multiplicities along the higher ADE strata.

The full comparison sufficient for primitive S-duality is
\begin{equation}\label{eq:quintic-residual-period}
             \eta^{55}h^{\mathrm{DT}}_1=\mathcal N_X .
\end{equation}
Proposition~\ref{prop:quintic-cubic-period} identifies the twisted-cubic terms in every point length, and \eqref{eq:quintic-two-boundary} integrates two DT coefficients over the full discriminant. To prove \eqref{eq:quintic-residual-period} one must still compare excess Hodge intersections with ambient virtual multiplicities, including those of noninvertible hulls and point defects on singular supports. The implication from that identity is
\[
 \eta^{55}h^{\mathrm{DT}}_1=\mathcal N_X
 \ \Longrightarrow\ h^{\mathrm{DT}}_1\text{ has weight }-3/2
 \ \Longrightarrow\ h^{\mathrm{DT}}_1=Z^{\mathrm{GSY}},
\]
by Proposition~\ref{prop:numeratorsuffices}, the principal part in Proposition~\ref{prop:polar-from-specializations} and negative-weight uniqueness.

Part~III \cite{QuinticThreePerfect} evaluates some of the singular-support terms. If $G$ is locally free along a perfect point defect $Z\subset T$, then the ambient Behrend weights satisfy $\nu(i_*I_ZG)=\nu(i_*G)$. At a free stalk of type $A_{r-1}$ the perfect punctual factor is $\prod_{j\ge1}(1-q^{rj})^{-1}$. At a nonfree stalk this independent factor is absent. For the trivial hull a second comparison applies to a possibly nonperfect ideal:
\[
 H^0(X,H)\twoheadrightarrow H^0(Z,H|_Z)
 \quad\Longrightarrow\quad
 \nu(i_* I_{Z/T})=(-1)^{4-n}\nu_{X^{[n]}}(Z),\qquad n=\length Z .
\]
The remaining weights enter \eqref{eq:quintic-residual-period}.

\subsection{The finite Fourier calculation}
For a finite quadratic module $D$, put $e(t)=e^{2\pi i t}$ and $B_D(x,y)=Q_D(x+y)-Q_D(x)-Q_D(y)$. Our Weil matrices are
\[
 \rho_D(T)e_x=e(Q_D(x))e_x,\qquad
 \rho_D(S)e_x=\frac{e(-\operatorname{sig}(D)/8)}{\sqrt{|D|}}
 \sum_y e(-B_D(x,y))e_y .
\]

Define
\begin{equation}\label{eq:gfinitefull}
 \mathfrak g:\CC[D_4]\otimes\CC[D_N]\to\CC[D_5],\qquad
 \mathfrak g(e_a\otimes e_{b,\mu})=\delta_{a,b}e_\mu.
\end{equation}

\begin{proposition}\label{prop:intertwinerfull}
The map $\mathfrak g$ intertwines $\rho_4\otimes\rho_N$ with $\rho_5$.  In the chosen bases its contraction coefficient is $1$.
\end{proposition}
\begin{proof}
On $e_a\otimes e_{b,\mu}$ both composites with $T$ equal $\delta_{ab}e(3\mu^2/5)e_\mu$, since the additional phase is $e((a^2-b^2)/8)$.

For $S$, the signature phases multiply to $e(-(1+3)/8)=-1$. The character sum on the contracted index gives
\[
 \sum_{x\bmod4}e((b-a)x/4)=4\delta_{ab},\qquad 4/\sqrt{80}=1/\sqrt5.
\]
Consequently $\mathfrak g(\rho_4(S)\otimes\rho_N(S))(e_a\otimes e_{b,\mu})$
equals $-\delta_{ab}\sum_{\nu\bmod5}e(-\mu\nu/5)e_\nu/\sqrt5
=\rho_5(S)\mathfrak g(e_a\otimes e_{b,\mu})$.
The two identities suffice since $S,T$ generate the metaplectic group with their standard lifts.
\end{proof}

For any $\Theta\in M_{31/2}(\rho_N)$, this gives
\begin{equation}\label{eq:arithcontractionfull}
 \mathfrak g(\widehat\Phi\otimes\Theta)_\mu
 =\sum_{a\bmod4}\widehat\Phi_a\Theta_{a,\mu}
 \in M_{26}(\rho_5),
\end{equation}
since $21/2+31/2=26$.

\subsection{The modular operator associated with a geometric kernel}
For a holomorphic kernel $\Theta\in M_{31/2}(\rho_4^\vee\otimes\rho_5)$ whose coefficients are independently realized by the relative correspondence, define
\begin{equation}\label{eq:operatorTfull}
 (T_\Theta u)_\mu=\eta^{-31}\sum_{a\bmod4}u_a\Theta_{a,\mu}.
\end{equation}
Since $\operatorname{wt}\eta^{-31}=-31/2$, $T_\Theta$ preserves weight. Equation~\eqref{eq:gfinitefull} contracts the two finite representation factors, and $\eta^{-31}$ contributes its scalar multiplier. The next section gives the geometric condition for applying this operator to the DT series.

\section{The two virtual correspondences and the modular implication}\label{sec:twodegoperators}
Let $\mathcal B$ denote pairing with the pencil complement and $\mathcal C$ pairing with the quintic complement. Both act on proper quartic virtual pushforwards.

\subsection{Boundary Chow classes and their two numerical evaluations}
For each nonzero boundary sector $\lambda$, an actual proper relative quartic moduli space gives a cycle
\[
 A_\nu=(\operatorname{ev}_Q)_*[M_{Q,\nu}^{\rel}]^{\vir}
 \in A_3(J_\lambda)_{\QQ}.
\]
Here $\nu$ retains the component charge, determinant data, and any insertions.  The complementary $B$ and $W$ moduli give corresponding cycles, and the two degree formulas pair them with $A_\nu$.

One may pass to $H^6(J_\lambda,\QQ)$ for the numerical pairing, since cycle class preserves degree of intersection on the smooth proper sixfold.  For a basis $\eta_1,\ldots,\eta_N$ write
\[
 A_\nu=\sum_i a_{\nu,i}\eta_i,\quad
 C_\gamma=\sum_j c_{\gamma,j}\eta_j,\quad
 B_\delta=\sum_j b_{\delta,j}\eta_j,\quad
 I_{ij}=\int_{J_\lambda}\eta_i\eta_j.
\]
The two evaluations are then the explicit linear functionals
\begin{equation}\label{eq:matrixpairingsfull}
 \langle A_\nu,C_\gamma\rangle
 =\sum_{i,j}a_{\nu,i}I_{ij}c_{\gamma,j},\qquad
 \langle A_\nu,B_\delta\rangle
 =\sum_{i,j}a_{\nu,i}I_{ij}b_{\delta,j}.
\end{equation}
For the elementary cycles of Proposition~\ref{prop:sfintersectionsfull}, $\langle f,f\rangle=0$ and $\langle f,s\rangle=1$. Thus pairing with $f$ alone loses information.

Let $V$ be the span of the common quartic virtual classes, completed in the point charge. The bounded charge decomposition makes the sums coefficientwise finite. Pairing with the $W$ and $B$ states, respectively, gives
\[
 \mathcal B:V\longrightarrow W_4,\qquad \mathcal C:V\longrightarrow W_5 ,
\]
where the subscripts count the normalized charge sectors. A common $V$ requires the relative stability problems on $Q$ to agree or their virtual classes to be related by wall crossing.

Assuming those identifications and the two complete virtual formulas, write
\begin{equation}\label{eq:twovirtualstatesfull}
 K^{\mathrm{raw}}=\mathcal B\mathcal A+R_Y,
 \qquad h^{\mathrm{DT}}_1=\mathcal C\mathcal A+R_X.
\end{equation}
Here $R_X$ and $R_Y$ sum the omitted geometric sectors in their respective degenerations, including the zero-boundary and outer-divisor terms.

\begin{figure}[htbp]
\centering
\begin{tikzpicture}[x=1cm,y=1cm,font=\small,>=Stealth]
\definecolor{qroot}{RGB}{219,232,242}
\definecolor{qpair}{RGB}{223,237,228}
\foreach \y in {0,1.65}{
  \filldraw[fill=qroot,draw=black!65] (0,\y-.32) rectangle (1.8,\y+.32);
  \node at (.88,\y) {$Q$};
  \filldraw[fill=qpair,draw=black!65] (1.8,\y-.32) rectangle (4.8,\y+.32);
  \draw[line width=1.3pt] (1.8,\y-.32)--(1.8,\y+.32);
  \node[anchor=south] at (1.8,\y+.30) {$S$};
}
\node[anchor=east] at (-.25,1.65) {quintic};
\node[anchor=east] at (-.25,0) {quartic pencil};
\node at (3.32,1.65) {$B=\mathrm{Bl}_{\Gamma}\PP^3$};
\node at (3.32,0) {$W=\mathrm{Bl}_{C_*}P$};
\draw[line width=1.4pt,dashed] (4.8,-.32)--(4.8,.32);
\node[anchor=west] at (4.88,0) {$B_0$};
\end{tikzpicture}
\caption{The two evaluations of a relative quartic class. Each solid seam is the K3 where the restrictions are matched. On the pencil side the complement $W$ also carries the outer relative divisor $B_0$ (dashed edge).}
\label{fig:quarticpairings}
\end{figure}
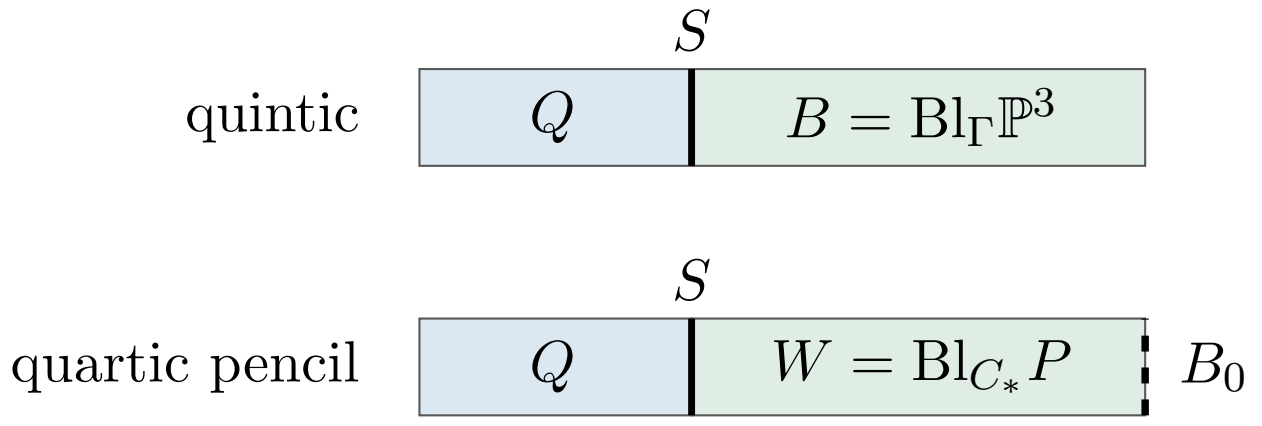

\subsection{The exact factorization condition}
The two evaluations in Figure~\ref{fig:quarticpairings} factor precisely under the following kernel condition.

\begin{lemma}\label{lem:kernelcriterionfull}
For linear maps $\mathcal B:V\to W_4$ and $\mathcal C:V\to W_5$ over a field, there exists a map $T:\operatorname{im}\mathcal B\to W_5$ with $\mathcal C=T\mathcal B$ if and only if
\begin{equation}\label{eq:kernelcriterionfull}
 \ker\mathcal B\subset\ker\mathcal C.
\end{equation}
Such a map is unique on $\operatorname{im}\mathcal B$.  For completed series, continuity, coefficientwise finiteness, and compatibility of the maps on all finite cutoffs are additional requirements.
\end{lemma}
\begin{proof}
If $\mathcal C=T\mathcal B$, then $\mathcal C(\ker\mathcal B)=0$. Conversely, under \eqref{eq:kernelcriterionfull}, the quotient $q:V\to V/\ker\mathcal B$ gives
\[
 \mathcal C=\overline{\mathcal C}\,q,\qquad
 \mathcal B=\overline{\mathcal B}\,q,\qquad
 T=\overline{\mathcal C}\,\overline{\mathcal B}^{-1}
       \quad\text{on }\operatorname{im}\mathcal B .
\]
Surjectivity onto the image gives uniqueness. For compatible finite cutoffs, $(T_n)_n$ induces $\varprojlim_nT_n$ on the inverse limit of their images.
\end{proof}

The diagram to prove geometrically is
\begin{equation}\label{diag:factorizationfull}
\begin{tikzcd}[column sep=huge,row sep=large]
 V\arrow[r,"\mathcal B"]\arrow[dr,"\mathcal C"']
 &W_4\arrow[d,"T_\Theta"]\\
 &W_5.
\end{tikzcd}
\end{equation}
The diagram must commute on every quartic virtual pushforward used in both degenerations.

To make its content coefficientwise explicit, write
\[
 \mathcal B(v)_a=\sum_n\beta_{a,n}(v)q^n,\quad
 \mathcal C(v)_\mu=\sum_m\gamma_{\mu,m}(v)q^m,\quad
 (T_\Theta)_{\mu a}=\sum_r t_{\mu a}(r)q^r.
\]
Then \eqref{diag:factorizationfull} means, for every $v\in V$ and every output degree $m$,
\begin{equation}\label{eq:coefficientfunctionalidentityfull}
 \gamma_{\mu,m}(v)=
 \sum_{a\bmod4}\ \sum_{r+n=m}t_{\mu a}(r)\beta_{a,n}(v).
\end{equation}
For fixed $m$ the sum is finite in the chosen charge completion. In Theorem~\ref{thm:twodegmodularityfull} this identity is required for every $v$.

\subsection{Recovering numerical boundary pairings with insertions}
Fix a smooth projective boundary space $J$ of dimension $2d$. Let
$V\subset N_d(J)_\QQ$ be the span of the quartic virtual
pushforwards in the chosen charges. A complementary class $b$ defines
$\ell_b(A)=\deg_J(A\cdot b)$. Given $w_1,\ldots,w_m$ from specified
insertions on the second complement, define $P(A)_j=\ell_{w_j}(A)$.
Then
\[
 \ker P\subseteq\ker(\ell_b|_V)
 \ \Longleftrightarrow\ 
 \exists\,t_1,\ldots,t_m\in\QQ:\quad
 \ell_b(A)=\sum_{j=1}^m t_jP(A)_j\quad(A\in V).
\]
Indeed $P(A)\mapsto\ell_b(A)$ descends to $\operatorname{im}P$ by
the kernel inclusion and extends linearly to $\QQ^m$. On the Abel
cycles of Proposition~\ref{prop:route19-abel} the line contact needs
one measurement, with coefficient $15$. For the length-two cycles of
Proposition~\ref{prop:mark20-twoprobe} the coefficients are $-50,10$.

For a general primitive sector with ordinary Jacobian evaluation, the
next two theorems give a list of measurements independent of the quartic
cycle. On boundary sectors of general perfect complexes the virtual
class of the proper matched locus in Section~\companionref{sec:boundaryisomproper}
continues to enter the formula.
\subsection{Markman's diagonal for the primitive K3 boundary}
\label{sec:mark20}
Let $\Pic(S)=\mathbb Zh$, $h^2=4$, and fix the ordinary boundary space
$J=J_\chi=M_h(0,h,\chi)$ of Proposition~\companionref{prop:Jacobian}.
We calculate the degree of its diagonal pairing from finitely many insertions
on each relative factor. We use singular cohomology with rational coefficients
and the cycle class of the virtual cycle.

The primitivity of $h$ makes
$v=(0,h,\chi)$ primitive. Its effective one-dimensional support cycle
has a single component of multiplicity one, so semistability and stability
coincide. In particular $h$ is $v$-generic in the sense of
\cite[Definition~5 and Condition~6]{MarkmanIntegral}.

For $x,y\in K^0_{\mathrm{top}}(S)$ write
$\kappa(x,y)=\int_S\ch(x^\vee)\ch(y)\operatorname{td}(S)$.
Choose a basis $x_1,\ldots,x_{24}$ and the topological universal class
$u_\chi$ of \cite[Definition~26 and Lemma~28]{MarkmanIntegral}.
The K3 has $K^1_{\mathrm{top}}(S)=0$, so
\[
 u_\chi=\sum_{a=1}^{24}x_a\boxtimes u_a,
 \qquad u_a\in K^0_{\mathrm{top}}(J),\qquad
 \kappa_{ab}=\kappa(x_a,x_b).
\]
In the nonfine case a topological trivialization of the twist of the universal
sheaf defines $u_\chi$, up to multiplication by a topological line bundle from
$J$. Write $z_{a,p}=\ch_p(u_a)$, with $z_{a,0}=\rk u_a$.

\begin{theorem}[Finite reconstruction of a numerical boundary pairing]
\label{thm:mark20-reconstruction}
Let $e_i:M_i\to J$ be proper evaluation maps with
$\alpha_i=[M_i]^{\vir}\in A_3(M_i)_{\QQ}$. Set
$\lambda_i(m)=\int_{\operatorname{cl}(\alpha_i)}e_i^*m$
for $m\in H^6(J,\QQ)$. The required monomial insertions belong to the three lists
\begin{equation}\label{eq:mark20-list}
 z_{a,3},\qquad z_{a,2}z_{b,1},\qquad z_{a,1}z_{b,1}z_{c,1}
 \qquad(1\le a,b,c\le24).
\end{equation}
Define
\begin{align*}
 \mathscr L(X,Y)&=
 \sum_{a,b=1}^{24}\kappa_{ab}
 \sum_{\substack{0\le p,q\le3\\p+q>0}}
   (-1)^q(p+q-1)!\,(z_{a,p}\boxtimes z_{b,q}) X^pY^q ,
 \\
 \mathscr D_{33}&=[X^3Y^3]\exp\mathscr L(X,Y).
\end{align*}
Then
\begin{equation}\label{eq:mark20-pairing}
 \deg_J(e_{1*}\alpha_1\cdot e_{2*}\alpha_2)
 =(\lambda_1\otimes\lambda_2)(\mathscr D_{33}).
\end{equation}
In particular the right side can be calculated from the two individual
lists of virtual integrals and the Euler form of the K3. For a marked gluing
satisfying Theorem~\companionref{thm:diagonal}, this is the degree of its virtual class.
\end{theorem}
\begin{proof}
On $J\times S\times J$ pull $u_\chi$ back from the first $J$ and $S$
to obtain $u_L$, and similarly form $u_R$ from $S$ and the second $J$.
Write $\pi$ for the projection to $J\times J$ and $\pi_!$ for the
complex-oriented topological $K$-theory pushforward. Markman's diagonal
formula \cite[Theorem~8, Proposition~24 and Corollary~25]{MarkmanIntegral}
gives
\[
 \operatorname{PD}[\Delta_J]=c_6(\mathcal K),\qquad
 \mathcal K=-\pi_!(u_L^\vee\cdot u_R)
           =-\sum_{a,b}\kappa_{ab}\,u_a^\vee\boxtimes u_b.
\]
The fine case also has the algebraic formula
of \cite[Theorem~1]{MarkmanDiagonal}. In either case the displayed
cohomological class is canonical. With $\ch_{p,q}$ denoting the
$H^{2p}(J)\otimes H^{2q}(J)$ part of the Chern character,
\[
 \ch_{p,q}(\mathcal K)
   =-(-1)^p\sum_{a,b}\kappa_{ab}\,z_{a,p}\boxtimes z_{b,q}.
\]
In the identity $\log c=\sum_{j\ge1}(-1)^{j-1}(j-1)!\ch_j$,
the summand of bidegree $(p,q)$ thus has sign
$(-1)^{p+q-1}(-1)^{p+1}=(-1)^q$.
Put $a_i=e_{i*}\operatorname{cl}(\alpha_i)$. The diagonal identity and the projection formula give
\[
 \deg(e_{1*}\alpha_1\cdot e_{2*}\alpha_2)
 =\langle\operatorname{PD}[\Delta_J],a_1\times a_2\rangle
 =(\lambda_1\otimes\lambda_2)\bigl([X^3Y^3]\exp\mathscr L\bigr).
\]
Only bidegree $(3,3)$ in complex codimension pairs with $a_1\times a_2$. The partitions $3,\ 2+1,\ 1+1+1$ give precisely \eqref{eq:mark20-list}. These monomials span $H^6(J,\QQ)$ by \cite[Theorem~1]{MarkmanIntegral} and the rational Newton identities.

\end{proof}

For comparison, $v^2=4$ gives $J$ of deformation type $S^{[3]}$ by \cite[Theorem~8.1]{YoshiokaDeformation}. G\"ottsche's formula \cite{GottscheBetti} gives
\begin{align*}
 b_6(J)&=[u^3t^6]\prod_{m=1}^3
 \frac{1}{(1-u^mt^{2m-2})(1-u^mt^{2m})^{22}(1-u^mt^{2m+2})}\\
 &=\left(\binom{24}{3}+22\right)+(22^2+2)+22=2554 .
\end{align*}
Fewer insertions can suffice on the subspace spanned by the quartic cycles.
Write $\mathscr D_{33}=\sum_{IJ}d_{IJ}\,m_I\boxtimes m_J$ in the monomials of \eqref{eq:mark20-list}. If $a,b,w$ record their integrals against the quartic, blowup and pencil cycles, the two contributions are
\[
       a^{\mathsf t}d\,b,\qquad a^{\mathsf t}d\,w .
\]
The kernel identity tests $a^{\mathsf t}d(b-Tw)=0$ for every relevant $a$, with $T$ acting on the charge index. This determines which K3-fibration series with insertions must be evaluated. For more general derived boundary objects one needs the virtual class of their proper matching space.

\subsection{Eight middle-degree insertions for a very general quartic}
\label{sec:mark20-eight}
There is a smaller list when the rational Hodge group of
$T=h^\perp\subset H^2(S,\QQ)$ is $\mathrm{SO}(T)$.
We impose this additional hypothesis in this subsection. It holds for a very
general polarized quartic: for a nonscalar rational endomorphism of $T$ the
period lines it preserves form a proper closed subset of the period domain
(a finite union of eigenspaces intersected with the open period quadric).

There are only countably many rational endomorphisms. Outside these loci
$\End_{\rm HS}(T)=\QQ$, and \cite[\S2]{ZarhinHG} gives the asserted Hodge group.
Let $p$ be the point class on $S$, put $w=(1,0,-2)$, and define
\[
 P_T=[\Delta_S]-1\boxtimes p-p\boxtimes1-\tfrac14h\boxtimes h
     =\sum_{i=1}^{21} t_i\boxtimes t^i
     \quad\text{in } H^4(S^2,\QQ)
\]
for any pair of dual bases of $T$.
For the Mukai pairing, let $s_z(x)=x-2(x,z)z/(z,z)$.
As $v^2=w^2=4$, an explicit rational Mukai Hodge isometry with $g(v)=w$ is
\[
 g=\begin{cases}
   s_{v-w},&(v-w)^2\ne0,\\
   -s_{v+w},&(v-w)^2=0.
   \end{cases}
\]
In the second case $(v+w)^2=16$. Let
$\gamma(g):H^*(J,\QQ)\xrightarrow{\sim}H^*(S^{[3]},\QQ)$
be the cohomology operator associated to $g$ by Markman \cite[\S3]{MarkmanDiagonal}.

\begin{theorem}[Eight insertions for numerical gluing]\label{thm:mark20-eight}
Write $q_j(a)$ for the unscaled Nakajima creation operator of length $j$,
let products act on the vacuum $|0\rangle$, and sum any repeated index $i$
over dual bases as above. On $S^{[3]}$ set, in the indicated order,
\[
\begin{array}{ll}
 b_1=q_1(h)^3,& b_2=\sum_iq_1(h)q_1(t_i)q_1(t^i),\\
 b_3=q_1(1)q_1(h)q_1(p),& b_4=q_2(h)q_1(h),\\
 b_5=\sum_iq_2(t_i)q_1(t^i),& b_6=q_2(1)q_1(p),\\
 b_7=q_2(p)q_1(1),& b_8=q_3(h).
\end{array}
\]
Their inverse images $\beta_i=\gamma(g)^{-1}b_i$ form a basis of
$H^6(J,\QQ)\cap H^{3,3}(J)$. For $A,B\in A_3(J)_{\QQ}$ put
$a_i=\int_{\operatorname{cl}(A)}\beta_i$ and
$b_i'=\int_{\operatorname{cl}(B)}\beta_i$. Then
\begin{equation}\label{eq:mark20-eight}
\begin{split}
 \deg(A\cdot B)={}&\frac{a_1 b_1'}{384}
 +\frac{a_2 b_2'}{168}+\frac{a_3 b_3'}4
 -\frac{a_4 b_4'}{32}-\frac{a_5 b_5'}{42}\\
 &-\frac{a_6 b_7'+a_7 b_6'}2+\frac{a_8 b_8'}{12}.
\end{split}
\end{equation}
\end{theorem}
\begin{proof}
Markman's formula preserves degree and the cup product \cite[Remark~2.1, Theorem~2.4, (2.3)]{OberdieckUni}. To check Hodge type, grade $H^*(S,\CC)$ by $V_d=\bigoplus_{p-q=d}H^{p,q}(S)$. Each reflection defining $g$ preserves $V_d$: its reflecting vector is in $\QQ 1\oplus\QQ h\oplus\QQ p$ and it acts by $\pm\id$ on $T$.

The normalized universal characters have total $p-q=0$, including in the nonfine case \cite[\S3]{MarkmanDiagonal}, \cite[\S2.1]{OberdieckUni}. Pullback, contraction over $S$, duality and products preserve this condition. The degree-twelve kernel thus has type $(6,6)$, giving $\gamma(g)(H^{p,q})\subset H^{p,q}$.

By the Nakajima basis theorem and the algebraicity of its incidence
correspondences, the middle cohomology, grouped by partitions
of $3$, has the following rational Hodge-group representations,
suppressing Tate twists. Write $V=\QQ h\oplus T$.
\[
\begin{array}{c|c|c}
 \text{partition}&\text{summands in }H^6&\dim(\text{invariants})\\ \hline
 1+1+1 &\operatorname{Sym}^3V\oplus V &(1+1)+1\\
 2+1   &V\otimes V\oplus\QQ^2 &(1+1)+2\\
 3     &V&1 .
\end{array}
\]
Indeed $q_j(x)$ has degree $\deg(x)+2j-2$.
In tensor degrees at most three the only invariant tensors of the
standard $21$-dimensional special orthogonal representation
besides scalars are multiples of the inverse form in $T\otimes T$.
This gives the displayed eight basis vectors, including $P_T$
in $b_2,b_5$. Their Gram matrix for the cup product is block diagonal:
\[
 G=\operatorname{diag}(384,168,4,-32,-42)
      \ \oplus\begin{pmatrix}0&-2\\-2&0\end{pmatrix}
      \ \oplus (12).
\]
For completeness, the adjoint commutator for positive creation indices is
$[q_j(x)^\dagger,q_k(y)]
 =\delta_{jk}(-1)^{j-1}j\int_Sxy$,
with $\dagger$ the ordinary Poincar\'e adjoint
\cite[Theorem~8.13 and (8.14)]{NakajimaLectures}.
It pairs distinct partitions to zero. Within partition $1+1+1$ the
three self-pairings are $3!\,4^3$, $4\cdot2\cdot21$, and $4$;
within $2+1$ each contraction has factor $-2$;
within partition $3$ the self-pairing is $3\cdot4$.

The isometry $\gamma(g)$ therefore gives the same $G$ on $J$.
The Poincar\'e duals of $A,B$ lie in this invariant subspace;
their coordinate vectors are $G^{-1}a,G^{-1}b'$. Thus
$\int_J\operatorname{PD}(A)\operatorname{PD}(B)=a^{\mathsf t}G^{-1}b'$.
\end{proof}
For a relative virtual pushforward the eight entries are the respective
integrals of $e^*\beta_i$ on the relative virtual class. Formula
\eqref{eq:mark20-eight} specifies the insertions required on the smooth
primitive boundary. Its coefficients are independent of the threefolds
and of the relative virtual multiplicities.

\begin{corollary}[Coprime multiplicities on the ordinary K3 boundary]
\label{cor:mark20-higher}
Suppose $\Pic(S)=\mathbb Zh$, $h^2=4$. Let $b\geq1$ and
$\gcd(b,\chi)=1$, and put
\[
 v_b=(0,bh,\chi),\qquad N=2b^2+1,\qquad J_{b,\chi}=M_h(v_b).
\]
The moduli of pure stable sheaves $J_{b,\chi}$ is a smooth projective
$2N$-fold. For any two proper evaluations to $J_{b,\chi}$ with
virtual classes of dimension $N$, formula \eqref{eq:mark20-pairing}
holds with $\mathscr D_{33}$ replaced by
\[
 [X^NY^N]\exp\!\left(
   \sum_{a,d=1}^{24}\kappa_{ad}
   \sum_{\substack{0\leq i,j\leq N\\i+j>0}}
    (-1)^j(i+j-1)!\,(z_{a,i}\boxtimes z_{d,j})X^iY^j
  \right).
\]
Here the $z_{a,i}$ are taken for $J_{b,\chi}$. The necessary insertions
are the monomials $\prod_{k=1}^l z_{a_k,i_k}$ with
$i_k\geq1$ and $\sum_k i_k=N$.
\end{corollary}
\begin{proof}
Write $P_E(m)=4bm+\chi$. If $0\ne E'\subsetneq E$ has the same reduced polynomial, write $c_1(E')=dh$, $1\le d\le b$. Then
\[
 b\chi(E')=d\chi,\quad \gcd(b,\chi)=1
 \ \Longrightarrow\ d=b,\quad P_{E/E'}=0 ,
\]
contrary to $E/E'\ne0$. Stability thus equals semistability, and projectivity follows from \cite[Theorem~1.21]{Simpson}. For a smooth $i:C\hookrightarrow S$ in $|bh|$ take $L\in\Pic^{\chi+g(C)-1}(C)$. Every nonzero proper subsheaf of $i_*L$ has quotient of positive length. This supplies a stable point.

At any such stable $E$, Serre duality identifies the obstruction trace with the dual of $\CC\xrightarrow{1\mapsto\id_E}\End(E)$. Hence
\[
 \Ext^2(E,E)_0=0,\qquad
 \ext^1(E,E)=2+\langle v_b,v_b\rangle=2+4b^2=2N .
\]
The same polynomial test shows that $h$ is $v_b$-generic. Markman's formula \cite[Theorem~8 and Proposition~24]{MarkmanIntegral} gives
$c_{2N}(-\pi_!(u_L^\vee u_R))=\operatorname{PD}[\Delta_{J_{b,\chi}}]$.
Apply the preceding diagonal calculation in bidegree $(N,N)$. For example, $b=2$ and odd $\chi$ give $\dim J=18$, with insertions of codimension nine on each branch. The splitting of the threefold character fixes the boundary multiplicity $b$.

\end{proof}

\subsection{Computing the insertions on a family of perfect restrictions}
Suppose a relative universal sheaf $\mathcal F$ on a threefold
restricts at the boundary to a perfect complex
$E=\mathbf L i^*\mathcal F$ on $S\times T$. One can form cohomology
insertions on $T$ before identifying a moduli space for $E$. We use
the projections $p:S\times T\to T$, $q:S\times T\to S$ and suppose that
the Mukai vector of each fiber of $E$ is $(0,bh,\chi)$, $b>0$.
\begin{proposition}\label{prop:mark20-normalized}
For $x\in K^0_{\rm top}(S)$ set
\[
 P_x(E)=p_!\big(q^*x^\vee\cdot[E]\big),\quad
 x_0=[\OO_S(-h)]-[\OO_S],\quad
 \ell(E)=\frac{c_1(P_{x_0}(E))}{4b}.
\]
For algebraic $x$ the index is the class of
$\mathbf Rp_*(q^*x^\vee\otimes^{\mathbf L}E)$. The classes
\begin{equation}\label{eq:mark20-relativeinput}
 Z_x(E)=e^{-\ell(E)}\ch(P_x(E))
  =p_*\!\left(q^*(\ch(x^\vee)\operatorname{td}(S))
          \,\ch(E)e^{-p^*\ell(E)}\right)
\end{equation}
are invariant under $E\mapsto E\otimes p^*L$ for a line bundle $L$
on $T$. On an ordinary stable boundary as in
Corollary~\ref{cor:mark20-higher}, take $x=y_a$ in the
basis Euler-dual to $x_a$, so $\kappa(y_a,x_d)=\delta_{ad}$.
In the numerical diagonal formula one may replace each $z_{a,j}$ by
$[Z_{y_a}]_{2j}$, where $[\ ]_{2j}$ denotes cohomological degree $2j$.
\end{proposition}
\begin{proof}
The Euler index and the projection formula give
\[
 \rk P_{x_0}(E)=\chi(E_t(h))-\chi(E_t)=4b,\qquad
 P_x(E\otimes p^*L)=P_x(E)\cdot [L].
\]
Thus $\ell(E\otimes p^*L)=\ell(E)+c_1(L)$.
Riemann--Roch for the smooth projection $p$ gives
\eqref{eq:mark20-relativeinput}. This calculation uses the class
$[E]$ and remains valid for complexes with cohomology in two degrees.
On an ordinary boundary $J=J_{b,\chi}$ use the topological universal
class of Theorem~\ref{thm:mark20-reconstruction}; the line-bundle
independence makes the expression independent of its trivialization.

The expansion in the Euler-dual basis gives $P_{y_a}(u)=u_a$.
For $m=\dim J$ write $K=-\pi_!(u_L^\vee u_R)$ and
$t=\operatorname{pr}_1^*\ell(u)-\operatorname{pr}_2^*\ell(u)$.
The normalized external-product expression replaces $\ch(K)$ by
$\ch(K)e^t$. By Riemann--Roch, $\rk K=m-2$. The splitting principle
for the degree-$m$ Chern class of a rational formal twist gives
\[
 \sum_{j=0}^m \binom{m-2-j}{m-j}t^{m-j}c_j(K)
 =c_m(K)-t c_{m-1}(K).
\]
By \cite[Lemma~9 and Proposition~24(1),(2)]{MarkmanIntegral},
$c_{m-1}(K)=0$ and $c_m(K)=\operatorname{PD}[\Delta_J]$.
Extracting the middle bidegree therefore gives the same numerical
pairing.
\end{proof}
For a proper relative moduli space $M$ with a virtual $N$-cycle,
the quantities to evaluate can now be written on the original family as
\[
   \int_{[M]^{\vir}}
   \prod_k [Z_{y_{a_k}}(\mathbf L i^*\mathcal F)]_{2j_k},
       \qquad \sum_k j_k=N .
\]
When the family only exists on a universal-sheaf gerbe, the expressions
can be pulled back from $J$ on the ordinary evaluation locus; on a chart
they are given by the displayed line-bundle-independent formula. For a
general derived boundary the same family classes can be tested on charts.

An identity with a gluing degree there requires its own proper matching
space and compatible obstruction theories. Indeed all the displayed
classes depend on $[E]$ and take the same values on $E$ and
$E\oplus P\oplus P[1]$ for any perfect $P$. Such complexes can
have different automorphisms and deformation spaces.

\subsection{The affine correction and the modular conclusion}
\begin{theorem}\label{thm:twodegmodularityfull}
Assume that the proper relative spaces and their obstruction theories give \eqref{eq:twovirtualstatesfull} with a common quartic class $\mathcal A$. Make the following assumptions on their evaluations.
\begin{enumerate}
\item A geometric correspondence gives \eqref{eq:coefficientfunctionalidentityfull} with $T=T_\Theta$ for an independently constructed $\Theta\in M_{31/2}(\rho_N)$.
\item In the very general pencil compute $c_Y$ from the point-ideal virtual degrees by \eqref{eq:pencil-neutral-reduction}, with zero components at nonzero residues, so that $K^{\rm raw}+c_Y=K^{\rm per}$ with the normalization \eqref{eq:Kcorrectednew}.
\item The remaining relative sectors and this pencil correction have the Fourier expansion of a holomorphic modular form $U$ in the following identity:
\begin{equation}\label{eq:Ugeometricfull}
 U=\eta^{55}\bigl(R_X-T_\Theta(R_Y+c_Y)\bigr)\in M_{26}(\rho_5).
\end{equation}
\end{enumerate}

Then the primitive quintic DT series satisfies S-duality.  More precisely,
\begin{equation}\label{eq:fullmodularidentity}
 \eta^{55}h^{\mathrm{DT}}_1
 =\mathfrak g(\widehat\Phi\otimes\Theta)+U,
\end{equation}
with the charge conventions above.
\end{theorem}
\begin{proof}
By \eqref{eq:twovirtualstatesfull} and the first two hypotheses,
\[
 h^{\mathrm{DT}}_1-R_X=T_\Theta\mathcal B\mathcal A
   =T_\Theta(K^{\rm raw}-R_Y)
   =T_\Theta(K^{\rm per}-c_Y-R_Y).
\]
Substitute $K^{\rm per}=\eta^{-24}\widehat\Phi$ in \eqref{eq:operatorTfull}. The eta powers cancel:
\[
 \eta^{55}h^{\mathrm{DT}}_1
 =\eta^{55-31-24}\mathfrak g(\widehat\Phi\otimes\Theta)+U
 \in M_{26}(\rho_5)
\]
by Proposition~\ref{prop:intertwinerfull} and \eqref{eq:Ugeometricfull}.
Since $\eta$ is nonvanishing on the upper half-plane and has cusp order $1/24$, division by $\eta^{55}$ gives the transformation laws and pole bound of Proposition~\ref{prop:numeratorsuffices}.
\end{proof}

\begin{proposition}[The primitive principal part]\label{prop:polar-from-specializations}
With $\operatorname{pp}$ denoting the sum of the terms of strictly negative Fourier exponent, the following identities hold:
\begin{align*}
 \operatorname{pp}(h^{\mathrm{DT}}_{1,0})&=5q^{-55/24}-800q^{-31/24}
                                     +58\,500q^{-7/24},\\
 \operatorname{pp}(h^{\mathrm{DT}}_{1,1})&=
 \operatorname{pp}(h^{\mathrm{DT}}_{1,4})=8625q^{-83/120},\\
 \operatorname{pp}(h^{\mathrm{DT}}_{1,2})&=
 \operatorname{pp}(h^{\mathrm{DT}}_{1,3})=0 .
\end{align*}
\end{proposition}
\begin{proof}
Use $x=25/12+(2e-5)^2/40-\chi$ and $\mu=-e\bmod5$. The hull estimates in Lemmas~\ref{lem:polar-degree-zero} and~\ref{lem:polar-small-negative-degree} give
\[
\begin{array}{c|rrrrr}
e&0&-1&1&-2&2\\ \hline
\max\chi\ \text{allowed by the estimates}&5&4&3&4&2\\
\min x\ \text{allowed by the estimates}&-55/24&-83/120&-83/120&13/120&13/120 .
\end{array}
\]
Reducing $\chi$ by one increases $x$ by one. The negative exponents are therefore exactly those displayed once the nonzero coefficients are evaluated.

For the moduli in Proposition~\ref{prop:polar-incidence-functors}, Theorems~\ref{thm:vacuumvirtual} and~\ref{thm:direct-line} give the structure-sheaf and line degrees. Propositions~\ref{prop:one-defect-positive-marked} and~\ref{prop:one-defect-zero-count} give $-320-360-120=-800$ for one point. For two points the full Chow specialization in Theorem~\ref{thm:two-points-direct} gives $58\,500$. This gives the principal part in the convention of \cite[\S3.3, (3.10)]{GSY}.
\end{proof}

\begin{corollary}[From the transformation law to the GSY identity]\label{prop:principalpart-criterion}
Suppose $h^{\mathrm{DT}}_1$ is weakly holomorphic of weight $-3/2$ with the multiplier in \eqref{eq:Tmodnew}--\eqref{eq:Smodnew}. This follows under the hypotheses of Theorem~\ref{thm:twodegmodularityfull}, or from \eqref{eq:quintic-residual-period}. Then $[q^x]h^{\mathrm{DT}}_{1,\mu}=[q^x]Z^{\mathrm{GSY}}_\mu$ for every $\mu\in\mathbb Z/5$ and every Fourier exponent $x$, with $Z^{\mathrm{GSY}}$ as in \cite[(3.10)--(3.11)]{GSY}.
\end{corollary}
\begin{proof}
Let $g=h^{\mathrm{DT}}_1-Z^{\mathrm{GSY}}$, $k=-3/2$, $y=\operatorname{Im}\tau$. By Proposition~\ref{prop:polar-from-specializations}, $g$ is bounded at the cusp. Both vectors have the same unitary multiplier by hypothesis, so $u=y^k\|g\|^2$ is invariant and tends to zero at the cusp. If $u$ is positive somewhere, it attains a positive maximum in the upper half-plane modulo the modular group.

At any lift of such a maximum, $g\ne0$. With $\Delta=4\partial_\tau\partial_{\bar\tau}$,
\[
 \Delta\log u=-\frac{k}{y^2}
    +4\,\frac{\|g\|^2\|g'\|^2-|\langle g',g\rangle|^2}{\|g\|^4}>0
\]
by Cauchy--Schwarz and $k<0$. This contradicts $\Delta\log u\le0$ at the maximum. Hence $g=0$.
\end{proof}

\subsection{Reconstruction of the primitive modular  from four marked degrees}\label{sec:step21-finite}
The four marked degrees determine a vector in the prescribed modular space. The two-point and mixed line degrees and the three smooth-quintic evaluations in Section~\ref{sec:smooth-evaluations} give five independent checks.

Use the charge-symmetric components indexed by $\mu=0,1,2$, extended by
$f_{5-\mu}=f_\mu$. For a vector satisfying
\eqref{eq:Tmodnew}--\eqref{eq:Smodnew} this symmetry follows by applying
the second transformation twice. Set $q=e^{2\pi i\tau}$,
$D=q\,d/dq$, $\nabla_w=D-wE_2/12$, and let $\nabla^{(j)}$ mean
$j$ successive Serre derivatives beginning in weight $3/2$. The unary
theta series of \cite[(5.2)--(5.5)]{AFKPS} in these conventions is
\[
 \theta_\mu=\tfrac12\sum_{s\equiv 2\mu+5\ (10)}
      (-1)^{(s-1)/2}s q^{s^2/40}.
\]
Its leading exponents and those of $\eta^{-70}\theta_\mu$ are
\[
 b=(25,9,1)/40,\qquad s=b-70/24=(-275,-323,-347)/120 .
\]
For any $f$ put $c_{\mu,n}(f)=[q^{s_\mu+n}]f_\mu$.
The hull bounds in Proposition~\ref{prop:polar-from-specializations}
give the five vanishings in the row list
\[
 I_0=((1,0),(1,1),(2,0),(2,1),(2,2))
\]
and exclude all $n<0$. Set $I_1=((0,0),(0,1),(1,2),(2,3))$.
The last entry has positive Fourier exponent $13/120$.

\begin{theorem}[Four marked degrees determine the primitive modular candidate]
\label{thm:step21-four-seeds}
There exists a unique weakly holomorphic modular vector $f$ with the
two stated primitive transformation laws, $c_{\mu,n}=0$ for $n<0$ and
$(c_I)_{I\in I_0}=0$, having any prescribed values
$(a,b',c,d)=(c_I)_{I\in I_1}$. In particular,
\begin{align}
 c_{0,2}(f)&=-1255a-129b'+(119c+d)/5,\label{eq:step21-neutral}\\
 c_{1,3}(f)&=109375a+3125b'+377c+2d.\label{eq:step21-mixed}
\end{align}
For the four marked degrees of the structure sheaf, one point, a
line and a conic, write $F^{\mathrm{DT}}$ for this vector. Then $F^{\mathrm{DT}}=Z^{\mathrm{GSY}}$, with $Z^{\mathrm{GSY}}$ as in \cite[(3.10)--(3.11)]{GSY}. The two formulae give
$58\,500$ and $-1\,138\,500$, respectively. If the geometric series
$h^{\mathrm{DT}}_1$ satisfies the stated modular laws and analytic conditions, then
$h^{\mathrm{DT}}_1=F^{\mathrm{DT}}$ in every charge.
\end{theorem}
\begin{proof}
Let $C$ be the $3\times3$ matrix with columns
$\theta,\nabla^{(1)}\theta,\nabla^{(2)}\theta$. These have weights
$3/2,7/2,11/2$ and the same multiplier. The leading coefficients of
the three rows of $\theta$ are $5,-3/2,1/2$, respectively. Factoring
out the row powers $q^{b_\mu}$ leaves a constant determinant
\[
 (5)(-3/2)(1/2)\prod_{\mu<\nu}(b_\nu-b_\mu)\ne0 .
\]
Indeed, the leading coefficients in the three columns are given by
monic polynomials of degrees $0,1,2$ in $b_\mu$. Hence
$\det(C)$ has weight $21/2$ and order $7/8$ at the cusp.
Its quotient by $\eta^{21}$ is holomorphic on the upper half-plane
and bounded at the cusp.

Its weight-zero character is trivial: the
ratio of the multipliers is an ordinary character, trivial on $T$
by the order, and $T$ generates the abelianization of
$\mathrm{SL}_2(\mathbb Z)$. Thus $\det(C)/\eta^{21}$ is a nonzero
constant. In particular $C$ is invertible throughout the upper half-plane.

For a modular $f$ under consideration, $\eta^{70}f$ has weight $67/2$
and the multiplier of $\theta$ by the unary theta transformation in
\cite[(5.2)--(5.5)]{AFKPS}. Cramer's rule gives
\[
 \eta^{70}f=\sum_{j=0}^2 P_j\nabla^{(j)}\theta,\qquad
                P_j\in M_{32-2j}(\mathrm{SL}_2(\mathbb Z)).
\]
In fact each determinant numerator has row orders at least $b_\mu$,
so the scalar quotient is holomorphic at the cusp as well.
Conversely every such triple produces a vector with the required
transformations and no term below $s_\mu$. For $w_j=32-2j$ put
$\epsilon_j=1$ if $w_j\equiv2\pmod4$, and $0$ otherwise.
With $\Delta=\eta^{24}$ define
\[
 U_{tj}=\eta^{-70}\Delta^t E_4^{(w_j-6\epsilon_j-12t)/4}
                  E_6^{\epsilon_j}\nabla^{(j)}\theta,
                   \qquad 0\le t,j\le2 .
\]
For each $j$ the identity $E_6^2=E_4^3-1728\Delta$ spans $M_{w_j}$ by the three displayed scalar monomials. Their cusp orders are $0,1,2$, so they are independent.

Let $K_{I,(t,j)}=c_I(U_{tj})$ with rows $I_0,I_1$ and columns
in lexicographic order $(t,j)$. Exact expansion through row index
$n=3$ gives
\begin{align*}
 \det K&=729/25000,\\
 \big(c_{0,2}(U_{tj})\big)_{tj} K^{-1}\big|_{I_1}
      &=(-1255,-129,119/5,1/5),\\
 \big(c_{1,3}(U_{tj})\big)_{tj} K^{-1}\big|_{I_1}
      &=(109375,3125,377,2).
\end{align*}
For explicit evaluation, use $E_2=1-24\sum_{n\ge1}\sigma_1(n)q^n$,
$E_4=1+240\sum\sigma_3(n)q^n$, $E_6=1-504\sum\sigma_5(n)q^n$,
and $\eta=q^{1/24}\prod_{m\ge1}(1-q^m)$. Each derivative acts on the
row $q^\gamma\sum u_nq^n$ by
\[
 u_n\longmapsto(\gamma+n)u_n-
        \tfrac{w}{12}\sum_{l=0}^n [q^l]E_2\,u_{n-l}.
\]
These equations specify the finite matrix without Fourier coefficients of the physics candidate. After removing each leading power, it suffices to truncate each product at $q^3$, form $K$ in the stated row and column orders, and solve $K u=(0,0,0,0,0,a,b',c,d)^{\mathsf t}$ over $\QQ$.

For $(a,b',c,d)=(5,-800,8625,-1218500)$, writing $F^{\mathrm{DT}}=\sum u_{tj}U_{tj}$
gives
\[
 (u_{tj})_{t,j=0}^2=
 \begin{pmatrix}
 -1/24&25/36&25/6\\
 1675/24&-10625/12&-15625/2\\
 -177095/12&1672375/18&7064875/3
 \end{pmatrix}.
\]
Using $1728\Delta=E_4^3-E_6^2$ gives \cite[(3.10)--(3.11)]{GSY}, equivalently \cite[(5.7)]{AFKPS}. Two further Fourier coefficients are evaluated by \eqref{eq:step21-neutral}--\eqref{eq:step21-mixed}. They agree with Theorems~\ref{thm:two-points-direct} and~\ref{thm:mixed-all-marked} in the indicated charges and their conjugates. If the DT vector belongs to the stated modular space, its four marked degrees and five vanishings give the same $u$ by invertibility of $K$.

\end{proof}
For a formal DT vector with the established support bound there is a
concrete alternative target for the geometric proof. For its three
representatives define
\[
 (P_0^{\rm DT},P_1^{\rm DT},P_2^{\rm DT})^{\mathsf t}
       = C^{-1}\eta^{70}(h^{\mathrm{DT}}_{1,0},h^{\mathrm{DT}}_{1,1},h^{\mathrm{DT}}_{1,2})^{\mathsf t}.
\]
After removing $q^{b_\mu}$ row by row the constant matrix is
invertible, so all three scalar series belong to $\QQ[[q]]$.
By Theorem~\ref{thm:rank-zero-charge-conjugation}, showing $P_j^{\rm DT}\in M_{32-2j}$ by the two degenerations would prove all-orders modularity. The theorem then
supplies every coefficient from the four marked degrees. Finite
agreement alone leaves a formal freedom: a term
$q^{s_0+N}$ in the neutral component for arbitrarily large $N$
preserves those comparisons.

For the fixed-polynomial Hall logarithm in \textup{(\companionref{eq:Hallloggeneral})},
$r=1$ and $r_i\ge1$ force $j=1$. For the positive surface classes
occurring there at any $r$, equality of slopes gives
$r_2d_1-r_1d_2=0$ in \textup{(\companionref{eq:generalEulerpairing})}. Thus this
Hall logarithm supplies no induction in the D0 coordinate of a
primitive surface sheaf.

Wall crossing through other stability
conditions introduces invariants of the destabilizing classes and
requires separate inputs for them. For example,
\cite[Theorem~1.1]{FeyzOSV} expresses rank-zero DT invariants using
rank-one DT and stable-pair invariants in its range of validity.
\begin{corollary}[Finite test under a depth bound]\label{cor:step21-depth}
Suppose a geometric calculation places each $P_j^{\rm DT}$ in the space
$M_{w_j}+E_2 M_{w_j-2}$ of homogeneous quasimodular forms of weight
$w_j$ and depth at most one. Then agreement with $F^{\mathrm{DT}}$ at the first
twelve Fourier positions allowed by the support bounds implies equality
in every charge. In the notation $[q^{s_\mu+n}]$ these are
\[
 \begin{array}{c|c}
 \mu&n\\\hline
 0&0,1,2,3,4\\
 1&2,3,4,5\\
 2&3,4,5.
 \end{array}
\]
\end{corollary}
\begin{proof}
Use the seventeen columns
\[
 \eta^{-70} E_2^p\Delta^tE_4^{(w_j-2p-6\epsilon-12t)/4}
 E_6^\epsilon\nabla^{(j)}\theta ,
\]
ordered lexicographically by $(j,p,t)$ with $j=0,1,2$, $p=0,1$,
$w_j-2p\equiv2\epsilon\pmod4$, $\epsilon\in\{0,1\}$, and
$t\ge0$ chosen to give a nonnegative exponent of $E_4$.
The three column counts are $6,6,5$. Let $A$ be the resulting $17\times17$ matrix with rows $I_0$ followed by the twelve stated positions in increasing order of $s_\mu+n$. The recurrence for the entries in Theorem~\ref{thm:step21-four-seeds} gives
\[
 \det A=\frac{2^{36}3^{35}7^7 11^6 13^5 17^4 19^4 23^3 29\cdot31}{5^2}
 \ne0 .
\]
For this determinant the same row recurrence and product expansions, truncated at index $n=5$, give every entry over $\QQ$. By the depth hypothesis, $h^{\mathrm{DT}}_1-F^{\mathrm{DT}}$ has a column-coordinate vector $u$ in the three representative residues. The five vanishings and the twelve equalities give $Au=0$, hence $u=0$. Theorem~\ref{thm:rank-zero-charge-conjugation} gives the other two residues.

\end{proof}
The geometric hypothesis here permits a nonzero $E_2$ anomaly prior
to the finite test. The fixed-K3 series of
Proposition~\ref{prop:mark20-k3bracket} is an example where an
independent quasimodular weight bound enables interpolation. A depth
and weight estimate for the three scalar series of the quintic would
require control of all relative virtual contributions.

Table~\ref{tab:twenty-rev27} gives the reconstructed modular candidate at twenty Fourier positions. Its first twelve rows give the finite test under the depth hypothesis above.

\section{Comparison with the quintic elliptic genus}\label{sec:verified-fourier-table}
Table~\ref{tab:intro-charge-comparison} counts each twist orbit once. The marked degrees at positions $1,2,3,5,6$ are the following sums:
\[
\begin{aligned}
  5 &=4+1, \qquad -800=-320-360-120,\\
  8625 &=3(1600)+(2+1)(1275),\\
  -1\,218\,500 &=-516\,400-326\,400-187\,850-187\,850,\\
  -1\,138\,500 &=1600(-158-238)+1275(-160-114-122).
\end{aligned}
\]
The line terms use the four-secant degrees of Proposition~\ref{prop:four-secant-sheaf-marked}. For the conic the first term is $-50(5016)+10(-26\,560)$ (Corollary~\ref{cor:mark20-conictest}); Theorem~\ref{thm:marked-conics-18} gives the remaining three. In the two mixed charges Theorem~\ref{thm:mixed-all-marked} identifies the quotient charts and restricted boundary diagonals, including the point modification in the opposite residue.

Theorem~\ref{thm:two-points-direct} gives $58\,500$ at position $4$ by Chow specialization. Section~\ref{sec:smooth-evaluations} gives positions $7$--$9$ from the full moduli on the smooth quintic, including the additional hulls of Propositions~\ref{prop:third-neutral-additional-hulls} and~\ref{prop:additional-conic-hulls}. The four geometric inputs to Theorem~\ref{thm:step21-four-seeds} also occur in \cite[text preceding (3.10)]{GSY}; the mixed value occurs in \cite[(3.5), (3.10)]{GSY}.

For all primitive charges Corollary~\ref{prop:principalpart-criterion} applies once the transformation law is proved. Theorem~\ref{thm:twodegmodularityfull} reduces this to a common proper quartic class, $\mathcal C=T_\Theta\mathcal B$ and \eqref{eq:Ugeometricfull}. The direct period comparison \eqref{eq:quintic-residual-period} is another sufficient condition. For $r>1$ the charge symmetry and the wall-crossing range are given in Theorems~\ref{thm:rank-zero-charge-conjugation} and~\ref{thm:hall22strip}.

\subsection{The first twenty Fourier positions}
Let $F^{\mathrm{DT}}$ be the modular vector determined by four marked degrees and the geometric support bound (Theorem~\ref{thm:step21-four-seeds}). Table~\ref{tab:twenty-rev27} compares its expansion with \cite[(3.10)--(3.11)]{GSY}, evaluated in the theta-derivative basis of \cite[(5.2), (5.7)]{AFKPS}. The individual surface-sheaf degrees proved here are listed in Table~\ref{tab:intro-charge-comparison}.

\begin{table}[!htbp]\centering\small
\caption{Fourier coefficients of the reconstructed modular form $F^{\mathrm{DT}}$ and the physics form $Z^{\mathrm{GSY}}$. Identification of the full geometric DT degree is proved for the rows of Table~\ref{tab:intro-charge-comparison}.}\label{tab:twenty-rev27}
\begin{tabular}{@{}rccrr@{}}\toprule
No. & $\mu$ & $x$ & $[q^x]F^{\mathrm{DT}}_\mu$ & $[q^x]Z^{\mathrm{GSY}}_\mu$\\\midrule
1 & $0$ & $-55/24$ & $5$ & $5$\\
2 & $0$ & $-31/24$ & $-800$ & $-800$\\
3 & $1,4$ & $-83/120$ & $8\,625$ & $8\,625$\\
4 & $0$ & $-7/24$ & $58\,500$ & $58\,500$\\
5 & $2,3$ & $13/120$ & $-1\,218\,500$ & $-1\,218\,500$\\
6 & $1,4$ & $37/120$ & $-1\,138\,500$ & $-1\,138\,500$\\
7 & $0$ & $17/24$ & $5\,817\,125$ & $5\,817\,125$\\
8 & $2,3$ & $133/120$ & $441\,969\,250$ & $441\,969\,250$\\
9 & $1,4$ & $157/120$ & $3\,777\,474\,000$ & $3\,777\,474\,000$\\
10 & $0$ & $41/24$ & $75\,474\,060\,100$ & $75\,474\,060\,100$\\
11 & $2,3$ & $253/120$ & $953\,712\,511\,250$ & $953\,712\,511\,250$\\
12 & $1,4$ & $277/120$ & $3\,102\,750\,380\,125$ & $3\,102\,750\,380\,125$\\
13 & $0$ & $65/24$ & $28\,096\,675\,153\,255$ & $28\,096\,675\,153\,255$\\
14 & $2,3$ & $373/120$ & $217\,571\,250\,023\,750$ & $217\,571\,250\,023\,750$\\
15 & $1,4$ & $397/120$ & $577\,727\,215\,123\,000$ & $577\,727\,215\,123\,000$\\
16 & $0$ & $89/24$ & $3\,756\,542\,229\,485\,475$ & $3\,756\,542\,229\,485\,475$\\
17 & $2,3$ & $493/120$ & $22\,258\,695\,264\,509\,625$ & $22\,258\,695\,264\,509\,625$\\
18 & $1,4$ & $517/120$ & $52\,559\,194\,851\,824\,125$ & $52\,559\,194\,851\,824\,125$\\
19 & $0$ & $113/24$ & $277\,591\,744\,202\,815\,875$ & $277\,591\,744\,202\,815\,875$\\
20 & $2,3$ & $613/120$ & $1\,374\,043\,315\,791\,020\,500$ & $1\,374\,043\,315\,791\,020\,500$\\
\bottomrule\end{tabular}\end{table}

Uniqueness within the modular space gives $F^{\mathrm{DT}}=Z^{\mathrm{GSY}}$. The computed complete sheaf invariants give nine entries of the proposed identity $h^{\mathrm{DT}}_1=Z^{\mathrm{GSY}}$. The transverse hulls in Corollary~\ref{cor:twoskew-all-lengths} and \eqref{eq:quintic-cubic-period-sector} give specified contributions at infinitely many exponents. 

\subsection*{AI use declaration}
AI assisted with proof readings, the figures and some of the numerical verifications.

\end{document}